\documentclass[12pt,letter,leqno]{article}
\usepackage{natbib}
\setcitestyle{authoryear,open={(},close={)},citesep={,}}
\usepackage{fancyheadings}
\usepackage[]{hyperref}
\usepackage{xcolor}
\hypersetup{colorlinks = true, citecolor = blue, linkcolor=blue, linkbordercolor = {white}}
\usepackage{amsmath}
\usepackage{amsfonts}
\usepackage{mathrsfs}
\usepackage{amsthm}
\usepackage{amssymb}
\usepackage{mathtools}
\usepackage{bm}
\usepackage[normalem]{ulem}
\usepackage{graphicx}
\usepackage{caption}

\allowdisplaybreaks

\newtheorem{theorem}{Theorem}

\newtheorem{proposition}{Proposition}

\newtheorem{example}{Example}
\newtheorem{remark}{Remark}
\newtheorem{condition}{Condition}

\newenvironment{customcon}[1]
{\renewcommand\thecondition{#1}\condition}
{\endcondition}

\newtheorem{innercustomgeneric}{\customgenericname}
\providecommand{\customgenericname}{}
\newcommand{\newcustomtheorem}[2]{%
  \newenvironment{#1}[1]
  {%
   \renewcommand\customgenericname{#2}%
   \renewcommand\theinnercustomgeneric{##1}%
   \innercustomgeneric
  }
  {\endinnercustomgeneric}
}
\newcustomtheorem{customlemma}{Lemma}
\newcustomtheorem{customproposition}{Proposition}
\newcustomtheorem{customtheorem}{Theorem}
\newcustomtheorem{customdefinition}{Definition}
\newcustomtheorem{customexample}{Example}
\newcustomtheorem{customremark}{Remark}

\newcommand{\op}{\oplus}
\newcommand{\om}{\ominus}
\newcommand{\od}{\odot}
\newcommand{\bzero}{\mathbf{0}}

\newcommand{\lng}{{\langle}}
\newcommand{\rng}{{\rangle}}

\newcommand{\mcH}{\mathcal{H}}
\newcommand{\mcM}{\mathcal{M}}

\newcommand{\mbR}{\mathbb{R}}
\newcommand{\mbH}{\mathbb{H}}
\newcommand{\mbB}{\mathbb{B}}
\newcommand{\mbY}{\mathbb{Y}}
\newcommand{\mbX}{\mathbb{X}}

\newcommand{\ma}{\mathbf{a}}
\newcommand{\mb}{\mathbf{b}}
\newcommand{\mc}{\mathbf{c}}
\newcommand{\me}{\mathbf{e}}
\newcommand{\mf}{\mathbf{f}}
\newcommand{\mg}{\mathbf{g}}
\newcommand{\mh}{\mathbf{h}}
\newcommand{\mk}{\mathbf{k}}
\newcommand{\mm}{\mathbf{m}}
\newcommand{\mr}{\mathbf{r}}

\newcommand{\mv}{\mathbf{v}}
\newcommand{\mw}{\mathbf{w}}
\newcommand{\mx}{\mathbf{x}}

\newcommand{\bu}{\mathbf{u}}

\newcommand{\mA}{\mathbf{A}}
\newcommand{\mG}{\mathbf{G}}

\newcommand{\mQ}{\mathbf{Q}}

\newcommand{\mS}{\mathbf{S}}

\newcommand{\mU}{\mathbf{U}}

\newcommand{\mW}{\mathbf{W}}
\newcommand{\mX}{\mathbf{X}}
\newcommand{\mY}{\mathbf{Y}}
\newcommand{\mZ}{\mathbf{Z}}
\newcommand{\mDelta}{\mathbf{\Delta}}
\newcommand{\mTheta}{\mathbf{\Theta}}

\newcommand{\mdelta}{\bm{\delta}}
\newcommand{\mepsilon}{\bm{\epsilon}}
\newcommand{\meta}{\bm{\eta}}

\newcommand{\mxi}{\bm{\xi}}

\newcommand{\mvarphi}{\bm{\varphi}}
\newcommand{\mpsi}{\bm{\psi}}
\newcommand{\mphi}{\bm{\phi}}

\newcommand{\mLambda}{\bm{\Lambda}}

\DeclareMathOperator*{\esssup}{ess\,sup}
\DeclareMathOperator*{\argmin}{\arg\min}

\newcommand{\oplusin}{\underset{i=1}{\overset{n}{\bigoplus}}}
\newcommand{\oplusjd}{\underset{j=1}{\overset{d}{\bigoplus}}}
\newcommand{\opluskj}{\underset{k\neq j}{\bigoplus}}
\newcommand{\sumin}{\underset{i=1}{\overset{n}{\sum}}}
\newcommand{\sumjd}{\underset{j=1}{\overset{d}{\sum}}}

\newcommand{\sumkd}{\underset{k=1}{\overset{d}{\sum}}}
\newcommand{\limn}{\underset{n\rightarrow\infty}{\lim}}
\newcommand{\allj}{1\leq j\leq d}

\newcommand{\tran}{^{\mathstrut\scriptscriptstyle{\top}}}

\newcommand{\hp}{\hat{p}}

\newcommand{\Var}{{\rm Var}}
\newcommand{\E}{{\rm E}}
\newcommand{\Cov}{{\rm Cov}}
\newcommand{\I}{{\rm I}}
\newcommand{\Prob}{{\rm P}}
\newcommand{\Leb}{{\rm Leb}}
\newcommand{\CV}{{\rm CV}}
\newcommand{\Exp}{{\rm Exp}}
\newcommand{\inj}{{\rm inj}}
\newcommand{\const}{{\rm (const.)}}

\newcommand{\Log}{{\rm Log}}
\newcommand{\HS}{{\rm HS}}

\def\ba#1\ea{\begin{align*}#1\end{align*}}

\def\thetable{\arabic{table}}
\renewcommand{\theequation}{\thesection.\arabic{equation}}

\renewcommand{\abstractname}{\sc Abstract}

\def\qed{\space$\square$ \par \vspace{.15in}}

\begin{document}

\Large
\centerline{\bf Additive Regression with General Imperfect Variables}
\large
\vspace{.2cm} \centerline{Jeong Min Jeon and Germain Van Bever}
\vspace{.2cm} \centerline{University of Namur, Belgium}
\vspace{.2cm}
\small
\begin{quotation}
\noindent{\it Abstract:}
In this paper, we study an additive model where the response variable is Hilbert-space-valued and predictors are multivariate Euclidean, and both are possibly imperfectly observed. Considering Hilbert-space-valued responses allows to cover Euclidean, compositional, functional and density-valued variables. By treating imperfect responses, we can cover functional variables taking values in a Riemannian manifold and the case where only a random sample from a density-valued response is available. This treatment can also be applied in semiparametric regression. Dealing with imperfect predictors allows us to cover various principal component and singular component scores obtained from Hilbert-space-valued variables. For the estimation of the additive model having such variables, we use the smooth backfitting method originated by Mammen et al. (1999). We provide full non-asymptotic and asymptotic properties of our regression estimator and present its wide applications via several simulation studies and real data applications. \par
\vspace{.2cm}
\end{quotation}
\normalsize

\section{Introduction}

\setcounter{equation}{0}
\setcounter{subsection}{0}

Analysis of non-Euclidean data has been an important topic in statistics. Examples of non-Euclidean data include compositional, functional, density-valued or manifold-valued data.
Since non-Euclidian objects are often infinite-dimensional or subject to certain geometric constraints, the analysis of non-Euclidian data is often much more challenging than that of Euclidian data. During the recent decades, there have been many attempts to analyze non-Euclidean data. We refer to \cite{Abdulaziz (2021)} and \cite{Wang et al. (2016)} for recent reviews on compositional data analysis and functional data analysis, respectively. Regression analysis for other types of non-Euclidean data include \cite{Talska et al. (2018)} and \cite{Maier et al. (2021)} for density-valued responses, \cite{Cheng and Wu (2013)} and \cite{Jeon et al. (2022)} for manifold-valued predictors, \cite{Lin and Yao (2019)} for manifold-valued functional predictors and \cite{Petersen and Muller (2019)} for metric-space-valued responses, among many others. However, the aforementioned works and most of other works studied parametric regression or full-dimensional nonparametric regression. Parametric model assumptions are somewhat strong and full-dimensional nonparametric regression suffers from the curse of dimensionality, that is, the estimation performance gets severely worse as the number of predictors increases. To overcome these issues, we study additive regression.

Classical additive models with $Y\in\mathbb{R}$ and $X_j\in\mathbb{R}$, $j=1,\dots,d$, take the form of
\ba
Y=f_0+\sum_{j=1}^df_j(X_j)+\epsilon,
\ea
where $f_0\in\mathbb{R}$ is an unknown constant and $f_j$ are unknown functions, called component functions. Among the estimation techniques for additive models, the smooth backfitting (SBF) method, originated by \cite{Mammen et al. (1999)}, is known to avoid the curse of dimensionality under weak conditions. This powerful method has been further studied in various structured nonparametric models for Euclidean data (e.g., \citet{Mammen and Park (2005), Mammen and Park (2006), Linton et al. (2008), Yu et al. (2008), Lee et al. (2010), Lee et al. (2012), Han and Park (2018)}). The SBF method is also applied to functional additive regression models (\citet{Han et al. (2018), Park et al. (2018)}). \cite{Han et al. (2018)} studied the case where there are multiple functional predictors and \cite{Park et al. (2018)} studied function-on-function regression for a single functional predictor. Recently, \cite{Han et al. (2020)} and \cite{Lin et al. (2022)} applied the SBF method to treat density-valued responses and manifold-valued responses, respectively. Also, \cite{Jeon and Park (2020)} extended the SBF technique to general Hilbert-space-valued (Hilbertian) responses, and \cite{Jeon et al. (2021b)} further extended \cite{Jeon and Park (2020)} to cover discrete-type predictors and censored or missing responses. However, all the aforementioned works on the SBF method only studied `univariate structured nonparametric models' where the associated component functions are defined on $\mathbb{R}$.

Recently, \citet{Jeon et al. (2021a)} extended \cite{Jeon and Park (2020)} to a `multivariate additive model' for which the associated component functions are defined on a finite-dimensional Hilbert space or a Riemannian manifold. However, \cite{Jeon et al. (2021a)} does not treat imperfectly observed predictors or responses. This limits the scope of their models as they do not cover, for examples, variables obtained from dimension reduction techniques such as principal component scores and singular component scores. Indeed, such variables are not directly observable in general. This prevents \cite{Jeon et al. (2021a)} from dealing with high/infinite-dimensional predictors. It also does not cover the case where the response variable is a random density and only a random sample from the random density is observable. Also, \citet{Jeon et al. (2021a)} assumes that the predictors in the model have compact supports. However, the support of a random element is not necessarily compact and it can be unbounded.

To introduce our model, we let $\mbH$ be a separable Hilbert space. We note that Euclidean spaces are examples of $\mbH$. Examples of non-Euclidean $\mbH$ can be found in Section \ref{examples of hilbert}. We also let $\op$, $\od$, $\bzero$, $\langle\cdot,\cdot\rangle$ and $\|\cdot\|$ denote a vector addition, a scalar multiplication, a zero vector, an inner product and a norm on $\mbH$, respectively. We note that $\op$, $\od$, $\bzero$, $\langle\cdot,\cdot\rangle$ and $\|\cdot\|$ for $\mbH=\mathbb{R}^k$ correspond to +, $\times$, $(0,\ldots,0)\in\mathbb{R}^k$, the dot product and the $\ell_2$-norm, respectively. Examples of those for other $\mbH$ can be found in the Supplementary Material \ref{examples}. In this paper, we consider the multivariate Hilbertian additive model
\begin{align}\label{model}
\mY=\mf_0\op\bigoplus_{j=1}^d\mf_j(\xi_j)\op\mepsilon,
\end{align}
where $\mY\in\mbH$ is a response variable satisfying $\E(\|\mY\|^2)<\infty$, $\xi_j=(\xi_{j1},\ldots,\xi_{jL_j})\in\mathbb{R}^{L_j}$ with $L_j\in\mathbb{N}$ are multivariate predictors, $\mepsilon\in\mbH$ is an error term satisfying $\E(\|\mepsilon\|^2)<\infty$ and $\E(\mepsilon|\xi_1,\cdots,\xi_d)=\bzero$, $\mf_0\in\mbH$ is an unknown constant and $\mf_j:\mathbb{R}^{L_j}\rightarrow\mbH$ are unknown component maps. Here, the conditional expectation $\E(\mepsilon|\xi_1,\cdots,\xi_d)$ is defined through Bochner integration, which is a generalization of Lebesgue integration to Banach-space-valued maps. The definition of Bochner integral can be found in the Supplementary Material \ref{bochner definition}.
Note that the multivariate additive model at (\ref{model}) fills the gap between the univariate additive model and the full-dimensional nonparametric model.

In model (\ref{model}), we do not assume that each $\xi_j$ has a compact support. Instead, we estimate each $\mf_j$ on an arbitrary compact set $D_j\subset\mathbb{R}^{L_j}$ of interest. Note that efficient nonparametric estimation of $\mf_j$ on the whole Euclidean space $\mathbb{R}^{L_j}$ is not feasible since the collection of observed data is bounded. This setting is new in a multivariate additive model, and while it gives additional complexity in the asymptotic theory compared to the case of compact support, it allows for much more flexible models.

We allow for the case where $\xi_{jl}$, for different $1\leq j\leq d$ and $1\leq l\leq L_j$, come from different sources. For example, $\xi_{jl}$ can be usual scalar predictors and the principal component or singular component scores of some high/infinite-dimensional variables. Since the usual scalar predictors are sometimes contaminated by vanishing measurement errors (e.g., \citet{Fan (1992), Delaigle (2008), van Es and Gugushvili (2010)}) and the component scores are unobservable in general, we consider the case where we cannot observe some $\xi_{jl}$ but we can only obtain a proxy $\tilde{\xi}_{jl}\in\mathbb{R}$ for such $\xi_{jl}$. We also allow for the case where we only have a proxy $\tilde{\mY}$ instead of observing $\mY$. The latter treatment covers vanishing measurement errors on the responses, imperfect density-valued responses and Riemannian-manifold-valued (Riemannian) functional responses. It can be also applied to the estimation of nonparametric part in profiling-based semiparametric regression. This general setting, dealing with various perfect/imperfect predictors and responses in one model, has not been studied in the literature even in parametric models. It has, however, wide applications including many new regression problems. We believe that this unified framework is a useful source which justifies how past and future Euclidean/non-Euclidean data analysis without regression setting (e.g., \citet{Yang et al. (2011), Petersen and Muller (2016), Dai and Muller (2018)}) can be applied to regression analysis.

This paper is organized as follows: we present examples of separable Hilbert space and imperfect variables in Section \ref{examples of hilbert and imperfect}. In the same section, we also provide instances of proxies $\tilde{\xi}_{jl}$ and $\tilde{\mY}$ and investigate their asymptotic properties. Our estimation method for model (\ref{model}) and its asymptotic properties are respectively given in Sections \ref{estimation method} and \ref{asymptotic properties sbf}. Section \ref{extension} is devoted to the case where $\mY$ is obtained from a Riemannian functional variable. We present two simulation studies and two real data applications in Section \ref{numerical study} and summarize our main contributions in Section \ref{conclusion}. All technical proofs and mathematical preliminaries are collected in the Supplementary Material.

\section{Examples of Hilbert spaces and imperfect variables}\label{examples of hilbert and imperfect}

\setcounter{equation}{0}
\setcounter{subsection}{0}

\subsection{Examples of Hilbert space}\label{examples of hilbert}

In this section, we introduce simplices, Bayes-Hilbert spaces on topological spaces, $L^2$-spaces of Hilbertian functions on topological spaces and product Hilbert spaces as non-Euclidean examples of separable Hilbert space.

\begin{example}\label{pc simplex}
(Simplices) Define
$\mathcal{S}_1^k=\{(v_1,\ldots,v_k)\in(0,1)^k:\sum_{j=1}^k v_j=1\}$ for $k\geq2$. It is well known that $\mathcal{S}_1^k$ equipped with the vector operations and inner product defined in the Supplementary Material \ref{examples} forms a $(k-1)$-dimensional Hilbert space. \qed
\end{example}

\begin{example}\label{Density valued}
(Bayes-Hilbert spaces on topological spaces) Define
\begin{align}\label{B-H space}
\mathfrak{B}^2(\mathcal{S},\mathcal{B}(\mathcal{S}),\nu)=\bigg\{g:\mathcal{S}\rightarrow(0,\infty)\bigg|\int_Sg(s)d\nu(s)=1\text{~and~}\int_{\mathcal{S}}(\log g(s))^2d\nu(s)<\infty\bigg\},
\end{align}
where $\mathcal{S}$ is a second countable topological space, $\mathcal{B}(\mathcal{S})$ is the Borel sigma-field of $\mathcal{S}$ and $\nu$ is a finite Borel measure on $\mathcal{S}$. From Proposition 2.2 in \cite{Maier et al. (2021)}, Proposition 3.4.5 in \cite{Cohn (2013)} and Proposition 3.1 in \cite{Preston (2008)}, it follows that $\mathfrak{B}^2(\mathcal{S},\mathcal{B}(\mathcal{S}),\nu)$ equipped with the vector operations and inner product defined in the Supplementary Material \ref{examples} forms a separable Hilbert space. The case $\mathcal{S}=[a,b]$ with $-\infty<a<b<\infty$ has been usually considered in density-valued data analysis (e.g., \citet{Hron et al. (2016), Talska et al. (2018)}). Here, we also consider a multi-dimensional Euclidean space or a Riemannian manifold for $\mathcal{S}$. Hence,  $\mathfrak{B}^2(\mathcal{S},\mathcal{B}(\mathcal{S}),\nu)$ covers a wide range of density-valued data. \qed
\end{example}

\begin{example}\label{Fully observed fd}
($L^2$-spaces of Hilbertian functions on topological spaces) Define
\ba
L^2((\mathcal{S},\mathcal{B}(\mathcal{S}),\nu),\mathcal{H})=\bigg\{\mg:\mathcal{S}\rightarrow\mathcal{H}\bigg|\int_{\mathcal{S}}\|\mg(s)\|^2_{\mathcal{H}}\,d\nu(s)<\infty\bigg\},
\ea
where the triple $(\mathcal{S},\mathcal{B}(\mathcal{S}),\nu)$ is defined as in Example \ref{Density valued}, $\mathcal{H}$ is a separable Hilbert space and $\|\cdot\|_{\mathcal{H}}$ is a norm on $\mathcal{H}$. Proposition 3.3 in \cite{Bunkure (2019)}, Proposition 3.4.5 in \cite{Cohn (2013)} and Proposition 3.1 in \cite{Preston (2008)} imply that $L^2((\mathcal{S},\mathcal{B}(\mathcal{S}),\nu),\mathcal{H})$, equipped with the vector operations and inner product defined in the Supplementary Material \ref{examples}, forms a separable Hilbert space by identifying Hilbertian functions that are equal almost everywhere on $\mathcal{S}$ with respect to $\nu$. We call this space a Lebesgue-Bochner-Hilbert space. The case $(\mathcal{S},\mathcal{H})=([a,b],\mathbb{R})$ with $-\infty<a<b<\infty$ has been usually considered in functional data analysis. Here, the general $(\mathcal{S},\mathcal{H})$ covers various domains and co-domains of random functions. \qed
\end{example}

\begin{example}
(Product of Hilbert spaces) Suppose that $\mathbb{H}_1,\cdots,\mathbb{H}_k$ with $k\geq2$ are separable Hilbert spaces. Then, the product space $\mbH=\prod_{m=1}^k\mathbb{H}_m$ equipped with the elementwise vector operations and the inner product $\lng\cdot,\cdot\rng$ defined by $\lng\mh_1,\mh_2\rng=\sum_{m=1}^k\lng \mh_{1m},\mh_{2m}\rng_m$ forms a separable Hilbert space, where $\lng\cdot,\cdot\rng_m$ is an inner product on $\mathbb{H}_m$. \qed
\end{example}

\begin{remark}
Examples \ref{Density valued} and \ref{Fully observed fd} are new and they generalize the examples in \cite{Jeon and Park (2020)}, \cite{Jeon et al. (2021a)} and \cite{Jeon et al. (2021b)}, since the latter works considered Bayes-Hilbert spaces on Euclidean $\mathcal{S}$ and $L^2$-spaces of real-valued functions on Euclidean $\mathcal{S}$, as examples of separable Hilbert space. This indicates that we can actually cover wider types of variables than those considered in the literature by treating Hilbertian variables. \qed
\end{remark}

In Section \ref{examples of imperfect}, we also introduce tensor Hilbert spaces, which are also separable Hilbert spaces originated by \cite{Lin and Yao (2019)}. They are used to deal with Riemannian functional variables; see Example \ref{Riemannian fd} for details.

\subsection{Examples of imperfect variables}\label{examples of imperfect}

In this section, we introduce several examples of imperfect variables with the ways to approximate them and their asymptotic properties. Let $(\mxi^i,\mY^i)$, $1\leq i\leq n$, be i.i.d. copies of $(\mxi,\mY)$, where $\mxi=(\xi_{jl}:1\leq j\leq d,1\leq l\leq L_j)$. Also, let $\tilde{\xi}^i_{jl}$ be random variables observed instead of $\xi^i_{jl}$, and $\tilde{\mY}^i$ be a $\mbH$-valued random element observed instead of $\mY^i$. We put $\tilde{\xi}^i_{jl}=\xi^i_{jl}$ (resp. $\tilde{\mY}^i=\mY^i$) if $\xi^i_{jl}$ (resp. $\mY^i$) is perfectly observed. For $\tilde{\mxi}^i=(\tilde{\xi}^i_{jl}:1\leq j\leq d,1\leq l\leq L_j)$, we do not assume that $(\tilde{\mxi}^i,\tilde{\mY}^i)$ are identically distributed nor independent across $1\leq i\leq n$ since some imperfect cases do not satisfy this. Instead, we consider the case where $\max_{1\leq i\leq n}|\tilde{\xi}^i_{jl}-\xi^i_{jl}|=O_p(a_{njl})$ for some sequence $a_{njl}=o(1)$ for all $1\leq j\leq d$ and $1\leq l\leq L_j$. In case $\xi_{jl}$ is a perfectly observable predictor, $\tilde{\xi}^i_{jl}=\xi^i_{jl}$ for all $i$ and hence $a_{njl}=0$. We also consider the case where $\max_{1\leq i\leq n}\|\tilde{\mY}^i\om\mY^i\|=O_p(b_n)$ for some sequence $b_n=o(1)$. In case $\mY$ is a perfectly observable response, $\tilde{\mY}^i=\mY^i$ for all $i$ and hence $b_n=0$. 
The vanishing speeds of $a_{njl}$ and $b_n$ affect the asymptotic properties of our estimator for model (\ref{model}). We first provide the examples of $\tilde{\xi}_{jl}$ and $a_{njl}$, and then provide the examples of $\tilde{\mY}$ and $b_n$. They contain some new propositions, which are of importance in their own right.

\begin{customexample}{A.1}\label{Measurement error case}
(Scalar predictors with vanishing measurement errors) We consider the case where $\tilde{\xi}^i_{jl}=\xi^i_{jl}+\sigma_n\cdot U^i$ for some $(j,l)$. Here, $\sigma_n$ is a sequence such that $\sigma_n=o(1)$ and $U^i$ are unobservable i.i.d. copies of a measurement error $U$. Such vanishing measurement error case is considered in density estimation problems such as \cite{Fan (1992)}, \cite{Delaigle (2008)} and \cite{van Es and Gugushvili (2010)}, with justification for the vanishing measurement error. In this case, the knowledge of the distribution of $U$ is not necessary unlike the classical measurement error problems where measurement errors do not vanish and one estimates its distribution or assumes that it is known (e.g., \citet{Fan (1991), Delaigle et al. (2008), Johannes (2009), Han and Park (2018), Bertrand et al. (2019)}).

\cite{Fan (1992)}, \cite{Delaigle (2008)} and \cite{van Es and Gugushvili (2010)} considered various vanishing rates for $\sigma_n$ including polynomial rates. In fact, the vanishing speed of $\sigma_n$ affects the vanishing speed of $a_{njl}$. If $\E(|U|^{\tau})<\infty$ for some $\tau>0$, then Lemma \ref{markov} in the Supplementary Material \ref{collection of lemmas} gives that $\max_{1\leq i\leq n}|U^i|=O_p(n^{1/\tau})$. In this case, $a_{njl}=\sigma_n\cdot n^{1/\tau}$. In case $U$ has an exponential moment, say $\E(\exp(c\cdot|U|))<\infty$ for some constant $c>0$, we get $\max_{1\leq i\leq n}|U^i|=O_p(\log{n})$ (Lemma \ref{markov}) and thus $a_{njl}=\sigma_n\cdot\log n$. In case $U$ is a sub-Gaussian random variable, it is well known that $\max_{1\leq i\leq n}|U^i|=O_p(\sqrt{\log n})$ and thus $a_{njl}=\sigma_n\sqrt{\log n}$. In case $U$ is a bounded random variable, we have $\max_{1\leq i\leq n}|U^i|=O_p(1)$ and thus $a_{njl}=\sigma_n$. Hence, a faster vanishing rate for $\sigma_n$ or a higher moment for $U$ gives a faster vanishing rate for $a_{njl}$. \qed
\end{customexample}

\begin{customexample}{A.2}\label{Fully observed hilbert}
(Principal component scores of Hilbertian variables) We consider the case where some $\xi_{jl}$ are the principal component (PC) scores of a Hilbertian random element. Let $\mathbb{H}_*$ be a separable Hilbert space, possibly different from $\mbH$, and equipped with a vector addition $\op_*$, a scalar multiplication $\od_*$, a zero vector $\bzero_*$ and an inner product $\langle\cdot,\cdot\rangle_*$. Let $\|\cdot\|_*$ denote the norm induced by $\langle\cdot,\cdot\rangle_*$, and $\mX$ be a random element taking values in $\mathbb{H}_*$ and satisfying $\E(\|\mX\|_*^2)<\infty$. Let $\mathbb{X}=\mX\om_*\E(\mX)$, where $\om_*$ is the vector subtraction on $\mathbb{H}_*$ defined by $\mh_1\om_*\mh_2=\mh_1\op_*(-1\od_*\mh_2)$ and the expectation is in the sense of Bochner integral. Define the covariance operator $C_{\mX}:\mathbb{H}_*\rightarrow\mathbb{H}_*$ of $\mX$ by $C_{\mX}(\mh)=\E(\lng \mathbb{X},\mh\rng_*\od_*\mathbb{X})$. Let $d_\mX=\dim(\overline{{\rm{Im}}(C_{\mX})})\in\mathbb{N}\cup\{+\infty\}$ be the dimension of the closure of the image of $C_{\mX}$. Then, $C_{\mX}$ admits the decomposition
\begin{align*}
C_{\mX}(\mh)=\sideset{}{_*}\bigoplus_{r=1}^{d_\mX}(\lng \lambda_r\od_*\mpsi_r,\mh\rng_*\od_*\mpsi_r)
\end{align*}
for all $\mh\in\mbH_*$, where $\{\lambda_r\}_{r=1}^{d_\mX}$ is the set of positive eigenvalues of $C_{\mX}$ with decreasing order $\lambda_1\geq\lambda_2\geq\cdots$ and $\{\mpsi_r\}_{r=1}^{d_\mX}$ is an orthonormal basis of $\overline{{\rm{Im}}(C_{\mX})}$ consisting of the corresponding eigenvectors of $C_{\mX}$ (Theorem 7.2.6 in \cite{Hsing and Eubank (2015)}). Also, it holds that
\begin{align}\label{optimal basis expansion}
\mathbb{X}=\sideset{}{_*}\bigoplus_{r=1}^{d_\mX}(\lng \mathbb{X},\mpsi_r\rng_*\od_*\mpsi_r)
\end{align}
almost surely, where $\eta_r:=\lng \mathbb{X},\mpsi_r\rng_*$ are uncorrelated random variables with $\E(\eta_r)=0$ and $\Var(\eta_r)=\lambda_r$ (Theorem 7.2.7 in \cite{Hsing and Eubank (2015)}). The random variables $\eta_r$ are called the PC scores of $\mX$.
We may take a finite number of $\eta_r$ as a subset of $\{\xi_{jl}:1\leq j\leq d,1\leq l\leq L_j\}$ when $\mX$ is a predictor. However, $\eta_r$ are not observable since $\E(\mX)$ and $\mpsi_r$ are unknown. Thus, we need to estimate them and we take the estimated $\eta_r$ as a subset of $\{\tilde{\xi}_{jl}:1\leq j\leq d,1\leq l\leq L_j\}$. This Hilbertian PCA generalizes Euclidean PCA, compositional PCA (e.g., \cite{Wang et al. (2015)}), functional PCA (FPCA) and multivariate FPCA (e.g., \cite{Ramsay and Silverman (2005)}) considering the examples of separable Hilbert space in Section \ref{examples of hilbert}.

Let $\{\mX^i:1\leq i\leq n\}$ be i.i.d. observations of $\mX$. Estimation of $\E(\mX)$ and the covariance operator $C_{\mX}$ is achieved using the sample mean $\bar{\mX}$ and the unbiased sample covariance operator $\hat{C}_{\mX}:\mbH_*\rightarrow\mbH_*$, respectively defined by
\begin{align}\label{sample mean X}
\bar{\mX}=\frac{1}{n}\od_*\sideset{}{_*}\bigoplus_{i=1}^n\mX^i\quad\textrm{ and }\quad\hat{C}_{\mX}(\mh)=\frac{1}{n-1}\od_*\sideset{}{_*}\bigoplus_{i=1}^n(\lng\mX^i\om_*\bar{\mX},\mh\rng_*\od_*(\mX^i\om_*\bar{\mX})).
\end{align}
As $\hat{C}_{\mX}$ is a nonnegative-definite compact self-adjoint operator, the spectral theorem for compact self-adjoint operators (e.g., Theorem 4.2.4 in \cite{Hsing and Eubank (2015)}) implies that
\begin{align*}
\hat{C}_{\mX}(\mh)=\sideset{}{_*}\bigoplus_{r=1}^{\tilde{d}_\mX}(\lng \hat{\lambda}_r\od_*\hat{\mpsi}_r,\mh\rng_*\od_*\hat{\mpsi}_r)
\end{align*}
for any $\mh\in\mbH_*$, where $\tilde{d}_\mX=\dim\big(\overline{{\rm{Im}}(\hat{C}_{\mX})}\big)$, $\{\hat{\lambda}_r\}_{r=1}^{\tilde{d}_\mX}$ is the set of positive eigenvalues of $\hat{C}_{\mX}$ with decreasing order and $\hat{\mpsi}_r\in\mathbb{H}_*$ are the associated orthonormal eigenvectors. We take $\hat{\mpsi}_r$ as estimators of $\mpsi_r$. Finally, we estimate $\eta^i_r=\lng\mX^i\om_*\E(\mX),\mpsi_r\rng_*$ by $\hat{\eta}^i_r=\lng\mX^i\om_*\bar{\mX},\hat{\mpsi}_r\rng_*$. The following proposition gives the vanishing speed of $a_{njl}$ for the case where $\xi^i_{jl}=\eta^i_r$ and $\tilde{\xi}^i_{jl}=\hat{\eta}^i_r$ for some $r$. Its proof is given in the Supplementary Material \ref{HPC proof}. Note that, if $\mh$ is an eigenvector of $C_{\mX}$, then $-1\od_*\mh$ is also an eigenvector of $C_{\mX}$. As is standard in functional data analysis, $\eta^i_r$ is defined using the unique eigenvector $\mpsi_r$ satisfying $\lng\hat{\mpsi}_r,\mpsi_r\rng_*>0$. 

\begin{proposition}\label{HPC}
Assume that the first $r$ eigenvalues of $C_{\mX}$ have multiplicity one, that is, $\lambda_1>\lambda_2>\cdots>\lambda_{r+1}$. Then, for $\mpsi_r$ with $\lng\hat{\mpsi}_r,\mpsi_r\rng_*>0$, $\max_{1\leq i\leq n}|\lng\mX^i\om_*\bar{\mX},\hat{\mpsi}_r\rng_*-\lng\mX^i\om_*\E(\mX),\mpsi_r\rng_*|$ achieves the following rates:
\ba
\begin{cases}
O_p(n^{-1/2+1/\tau}), & \text{if}\ \E(\|\mX\|_*^{\tau})<\infty\ \text{for some}\ \tau\geq4, \\
O_p(n^{-1/2}\cdot\log{n}), & \text{if}\ \E\big(\exp(c\cdot\|\mX\|_*)\big)<\infty\ \text{for some}\ c>0, \\
O_p(n^{-1/2}), & \text{if}\ \|\mX\|_*<C\ \text{almost surely for some}\ C>0.
\end{cases}
\ea
\end{proposition}

Note that a version of the eigenvalue condition in Proposition \ref{HPC} is typically assumed in the literature of functional data analysis; see e.g., \cite{Cai and Hall (2006)} and \cite{Hall and Horowitz (2007)}. \qed

\end{customexample}

\begin{customexample}{A.3}\label{singular hilbert}
(Singular component scores of Hilbertian variables) Recall the notations and setting in Example \ref{Fully observed hilbert}. We introduce the notion of singular component (SC) score and consider the case where some $\xi_{jl}$ are SC scores obtained from $\mX$ and $\mY$. Let $\mbY=\mY\om\E(\mY)$, where $\om$ is the vector subtraction on $\mbH$ defined by $\mh_1\om\mh_2=\mh_1\op(-1\od\mh_2)$. Define the cross-covariance operators $C_{\mX\mY}:\mathbb{H}\rightarrow\mathbb{H}_*$ and $C_{\mY\mX}:\mathbb{H}_*\rightarrow\mathbb{H}$ of $\mX$ and $\mY$ by $C_{\mX\mY}(\mh)=\E(\lng \mbY,\mh\rng\od_*\mbX)$ and $C_{\mY\mX}(\mh_*)=\E(\lng \mbX,\mh_*\rng_*\od\mbY)$, respectively. Also, define $C_{\mX\mY\mX}=C_{\mX\mY}\circ C_{\mY\mX}:\mathbb{H}_*\rightarrow\mathbb{H}_*$ and $C_{\mY\mX\mY}=C_{\mY\mX}\circ C_{\mX\mY}:\mathbb{H}\rightarrow\mathbb{H}$. Note that $C_{\mX\mY\mX}$ and $C_{\mY\mX\mY}$ are compact operators since they are compositions of two compact operators. Since the adjoint of $C_{\mX\mY}$ is $C_{\mY\mX}$ by Theorem 7.2.9 in \cite{Hsing and Eubank (2015)}, $C_{\mX\mY\mX}$ and $C_{\mY\mX\mY}$ are nonnegative-definite self-adjoint operators. 
Let $d_{\mX\mY}=\dim(\overline{{\rm{Im}}(C_{\mX\mY})})$, which can be finite or infinite. For any bounded linear operator $A:\mbH\rightarrow\mbH_*$, let $A^*:\mbH_*\rightarrow\mbH$ denote its adjoint. Since $\overline{{\rm{Im}}(A\circ A^*)}=\overline{{\rm{Im}}(A)}$ (e.g., Theorem 3.3.7 in \cite{Hsing and Eubank (2015)}), it holds that $\dim(\overline{{\rm{Im}}(C_{\mX\mY\mX})})=d_{\mX\mY}$. Then, by the spectral theorem for compact self-adjoint operators, $C_{\mX\mY\mX}$ admits the decomposition
\begin{align*}
C_{\mX\mY\mX}(\mh_*)=\sideset{}{_*}\bigoplus_{r=1}^{d_{\mX\mY}}(\lng \sigma_r^2\od_*\mphi_r,\mh_*\rng_*\od_*\mphi_r)
\end{align*}
for all $\mh_*\in\mbH_*$, where $\{\sigma_r\}_{r=1}^{d_{\mX\mY}}$ is the set of positive eigenvalues of $C_{\mX\mY\mX}$ with decreasing order $\sigma_1\geq\sigma_2\geq\cdots$ and $\{\mphi_r\}_{r=1}^{d_{\mX\mY}}\subset\mbH_*$ is an orthonormal basis of $\overline{{\rm{Im}}(C_{\mX\mY\mX})}$ consisting of the corresponding eigenvectors of $C_{\mX\mY\mX}$. Since $\dim({\rm{Im}}(A))=\dim({\rm{Im}}(A^*))$ for any bounded linear operator $A:\mbH\rightarrow\mbH_*$ (e.g., Theorem 3.3.7 in \cite{Hsing and Eubank (2015)}) and ${\rm{Im}}(A)=\overline{{\rm{Im}}(A)}$ for any $A$ with $\dim({\rm{Im}}(A))<\infty$, we may prove that $\dim(\overline{{\rm{Im}}(C_{\mY\mX})})=d_{\mX\mY}$, which implies that $\dim(\overline{{\rm{Im}}(C_{\mY\mX\mY})})=d_{\mX\mY}$. Also, since $A\circ A^*$ and $A^*\circ A$ have the same eigenvalues for any bounded linear operator $A:\mbH\rightarrow\mbH_*$ (e.g., Theorem 4.3.1 in \cite{Hsing and Eubank (2015)}), $C_{\mX\mY\mX}$ and $C_{\mY\mX\mY}$ have the common eigenvalues $\sigma_r^2$. Hence, the spectral theorem implies that
\ba
C_{\mY\mX\mY}(\mh)=\bigoplus_{r=1}^{d_{\mX\mY}}(\lng \sigma_r^2\od\mvarphi_r,\mh\rng\od\mvarphi_r)
\ea
for all $\mh\in\mbH$, where $\{\mvarphi_r\}_{r=1}^{d_{\mX\mY}}\subset\mathbb{H}$ is an orthonormal basis of $\overline{{\rm{Im}}(C_{\mY\mX\mY})}$ consisting of the eigenvectors of $C_{\mY\mX\mY}$. We call the values $\sigma_r$ singular values and the vectors $\mphi_r$ and $\mvarphi_r$ singular vectors. Also, the orthogonal decomposition theorem of Hilbert spaces (e.g., Theorem 3.3.7 in \cite{Hsing and Eubank (2015)}) implies that
\ba
\mbX=\sideset{}{_*}\bigoplus_{r=1}^{d_{\mX\mY}}(\lng \mathbb{X},\mphi_r\rng_*\od_*\mphi_r)\op_*\mathbf{R}_{\mX},\quad\mbY=\bigoplus_{r=1}^{d_{\mX\mY}}(\lng \mathbb{Y},\mvarphi_r\rng\od\mvarphi_r)\op\mathbf{R}_{\mY},
\ea
where $\mathbf{R}_{\mX}\in\mbH_*$ and $\mathbf{R}_{\mY}\in\mbH$ are stochastic remainders such that $C_{\mX\mY\mX}(\mathbf{R}_{\mX})=\bzero_*$ and $C_{\mY\mX\mY}(\mathbf{R}_{\mY})=\bzero$. Here, we call the mean 0 random variables $\theta_r:=\lng \mathbb{X},\mphi_r\rng_*$ and $\vartheta_r:=\lng \mathbb{Y},\mvarphi_r\rng$ the SC scores of $\mX$ and $\mY$, respectively. Note that $\theta_r$ are constructed from both $\mX$ and $\mY$ unlike the PC scores $\eta_r$ constructed from $\mX$ only. Hence, $\theta_r$ may contain more significant information than $\eta_r$ for assessing the relationship between $\mX$ and $\mY$. Thus, we may take a finite number of $\theta_r$ as a subset of $\{\xi_{jl}:1\leq j\leq d,1\leq l\leq L_j\}$ when $\mX$ is a predictor. However, $\theta_r$ are not observable since $\E(\mX)$, $\E(\mY)$ and $\mphi_r$ are unknown. Therefore, we need to estimate them and we take the estimated $\theta_r$ as a subset of $\{\tilde{\xi}_{jl}:1\leq j\leq d,1\leq l\leq L_j\}$. This Hilbertian singular component analysis is a generalization of partial least squares singular value decomposition (e.g., \cite{Lafaye de Micheaux et al. (2019)}), which is an alternative dimension reduction technique to Euclidean PCA when both $\mbH$ and $\mbH_*$ are Euclidean spaces. It also generalizes functional singular component analysis developed by \cite{Yang et al. (2011)}, which is an alternative dimension reduction technique to FPCA when both $\mbH$ and $\mbH_*$ are $L^2$-spaces of real-valued functions on $\mathbb{R}$. Our general Hilbertian singular component analysis can be also applied to various $\mX$ and $\mY$ of possibly different types.

We now study the problem of estimating SC scores. For more general results, we remain the case where we only have $\{\tilde{\mY}^i:1\leq i\leq n\}$ instead of $\{\mY^i:1\leq i\leq n\}$. For the estimation of $\E(\mX)$ and $\E(\mY)$, we may use $\bar{\mX}$ defined at (\ref{sample mean X}) and $\bar{\tilde{\mY}}=n^{-1}\od\bigoplus_{i=1}^n\tilde{\mY}^i$, respectively. For the estimation of $\mphi_r$, we estimate $C_{\mX\mY\mX}$ by $\hat{C}_{\mX\mY\mX}=\hat{C}_{\mX\mY}\circ\hat{C}_{\mY\mX}$, where $\hat{C}_{\mX\mY}:\mbH\rightarrow\mbH_*$ and $\hat{C}_{\mY\mX}:\mbH_*\rightarrow\mbH$ are defined by
\ba
\hat{C}_{\mX\mY}(\mh)&=\frac{1}{n-1}\od_*\sideset{}{_*}\bigoplus_{i=1}^n(\lng\tilde{\mY}^i\om\bar{\tilde{\mY}},\mh\rng\od_*(\mX^i\om_*\bar{\mX})),\\
\hat{C}_{\mY\mX}(\mh_*)&=\frac{1}{n-1}\od\bigoplus_{i=1}^n(\lng\mX^i\om_*\bar{\mX},\mh_*\rng_*\od(\tilde{\mY}^i\om\bar{\tilde{\mY}})).
\ea
We then take the eigenvectors $\hat{\mphi}_r$ of $\hat{C}_{\mX\mY\mX}$ as estimators of $\mphi_r$. Finally, we estimate $\theta^i_r=\lng\mX^i\om_*\E(\mX),\mphi_r\rng_*$ by $\hat{\theta}^i_r=\lng\mX^i\om_*\bar{\mX},\hat{\mphi}_r\rng_*$. The following proposition, whose proof is given in the Supplementary Material \ref{HSC proof}, gives the vanishing speed of $a_{njl}$ for the case where $\xi^i_{jl}=\theta^i_r$ and $\tilde{\xi}^i_{jl}=\hat{\theta}^i_r$ for some $r$. Let $\|\cdot\|_{\HS}$ denote the Hilbert-Schmidt norm.

\begin{proposition}\label{HSC}
Assume that the first $r$ eigenvalues of $C_{\mX\mY\mX}$ have multiplicity one and that $\E(\|\mX\|_*^2\|\mY\|^2)<\infty$. Then, it holds that $\|\hat{C}_{\mX\mY\mX}-C_{\mX\mY\mX}\|_{\rm HS}=O_p(n^{-1/2}+b_n)$. Also, for $\mphi_r$ with $\lng\hat{\mphi}_r,\mphi_r\rng_*>0$, it holds that $\|\hat{\mphi}_r\om_*\mphi_r\|_*=O_p(n^{-1/2}+b_n)$. Moreover, $\max_{1\leq i\leq n}|\lng\mX^i\om_*\bar{\mX},\hat{\mphi}_r\rng_*-\lng\mX^i\om_*\E(\mX),\mphi_r\rng_*|$ achieves the following rates:
\ba
\begin{cases}
O_p((n^{-1/2}+b_n)\cdot n^{1/\tau}), & \text{if}\ \E(\|\mX\|_*^{\tau})<\infty\ \text{for some}\ \tau\geq2, \\
O_p((n^{-1/2}+b_n)\cdot \log{n}), & \text{if}\ \E\big(\exp(c\cdot\|\mX\|_*)\big)<\infty\ \text{for some}\ c>0, \\
O_p(n^{-1/2}+b_n), & \text{if}\ \|\mX\|_*<C\ \text{almost surely for some}\ C>0.
\end{cases}
\ea
\end{proposition}

We note that the above eigen-analysis for general Hilbertian variables does not exist in the literature even for the case $b_n=0$, to the best of our knowledge. \qed
\end{customexample}


\begin{customexample}{A.4}\label{Riemannian fd}
(Intrinsic Riemannian functional principal component scores) We consider the case where some $\xi_{jl}$ are functional PC (FPC) scores obtained from a Riemannian functional variable. This case is somewhat similar to Example \ref{Fully observed hilbert} but the procedure for obtaining $\tilde{\xi}_{jl}$ and the related asymptotic properties of $\tilde{\xi}_{jl}$ are distinct.

We first introduce tensor Hilbert spaces, as in \cite{Lin and Yao (2019)}. We defer the definitions of associated mathematical notions to the Supplementary Material \ref{Riemannian definition}. Let $\mu:\mathcal{T}\rightarrow\mathcal{M}$ be a measurable function, where $\mathcal{T}\subset\mathbb{R}$ is a compact set equipped with a finite Borel measure $\nu$ and $\mathcal{M}$ is a complete and connected Rimannian manifold. Let $T_{\mu(t)}\mathcal{M}$ denote the tangent space of $\mathcal{M}$ at $\mu(t)$ for $t\in \mathcal{T}$ and $\lng\cdot,\cdot\rng_{\mu(t)}$ denote the inner product on $T_{\mu(t)}\mathcal{M}$. Define
$\mathfrak{T}(\mu)$ be the collection measurable functions $V:\mathcal{T}\rightarrow \bigcup_{t\in\mathcal{T}}T_{\mu(t)}\mathcal{M}$ satisfying $V(t)\in T_{\mu(t)}\mathcal{M}$ for all $t\in \mathcal{T}$ and $\int_{\mathcal{T}}\lng V(t),V(t)\rng_{\mu(t)}d\nu(t)<\infty$. It has the vector space structure with vector addition $V_1+V_2$ defined by $(V_1+V_2)(t)=V_1(t)+V_2(t)$ and scalar multiplication $c\cdot V$ defined by $(c\cdot V)(t)=c\cdot V(t)$, where the latter $+$ and $\cdot$ are the vector addition and scalar multiplication in $T_{\mu(t)}\mathcal{M}$. Endow an inner product $\lng \cdot,\cdot\rng_{\mathfrak{T}(\mu)}$ on $\mathfrak{T}(\mu)$ by $\lng V_1,V_2\rng_{\mathfrak{T}(\mu)}=\int_{\mathcal{T}}\lng V_1(t),V_2(t)\rng_{\mu(t)}d\nu(t)$. Then, the inner product space $\mathfrak{T}(\mu)$ forms a separable Hilbert space by identifying functions that are equal almost everywhere on $\mathcal{T}$ with respect to $\nu$ (Theorem 1 in \cite{Lin and Yao (2019)}) and is called the tensor Hilbert space along $\mu$.

Now, let $Z$ be a $\mathcal{M}$-valued random function defined on $\mathcal{T}$ such that the intrinsic mean $\mu_{Z}(t)$ defined by $\mu_{Z}(t)=\argmin_{p\in\mathcal{M}}\E((d_{\mathcal{M}}(p,Z(t)))^2)$ exists for all $t\in\mathcal{T}$ and $\Prob(Z(t)\in\mathcal{M}\setminus{{\rm Cut}}(\mu_{Z}(t))\text{~for all~}t\in\mathcal{T})=1$, where $d_{\mathcal{M}}$ is the Riemannian distance function on $\mathcal{M}$ 
and ${{\rm Cut}}(\mu_{Z}(t))$ is the cut locus of $\mu_{Z}(t)$. Let $\Log_{\mu_{Z}(t)}:\mathcal{M}\setminus{{\rm Cut}}(\mu_{Z}(t))\rightarrow T_{\mu_{Z}(t)}\mathcal{M}$ denote the Riemannian logarithm map at $\mu_{Z}(t)$. The specific form of $\Log_{\mu_{Z}(t)}$ for several $\mathcal{M}$ can be found in \cite{Dai and Muller (2018)} and \cite{Lin and Yao (2019)}, for examples. Then, the random process $X:=\Log_{\mu_{Z}}Z$ that maps $t$ to $\Log_{\mu_{Z}(t)}Z(t)$ is well defined almost surely and it satisfies that $\E(X(t))=0_{{\mu_{Z}(t)}}$ for all $t\in\mathcal{T}$, where $0_{{\mu_{Z}(t)}}$ is the zero vector of $T_{\mu_{Z}(t)}\mathcal{M}$. Also, $X$ can be viewed as a random element taking values in $\mathfrak{T}(\mu_{Z})$ provided that $Z$ has continuous sample paths and $\E(\lng X,X\rng_{\mathfrak{T}(\mu_Z)})<\infty$, as described in \cite{Lin and Yao (2019)}. In this case, the PC scores of $X$ are given by the mean 0 uncorrelated iRFPC scores $\lng X,\mpsi_r\rng_{\mathfrak{T}(\mu_{Z})}$, where $\mpsi_r$ are the orthonormal eigenfunctions of the covariance operator $C_{X}:\mathfrak{T}(\mu_{Z})\rightarrow\mathfrak{T}(\mu_{Z})$ of $X$. The intrinsic Riemannian functional linear model with scalar response $Y$ and predictor $Z$ introduced in \cite{Lin and Yao (2019)} takes the form of
\begin{align}\label{lin yao model}
Y=\alpha+\lng X,\beta\rng_{\mathfrak{T}(\mu_{Z})}+\epsilon,
\end{align}
which is a natural extension of the standard functional linear model. Similarly to the existing nonparametric extensions of the standard functional linear model, we may also consider a nonparametric version of (\ref{lin yao model}) that uses a finite number of iRFPC scores. 
For this, we need to estimate $\mu_{Z}$ and $\mpsi_r$. Note that $X$ can be understood as a $\mathcal{M}$-valued random variable when $\mathcal{T}$ is a singleton and $\nu$ is the counting measure. In the latter case, the above PCA is useful when the dimension of $\mathcal{M}$ is large.

Let $\{Z^i:1\leq i\leq n\}$ be i.i.d. observations of $Z$. For the estimation of $\mu_{Z}$, we may use the sample intrinsic mean function $\hat{\mu}_{Z}$ defined by
\begin{align}\label{mu z estimator}
\hat{\mu}_{Z}(t)=\argmin_{p\in\mathcal{M}}\sum_{i=1}^n(d_{\mathcal{M}}(p,Z^i(t)))^2.
\end{align}
As estimators of $\mpsi_r$, we take the orthonormal eigenfunctions $\hat{\mpsi}_r\in\mathfrak{T}(\hat{\mu}_{Z})$ of $\hat{C}_{X}:\mathfrak{T}(\hat{\mu}_{Z})\rightarrow\mathfrak{T}(\hat{\mu}_{Z})$ defined by $\hat{C}_{X}(\hat{V})=n^{-1}\sum_{i=1}^n(\lng\Log_{\hat{\mu}_{Z}}Z^i,\hat{V}\rng_{\mathfrak{T}(\hat{\mu}_{Z})}\cdot\Log_{\hat{\mu}_{Z}}Z^i)$.
Finally, we estimate $\lng\Log_{\mu_{Z}}Z^i,\mpsi_r\rng_{\mathfrak{T}(\mu_{Z})}$ by $\lng\Log_{\hat{\mu}_{Z}}Z^i,\hat{\mpsi}_r\rng_{\mathfrak{T}(\hat{\mu}_{Z})}$. We provide a proposition that gives the vanishing speed of $a_{njl}$ for the case where $\xi^i_{jl}=\lng\Log_{\mu_{Z}}Z^i,\mpsi_r\rng_{\mathfrak{T}(\mu_{Z})}$ and $\tilde{\xi}^i_{jl}=\lng\Log_{\hat{\mu}_{Z}}Z^i,\hat{\mpsi}_r\rng_{\mathfrak{T}(\hat{\mu}_{Z})}$, for some $r$. For this, we introduce more notions on Riemannian manifolds. Let $\|\cdot\|_{\mathfrak{T}(\mu_{Z})}$ denote the norm induced by $\lng\cdot,\cdot\rng_{\mathfrak{T}(\mu_{Z})}$. For any two measurable functions $f,g:\mathcal{T}\rightarrow\mathcal{M}$, we also let $\Gamma_{f,g}:\mathfrak{T}(f)\rightarrow\mathfrak{T}(g)$ denote the parallel transport of vector fields (\cite{Lin and Yao (2019)}) induced by the Levi-Civita connection and a family $\{\gamma(t,\cdot):t\in\mathcal{T}\}$ of smooth curves $\gamma(t,\cdot):[0,1]\rightarrow\mathcal{M}$ such that $\gamma(\cdot,u)$ is measurable for each $u\in[0,1]$ and $\gamma(t,\cdot)$ is the minimizing geodesic between $\gamma(t,0)=f(t)$ and $\gamma(t,1)=g(t)$ for each $t\in\mathcal{T}$. The following new proposition is based on condition \ref{L} given in the Supplementary Material \ref{Lin condition}. Condition \ref{L} is from \cite{Lin and Yao (2019)} and \citet{Lin et al. (2022)} and it is on a Riemannian manifold $\mathcal{M}^*$ and a $\mathcal{M}^*$-valued random function $Z^*$, which is satisfied with various Riemannian manifolds such as compact and connected Riemannian manifolds.

\begin{proposition}\label{iRFPC}
Assume that the first $r$ eigenvalues of $C_X$ have multiplicity one and that condition \ref{L} holds for $\mathcal{M}^*=\mathcal{M}$ and $Z^*=Z$. Then, for $\mpsi_r$ with $\lng\Gamma_{\hat{\mu}_{Z},\mu_{Z}}(\hat{\mpsi}_r),\mpsi_r\rng_{\mathfrak{T}(\mu_{Z})}>0$, $\max_{1\leq i\leq n}|\lng \Log_{\hat{\mu}_{Z}}Z^i,\hat{\mpsi}_r\rng_{\mathfrak{T}(\hat{\mu}_{Z})}-\lng\Log_{\mu_{Z}}Z^i,\mpsi_r\rng_{\mathfrak{T}(\mu_{Z})}|$ achieves the rates given in Proposition \ref{HPC} under the same moment conditions with $\|\mX\|_*$ being replaced by $\|X\|_{\mathfrak{T}(\mu_{Z})}$.
\end{proposition}

We note that the proof of Proposition \ref{iRFPC} is quite different to the one of Proposition \ref{HPC} since the former involves an additional technique in quantifying the difference between $\Gamma_{\hat{\mu}_{Z},\mu_{Z}}(\Log_{\hat{\mu}_{Z}}Z^i)$ and $\Log_{\mu_{Z}}Z^i$. We refer to the Supplementary Material \ref{iRFPC proof} for details. Proposition \ref{iRFPC} also holds for Riemannian non-functional predictors by considering a singleton for $\mathcal{T}$ and the counting measure for $\nu$ or equivalently by treating them as random constant functions defined on $\mathcal{T}$. \qed
\end{customexample}

\begin{customexample}{A.5}\label{Riemannian Hilbertian fd}
(Intrinsic Riemannian Hilbertian singular component scores) Recall the notations and setting in Example \ref{Riemannian fd}. We consider the case where some $\xi_{jl}$ are SC scores obtained from $X$ and $\mY$. This case is differentiated from Example \ref{singular hilbert} in terms of estimation procedure and theory. Define the cross-covariance operators $C_{X\mY}$ and $C_{\mY X}$ and the operator $C_{X\mY X}$ in the same way as in Example \ref{singular hilbert} with $\mbH_*=\mathfrak{T}(\mu_{Z})$. For the orthonormal eigenfunctions $\mphi_r$ of $C_{X\mY X}$ induced from the spectral theorem on $C_{X\mY X}$, we call the mean 0 random variables $\lng X,\mphi_r\rng_{\mathfrak{T}(\mu_{Z})}$ the intrinsic Riemannian Hilbertian SC (iRHSC) scores of $X$. To the best of our knowledge, the notion of iRHSC score has not been introduced in the literature.

To estimate $\lng X,\mphi_r\rng_{\mathfrak{T}(\mu_{Z})}$, we may estimate $\mu_{Z}$ by $\hat{\mu}_{Z}$ defined at (\ref{mu z estimator}) and $\E(\mY)$ by $\bar{\tilde{\mY}}=n^{-1}\od\bigoplus_{i=1}^n\tilde{\mY}^i$. We take the orthonormal eigenfunctions $\hat{\mphi}_r\in\mathfrak{T}(\hat{\mu}_{Z})$ of $\hat{C}_{X\mY X}=\hat{C}_{X\mY}\circ\hat{C}_{\mY X}:\mathfrak{T}(\hat{\mu}_{Z})\rightarrow\mathfrak{T}(\hat{\mu}_{Z})$ as estimators of $\mphi_r$, where $\hat{C}_{X\mY}:\mbH\rightarrow\mathfrak{T}(\hat{\mu}_{Z})$ and $\hat{C}_{\mY X}:\mathfrak{T}(\hat{\mu}_{Z})\rightarrow\mbH$ are defined by
\ba
\hat{C}_{X\mY}(\mh)&=\frac{1}{n-1}\sum_{i=1}^n(\lng\tilde{\mY}^i\om\bar{\tilde{\mY}},\mh\rng\cdot\Log_{\hat{\mu}_{Z}}Z^i),\\
\hat{C}_{\mY X}(\hat{V})&=\frac{1}{n-1}\od\bigoplus_{i=1}^n(\lng\Log_{\hat{\mu}_{Z}}Z^i,\hat{V}\rng_{\mathfrak{T}(\hat{\mu}_{Z})}\od(\tilde{\mY}^i\om\bar{\tilde{\mY}})).
\ea
Finally, we estimate $\lng\Log_{\mu_{Z}}Z^i,\mphi_r\rng_{\mathfrak{T}(\mu_{Z})}$ by $\lng\Log_{\hat{\mu}_{Z}}Z^i,\hat{\mphi}_r\rng_{\mathfrak{T}(\hat{\mu}_{Z})}$. The following proposition gives the vanishing speed of $a_{njl}$ for the case where $\xi^i_{jl}=\lng\Log_{\mu_{Z}}Z^i,\mphi_r\rng_{\mathfrak{T}(\mu_{Z})}$ and $\tilde{\xi}^i_{jl}=\lng\Log_{\hat{\mu}_{Z}}Z^i,\hat{\mphi}_r\rng_{\mathfrak{T}(\hat{\mu}_{Z})}$ for some $r$. To state the proposition, let $\HS(\mathfrak{T}(f))$ denote the set of all Hilbert-Schmidt operators on $\mathfrak{T}(f)$ for each measurable function $f:\mathcal{T}\rightarrow\mathcal{M}$. For any two measurable functions $f,g:\mathcal{T}\rightarrow\mathcal{M}$, we define $\Phi_{f,g}:\HS(\mathfrak{T}(f))\rightarrow \HS(\mathfrak{T}(g))$ by $\Phi_{f,g}(\mathcal{C})(V)=\Gamma_{f,g}(\mathcal{C}(\Gamma_{g,f}(V)))$ for $\mathcal{C}\in\HS(\mathfrak{T}(f))$ and $V\in\mathfrak{T}(g)$, where the definitions of $\Gamma_{f,g}:\mathfrak{T}(f)\rightarrow\mathfrak{T}(g)$ and $\Gamma_{g,f}:\mathfrak{T}(g)\rightarrow\mathfrak{T}(f)$ are given in Example \ref{Riemannian fd}. The map $\Phi_{f,g}$ was introduced in \cite{Lin and Yao (2019)}. 

\begin{proposition}\label{iRHSC}
Assume that the first $r$ eigenvalues of $C_{X\mY X}$ have multiplicity one, that condition \ref{L} holds for $\mathcal{M}^*=\mathcal{M}$ and $Z^*=Z$, and that $\E(\|X\|_{\mathfrak{T}(\mu_{Z})}^2\cdot\|\mY\|^2)<\infty$. Then, it holds that $\|\Phi_{\hat{\mu}_{Z},\mu_{Z}}(\hat{C}_{X\mY X})-C_{X\mY X}\|_{\rm HS}=O_p(n^{-1/2}+b_n)$. Also, for $\mphi_r$ with $\lng\Gamma_{\hat{\mu}_{Z},\mu_{Z}}(\hat{\mphi}_r),\mphi_r\rng_{\mathfrak{T}(\mu_{Z})}>0$, it holds that $\|\Gamma_{\hat{\mu}_{Z},\mu_{Z}}(\hat{\mphi}_r)-\mphi_r\|_{\mathfrak{T}(\mu_{Z})}=O_p(n^{-1/2}+b_n)$. Moreover, $\max_{1\leq i\leq n}|\lng \Log_{\hat{\mu}_{Z}}Z^i,\hat{\mphi}_r\rng_{\mathfrak{T}(\hat{\mu}_{Z})}-\lng\Log_{\mu_{Z}}Z^i,\mphi_r\rng_{\mathfrak{T}(\mu_{Z})}|$ achieves the rates given in Proposition \ref{HSC} under the same moment conditions with $\|\mX\|_*$ being replaced by $\|X\|_{\mathfrak{T}(\mu_{Z})}$.
\end{proposition}

We note that Proposition \ref{iRHSC} is not immediate from Propositions \ref{HSC} and \ref{iRFPC} since deriving the error rates for the cross-covariance operators $C_{X\mY}$ and  $C_{\mY X}$ has not been studied in the literature. We refer to the Supplementary Material \ref{iRHSC proof} for details. \qed
\end{customexample}

Since we allow for different types of perfect/imperfect $\xi_{jl}$ in one model, the above examples indicate that we can treat Euclidean predictors possibly subject to vanishing measurement errors, compositional predictors, density-valued predictors with flexible density supports, Hilbertian functional predictors with flexible function domains and Riemannian functional/non-functional predictors simultaneously. Also, PC and SC scores from different non-Euclidean spaces, and other predictor types induced from different dimension reduction techniques or other imperfect scenarios can be treated.
Our general coverage in the predictor side together with the Hilbertian response $\mY$ covers not only scalar-on-function regression, function-on-scalar regression and function-on-function regression but also various other non-Euclidean regression beyond functional regression. Many of such non-Euclidean regression problems are new even for perfectly observed scalar responses since most of the literature focused on a single predictor or multiple predictors of the same type. 
We now provide examples of $\tilde{\mY}$ and $b_n$. As for $\mX$ above, some of the regression settings below are new even for perfectly observed scalar predictors.

\begin{customexample}{B.1}\label{B.1}
(Responses with vanishing measurement errors) We may consider a version of Example \ref{Measurement error case} for the response variable. In this case, $b_n$ is determined by the moment of a Hilbertian measurement error and the associated vanishing speed. \qed
\end{customexample}

\begin{customexample}{B.2}\label{B.2}
(Imperfect density-valued responses supported on Euclidean spaces) We may consider the case where $\mY^i=Y^i(\cdot)$ is a random density taking values in a Bayes-Hilbert space $\mathfrak{B}^2(\mathcal{S},\mathcal{B}(\mathcal{S}),\nu)$ defined at (\ref{B-H space}) with $\mathcal{S}\subset\mathbb{R}^q$ being a compact and convex set for $q\in\mathbb{N}$ and $\nu$ being the corresponding Lebesgue measure. Suppose that we only observe a random sample $\{Y^*_{ij}:1\leq j\leq n^*_i\}\subset\mathcal{S}$ from $Y^i(\cdot)$ instead of $Y^i(\cdot)$ itself for each $1\leq i\leq n$. The latter scenario is also considered in \citet{Petersen and Muller (2016)} and \citet{Han et al. (2020)} for a class of random densities supported on $[0,1]$. In our case, we take the multivariate density estimator $\tilde{Y}^i(\cdot):=\hat{g}_i(\cdot)/\int_{\mathcal{S}}\hat{g}_i(s)ds$ as $\tilde{\mY}^i$, where
\ba
\hat{g}_i(s)=\frac{w(s,h_i^*)}{n_i^*\cdot(h_i^*)^{q}}\sum_{j=1}^{n_i^*}K^*\bigg(\frac{\|s-Y^*_{ij}\|_{\mathbb{R}^{q}}}{h_i^*}\bigg).
\ea
Here, $h_i^*>0$ is a bandwidth, $K^*:[0,\infty)\rightarrow[0,\infty)$ is a Lipschitz continuous function that vanishes on $(1,\infty)$ and satisfies that $K^*>0$ on $[0,\zeta]$ for some $\zeta\in(0,1)$, $\|\cdot\|_{\mathbb{R}^q}$ is the Euclidean norm on $\mathbb{R}^q$ and $w(s,h_i^*)=(\int_{\mathcal{S}}(h_i^*)^{-q}K^*(\|s-t\|_{\mathbb{R}^{q}}/h_i^*)dt)^{-1}$. The term $w(s,h_i^*)$ is for boundary correction. This new multivariate density estimator $\tilde{Y}^i(\cdot)$ is a bona-fide density on $\mathcal{S}$. \citet{Petersen and Muller (2016)} and \cite{Han et al. (2020)} used a similar density estimator to reconstruct random densities supported on $[0,1]$ but applied a common bandwidth across $1\leq i\leq n$. The following proposition gives a maximum uniform rate of $\tilde{Y}^i(\cdot)$ over $1\leq i\leq n$ and the vanishing speed of $b_n$.

\begin{proposition}\label{imperfect density prop}
Assume that $\log{n}/\min_{1\leq i\leq n}(n_i^*\cdot(h_i^*)^q)=o(1)$, that $\max_{1\leq i\leq n}h_i^*=o(1)$, that $n^c\cdot (\min_{1\leq i\leq n}h_i^*)^{-1}=O(1)$ for some $c<0$ and that $Y(\cdot)$ takes values in $\mathcal{G}(a,b,C)=\{g\in\mathfrak{B}(\mathcal{S},\mathcal{B}(\mathcal{S}),\nu):a\leq\inf_{s\in \mathcal{S}}g(s)\leq\sup_{s\in \mathcal{S}}g(s)\leq b\text{~and~}|g(s_1)-g(s_2)|\leq C\|s_1-s_2\|_{\mathbb{R}^q}\text{~for all~}s_1,s_2\in \mathcal{S}\}$ for some $a,b,C>0$.
Then, $\max_{1\leq i\leq n}\sup_{s\in\mathcal{\mathcal{S}}}|\tilde{Y}^i(s)-Y^i(s)|$ and $\max_{1\leq i\leq n}\|\tilde{\mY}^i\om\mY^i\|$ achieve the rate $O_p(\max_{1\leq i\leq n}h_i^*+\sqrt{\log{n}/\min_{1\leq i\leq n}(n_i^*\cdot(h_i^*)^q)})$.
\end{proposition}

We note that similar assumptions to those in Proposition \ref{imperfect density prop} are adopted in \cite{Han et al. (2020)}. However, the assumptions and conclusion in Proposition \ref{imperfect density prop} are somewhat different to those in \cite{Han et al. (2020)} even when $\mathcal{S}=[0,1]$ and $h_i^*\equiv h^*$ for a common bandwidth $h^*$, as the latter work did not consider Bayes-Hilbert spaces. Proposition \ref{imperfect density prop} is the first asymptotic result on imperfectly observed random densities in the literature of Bayes-Hilbert space. It is also the first asymptotic result in the analysis of random densities supported on a multi-dimensional space, to the best of our knowledge. \qed
\end{customexample}

\begin{customexample}{B.3}\label{B.3}
(Imperfect density-valued responses supported on Riemannian manifolds) The definitions of the mathematical notions in this example can be found in the Supplementary Material \ref{Riemannian definition}. Let $\mathcal{S}$ be a $q$-dimensional compact and connected Riemannian manifold for $q\in\mathbb{N}$ and $\nu$ be the Riemannian volume measure on $\mathcal{S}$. Riemannian volume measures are commonly used regular Borel measures on Riemannian manifolds. We may consider the case where $\mY^i=Y^i(\cdot)$ is a random density taking values in $\mathfrak{B}^2(\mathcal{S},\mathcal{B}(\mathcal{S}),\nu)$ and only a random sample $\{Y^*_{ij}:1\leq j\leq n^*_i\}\subset\mathcal{S}$ from $Y^i(\cdot)$ is observable instead of $Y^i(\cdot)$. In this case, we take the density estimator
\begin{align}\label{Pelletier density}
\tilde{Y}^i(\cdot)=\frac{1}{n_i^*\cdot(h_i^*)^{q}}\sum_{j=1}^{n_i^*}K^*\bigg(\frac{d_{\mathcal{S}}(\cdot,Y^*_{ij})}{h_i^*}\bigg)\frac{1}{\theta_{\mathcal{S}}(\cdot;Y^*_{ij})}
\end{align}
originated by \cite{Pelletier (2005)} as $\tilde{\mY}^i$, where $h_i^*>0$ is a bandwidth, $K^*:[0,\infty)\rightarrow[0,\infty)$ is a Lipschitz continuous function that vanishes on $(1,\infty)$ and satisfies that $\int_{\mathbb{R}^q}K^*(\|u\|_{\mathbb{R}^q})du=1$, $d_{\mathcal{S}}$ is the Riemannian distance function on $\mathcal{S}$ and $\theta_{\mathcal{S}}(\cdot;Y^*_{ij})>0$ is the volume density function of $\mathcal{S}$ at $Y^*_{ij}$. The volume density function is a normalization term such that $\int_{\mathcal{S}}\tilde{Y}^i(s)d\nu(s)=\int_{\mathbb{R}^q}K^*(\|u\|_{\mathbb{R}^q})du=1$. Suppose that there exist positive constants $R$ and $L$ such that, for all $y\in \mathcal{S}$ and $s_1,s_2\in \mathcal{S}$ with $\max\{d_{\mathcal{S}}(s_1,y),d_{\mathcal{S}}(s_2,y)\}<R$, $|(\theta_{\mathcal{S}}(s_1;y))^{-1}-(\theta_{\mathcal{S}}(s_2;y))^{-1}|\leq L\cdot d_{\mathcal{S}}(s_1,s_2)$. 
This condition is satisfied for various $\mathcal{S}$ including finite-dimensional toruses, finite-dimensional spheres and planar shape spaces, as illustrated in \cite{Jeon et al. (2021a)}. The following proposition gives a similar result to Proposition \ref{imperfect density prop}. 

\begin{proposition}\label{imperfect density prop 2}
Under the same assumptions on $n_i^*$, $h_i^*$ and $Y^i(\cdot)$ as in Proposition \ref{imperfect density prop} with $\|s_1-s_2\|_{\mathbb{R}^q}$ being replaced by $d_\mathcal{S}(s_1,s_2)$, the same conclusion holds.
\end{proposition}

We note that dealing with random densities supported on a manifold has not been studied in the literature of density-valued data analysis, to the best of our knowledge. \qed
\end{customexample}

\begin{remark}
We note that \cite{Han et al. (2020)} studied univariate additive regression where the response $Y(\cdot)$ is a random density supported on $[0,1]$ and the scalar predictors are perfectly observed. In the latter work, an additive model is applied to $\Psi(Y(\cdot))(s)\in\mathbb{R}$ for each $s\in[0,1]$, where $\Psi$ is a transformation that maps $Y(\cdot)$ into $L^2([0,1])$. Since it is pointwise regression performed for each $s$, an additional smoothing term is included in its estimation algorithm to ensure that the values of the regression estimator are smooth densities. Such additional smoothing term is not required in our estimation algorithm to be presented, since we deal with infinite-dimensional random densities directly. In fact, if $K^*$ is smooth, then the estimated densities $\tilde{Y}^i(\cdot)$ in Examples \ref{B.2} and \ref{B.3} are smooth, and our algorithm at (\ref{simple algorithm}) in the Supplementary Material \ref{suppl-imple} guarantees that the values of our regression estimator are smooth densities. \qed 
\end{remark}

\begin{customexample}{B.4}\label{B.4}
(Semiparametric regression) Consider the semiparametric regression model
\begin{align}\label{original PLAM}
W=\sum_{k=1}^q\beta_k V_k+\sum_{j=1}^df_j(\xi_j)+\epsilon,
\end{align}
where $W$ and $V_k$ are random variables, $\beta_k\in\mathbb{R}$ are unknown coefficients and $q\in\mathbb{N}$. We may consider (\ref{original PLAM}) as a version of (\ref{model}) for which $\mbH=\mathbb{R}$ and $W-\sum_{k=1}^q\beta_k V_k$ takes the role of $\mY$. From (\ref{original PLAM}), we have
\begin{align}\label{new PLAM}
\E(W|\mxi)=\sum_{k=1}^q\beta_k\E(V_k|\mxi)+\sum_{j=1}^df_j(\xi_j).
\end{align}
By subtracting (\ref{new PLAM}) from (\ref{original PLAM}), we get $W-\E(W|\mxi)=\sum_{k=1}^q\beta_k (V_k-\E(V_k|\mxi))+\epsilon$. Since the latter model takes the form of the standard linear model with response variable $W-\E(W|\mxi)$ and predictors $V_k-\E(V_k|\mxi)$, we may estimate $\beta=(\beta_1,\ldots,\beta_q)\tran$ by $\hat{\beta}=(\sum_{i=1}^n\tilde{V}_i\tilde{V}_i\tran)^{-1}\cdot(\sum_{i=1}^n\tilde{V}_i \tilde{W}_i)$, where $\tilde{V}_i=(V_{i1}-\hat{\E}(V_1|\mxi=\mxi_i),\ldots,V_{iq}-\hat{\E}(V_q|\mxi=\mxi_i))\tran$ and $\tilde{W}_i=W_i-\hat{\E}(W|\mxi=\mxi_i)$. Here, $\hat{\E}(V_k|\mxi)$ and $\hat{\E}(W|\mxi)$ are some estimators of $\E(V_k|\mxi)$ and $\E(W|\mxi)$, respectively. Then, we may estimate $f_j$ based on $W-\sum_{k=1}^q\hat{\beta}_k V_k$ instead of the unobservable $W-\sum_{k=1}^q\beta_k V_k$. Hence, $W-\sum_{k=1}^q\hat{\beta}_k V_k$ takes the role of $\tilde{\mY}$. Note that
$\max_{1\leq i\leq n}\big|W_i-\sum_{k=1}^q\hat{\beta}_k V_{ik}-\big(W_i-\sum_{k=1}^q\beta_k V_{ik}\big)\big|\leq\sum_{k=1}^q\big(|\hat{\beta}_k-\beta_k|\cdot\max_{1\leq i\leq n}|V_{ik}|\big)$.
With suitable estimators $\hat{\E}(V_k|\mxi)$ and $\hat{\E}(W|\mxi)$, one may obtain $\hat{\beta}_k-\beta_k=O_p(n^{-1/2})$. For example, the latter rate is achieved by taking the SBF estimators for $\hat{\E}(V_k|\mxi)$ and $\hat{\E}(W|\mxi)$ as shown in \cite{Yu et al. (2011)} in the case of $L_j\equiv1$ and perfectly observed $\xi_j$. Also, $\max_{1\leq i\leq n}|V_{ik}|$ achieves a rate, say $O_p(v_n)$, depending on the moment of $V_k$, as illustrated in Example \ref{Measurement error case}. In this case, $n^{-1/2}\cdot v_n$ takes the role of $b_n$, and using this and our asymptotic theory presented in Section \ref{asymptotic properties sbf} we may obtain asymptotic properties of an estimator of $f_j$. Similar procedures occur in other semiparametric regression models, and our general framework may be also useful for such models. \qed
\end{customexample}

In Section \ref{extension}, we also treat the case where $\mY$ takes values in a tensor Hilbert space induced from the intrinsic mean function of a Riemannian functional variable and $\tilde{\mY}$ takes values in a tensor Hilbert space induced from an estimated intrinsic mean function. In that case, $\mY$ and $\tilde{\mY}$ take values in different tensor Hilbert spaces.

\begin{remark}
We note that \citet{Jeon et al. (2021b)} also studied additive regression for possibly imperfect Hilbertian responses. However, they considered the case where there exists a transformation $\Psi$ such that $\E(\Psi(\mW,\mY^*)|\mW)=\E(\mY|\mW)$, where $\mW$ is a set of perfectly observed scalar predictors, $\mY$ is an unobservable response variable and $\mY^*$ is an observed variable instead of $\mY$. They showed that such case occurs in regression with missing responses or censored responses. This setting is different to our setting since we do not require the existence of such transformation. In fact, Examples \ref{B.1}--\ref{B.4} cannot be covered by the setting of \citet{Jeon et al. (2021b)}. \qed
\end{remark}

We close this section by giving a detailed comparison between our model and those of \cite{Han et al. (2018)} and \cite{Park et al. (2018)}, both works also covering some functional predictors in the literature of SBF. Some theoretical differences between them can be found in Sections \ref{estimation method} and \ref{asymptotic properties sbf}.

\cite{Han et al. (2018)} considered perfectly observed scalar responses and vanishing errors-in-predictors with its theory being designed for the case where the predictors are the FPC scores obtained from $Q<\infty$ real-valued random functions defined on a compact interval. In fact, it focused on the univariate additive model
\begin{align}\label{Han model}
Y=f_0+\sum_{q=1}^Q\sum_{r=1}^{R_q}f_{qr}(\eta_{qr})+\epsilon,
\end{align}
where $R_q<\infty$, $1\leq q\leq Q$, are fixed numbers and $(\eta_{q1},\ldots,\eta_{qR_q})$ is the tuple of the first $R_q$ FPC scores of the $q$th functional predictor. Note that (\ref{Han model}) is a special case of (\ref{model}) for $\mbH=\mathbb{R}$, $d=\sum_{q=1}^QR_q$, $L_j\equiv1$ and when all the $\xi_{jl}$ are the FPC scores. Model (\ref{model}) also extends to the case where the FPC scores obtained from the same functional predictor consist of the same multivariate predictor, that is, $\xi_q=(\eta_{q1},\ldots,\eta_{qR_q})$ for all $1\leq q\leq Q$. In addition, our work allows for possibly imperfect $\mY$. Hence, our model is more general than that of \cite{Han et al. (2018)}

Meanwhile, \cite{Park et al. (2018)} considered the case where there is only one functional predictor and both the response $Y(\cdot)$ and predictor $X(\cdot)$ are real-valued random functions defined on a compact interval. For the set $\{\theta_r:1\leq r\leq R\}$ of the SC scores of $X(\cdot)$ and the set $\{\varphi_q(\cdot):1\leq q\leq R\}$ of the singular functions of $Y(\cdot)$, it studied the univariate additive model
\begin{align}\label{Park model}
\begin{split}
\E(Y(t)|\theta_1,\ldots,\theta_R)=\E(Y(t))+\sum_{q=1}^R\sum_{r=1}^Rf_{qr}(\theta_r)\varphi_q(t),
\end{split}
\end{align}
where $R<\infty$ is a fixed number. The above model may lose some information of $Y(\cdot)$ as it only considers the finite singular functions of $Y(\cdot)$. Also, it is pointwise regression at each $t$. Our approach differs, as it treats infinite-dimensional responses directly. In addition, (\ref{Park model}) does not cover multiple functional predictors and other types of perfect/imperfect predictors. Moreover, it does not cover other types of perfect/imperfect Hilbertian responses. Hence, our model is more general than that of \cite{Park et al. (2018)}.

\section{Estimation method}\label{estimation method}

In this section, we introduce an estimation method for model (\ref{model}) and provide its non-asymptotic properties. We estimate $\mf_j$ on an arbitrary compact domain $D_j\subset\mathbb{R}^{L_j}$ of interest for $1\leq j\leq d$. Our method and properties generalize those in the case where the support of $\mxi$ is compact, since the former reduces to the latter by setting $\prod_{j=1}^dD_j$ to the support of $\mxi$. Let $p$ denote the density of $\mxi=(\xi_1,\ldots,\xi_d)$ and let $p_0^D=\int_D p(\mx)d\mx>0$, where $D=\prod_{j=1}^dD_j$. 
Define $p^D(\mx)=p(\mx)/p_0^D$ for \mbox{$\mx=(x_1,\ldots,x_d)\in D$}. Note that $p^D$ is a density function on $D$. Also, define its marginal density functions $p^D_j(x_j)=\int_{D_{-j}}p^D(\mx)d\mx_{-j}$ and $p^D_{jk}(x_j,x_k)=\int_{D_{-jk}}p^D(\mx)d\mx_{-jk}$, where $D_{-j}=\prod_{m\neq j}D_m$, $D_{-jk}=\prod_{m\neq j,k}D_m$, and $\mx_{-j}$ and $\mx_{-jk}$ denote the respective $(d-1)$- and $(d-2)$-vector obtained by omitting $x_j$ and $(x_j,x_k)$ in $\mx$.

Note that $\mf_j$ are not identifiable in model (\ref{model}) since $\bigoplus_{j=0}^d\mf_j=\bigoplus_{j=0}^d(\mf_j\op\mc_j)$ for any constants $\mc_j\in\mbH$ satisfying $\bigoplus_{j=0}^d\mc_j=\bzero$. To ensure identifiability, we further impose
\begin{align}\label{constraint 1}
\int_{D_j}\mf_j(x_j)\odot p^D_j(x_j)dx_j=\bzero,\quad1\leq j\leq d.
\end{align}
Note that the constraints (\ref{constraint 1}) determine $\mf_0$ as
\begin{align}\label{f_0}
\begin{split}
\mf_0=\int_D\E(\mY|\mxi=\mx)\od p^D(\mx)d\mx=(p_0^D)^{-1}\od\E(\mY\od\I(\mxi\in D)),
\end{split}
\end{align}
where $\I(\cdot)$ is the indicator function. To introduce our estimation method, we multiply $p^D(\mx)$ on both sides of $\E(\mY|\mxi=\mx)=\mf_0\op\bigoplus_{j=1}^d\mf_j(x_j)$ and integrate them over $D_{-j}$, so that we have
\ba
\int_{D_{-j}}\E(\mY|\mxi=\mx)\od p^D(\mx)d\mx_{-j}=&\,\mf_0\od p_j^D(x_j)\op\mf_j(x_j)\od p_j^D(x_j)\\
&\op\bigoplus_{k\neq j}\int_{D_k}\mf_k(x_k)\od p^D_{jk}(x_j,x_k)dx_k,\quad1\leq j\leq d.
\ea
From this, we get
\begin{align}\label{Bochner equation}
\mf_j(x_j)=\mm_j(x_j)\om\mf_0\om\bigoplus_{k\neq j}\int_{D_k}\mf_k(x_k)\od\frac{p^D_{jk}(x_j,x_k)}{p^D_j(x_j)}dx_k,\quad1\leq j\leq d,
\end{align}
where
\begin{align}\label{mj}
\mm_j(x_j)=(p_j^D(x_j))^{-1}\od\int_{D_{-j}}\E(\mY|\mxi=\mx)\od p^D(\mx)d\mx_{-j}
\end{align}
and $\om$ is the vector subtraction on $\mbH$. Note that the Bochner integrals in (\ref{Bochner equation}) are well defined under mild conditions on $p_j^D$, $p_{jk}^D$ and $\mf_k$. For example, H\"{o}lder's inequality implies that they are well defined if $p^D_j(x_j)>0$, $\int_{D_k}(p^D_{jk}(x_j,x_k))^2/p^D_k(x_k)dx_k<\infty$ and
\begin{align}\label{constraint 2}
\int_{D_k}\|\mf_k(x_k)\|^2p^D_k(x_k)dx_k<\infty
\end{align}
for all $1 \le j \neq k \le d$ and $x_j\in D_j$. 

We estimate the system of Bochner integral equations (\ref{Bochner equation}) based on observations $\{(\tilde{\mxi}^i,\tilde{\mY}^i):1\leq i\leq n\}$, and then solve the estimated system of equations to obtain an estimator of $(\mf_j:1\leq j\leq d)$. Throughout this section, we assume that there exists at least one observation such that $\tilde{\mxi}^i\in D$. For the estimation of (\ref{Bochner equation}), we estimate (\ref{f_0}) by $\hat{\mf}_0=(\hat{p}_0^D\cdot n)^{-1}\od\bigoplus_{i=1}^n\big(\tilde{\mY}^i\od\I(\tilde{\mxi}^i\in D)\big)$, where $\hat{p}_0^D=n^{-1}\sum_{i=1}^n\I(\tilde{\mxi}^i\in D)$. We also estimate $p^D_j(x_j)$, $p^D_{jk}(x_j,x_k)$ and $\mm_j(x_j)$ by kernel smoothing estimators. For this, let $K_j:[0,\infty)\rightarrow[0,\infty)$ be a continuous function such that $K_j$ is positive on $[0,1)$ and is zero on $[1,\infty)$. Also, let $h_j>0$ be a bandwidth and $\|\cdot\|_j$ be the $\ell_2$-norm on $\mathbb{R}^{L_j}$. We take a kernel $K_{h_j}:\mathbb{R}^{L_j}\times\mathbb{R}^{L_j}\rightarrow[0,\infty)$ defined by
\begin{align}\label{kernel def}
K_{h_j}(x_j,u_j)=\frac{h_j^{-L_j}K_j(\|x_j-u_j\|_j/h_j)}{\int_{D_j}h_j^{-L_j}K_j(\|t_j-u_j\|_j/h_j)dt_j}
\end{align}
whenever the denominator is nonzero and $K_{h_j}(x_j,u_j)\equiv|D_j|^{-1}$ otherwise, where $|D_j|$ is the hypervolume of $D_j\subset\mathbb{R}^{L_j}$. The kernel $K_{h_j}$ satisfies
\begin{align}\label{normalization}
\int_{D_j}K_{h_j}(x_j,u_j)dx_j=1~\text{for all}~u_j\in\mathbb{R}^{L_j}.
\end{align}
Using this, we take
\begin{align*}
\hat{p}^D(\mx)&=(\hat{p}_0^D\cdot n)^{-1}\sum_{i=1}^n\prod_{j=1}^dK_{h_j}(x_j,\tilde{\xi}_j^i)\I(\tilde{\mxi}^i\in D),\\
\hat{p}^D_j(x_j)&=\int_{D_{-j}}\hat{p}^D(\mx)d\mx_{-j}=(\hat{p}_0^D\cdot n)^{-1}\sum_{i=1}^nK_{h_j}(x_j,\tilde{\xi}_j^i)\I(\tilde{\mxi}^i\in D),\\
\hat{p}^D_{jk}(x_j,x_k)&=\int_{D_{-jk}}\hat{p}^D(\mx)d\mx_{-jk}=(\hat{p}_0^D\cdot n)^{-1}\sum_{i=1}^nK_{h_j}(x_j,\tilde{\xi}_j^i)K_{h_k}(x_k,\tilde{\xi}_k^i)\I(\tilde{\mxi}^i\in D).
\end{align*}
For the estimation of $\mm_j(x_j)$, we estimate $\E(\mY|\mxi=\mx)\od p^D(\mx)$ in (\ref{mj}) by $n^{-1}\od\bigoplus_{i=1}^n\big(\tilde{\mY}^i\od(\prod_{j=1}^dK_{h_j}(x_j,\tilde{\xi}_j^i)\I(\tilde{\mxi}^i\in D))\big)$. This together with (\ref{normalization}) gives the estimator $\hat{\mm}_j(x_j)=(\hat{p}_j^D(x_j)\cdot\hat{p}_0^D\cdot n)^{-1}\od\bigoplus_{i=1}^n\big(\tilde{\mY}^i\od(K_{h_j}(x_j,\tilde{\xi}_j^i)\I(\tilde{\mxi}^i\in D))\big)$. The SBF estimator of $(\mf_j:1\leq j\leq d)$ on $D$ is then
defined as a solution $(\hat{\mf}_j:1\leq j\leq d)$ of the system of equations
\begin{align}\label{estimated Bochner equation}
\hat{\mf}_j(x_j)=\hat\mm_j(x_j)\ominus\hat{\mf}_0\ominus\bigoplus_{k\neq j}
\int_{D_k}\hat{\mf}_k(x_k)\odot\frac{\hat{p}^D_{jk}(x_j,x_k)}{\hat{p}^D_j(x_j)}dx_k,\quad1\leq j\leq d
\end{align}
subject to the constraints
\begin{align}\label{estimated constraint}
\int_{D_j}\hat{\mf}_j(x_j)\odot\hat{p}^D_j(x_j)dx_j=\bzero,\quad\int_{D_j}\|\hat{\mf}_j(x_j)\|^2\hat{p}^D_j(x_j)dx_j<\infty,\quad1\leq j\leq d.
\end{align}
The above constraints are empirical versions of the constraints (\ref{constraint 1}) and (\ref{constraint 2}). Note that the SBF estimator is not guaranteed to exist nor to be unique since the system of equations (\ref{estimated Bochner equation}) can have no solution or multiple solutions. Hence, it is essential to check the existence and uniqueness of the estimator. For this, we need a condition on $h_j$. Let $B_j(x_j,R)$ denote the open ball in $\mathbb{R}^{L_j}$ centered at $x_j$ with radius $R>0$. The following condition is not an asymptotic condition but a condition on a given dataset.

\begin{customcon}{(A)}\label{A}\leavevmode
For each $1\leq j\leq d$ and $x_j\in D_j$, there exists at least one observation $\tilde{\mxi}^i\in D$ such that $\tilde{\xi}^i_j\in B_j(x_j,h_j)$.
\end{customcon}

Condition \ref{A} is a minimal requirement such that $\hat{p}^D_j(x_j)>0$ for all $x_j\in D_j$ since $K_j$ is assumed to vanish on $[1,\infty)$. The latter property on $K_j$ is essential for the SBF technique. Note that the Bochner integrals in (\ref{estimated Bochner equation}) are well defined under condition \ref{A} since it implies that $\inf_{x_j\in D_j}\hat{p}^D_j(x_j)>0$ and $\sup_{(x_j,x_k)\in D_j\times D_k}\hat{p}^D_{jk}(x_j,x_k)<\infty$, as proved in the Supplementary Material \ref{proof existence}. To state the existence and uniqueness of the SBF estimator, we define $\mathcal{B}(D_j)=\{B\cap D_j:B\text{~is a Borel set in~}\mathbb{R}^{L_j}\}$. Note that $\mathcal{B}(D_j)$ forms a sigma-field on $D_j$. We also define a probability measure $\hat{P}_{\mxi}^D$ on the product sigma-field $\bigotimes_{j=1}^{d}\mathcal{B}(D_j)$ by $\hat{P}_{\mxi}^D(A)=\int_A\hat{p}^D(\mx)d\mx$. Let $\Leb_j$ denote the $L_j$-dimensional Lebesgue measure. The following proposition holds for a given dataset, and thus it is not required that $\{(\mxi^i,\mY^i):1\leq i\leq n\}$ are i.i.d. observations of $(\mxi,\mY)$. Also, it is not required that $\max_{1\leq i\leq n}|\tilde{\xi}^i_{jl}-\xi^i_{jl}|=O_p(a_{njl})$ and $\max_{1\leq i\leq n}\|\tilde{\mY}^i\om\mY^i\|=O_p(b_n)$ for some $a_{njl}=o(1)$ and $b_n=o(1)$. This kind of non-asymptotic property has not been covered in \cite{Han et al. (2018)}, \cite{Park et al. (2018)}, \cite{Han et al. (2020)} and \cite{Lin et al. (2022)}.

\begin{proposition}\label{existence}
Assume that condition \ref{A} holds. Then, there exists a solution of (\ref{estimated Bochner equation}) satisfying (\ref{estimated constraint}). In addition, if $(\hat{\mf}_j:1\leq j\leq d)$ and $(\hat{\mf}^\star_j:1\leq j\leq d)$ are solutions of (\ref{estimated Bochner equation}), then $\bigoplus_{j=1}^d\hat{\mf}_j(x_j)=\bigoplus_{j=1}^d\hat{\mf}^\star_j(x_j)$ almost everywhere on $D$ with respect to $\hat{P}_{\mxi}^D$. Moreover, if the solutions $(\hat{\mf}_j:1\leq j\leq d)$ and $(\hat{\mf}^\star_j:1\leq j\leq d)$ satisfy (\ref{estimated constraint}), and there exists at least one observation $\tilde{\mxi}^i\in D\cap\prod_{j=1}^dB_j(x_j,h_j)$ for each $\mx\in D$, then
$\hat{\mf}_j(x_j)=\hat{\mf}_j^\star(x_j)$ almost everywhere on $D_j$ with respect to $\Leb_j$ for all $1\leq j\leq d$.
\end{proposition}

Note that the uniqueness of the individual SBF estimators $\hat{\mf}_j$ is based on the stronger condition that there exists at least one observation $\tilde{\mxi}^i\in D\cap\prod_{j=1}^dB_j(x_j,h_j)$ for each $\mx\in D$.
In Section \ref{asymptotic properties sbf}, we show that such condition is not required when we consider the uniqueness of the individual SBF estimators in an asymptotic sense. The proposition above only tells about the existence and uniqueness, and it does not provide the closed form of the SBF estimator. Hence, we obtain it via a numerical algorithm. For this, we take any initial estimator $(\hat{\mf}^{[0]}_j:1\leq j\leq d)$ satisfying the constraints (\ref{estimated constraint}).
For example, we may take $\hat{\mf}^{[0]}_j\equiv\bzero$ for all $1\leq j\leq d$. 
For the $r$th iteration with $r\geq1$, we subsequently update the estimator as follows:
\begin{align}\label{algorithm}
\begin{split}
\hat{\mf}^{[r]}_j(x_j)=\hat\mm_j(x_j)&\ominus\hat{\mf}_0\ominus\bigoplus_{k<j}
\int_{D_k}\hat{\mf}^{[r]}_k(x_k)\odot\frac{\hat{p}^D_{jk}(x_j,x_k)}{\hat{p}^D_j(x_j)}dx_k\\
&\ominus\bigoplus_{k>j}
\int_{D_k}\hat{\mf}^{[r-1]}_k(x_k)\odot\frac{\hat{p}^D_{jk}(x_j,x_k)}{\hat{p}^D_j(x_j)}dx_k,\quad1\leq j\leq d.
\end{split}
\end{align}
We call the above algorithm the SBF algorithm. One can show that $(\hat{\mf}^{[r]}_j:1\leq j\leq d)$ satisfy (\ref{estimated constraint}) under condition \ref{A}. The SBF algorithm involving the Bochner integrals can be easily implemented using Lebesgue integrals in case each $\hat{\mf}^{[0]}_j(x_j)$ takes the form of $\bigoplus_{i=1}^n(w^{i,[0]}_j(x_j)\od\tilde{\mY}^i)$ for some weights $w^{i,[0]}_j(x_j)\in\mathbb{R}$; see the Supplementary Material \ref{suppl-imple} for details.
The next two propositions demonstrate that the sum $\hat\mf^{[r]}(\mx):=\hat{\mf}_0\op\bigoplus_{j=1}^{d}\hat\mf^{[r]}_j(x_j)$ and the related individual maps $\hat{\mf}^{[r]}_j(x_j)$, obtained at the $r$th iteration in (\ref{algorithm}), converge respectively to the regression estimator $\hat\mf(\mx):=\hat{\mf}_0\op\bigoplus_{j=1}^{d}\hat\mf_j(x_j)$ and to $\hat{\mf}_j(x_j)$ as $r\rightarrow\infty$. Both propositions also hold for a given dataset. This kind of non-asymptotic property has not been studied in \cite{Han et al. (2018)}, \cite{Park et al. (2018)}, \cite{Han et al. (2020)} and \cite{Lin et al. (2022)}. 

\begin{proposition}\label{convergence1}
Assume that condition \ref{A} holds. Then, it holds that (i) there exist constants $\tilde{c}_1>0$ and $\tilde{\varrho}_1\in(0,1)$ such that
$\int_{D}\|\hat{\mf}^{[r]}(\mx)\ominus\hat{\mf}(\mx)\|^2d\hat{P}_{\mxi}^D(\mx)\leq\tilde{c}_1\cdot\tilde{\varrho}_1^r$ for all $r\geq0$;
(ii) $\lim_{r\rightarrow\infty}\|\hat{\mf}^{[r]}(\mx)\ominus\hat{\mf}(\mx)\|=0$ almost everywhere on $D$ with respect to $\hat{P}_{\mxi}^D$;
(iii) for each $\varepsilon>0$, there exists a Borel set $S(\varepsilon)\subset D$ with $\hat{P}_{\mxi}^D(D\setminus S(\varepsilon))<\varepsilon$ such that $\lim_{r\rightarrow\infty}\sup_{\mx\in S(\varepsilon)}\|\hat{\mf}^{[r]}(\mx)\ominus\hat{\mf}(\mx)\|=0$.
\end{proposition}

Note that $\int_{D}\|\hat{\mf}^{[r]}(\mx)\ominus\hat{\mf}(\mx)\|^2d\hat{P}_{\mxi}^D(\mx)=\int_{D}\|\hat{\mf}^{[r]}(\mx)\ominus\hat{\mf}(\mx)\|^2\hat{p}^D(\mx)d\mx$. The first result in Proposition \ref{convergence1} tells that the $L^2$ convergence speed of $\hat{\mf}^{[r]}$ to $\hat{\mf}$ is geometric. The third result in Proposition \ref{convergence1} shows that $\hat{\mf}^{[r]}$ converges to $\hat{\mf}$ almost uniformly.

\begin{proposition}\label{convergence1 individual}
Assume that there exists at least one observation $\tilde{\mxi}^i\in D\cap\prod_{j=1}^dB_j(x_j,h_j)$ for each $\mx\in D$. Then, it holds that (i) there exist constants $\tilde{c}_2>0$ and $\tilde{\varrho}_2\in(0,1)$ such that, for all $1\leq j\leq d$ and $r\geq0$,
$\int_{D_j}\|\hat{\mf}_j^{[r]}(x_j)\ominus\hat{\mf}_j(x_j)\|^2dx_j\leq\tilde{c}_2\cdot\tilde{\varrho}_2^r$;
(ii) for all $1\leq j\leq d$, $\lim_{r\rightarrow\infty}\|\hat{\mf}_j^{[r]}(x_j)\ominus\hat{\mf}_j(x_j)\|=0$ almost everywhere on $D_j$ with respect to $\Leb_j$;
(iii) for each $1\leq j\leq d$ and $\varepsilon>0$, there exists a Borel set $S_j(\varepsilon)\subset D_j$ with $\Leb_j(D_j\setminus S_j(\varepsilon))<\varepsilon$ such that
$\lim_{r\rightarrow\infty}\sup_{x_j\in S_j(\varepsilon)}\|\hat{\mf}_j^{[r]}(x_j)\ominus\hat{\mf}_j(x_j)\|=0$.
\end{proposition}


\section{Asymptotic properties of SBF estimator}\label{asymptotic properties sbf}

\setcounter{equation}{0}
\setcounter{subsection}{0}

In this section we show that the individual SBF estimators $\hat{\mf}_j$ are uniquely determined with probability tending to one and that $\hat{\mf}^{[r]}_j$ converges to $\hat{\mf}_j$ as $r\rightarrow\infty$ with probability tending to one for each $1\leq j\leq d$. We also study the error rates for $\|\hat{\mf}_j\om\mf_j\|$. In addition, we derive the asymptotic distributions of $\hat{\mf}_j$ and $\hat{\mf}$.

\subsection{Asymptotic uniqueness and algorithm convergence for individual SBF estimators}\label{individual sbf estimator}

We first study the uniqueness of $\hat{\mf}_j$ and the convergence of $\hat{\mf}^{[r]}_j$ to $\hat{\mf}_j$ in an asymptotic sense. For this, we define $a_{nj}=\max_{1\leq l\leq L_j}a_{njl}$. Note that $\max_{1\leq i\leq n}\|\tilde{\xi}_j^i-\xi_j^i\|_j=O_p(a_{nj})$. We also define $D_j^+(R)=\bigcup_{x_j\in D_j}\bar{B}_j(x_j,R)$ and $D^+(R)=\prod_{j=1}^dD_j^+(R)$, where $\bar{B}_j(x_j,R)$ is the closed ball in $\mathbb{R}^{L_j}$ centered at $x_j\in D_j$ with radius $R>0$. We introduce the following conditions.

\begin{customcon}{(B)}\label{B}\leavevmode
For all $1\leq j\neq k\leq d$, the following holds.
\begin{itemize}
\item[{\rm(B1)}] $\E(\|\mY\|^{\alpha})<\infty$ for some $\alpha>2$ and $\E(\|\mY\|^{\alpha}|\mxi=\cdot)$ is bounded on $D^+(\varepsilon)$ for some $\varepsilon>0$.
\item[{\rm(B2)}] $p$ is continuous and positive on $D$ and is bounded on $D^+(\varepsilon)$ for some $\varepsilon>0$.
\item[{\rm(B3)}] $D_j$ is convex.
\item[{\rm(B4)}] $K_j(\|\cdot\|_j)$ is differentiable on $\mathbb{R}^{L_j}$ and all the partial derivatives of $K_j(\|\cdot\|_j)$ are Lipschitz continuous on $\mathbb{R}^{L_j}$.
\item[{\rm(B5)}] (i) $h_j=o(1)$, (ii) $n^{-1}\cdot(\log{n})\cdot h_j^{-L_j}\cdot h_k^{-L_k}=o(1)$ and (iii) $n^{-1+2\beta_j+2/\alpha}\cdot h_j^{-L_j}=o(1)$ for $\alpha$ in (B1) and some $\beta_j>0$.
\item[{\rm(B6)}] (i) $h_j^{-L_j}\cdot a_{nj}=o(1)$ and (ii) $h_j^{-L_j-2}\cdot a_{nj}^2=o(1)$.
\item[{\rm(B7)}] $\lim_{n\rightarrow\infty}\Prob\big(\int_{D_j}\|\hat{\mf}^{[0]}_j(x_j)\|^2dx_j<C\big)=1$ for some $C>0$.
\end{itemize}
\end{customcon}

Conditions (B1) and (B2) are versions of standard conditions in nonparametric regression. Let $\partial D_j$ denote the boundary of $D_j$. Condition (B3) guarantees that $\Leb_j(\partial D_j)=0$ and $\Leb_j(\bigcup_{x_j\in \partial D_j}\bar{B}_j(x_j,R))=O(R)$, which are necessary in our asymptotic analysis. The smoothness conditions on $K_j$ in (B4) are required to deal with the errors-in-variables efficiently in our asymptotic analysis. 
Examples of $K_j$ satisfying (B4) and the underlying conditions on $K_j$ given immediately above (\ref{kernel def}) include the biweight-type kernel
\begin{align}\label{biweight kernel}
K_j(t)=(1-t^2)^2\cdot\I(t\in[0,1]).
\end{align}
Condition (B5) is a bandwidth condition. 
When one takes $h_j$ of the form $h_j\asymp n^{-c_j}$ for $c_j>0$, (B5) holds provided that $L_jc_j+L_kc_k<1$ and $L_jc_j<(1-2/\alpha)$. If one takes $c_j=1/(L_j+4)$, which is a popular choice for a $L_j$-dimensional predictor, (B5) holds provided that $L_jL_k<16$ and $\alpha>(L_j+4)/2$. Note that we may not use large $L_j$ due to the curse of dimensionality and thus the sufficient conditions are not restrictive. We emphasize that, in case there is a high/infinite-dimensional predictor, the number of its component scores to be used in the model is not necessarily one of the $L_j$'s since they can be divided and placed in different component maps $\mf_j$. Condition (B6) tells that the vanishing speed of $a_{nj}$ should be faster than those of $h_j^{L_j}$ and $h_j^{L_j/2+1}$. Note that (B6)-(ii) is implied by (B5)-(i) and (B6)-(i) if $L_j\geq2$. In case one takes $h_j\asymp n^{-c_j}$ and $a_{nj}=n^{-r_j}$, (B6)-(i) is satisfied if $r_j>L_jc_j$, and (B6)-(ii) with $L_j=1$ is satisfied if $r_j>(3/2)c_j$. Hence, in case $c_j=1/(L_j+4)$ and $r_j=1/2$, for example, (B6) is satisfied if $L_j<4$. We note that \cite{Han et al. (2018)} and \cite{Park et al. (2018)} fixed their bandwidth rates to $n^{-1/5}$ and the error rates for the functional PC or SC scores to $O_p(n^{-(\beta-1)/(2\beta)})$ for some $\beta>5$, to study their models (\ref{Han model}) and (\ref{Park model}). Our general rates for $h_j$ and $a_{nj}$ in (B5) and (B6) cover the specific rates in \cite{Han et al. (2018)} and \cite{Park et al. (2018)}. We emphasize that obtaining asymptotic results with flexible $h_j$ and $a_{nj}$ in multivariate Hilbertian additive models is much more complex than the one with fixed $h_j$ and $a_{nj}$ in univariate additive models on $\mathbb{R}$, since one needs the exact evaluation of more complex stochastic terms in asymptotic expansions. Condition (B7) is a weak condition on the choice of an initial estimator. For example, the choice $\hat{\mf}_j^{[0]}\equiv\bzero$ satisfies (B7).

\begin{theorem}\label{convergence2}
Assume that condition \ref{B} holds. Then, with probability tending to one, the individual SBF estimators are unique in the sense that if $(\hat{\mf}_j:1\leq j\leq d)$ and $(\hat{\mf}^\star_j:1\leq j\leq d)$ are solutions of (\ref{estimated Bochner equation}) satisfying the constraints (\ref{estimated constraint}), then $\hat{\mf}_j=\hat{\mf}^\star_j$ almost everywhere on $D_j$ with respect to $\Leb_j$ for all $1\leq j\leq d$. Also, the conclusion in Proposition \ref{convergence1 individual} holds with probability tending to one.
\end{theorem}

Note that, for the proof of Theorem \ref{convergence2}, we prove the uniform and $L^2$ consistency of the marginal density estimators (Lemmas \ref{marginal density convergence 1} and \ref{marginal density convergence 2} in the Supplementary Material \ref{collection of lemmas}). Those results are of importance in their own right since density estimation with general errors-in-variables is also an important topic that has not been studied well. 
Proving such consistency involves much more complex analysis than the one for the error-free case; see the Supplementary Material for details. 

\subsection{Error rates}

In this section, we study the rates of convergence of the SBF estimator. For this, we consider two types of smoothness on the maps $\mf_j$: H\"{o}lder continuity or differentiability. We first consider the case where $\mf_j$ are H\"{o}lder continuous. This type of $\mf_j$ has not been treated in \cite{Han et al. (2018)}, \cite{Park et al. (2018)}, \cite{Han et al. (2020)} and \cite{Lin et al. (2022)}.
To do so, we introduce the next condition.
\begin{customcon}{(C)}\label{C}\leavevmode
For all $1\leq j\leq d$, the following holds.
\begin{itemize}
\item[{\rm(C1)}] $\mf_j$ is H\"{o}lder continuous on $D_j^+(\varepsilon)$ with exponent $\nu_j>0$ for some $\varepsilon>0$.
\item[{\rm(C2)}] $\E(\mepsilon^i|\mxi^1,\ldots\mxi^n,\tilde{\mxi}^1,\ldots,\tilde{\mxi}^n)=\bzero$ for all $1\leq i\leq n$, $\mepsilon^i$ are conditionally independent given $(\mxi^1,\ldots\mxi^n,\tilde{\mxi}^1,\ldots,\tilde{\mxi}^n)$ and $\max_{1\leq i\leq n}\E(\|\mepsilon^i\|^2|\mxi^1,\ldots\mxi^n,\tilde{\mxi}^1,\ldots,\tilde{\mxi}^n)<C$ almost surely for some $C>0$.
\item[{\rm(C3)}] $h_j^{-L_j}\cdot a_{nj}=O((\log{n})^{-1})$.
\item[{\rm(C4)}] $n^{-2\beta_j}\cdot(\log{n})\cdot(\sum_{j=1}^da_{nj})^{2/\alpha}\cdot(h_j^{-L_j}\cdot a_{nj}+\sum_{k\neq j}a_{nk})^{-1}=O(1)$ for $\alpha$ in (B1) and $\beta_j$ in (B5).
\end{itemize}
\end{customcon}

Condition (C2) is a weak condition on the error term $\mepsilon$ of model (\ref{model}). A simple sufficient condition for (C2) is that $(\mepsilon^i:1\leq i\leq n)$ is independent of $(\mxi^1,\ldots\mxi^n,\tilde{\mxi}^1,\ldots,\tilde{\mxi}^n)$. In this case, (C2) easily holds since $\mepsilon^i$ are i.i.d. and $\mepsilon$ satisfies that $\E(\mepsilon)=\bzero$ and $\E(\|\mepsilon\|^2)<\infty$.
Condition (C2) is required to get conditional large deviation inequalities for stochastic terms involving $\mepsilon^i$. The latter large deviation inequalities are required since $\tilde{\mxi}^i$ are not necessarily independent across $1\leq i\leq n$.
Condition (C3) is slightly stronger than (B6)-(i), but it still holds with various combinations of $(h_j,a_{nj})$. Condition (C4) is also a weak condition required for a fast uniform error rate of a stochastic term. Note that larger $\beta_j$ is required for (C4) to hold while smaller $\beta_j$ is required for (B5)-(iii) to hold. Both conditions can hold simultaneously with various combinations of $(h_j,a_{nj})$. As a simple example, consider the case where $h_j\asymp n^{-c_j}$ and $a_{nj}\asymp n^{-r}$ for all $1\leq j\leq d$. Then, the conditions hold if $L_jc_j<(1-2/\alpha)$ and $r<1$.

\begin{theorem}\label{holder rate}
Assume that conditions (B1)--(B5), (B6)-(ii) and \ref{C} hold. Then, for all $1\leq j\leq d$, it holds that (i) $\big(\int_{D_j}\|\hat{\mf}_j(x_j)\om\mf_j(x_j)\|^2dx_j\big)^{1/2}=O_p\big(\sum_{j=1}^dh_j^{\nu_j}+\sum_{j=1}^d(nh_j^{L_j})^{-1/2}+b_n\big)$; (ii) $\|\hat{\mf}_j(x_j)\om\mf_j(x_j)\|=O_p\big(\sum_{j=1}^dh_j^{\nu_j}+\sum_{j=1}^d(nh_j^{L_j})^{-1/2}+b_n\big)$ for all $x_j\in D_j$; (iii)
$\sup_{x_j\in D_j}\|\hat{\mf}_j(x_j)\om\mf_j(x_j)\|=O_p\big(\sum_{j=1}^dh_j^{\nu_j}+\sum_{k\neq j}(nh_k^{L_k})^{-1/2}+(nh_j^{L_j}/\log{n})^{-1/2}+b_n\big)$.
\end{theorem}

Theorem \ref{holder rate} tells that our estimator avoids the curse of dimensionality since the error rates do not involve $(n\prod_{j=1}^dh_j^{L_j})^{-1/2}$, which may appear in full-dimensional kernel regression. When $h_j\asymp h$, $\nu_j=\nu$ and $L_j=L$ for all $1\leq j\leq d$, the $L^2$ and pointwise error rates are simplified to $h^\nu+(nh^L)^{-1/2}+b_n$ and the uniform error rate is simplified to $h^\nu+(nh^L/\log{n})^{-1/2}+b_n$. When the response variable is perfectly observable and there is a single perfectly observable predictor of $L$-dimension, a standard kernel regression estimator achieves the above error rates without the term $b_n$. Hence, our estimator achieves the optimal error rates subject to the possibly negligible rate $b_n$. Recall that $b_n=0$ in case $\mY$ is perfectly observable. We now investigate the rates of convergence for differentiable $\mf_j$. For this, we make additional conditions. Let $\bzero_j$ denote the zero vector in $\mathbb{R}^{L_j}$.

\begin{customcon}{(D)}\label{D}\leavevmode
For all $1\leq j\leq d$, the following holds.
\begin{itemize}
\item[{\rm(D1)}] $\mf_j$ is Lipschitz continuous on $D_j^+(\varepsilon)$ for some $\varepsilon>0$ and is twice continuously Fr\'{e}chet differentiable on $D_j$. 
\item[{\rm(D2)}] $p$ is continuously differentiable on $D$.
\item[{\rm(D3)}] $\int_{B_j(\bzero_j,1)}t_j K_j(\|t_j\|_j)dt_j=\bzero_j$.
\end{itemize}
\end{customcon}

Conditions (D1) and (D2) are versions of standard conditions in nonparametric regression.
Condition (D3) is satisfied with many kernel functions. For example, it is satisfied with the kernel function defined at (\ref{biweight kernel}). To state the rates of convergence, we define $D_j^-(R)=D_j\setminus\bigcup_{x_j\in \partial D_j}B_j(x_j,R)$. Note that $D_j^-(2h_j)$ takes the role of the interior region of $D_j$ and $D_j\setminus D_j^-(2h_j)$ takes the role of the boundary region of $D_j$. We also define $\mathcal{A}_{nj}=(h_j^{-L_j}\cdot a_{nj}+h_j^{-L_j-2}\cdot a^2_{nj})\cdot\sum_{k=1}^dh_k+\sum_{k\neq j}a_{nk}$.
Clearly, $\mathcal{A}_{nj}=o(1)$ under (B5) and (B6). Define $\delta^{(1)}_{nj}=\sum_{k\neq j}h_k^2+\sum_{j=1}^d (nh_j^{L_j})^{-1/2}+\sum_{j=1}^d\mathcal{A}_{nj}+b_n$ and $\delta^{(2)}_{nj}=\sum_{k\neq j}h_k\cdot(nh_k^{L_k}/\log{n})^{-1/2}+(nh_j^{L_j}/\log{n})^{-1/2}$.

\begin{theorem}\label{differentiable rate}
Assume that conditions (B1)--(B5), (B6)-(ii), (C2)--(C4) and \ref{D} hold. Then, for all $1\leq j\leq d$, it holds that (i) $\big(\int_{D_j^-(2h_j)}\|\hat{\mf}_j(x_j)\om\mf_j(x_j)\|^2dx_j\big)^{1/2}=O_p\big(h_j^2+\delta^{(1)}_{nj}\big)$ and $\big(\int_{D_j}\|\hat{\mf}_j(x_j)\om\mf_j(x_j)\|^2dx_j\big)^{1/2}=O_p\big(h_j^{3/2}+\delta^{(1)}_{nj}\big)$; (ii) $\|\hat{\mf}_j(x_j)\om\mf_j(x_j)\|=O_p\big(h_j^2+\delta^{(1)}_{nj}\big)$ for all $x_j\in D_j^-(2h_j)$ and $\|\hat{\mf}_j(x_j)\om\mf_j(x_j)\|=O_p\big(h_j+\delta^{(1)}_{nj}\big)$ for all $x_j\in D_j\setminus D_j^-(2h_j)$; (iii) $\sup_{x_j\in D_j^-(2h_j)}\|\hat{\mf}_j(x_j)\om\mf_j(x_j)\|=O_p\big(h_j^2+\delta^{(1)}_{nj}+\delta^{(2)}_{nj}\big)$ and $\sup_{x_j\in D_j}\|\hat{\mf}_j(x_j)\om\mf_j(x_j)\|=O_p\big(h_j+\delta^{(1)}_{nj}+\delta^{(2)}_{nj}\big)$.
\end{theorem}

Theorem \ref{differentiable rate} shows that our estimator with differentiable type $\mf_j$ again avoids the curse of dimensionality. The errors-in-variables induce the additional rate $\sum_{j=1}^d\mathcal{A}_{nj}+b_n$. In many cases, however, the rate $\sum_{j=1}^d\mathcal{A}_{nj}$ is negligible. As a simple example, consider the case where $L_j\equiv L$, $h_j\asymp n^{-c}$ and $a_{nj}\asymp n^{-r}$ with $c, r>0$ for all $1\leq j\leq d$. In this case, $\mathcal{A}_{nj}=O(h_j^2+(nh_j^{L_j})^{-1/2})$ provided that $r\geq \max\{c(L+1),c(L/2-1)+1/2\}$. In particular, when $c=1/(L+4)$ and $r=1/2$, it holds provided that $L<3$. This indicates that our estimator can achieve the optimal error rates not only in univariate additive models but also in bivariate additive models when $\mf_j$ are of differential type. One may find other sufficient conditions with some $L_j\geq3$ for which $\sum_{j=1}^d\mathcal{A}_{nj}$ is negligible. Even when $\sum_{j=1}^d\mathcal{A}_{nj}$ is not negligible, Theorem \ref{differentiable rate} still gives useful information on how errors-in-variables with flexible $L_j$, $h_j$ and $a_{nj}$ affect error rates in additive regression. 

\subsection{Asymptotic distribution}\label{asymptotic distribution section}

In this section, we investigate the asymptotic distributions of $\hat{\mf}_j$ and $\hat{\mf}$. We first collect some conditions. To state the conditions, let $\{\me_m:m\geq1\}$ be an orthonormal basis of $\mbH$. Define $L_{\rm max}=\max_{1\leq j\leq d}L_j$ and $\mathcal{R}_n=\sum_{j=1}^d\big((nh^{L_j}/\log{n})^{-1/2}\cdot\big(h_j^{-L_j}\cdot a_{nj}+\sum_{k\neq j}a_{nk}\big)^{1/2}\big)+\sum_{j=1}^d\mathcal{A}_{nj}+b_n$.

\begin{customcon}{(E)}\label{E}\leavevmode
For all $m, m'\geq 1$ and $1\leq j\neq k\leq d$, the following holds.
\begin{itemize}
\item[{\rm(E1)}] $\E(\lng\mepsilon,\me_m\rng\lng\mepsilon,\me_{m'}\rng|\xi_j=\cdot,\mxi\in D)$ is continuous on $D_j$ and $\E(\lng\mepsilon,\me_m\rng\lng\mepsilon,\me_{m'}\rng|\xi_j=\cdot,\xi_k=\cdot,\mxi\in D)$ is bounded on $D_j\times D_k$.
\item[{\rm(E2)}] $h_j\asymp n^{-1/(L_j+4)}$ and $\alpha_j:=\lim_{n\rightarrow\infty}h_j\cdot n^{1/(L_{\rm max}+4)}$ exists.
\item[{\rm(E3)}] $n^{2/(L_{\rm max}+4)}\cdot \mathcal{R}_n=o(1)$.
\end{itemize}
\end{customcon}

Condition (E1) is a mild condition. When $\mbH=\mathbb{R}$, it reduces to the conditions on the conditional variances of $\mepsilon$, which are typical assumptions in nonparametric regression. Versions of (E1) are also used in \cite{Jeon and Park (2020)}, \cite{Jeon et al. (2021a)} and \cite{Jeon et al. (2021b)}. Here, we consider the optimal bandwidth rates as described in condition (E2). Note that $\alpha_j=0$ for $j$ with $L_j<L_{\rm max}$ under (E2). Condition (E3) is satisfied in many cases. As a simple example, consider the case where $L_j\equiv L$, $h_j\asymp n^{-1/(L+4)}$ and $a_{nj}\asymp n^{-r}$ with $r>0$ for all $1\leq j\leq d$. Then, (E3) holds provided that $r>(L+1)/(L+4)$ and $n^{2/(L+4)}\cdot b_n=o(1)$.

We define several maps to describe the asymptotic bias and variance of our estimator. Below, let $\mathfrak{D}^1$ denote the Fr\'{e}chet differential operator. For a map $\mg_j:D_j\rightarrow\mbH$, $\mathfrak{D}^1\mg_j$ is defined as a map from $D_j$ to $\mathcal{L}(\mathbb{R}^{L_j},\mbH)$, where $\mathcal{L}(\mathbb{R}^{L_j},\mbH)$ is the space of all bounded linear operators from $\mathbb{R}^{L_j}$ to $\mbH$. 
For $1\leq j\neq k\leq d$, define
\begin{align*}
\mdelta_j(x_j)&=(p^D_j(x_j))^{-1}\od\mathfrak{D}^1\mf_j(x_j)\bigg(\int_{B_j(\bzero_j,1)}t_j\cdot\mathfrak{D}^1p^D_j(x_j)(t_j)\cdot K_j(\|t_j\|_j)dt_j\bigg),\\
\mdelta_{jk}(x_j,x_k)
&=(p^D_{jk}(x_j,x_k))^{-1}\od\mathfrak{D}^1\mf_k(x_k)\bigg(\int_{B_k(\bzero_k,1)}t_k\cdot\mathfrak{D}^1_2\,p^D_{jk}(x_j,x_k)(t_k)\cdot K_k(\|t_k\|_k)dt_k\bigg),
\end{align*}
where $\mathfrak{D}^1_2\,p^D_{jk}(x_j,x_k)\in\mathcal{L}(\mathbb{R}^{L_j},\mathbb{R})$ denotes the partial Fr\'{e}chet derivative of $p^D_{jk}$ at $(x_j,x_k)$ with respect to the second argument.
We also define
\[
\tilde\mDelta_j(x_j)=\alpha^2_j\od\mdelta_j(x_j)\oplus \bigoplus_{k\neq j} \int_{D_k}
\mdelta_{jk}(x_j,x_k) \od \bigg(\alpha^2_k\cdot\frac{p^D_{jk}(x_j,x_k)}{p^D_j(x_j)}\bigg)dx_k,
\]
where $\alpha_j$ is the constant in (E2). Let $(\mDelta_j:1\leq j\leq d)$ be a solution of the system of equations
\begin{align}\label{Delta}
\mDelta_j(x_j)=\tilde{\mDelta}_j(x_j)\om\bigoplus_{k\neq j}\int_{D_k}
\mDelta_k(x_k)\od\frac{p^D_{jk}(x_j,x_k)}{p^D_j(x_j)}dx_k,~1\leq j\leq d
\end{align}
subject to the constraints
\begin{align}\label{Delta constraint}
\begin{split}
\int_{D_j}\mDelta_j(x_j)\od p^D_j(x_j)dx_j
&=\int_{D_j}\mdelta_j(x_j)\od(\alpha^2_j\cdot p^D_j(x_j))\,dx_j,~1\leq j\leq d,\\
\int_{D_j}\|\mDelta_j(x_j)\|^2p^D_j(x_j)dx_j&<\infty,~1\leq j\leq d.
\end{split}
\end{align}
We define $\mc_j(x_j)=\int_{B(\bzero_j,1)}\mathfrak{D}^2\mf_j(x_j)(t_j)(t_j)\od(K_j(\|t_j\|_j)/2)dt_j$,
where $\mathfrak{D}^2\mf_j(x_j)$ denotes the second order Fr\'{e}chet derivative of $\mf_j$ at $x_j$. $\mathfrak{D}^2\mf_j$ is a bounded linear operator from $D_j$ to $\mathcal{L}(\mathbb{R}^{L_j},\mathcal{L}(\mathbb{R}^{L_j},\mbH))$. It turns out that $\mTheta_j(x_j):=\alpha^2_j\od\mc_j(x_j)\op\mDelta_j(x_j)$ is the asymptotic bias of $\hat\mf_j(x_j)$. To describe the asymptotic variance, we define the outer product $\mepsilon\otimes\mepsilon:\mbH\rightarrow\mbH$ by $(\mepsilon\otimes\mepsilon)(\mh)=\lng\mepsilon,\mh\rng\od\mepsilon$. We also define the covariance operator $C_{j,x_j}:\mbH\rightarrow\mbH$ by
\begin{align*}
C_{j,x_j}(\mh)=\alpha^{-L_{\rm max}}_j\cdot (p^D_j(x_j)\cdot p_0^D)^{-1}\cdot\E((\mepsilon\otimes\mepsilon)(\mh)|\xi_j=x_j,\mxi\in D)\cdot
\int_{B_j(\bzero_j,1)}K_j^2(\|t_j\|_j)dt_j
\end{align*}
whenever $L_j=L_{\rm max}$ and $C_{j,x_j}(\mh)\equiv\bzero$ otherwise. We denote by $\mG(\mathbf{0},C_{j,x_j})$ the Hilbertian Gaussian random element with mean $\bzero$ and covariance operator $C_{j,x_j}$. It is a random element such that $\lng\mG(\mathbf{0},C_{j,x_j}),\mh\rng$ is the normal random variable with mean 0 and variance $\lng C_{j,x_j}(\mh),\mh\rng$.
When $\mbH=\mathbb{R}$, it reduces to a normal random variable. We are now ready to state the theorem. We emphasize that the derivation of asymptotic distribution is more involved than that of error rates since the former needs finding the exact forms of leading terms in asymptotic expansions. We denote by $\Leb$ the $\sum_{j=1}^dL_j$-dimensional Lebesgue measure.

\begin{theorem}\label{Asymptotic distribution}
Assume that conditions (B1)--(B6), (C2), (C4), \ref{D} and \ref{E} hold. Then, {\rm (i)} the solution of (\ref{Delta}) subject to (\ref{Delta constraint}) is unique in the sense that, if $(\mDelta^{\star}_j:1\leq j\leq d)$ is another solution, then $\mDelta_j=\mDelta_j^{\star}$ almost everywhere with respect to $\Leb_j$ for all $1\leq j\leq d$; {\rm (ii)} for almost everywhere $\mx\in\prod_{j=1}^d(D_j\setminus\partial D_j)$ with respect to $\Leb$, it holds that $(n^{2/(L_{\rm max}+4)}\odot(\hat{\mf}_j(x_j)\ominus\mf_j(x_j)):1\leq j\leq d)\overset{d}{\rightarrow}(\mathbf{\Theta}_j(x_j)\oplus\mG(\mathbf{0},C_{j,x_j}):1\leq j\leq d)$ jointly and that
$n^{2/(L_{\rm max}+4)}\odot(\hat{\mf}(\mx)\ominus\mf(\mx))\overset{d}{\rightarrow}\oplusjd\mathbf{\Theta}_j(x_j)\oplus\mG\big(\mathbf{0},\sum_{j=1}^dC_{j,x_j}\big)$,
where $\mathbf{\Theta}_1(x_1)\oplus\mG(\mathbf{0},C_{1,x_1}),\cdots,\mathbf{\Theta}_d(x_d)\oplus\mG(\mathbf{0},C_{d,x_d})$ are independent.
\end{theorem}

For $j$ with $L_j=L_{\rm max}$, let $\hat{\mf}^{\rm ora}_j$ be the Nadaraya-Watson-type `oracle' estimator of $\mf_j$ obtained by using the knowledge of all other component maps $(\mf_k: k\neq j)$.  Using Lemma \ref{A.I.} in the Supplementary Material \ref{collection of lemmas}, we may prove that $\hat{\mf}_j(x_j)$ and $\hat{\mf}^{\rm ora}_j(x_j)$ share the same asymptotic covariance operator and the difference $\mathbf{B}_j(x_j)$ of their asymptotic biases satisfies $\int_{D_j}\mathbf{B}_j(x_j)\odot p^D_j(x_j)dx_j=\mathbf{0}$. The latter equality tells that the bias difference is expected to be zero.


\begin{remark}
In the presence of a density-valued/functional predictor, we do not necessarily assume that the predictor is determined by a finite number of PC or SC scores, although some such variables are indeed finite-dimensional as illustrated in, e.g., \cite{Hall and Vial (2006)}. Instead, we aim to find the conditional expectation of $\mY$ given a set of selected component scores or target to predict the response values. Also, as in \cite{Han et al. (2018)} and \cite{Park et al. (2018)}, we do not include a diverging number of component scores in the model due to several technical obstacles. However, since the asymptotic properties of our SBF estimator do not depend on the size of $d$, one can take a large enough number of component scores in practices using, for example, the fraction of variance explained criterion. In fact, the usual number of component scores required to attain a certain level of prediction accuracy for high/infinite-dimensional data is not very large since the first few component scores tend to contain most of important information. 
\end{remark}

\section{Case of Riemannian functional responses}\label{extension}

\setcounter{equation}{0}
\setcounter{subsection}{0}


In this section, we study the case where the response variable is of Riemannian nature. Recently, \cite{Lin et al. (2022)} pioneered additive regression for responses taking values in a general Riemannian manifold. They studied univariate additive regression with perfectly observed scalar predictors. In this section, we aim to extend such additive regression to multivariate additive regression that covers Riemannian functional responses and various perfect/imperfect predictors.

Throughout this section, let $\mathcal{N}$ be a complete and connected Riemannian manifold and $d_{\mathcal{N}}$ be the Riemannian distance function on $\mathcal{N}$. Also, let $\mathcal{T}\subset\mathbb{R}$ be a compact set equipped with a finite Borel measure $\nu$. Let $W$ be a $\mathcal{N}$-valued random function defined on $\mathcal{T}$ such that the intrinsic mean function $\mu_{W}$ of $W$ exists. In case $\mathcal{T}$ is a singleton and $\nu$ is the counting measure, $W$ can be understood as a $\mathcal{N}$-valued random variable, which is the case considered in \cite{Lin et al. (2022)}. Let $\Log_{\mu_{W}(t)}$ be the Riemannian logarithm map at $\mu_{W}(t)$ for $t\in\mathcal{T}$.

Consider model (\ref{model}) with $\mY$ being the Hilbertian variable \mbox{${\rm{Log}}_{\mu_W}W\in\mathfrak{T}(\mu_W)$} defined by ${\rm{Log}}_{\mu_W}W:t\mapsto{\rm{Log}}_{\mu_W(t)}W(t)$ and satisfying $\E(\|{\rm{Log}}_{\mu_W}W\|^2_{\mathfrak{T}(\mu_W)})<\infty$. Note that applying an additive model directly to $W$ is not possible since there is no vector space structure on $\mathcal{N}$ in general. However, when $\mathcal{N}$ is a Lie group, the above additive model can be understood as a group additive model defined through the group structure of $\mathcal{N}$.
More details on this special case can be found in the Supplementary Material \ref{Lie case}. For i.i.d. copies $W^1,\ldots,W^n$ of $W$, we estimate $\mu_W(t)$ by $\hat{\mu}_W(t)=\argmin_{q\in\mathcal{N}}\sum_{i=1}^n(d_{\mathcal{N}}(q,W^i(t)))^2$.
For the estimation of the additive model, we use $\Log_{\hat{\mu}_{W}}W^i\in\mathfrak{T}(\hat{\mu}_{W})$ instead of $\Log_{\mu_{W}}W^i$ since $\mu_W$ is unknown. Hence, $\Log_{\hat{\mu}_{W}}W^i$ takes the role of $\tilde{\mY}^i$. Let $\Gamma_{\hat{\mu}_W,\mu_W}:\mathfrak{T}(\hat{\mu}_W)\rightarrow\mathfrak{T}(\mu_W)$ denote the parallel transport of vector fields as defined in Example \ref{Riemannian fd}. The following proposition quantifies the maximal discrepancy between $\Gamma_{\hat{\mu}_{W},\mu_{W}}(\Log_{\hat{\mu}_{W}}W^i)$ and $\Log_{\mu_{W}}W^i$, which takes the role of $b_n$.

\begin{proposition}\label{iRSC response}
Assume that condition \ref{L} holds for $\mathcal{M}^*=\mathcal{N}$ and $Z^*=W$. Also, assume that $\E(\|{\rm{Log}}_{\mu_W}W\|_{\mathfrak{T}(\mu_{W})}^{\tau})<\infty$ for some $\tau\geq2$. Then, $\max_{1\leq i\leq n}\|\Gamma_{\hat{\mu}_{W},\mu_{W}}(\Log_{\hat{\mu}_{W}}W^i)-\Log_{\mu_{W}}W^i\|_{\mathfrak{T}(\mu_{W})}$ achieves the rates given in Proposition \ref{HPC} under the same moment conditions with $\|\mX\|_*$ and $\tau\geq4$ respectively being replaced by $\|{\rm{Log}}_{\mu_W}W\|_{\mathfrak{T}(\mu_{W})}$ and $\tau\geq2$.
\end{proposition}

We define $\tilde{\mf}_0$, $\tilde{\mm}_j(x_j)$, $\tilde{\mf}_j(x_j)$, $\tilde\mf(\mx)$, $\tilde{\mf}^{[r]}_j(x_j)$ and $\tilde\mf^{[r]}(\mx)$ as in Section \ref{estimation method} with $\tilde{\mY}^i$ being replaced by $\Log_{\hat{\mu}_{W}}W^i$. Note that the estimators take values in $\mathfrak{T}(\hat{\mu}_{W})$. We also define 
$\hat{\mf}_j(x_j)=\Gamma_{\hat{\mu}_{W},\mu_{W}}(\tilde{\mf}_j(x_j))$, $\hat\mf(\mx)=\Gamma_{\hat{\mu}_{W},\mu_{W}}(\tilde\mf(\mx))$, $\hat{\mf}^{[r]}_j(x_j)=\Gamma_{\hat{\mu}_{W},\mu_{W}}(\tilde{\mf}^{[r]}_j(x_j))$ and $\hat\mf^{[r]}(\mx)=\Gamma_{\hat{\mu}_{W},\mu_{W}}(\tilde\mf^{[r]}(\mx))$. The latter quantities take values in $\mathfrak{T}(\mu_{W})$. The following theorem tells that the non-asymptotic properties in Section \ref{estimation method} hold for the new estimators $\tilde{\mf}_j(x_j)$, $\tilde\mf(\mx)$, $\tilde{\mf}^{[r]}_j(x_j)$ and $\tilde\mf^{[r]}(\mx)$, and the estimators are close to their targets when $\hat{\mu}_{W}$ is close to $\mu_{W}$.

\begin{theorem}\label{extension theorem}
Propositions \ref{existence}-\ref{convergence1 individual} remain valid for the new estimators $\tilde{\mf}_j(x_j)$, $\tilde\mf(\mx)$, $\tilde{\mf}^{[r]}_j(x_j)$ and $\tilde\mf^{[r]}(\mx)$. Also,
Theorems \ref{convergence2}-\ref{Asymptotic distribution} with $b_n$ being the rate of $\max_{1\leq i\leq n}\|\Gamma_{\hat{\mu}_{W},\mu_{W}}(\Log_{\hat{\mu}_{W}}W^i)-\Log_{\mu_{W}}W^i\|_{\mathfrak{T}(\mu_{W})}$ remain valid for the newly defined $\hat{\mf}_j(x_j)$, $\hat\mf(\mx)$, $\hat{\mf}^{[r]}_j(x_j)$ and $\hat\mf^{[r]}(\mx)$.
\end{theorem}

Note that $\tilde\mf(\mx)(t)\in T_{\hat{\mu}_W(t)}\mathcal{N}$ for each $t\in\mathcal{T}$, where $T_{\hat{\mu}_W(t)}\mathcal{N}$ is the tangent space of $\mathcal{N}$ at $\hat{\mu}_W(t)$. To predict the realization of the Riemannian functional variable $W$ given $\mxi=\mx$, we use $\Exp_{\hat{\mu}_W}\tilde\mf(\mx):\mathcal{T}\rightarrow\mathcal{N}$ defined by $\Exp_{\hat{\mu}_W}\tilde\mf(\mx):t\mapsto\Exp_{\hat{\mu}_W(t)}\tilde\mf(\mx)(t)$, where $\Exp_{\hat{\mu}_W(t)}:T_{\hat{\mu}_W(t)}\mathcal{N}\rightarrow\mathcal{N}$ is the Riemannian exponential map at $\hat{\mu}_W(t)$. The specific form of $\Exp_{\hat{\mu}_W(t)}$ for several $\mathcal{N}$ can be found in \citet{Dai and Muller (2018)} and \cite{Lin and Yao (2019)}. We may also consider SC scores between ${\rm{Log}}_{\mu_{W}}W$ and a Hilbertian functional predictor or between ${\rm{Log}}_{\mu_{W}}W$ and a Riemannian functional predictor. We call them Hilbertian intrinsic Riemannian singular component score and intrinsic Riemannian singular component score, respectively. We defer their formal definitions and the corresponding asymptotic properties to Examples \ref{Hilbertian Riemannian singular fd} and \ref{Riemannian singular fd} in the Supplementary Material \ref{SC scores Rimannian Y}.

\section{Numerical study}\label{numerical study}

\setcounter{equation}{0}
\setcounter{subsection}{0}

In this section, we present simulation studies and real data applications. Implementation details including the bandwidth selection can be found in the Supplementary Material \ref{suppl-imple}.

\subsection{Simulation study}\label{simulation study}

We conduct two simulation studies. Let $\mathbb{S}^k=\{u\in\mathbb{R}^{k+1}:\|u\|_{\mathbb{R}^{k+1}}=1\}$ for $k\in\mathbb{N}$.
The first simulation study considers the case where $\mY$ is a random density supported on the unit circle $\mathbb{S}^1$ and $\xi_j$ are perfectly observed. The second simulation study deals with the case where $\mY$ is a scalar random variable and $\xi_j$ are the iRFPC scores obtained from multiple $\mathbb{S}^2$-valued random functions. Such regression settings have not been covered in the literature. In particular, the second regression setting can be viewed as a spherical version of the standard functional regression with scalar responses and multiple functional predictors.

The first simulation study uses $d=2$ and $L_1=L_2=1$. We generate $\mxi=(\xi_1,\xi_2)\in\mathbb{R}^2$ from the bivariate normal distribution with $\E(\xi_1)=\E(\xi_2)=0.5$, $\Var(\xi_1)=\Var(\xi_2)=0.25^2$ and $\Cov(\xi_1,\xi_2)=0.25^3$. Let $g_1(\cdot)$, $g_2(\cdot)$ and $g_3(\cdot)$ be the densities of the von Mises distributions on $\mathbb{S}^1$ with concentration parameter 1 and respective means $-\pi/2$, $\pi/2$ and 0. Also, let $\delta$ be a $N(0,1)$ random variable independent of $\mxi$. We generate a random density $\mY=Y(\cdot)$ from the model
\begin{align}\label{simulation model 1}
Y(\cdot)=f_1(\xi_1)(\cdot)\cdot f_2(\xi_2)(\cdot)\cdot \epsilon(\cdot)\cdot\bigg(\int_{\mathbb{S}^1}f_1(\xi_1)(s)\cdot f_2(\xi_2)(s)\cdot \epsilon(s)d\nu(s)\bigg)^{-1},
\end{align}
where $f_1(x_1)(s)=(g_1(s))^{\cos(\pi x_1)}$, $f_2(x_2)(s)=(g_2(s))^{\sin(2\pi x_2)}$, $\epsilon(s)=(g_3(s))^{\delta}$ and $\nu$ is the Riemannian volume measure on $\mathbb{S}^1$. Note that the above model is an additive model in the geometry of the Bayes-Hilbert space $\mathfrak{B}^2(\mathbb{S}^1,\mathcal{B}(\mathbb{S}^1),\nu)$. We also generate a random sample $\{Y^*_{ik}:1\leq k\leq n^*_i\}$ from $Y^i(\cdot)$ and estimate $Y^i(\cdot)$ by the estimator $\tilde{Y}^i(\cdot)$ defined at (\ref{Pelletier density}). We consider $n=100,400$ and $n^*_i\equiv n^*=100,400$ to see the effects of $n$ and $n^*$. 
We then estimate the component maps $\mf_1$ and $\mf_2$, defined by $\mf_1(x_1)=f_1(x_1)(\cdot)$ and $\mf_2(x_2)=f_2(x_2)(\cdot)$, on the estimation domains $D_1=D_2=[0,1]$. Note that $\int_{D_1}\mf_1(x_1)\od p_1^D(x_1)dx_1=\int_{D_2}\mf_2(x_2)\od p_2^D(x_2)dx_2=\bzero$ since $\int_{D_1}\cos(\pi x_1)\cdot p_1^D(x_1)dx_1=\int_{D_2}\sin(2\pi x_2)\cdot p_2^D(x_2)dx_2=0$, where $\bzero$ is the uniform density on $\mathbb{S}^1$. For the estimation, we use the two sets of observations $\{(\mxi^i,Y^i(\cdot)):1\leq i\leq n\}$ and $\{(\mxi^i,\tilde{Y}^i(\cdot)):1\leq i\leq n\}$, and compare their estimation performance. The simulation is repeated $M=100$ times. As a measure of performance, we use the integrated squared bias (ISB), integrated variance (IV) and integrated mean squared error (IMSE) defined by
\begin{align}\label{measure of performace}
\begin{split}
{\rm IMSE}_j &= \int_{D_j} M^{-1} \underset{m=1}{\overset{M}{\sum}}\big\|\mf_j(x_j)
\ominus\hat{\mf}^{(m)}_j(x_j)\big\|^2 dx_j = {\rm ISB}_j
+{\rm IV}_j,\\
{\rm ISB}_j &=\int_{D_j} \bigg\|\mf_j(x_j) \ominus \bigg(M^{-1} \odot\bigoplus_{m=1}^M\hat{\mf}^{(m)}_j(x_j)\bigg)
\bigg\|^2 dx_j,\textrm{ and }\\
{\rm IV}_j &= \int_{D_j} M^{-1} \underset{m=1}{\overset{M}{\sum}}\bigg\|\hat{\mf}^{(m)}_j(x_j)
\ominus \bigg(M^{-1} \odot\bigoplus_{m=1}^M\hat{\mf}^{(m)}_j(x_j)\bigg)
\bigg\|^2 dx_j,
\end{split}
\end{align}
for $j\in\{1,2\}$, where $\hat{\mf}^{(m)}_j$ is the estimator of $\mf_j$ obtained from the $m$th Monte-Carlo sample.

Table \ref{S1} shows the IMSE, ISB and IV values. It shows that the estimation performance of our method gets better as $n$ increases. Also, the performance of our method based on $\{(\mxi^i,\tilde{Y}^i(\cdot)):1\leq i\leq n\}$ gets closer to the one based on $\{(\mxi^i,Y^i(\cdot)):1\leq i\leq n\}$ as $n^*$ increases. This indicates that we can achieve good regression performance for density-valued responses if we have an enough sample from each density even though we do not observe the densities directly. This confirms our theoretical developments in Example \ref{B.3}.

\begin{table*}[!t]
\begin{center}
\caption{\footnotesize Comparison of the IMSE, ISB and IV values
of the proposed estimator for model (\ref{simulation model 1}) with $M=100$ Monte-Carlo samples.}
\label{S1}
{\small
\begin{tabular}{ccccccccc}
\hline
&  & \multicolumn{3}{c}{$j=1$} & & \multicolumn{3}{c}{$j=2$}
\\ \cline{3-5} \cline{7-9}
&  & \multicolumn{2}{c}{$\tilde{Y}^i(\cdot)$} & $Y^i(\cdot)$ & & \multicolumn{2}{c}{$\tilde{Y}^i(\cdot)$} & $Y^i(\cdot)$
\\ \cline{3-4} \cline{7-8}
$n$ & Criterion & $n^*=100$ & $n^*=400$ & & & $n^*=100$ & $n^*=400$ &  \\
\hline
    & IMSE & 0.219 & 0.213 & 0.192 &  & 0.370 & 0.351 & 0.344 \\ \cline{2-9}
 100& ISB  & 0.052 & 0.038 & 0.031 &  & 0.087 & 0.071 & 0.061 \\ \cline{2-9}
    & IV   & 0.167 & 0.175 & 0.161 &  & 0.283 & 0.280 & 0.283 \\
\hline
    & IMSE & 0.087 & 0.076 & 0.071 &  & 0.123 & 0.116 & 0.108 \\ \cline{2-9}
 400& ISB  & 0.027 & 0.017 & 0.011 &  & 0.037 & 0.025 & 0.019 \\ \cline{2-9}
    & IV   & 0.060 & 0.059 & 0.060 &  & 0.086 & 0.091 & 0.089 \\
\hline
\end{tabular}
}
\end{center}
\end{table*}

For the second simulation study, we generate $\mathbb{S}^2$-valued random functions $Z_1(\cdot)$ and $Z_2(\cdot)$ defined on $[0,1]$. They are generated by $Z^i_k(t)=\Exp_{\mu_{Z_k}(t)}\sum_{r=1}^2(\eta^i_{kr}\cdot\psi_r(t))$ with intrinsic mean functions $\mu_{Z_k}(t)=(\cos(\theta_k(t))\sin(\phi_k(t)),\sin(\theta_k(t))\sin(\phi_k(t)),\cos(\phi_k(t)))\in\mathbb{S}^2$ for $k\in\{1,2\}$, where $\theta_1(t)=\pi/2\cdot t$, $\theta_2(t)=\pi/2\cdot t^2$, $\phi_1(t)=\pi\cdot t$, $\phi_2(t)=\pi\cdot t^2$ and $\Exp_{\mu_{Z_k}(t)}$ is the Riemannian exponential map defined in Section \ref{extension}. We take $n=100$ and 400. The iRFPC scores are generated by $\eta^i_{kr}=c_{kr}\cdot\varrho^i_{kr}$, where $c_{11}=1$, $c_{12}=0.75$, $c_{21}=0.75$, $c_{22}=0.5$ and $\varrho^i_{kr}$ are independently sampled from the truncated normal distribution $tN(0,0.5^2)$ truncated on $[-1,1]$. We choose the Fourier basis functions for $\psi_r$. We consider the \emph{spherical functional multivariate additive model}
\begin{align}\label{simulation model 2}
Y=f_1(\eta_{11})+f_2(\eta_{12})+f_3(\eta_{21},\eta_{22})+\epsilon,
\end{align}
where $f_1(x_1)=\sin(\pi x_1)/5$, $f_2(x_2)=-x_2^3$, $f_3(x_{31},x_{32})=x_{31}\cdot\tan(x_{32})$ and $\epsilon$ is the $N(0,0.1^2)$ random variable. We estimate $f_1$ on $D_1=[-1,1]$, $f_2$ on $D_2=[-0.75,0.75]$ and $f_3$ on $D_3=[-0.75,0.75]\times[-0.5,0.5]$ based on the estimated iRFPC scores and based on the true iRFPC scores under the multivariate additive model at (\ref{simulation model 2}). We also compare them with the univariate additive modelling that assumes $f_3$ as the sum of a function of $\eta_{21}$ and a function of $\eta_{22}$. The latter comparison is to see how the multivariate additive modelling is useful when the true model is not an univariate additive model. We used the R package `RFPCA' (\url{https://rdrr.io/github/CrossD/RFPCA/}) to generate $Z_k(\cdot)$ and estimate $\eta_{kr}$. As a measure of performance, we use (\ref{measure of performace}) with the norm and vector operations being replaced by those for $\mbH=\mathbb{R}$. We take $M=100$.

Table \ref{S2} shows the values of IMSE, ISB and IV. It shows that the estimation performance gets better as $n$ increases. The multivariate additive regression based on the estimated iRFPC scores works quite well compared to the one based on the true iRFPC scores. Also, the multivariate additive modelling outperforms the univariate counterpart. In particular, the biases for $f_3$ are much smaller in the multivariate modelling than those in the univariate modelling. This implies that our flexible model can be useful in practices.

\begin{table*}[!t]
\begin{center}
\caption{\footnotesize Comparison of the IMSE, ISB and IV values, multiplied by $10^2$,
of the proposed estimator for model (\ref{simulation model 2}) with $M=100$ Monte-Carlo samples.}
\label{S2}
{\small
\begin{tabular}{ccccccccccccc}
\hline
&  & \multicolumn{3}{c}{$j=1$} & & \multicolumn{3}{c}{$j=2$} & & \multicolumn{3}{c}{$j=3$}
\\ \cline{3-5} \cline{7-9} \cline{11-13}
& iRFPC & \multicolumn{2}{c}{Estimated} & True & & \multicolumn{2}{c}{Estimated} & True & & \multicolumn{2}{c}{Estimated} & True
\\ \cline{3-4} \cline{7-8} \cline{11-12}
$n$ & Criterion & Uni & Multi & Multi & & Uni & Multi & Multi & & Uni & Multi & Multi \\
\hline
    & IMSE & 0.472 & 0.370 & 0.301 &  & 0.522 & 0.440 & 0.357 &  & 2.765 & 0.684 & 0.625\\ \cline{2-13}
 100& ISB  & 0.108 & 0.094 & 0.067 &  & 0.244 & 0.208 & 0.160 &  & 2.454 & 0.418 & 0.371\\ \cline{2-13}
    & IV   & 0.364 & 0.276 & 0.234 &  & 0.278 & 0.232 & 0.197 &  & 0.311 & 0.266 & 0.254\\
\hline
    & IMSE & 0.166 & 0.128 & 0.103 &  & 0.152 & 0.121 & 0.096 &  & 2.559 & 0.286 & 0.268\\ \cline{2-13}
 400& ISB  & 0.037 & 0.030 & 0.026 &  & 0.061 & 0.051 & 0.037 &  & 2.459 & 0.169 & 0.150\\ \cline{2-13}
    & IV   & 0.129 & 0.098 & 0.077 &  & 0.091 & 0.070 & 0.059 &  & 0.100 & 0.117 & 0.118\\
\hline
\end{tabular}
}
\end{center}
\end{table*}

\subsection{Real data analysis}\label{real data analysis}

We present two real data applications. The first dataset contains a compositional response and mixed scalar predictors and compositional predictors. The second application is for a spherical functional response and mixed scalar predictors and Riemannian non-functional predictors.

\subsubsection{US presidential election}

It is believed that the population characteristics and underlying political orientation of a region are important factors determining election results in that region. To see how such factors affect election results in the United States, we analyze the 2020 US presidential election data. We collected the proportion of people who have a bachelor or a higher degree ($\xi_1$), the per capita income ($\xi_2$) and the median age ($\xi_3$) for each US state. For each state, we also obtained the compositional vector $\mX=(X_1,X_2,X_3)\in\mathcal{S}_1^3$, where $X_1$, $X_2$ and $X_3$ are the proportions of the white race, African American race and remaining races, respectively. We extracted the above population characteristics from \url{https://www.census.gov/acs/www/data/data-tables-and-tools/data-profiles/2019/}. As a measure of the underlying political orientation of each state, we obtained the compositional vector $\mZ=(Z_1,Z_2,Z_3,Z_4)\in\mathcal{S}_1^4$, where $Z_1$, $Z_2$, $Z_3$ and $Z_4$ are the proportions of votes earned by the Democratic party, Republican party, Libertarian party and remaining parties, respectively, in the 2016 US presidential election. We note that the first three parties are the most popular parties in US and only the three parties enrolled to all the states for the election. We took the three scalars $\xi_1$, $\xi_2$ and $\xi_3$ and the two compositional vectors $\mX$ and $\mZ$ as predictors. We took the compositional vector $\mY=(Y_1,Y_2,Y_3)\in\mathcal{S}_1^3$ obtained from the 2020 US presidential election as the response variable, where $Y_1$, $Y_2$ and $Y_3$ are the proportions of votes earned by the Democratic party, Republican party and remaining parties, respectively. We only took the three proportions for $\mY$ since it provides a clear visualization and the main interest of the election was on the first two parties. Since $Z_4=0$ in the Oklahoma state due to the absence of the corresponding parties, we analyzed the dataset $\{(\xi_1^i,\xi_2^i,\xi_3^i,\mX^i,\mZ^i,\mY^i):1\leq i\leq n\}$ with $n=51-1=50$.

To the best our knowledge, our method and the method of \cite{Jeon et al. (2021a)} are the only methods designed to cover regression with compositional responses and mixed multiple scalar predictors and multiple compositional predictors. Hence, we compare our method with the method of \cite{Jeon et al. (2021a)}. The latter treats the compositional predictors as multivariate predictors, so that the corresponding component maps are defined on simplices. For our method, we apply both PCA and singular component analysis to the compositional predictors for dimension reduction. Since the actual dimension of $\mX$ is 2 due to the constraint on $\mX$, we only consider the first component score of $\mX$ and take it as $\xi_4$. Similarly, we consider one or two component scores for $\mZ$, respectively as $\xi_5$ (univariate or bivariate component scores) or  $\xi_5$ and $\xi_6$ (two scalar predictors). The resulting 12 models for our method are illustrated in Table \ref{R2}. We compare all the models via the leave-one-out averaged squared prediction error (ASPE) defined by $n^{-1}\sum_{i=1}^n\|\mY^i\ominus \hat{\mY}^{i,(-i)}\|^2$, where $\hat{\mY}^{i,(-i)}$ is the prediction of $\mY^i$ based on the sample without the $i$th observation.

\begin{table*}[!t]
\begin{center}
\caption{\footnotesize ASPE values for US presidential election data}
\label{R2}
{\small
\begin{tabular}{c|cccccccccccc}
\hline
\# of scores ($\mZ$) & \multicolumn{4}{c}{1} & \multicolumn{8}{c}{2} \\
\hline
$\mf_j$ for $\mZ$ & \multicolumn{4}{c}{One univariate} & \multicolumn{4}{c}{Two univariate} & \multicolumn{4}{c}{One bivariate} \\
\hline
Score type ($\mZ$) & \multicolumn{2}{c}{PC} & \multicolumn{2}{c}{SC} & \multicolumn{2}{c}{PC} & \multicolumn{2}{c}{SC} & \multicolumn{2}{c}{PC} & \multicolumn{2}{c}{SC} \\
\hline
Score type ($\mX$) & PC & SC & PC & SC & PC & SC & PC & SC & PC & SC & PC & SC \\
\hline
ASPE & .160 & .154 & .151 & .148 & .137 & .134 & .096 & \textbf{.095} & .112 & .111 & .109 & .108 \\
\hline
\end{tabular}
}
\end{center}
\end{table*}

We obtained the ASPE value 0.128 for the method of \cite{Jeon et al. (2021a)}. Table \ref{R2} shows that most of our approaches based on 2 scores for $\mZ$ are superior to \cite{Jeon et al. (2021a)}. This indicates that an appropriate dimension reduction is indeed beneficial. Table \ref{R2} also shows that our methods based on the SC scores give better performance than those based on the PC scores. This is natural since the former component scores are constructed using both response and predictor variables, while the latter component scores exclude the response variable in their construction. Since the univariate additive model with two SC scores for $\mZ$ and one SC score for $\mX$ achieves the best performance, we fitted the dataset with this model to see how the predictors affect the response variable.

\begin{figure}[!t]
\center
\includegraphics[width=0.495\textwidth]{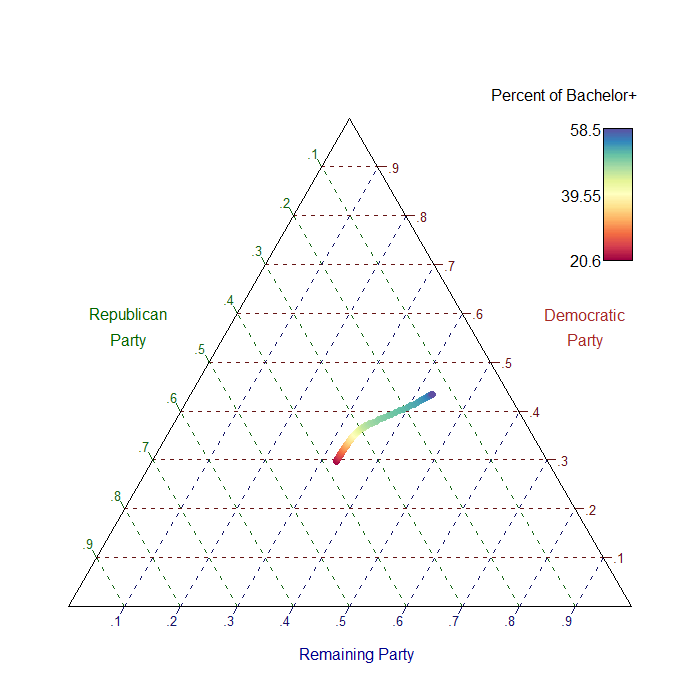}
\includegraphics[width=0.495\textwidth]{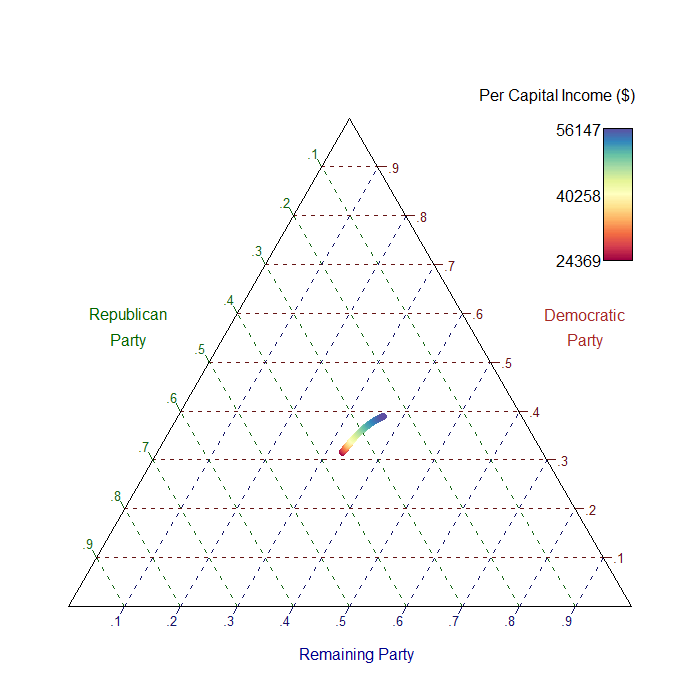}
\includegraphics[width=0.495\textwidth]{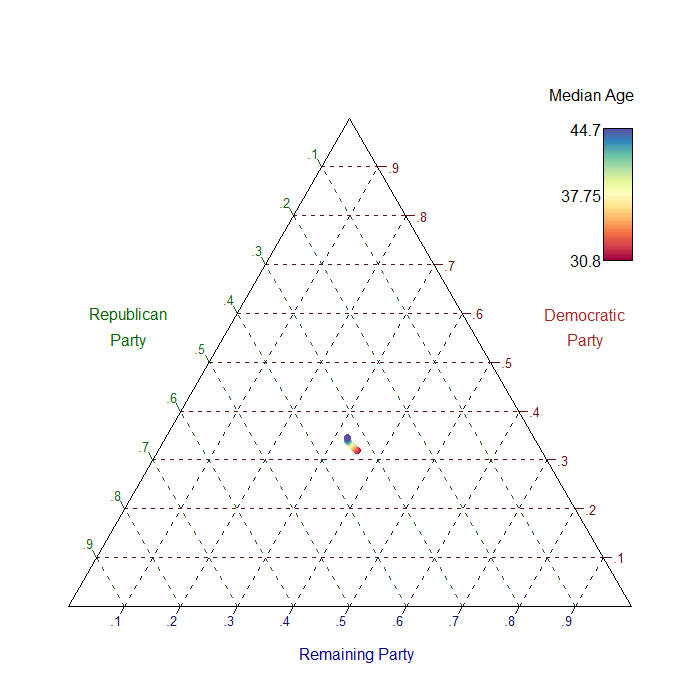}
\includegraphics[width=0.495\textwidth]{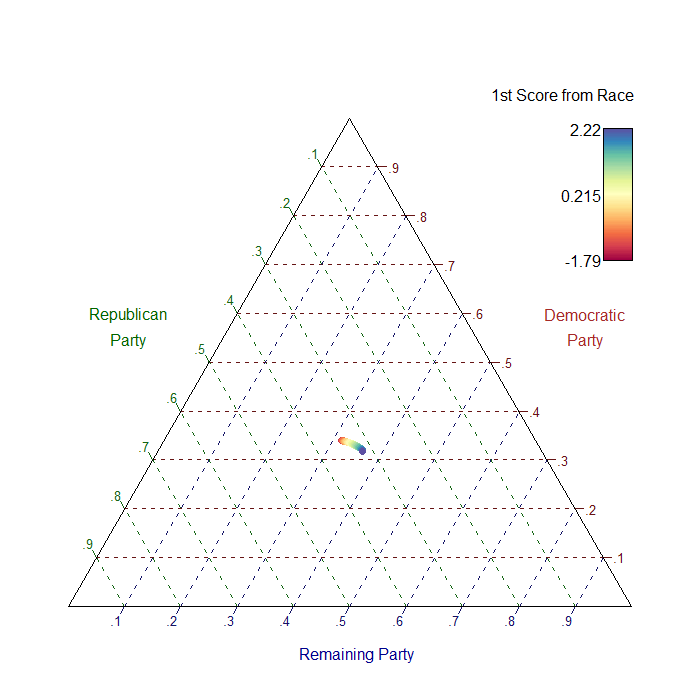}
\caption{\small Estimated component maps $\hat\mf_j$ for the population characteristics based on the proposed method applied to the US election data. For each point in the ternary plots, the proportions for the Democratic party, Republican party and remaining parties are respectively obtained by following the red, green and blue dashed lines.}
\label{fitted plots simplex 1}
\end{figure}

\begin{figure}[!t]
\center
\includegraphics[width=0.495\textwidth]{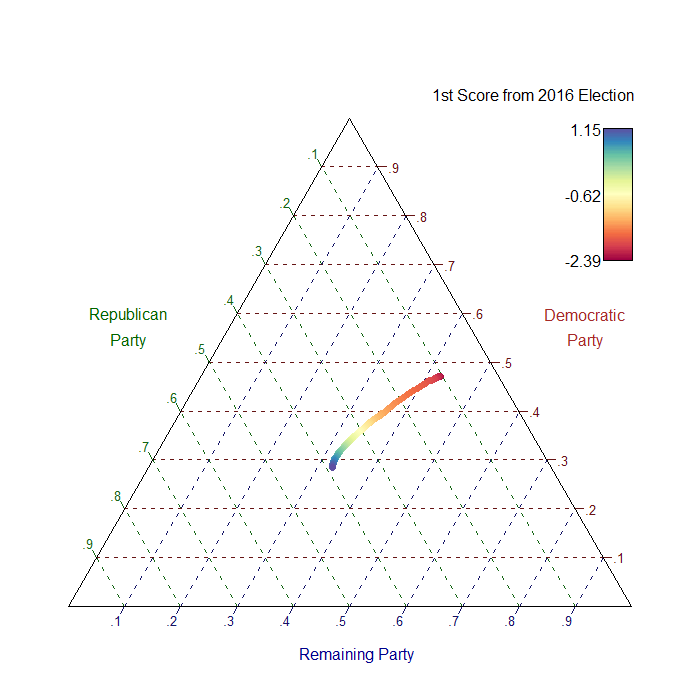}
\includegraphics[width=0.495\textwidth]{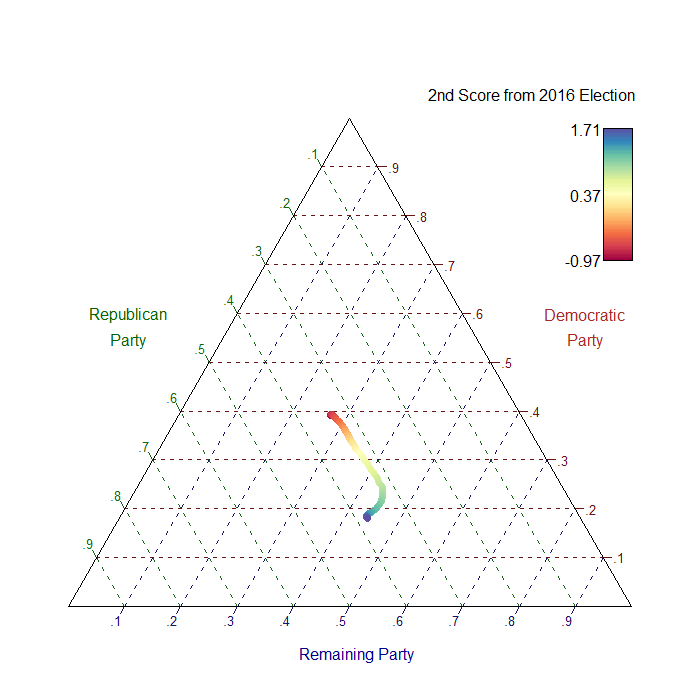}
\caption{\small The estimated component maps $\hat\mf_j$ for the underlying political orientation based on the proposed method applied to the US election data.}
\label{fitted plots simplex 2}
\end{figure}

Figures \ref{fitted plots simplex 1} and \ref{fitted plots simplex 2} depict the estimated component maps $(\hat{\mf}_j:1\leq j\leq 6)$. The top corners and bottom left corner of Figure \ref{fitted plots simplex 1} tell that the proportion of votes earned by the Democratic party increases as $\xi_1$, $\xi_2$ or $\xi_3$ increases. This indicates that more educated, richer or older people tend to prefer the Democratic party. The three plots also demonstrate that higher education level and higher income have negative effects on the votes for the Republican party and the age factor has no influence on the Republican party. Considering the lengths of the curves in the three plots, we may interpret that the education factor is more influential to the response than the other two factors. Regarding the race factor $\mX$, its estimated first singular vector $\hat{\mphi}_1=(\hat{\mphi}_{11},\hat{\mphi}_{12},\hat{\mphi}_{13})$ is approximately $(0.48,0.13,0.39)\in\mathcal{S}_1^3$. For two compositional vectors $\mv=(v_1,\ldots,v_k)$ and $\mw=(w_1,\ldots,w_k)$, their inner product in the geometry of $\mathcal{S}_1^k$ is given by $k^{-1}\sum_{l=1}^k\big((k-1)\log{w_l}-\sum_{l'\neq l}\log{w_{l'}}\big)\log{v_l}$.
Since $\log{v_l}<0$ with $\mv=\mX\om_*\bar{\mX}$ and the values $(k-1)\log{w_l}-\sum_{l'\neq l}\log{w_{l'}}$ with $\mw=\hat{\mphi}_1$ for $l\in\{1,2,3\}$ are respectively given by $2\log{\hat{\mphi}_{11}}-\log{\hat{\mphi}_{12}}-\log{\hat{\mphi}_{13}}\approx 1.5>0$, $2\log{\hat{\mphi}_{12}}-\log{\hat{\mphi}_{11}}-\log{\hat{\mphi}_{13}}\approx -2.4<0$ and $2\log{\hat{\mphi}_{13}}-\log{\hat{\mphi}_{11}}-\log{\hat{\mphi}_{12}}\approx 0.87>0$, the estimated SC score $\lng\mX\om_*\bar{\mX},\hat{\mphi}_1\rng$ increases as the proportion of the white race or the proportion of the remaining races increases, while it decreases as the proportion of the African American race increases. Having this in mind, the bottom right corner of Figure \ref{fitted plots simplex 1} tells that the increase of white or remaining races tends to increase the votes to the remaining parties but it tends to decrease the votes to the Democratic and Republican parties. Regarding the underlying political orientation factor $\mZ$, its estimated first singular vector $\hat{\mphi}_1$ and the second singular vector $\hat{\mphi}_2$ are approximately $(0.15,0.48,0.23,0.14)\in\mathcal{S}_1^4$ and $(0.10,0.22,0.27,0.41)\in\mathcal{S}_1^4$, respectively. Since the values $(k-1)\log{w_l}-\sum_{l'\neq l}\log{w_{l'}}$ with $\mw=\hat{\mphi}_1$ are approximately given by $-1.44$, $3.18$, $0.21$, $-1.95$ for $l\in\{1,2,3,4\}$, respectively, the decrease (increase) of the votes earned by the Democratic (Republican) party in the 2016 election increases the first SC score. Having this in mind, we can interpret the left plot of Figure \ref{fitted plots simplex 2} as that the political preference of each state is not well changed. The right plot of Figure \ref{fitted plots simplex 2} gives the same interpretation since the values $(k-1)\log{w_l}-\sum_{l'\neq l}\log{w_{l'}}$ with $\mw=\hat{\mphi}_2$ are approximately given by $-3.12$, $0.01$, $0.72$, $2.39$ for $l\in\{1,2,3,4\}$, respectively. Considering the lengths of the curves in the ternary plots for $\mX$ and $\mZ$, we may interpret that the political orientation factor is more influential to the response than the race factor.

\subsection{Typhoons in east Asia}

Typhoons equipped with strong wind and large amounts of rainfall damage many countries every year. Hence, exact prediction of a track of typhoon in an initial stage is very important to prevent such damage. The database of the Korea meteorological administration at \url{https://data.kma.go.kr/data/typhoonData/typInfoTYList.do?pgmNo=689} contains information of typhoons that appeared in east Asia from the year 2001. We collected information on each typhoon 12 hours after its emergence since the information from that time, say $T_0$, reflects the characteristics of typhoons better than the earlier times. We obtained the pre-smoothed 3-day trajectory ($W$) of each typhoon from time $T_0$ using the R package `tFrechet' (\url{https://rdrr.io/github/functionaldata/tFrechet/}). We also obtained the central air pressure ($\xi_1$), central maximum wind speed ($\xi_2$), strong wind range ($\xi_3$), moving speed ($\xi_4$), moving direction and location of each typhoon at time $T_0$. We created a `date' variable consisting of the month and day of $T_0$. Note that the moving direction and date are circular variables and the location is a spherical variable. Hence, we consider the iRFPC score ($\xi_5$) of the moving direction, the iRFPC score ($\xi_6$) of the date and the iRFPC scores ($\xi_7,\xi_8$) of the location with singleton time domain $\mathcal{T}=\{T_0\}$. We apply the model in Section \ref{extension} with response $\Log_{\mu_W}W$ and predictors $(\xi_j:1\leq j\leq 8)$, where $\mu_W$ is the intrinsic mean function of $W$. For simplicity, we estimate the univariate additive model with $n=265$ typhoons observed from 2001 to 2021 after excluding typhoons having missing observations. We then predict the trajectories of typhoons appeared in 2022 by applying $\Exp_{\hat{\mu}_W}$, defined in Section \ref{extension}, to the predicted $\Log_{\hat{\mu}_W}W$.

\begin{figure}[!t]
\center
\includegraphics[width=0.75\textwidth]{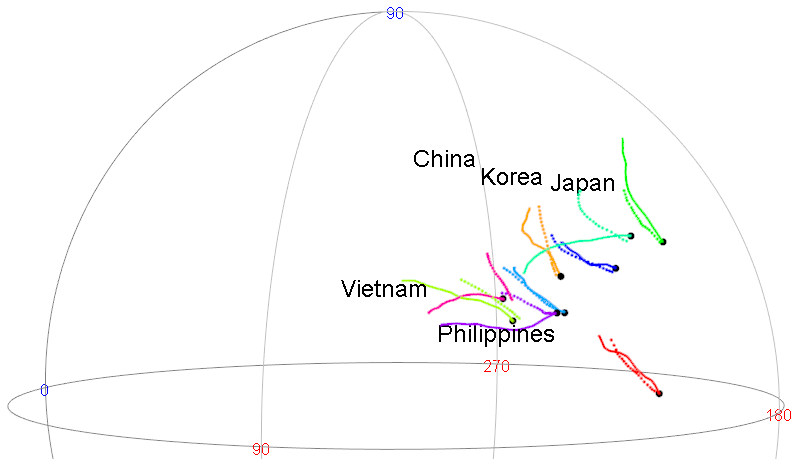}
\caption{\small The true trajectories of $W$ (solid line) and the predicted trajectories of $W$ (dotted line) for typhoons appeared in 2022. Black points are the starting points of the solid lines and the names of selected countries are placed on the locations of their capital cities.}
\label{predicted plot typhoon}
\end{figure}

Figure \ref{predicted plot typhoon} depicts the prediction result on the upper hemisphere. It shows reasonable performance although some predictions are less accurate due to the large variability of typhoons. We believe that our method combined with more significant predictors can be a useful tool for predicting the tracks of typhoons.

\section{Conclusion}\label{conclusion}

In this paper, we studied multivariate additive regression for Hilbertian responses where the variables are possibly subject to errors. We extended the scope of Hilbertian variables to a random density supported on a topological space and a random Hilbertian or Riemannian function defined on a topological space. We covered various scenarios of imperfect variables such as vanishing measurement errors on the variables, random samples from density-valued responses, semiparametric regression, Riemannian functional/non-functional responses and PCA for Hilbertian predictors and Riemannian functional/non-functional predictors with suitable theory for each scenario. We also introduced singular component analysis for various combinations of the response and predictor variables with several new notions and eigen-analysis. Our unified framework is not limited to such scenarios and has huge potential to many other instances. We also studied full non-asymptotic properties of our SBF estimator and established a perfect asymptotic theory demonstrating that our SBF estimator can avoid the curse of dimensionality in the general setting where the predictors have possibly non-compact supports. 
We believe that our work will serve as an useful tool to extend past and future dimension reduction techniques or imperfect data analysis to regression problems.

\newpage


\renewcommand{\theequation}{S.\arabic{equation}}
\renewcommand{\thesubsection}{S.\arabic{subsection}}
\setcounter{equation}{0}
\setcounter{subsection}{0}
\setcounter{page}{1}

\bigskip

\centerline{\bf \large Supplementary Material to}
\centerline{\bf \large `Additive Regression with General Imperfect Variables'}
\centerline{\bf \large by Jeong Min Jeon and Germain Van Bever}

\bigskip

In the Supplementary Material, we introduce the vector operations and inner products that define the separable Hilbert spaces in Section \ref{examples of hilbert} and the definition of Bochner integral. We also give the definitions of various notions on Riemannian manifolds and condition \ref{L} on Riemannian manifolds. We describe the group additive model introduced in Section \ref{extension} and implementation details used in Section \ref{numerical study}. The Supplementary Material also contains all technical proofs with a number of general lemmas. It also presents the definitions of Hilbertian intrinsic Riemannian singular component scores and intrinsic Riemannian singular component scores with their asymptotic properties. From now on, $\const$ denotes a generic positive constant.

\subsection{Examples of vector operations, inner product and norm}\label{examples}

\medskip
(i) {\it Simplex}. Recall the definition of $\mathcal{S}_1^k$ in Example \ref{pc simplex}. For this space, $\bzero=(1/k,\ldots,1/k)$.
For $c \in \mbR$ and $\mv,\mw\in\mathcal{S}_1^k$,
the vector addition $\mv \oplus \mw$ and scalar multiplication $c \odot \mv$ are defined by
\begin{eqnarray*}\begin{split}
\mv \oplus \mw &= \left(\frac{v_1\cdot w_1}{v_1\cdot w_1+\cdots+v_k \cdot w_k},\ldots,\frac{v_k \cdot w_k}{v_1 \cdot w_1+\cdots+v_k \cdot w_k}\right),\\
c \odot \mv &= \left(\frac{v_1^c}{v_1^c+\cdots+v_k^c},\ldots,\frac{v_k^c}{v_1^c+\cdots+v_k^c}\right).
\end{split}
\end{eqnarray*}
The inner product and norm are defined by
\[
\lng \mv,\mw \rng=\frac{1}{2k}\sum_{j=1}^k \sum_{l=1}^k \log(v_j/v_l)\log(w_j/w_l), \quad \|\mv\|=\left(\frac{1}{2k}\sum_{j=1}^k \sum_{l=1}^k \big(\log(v_j/v_l)\big)^2\right)^{1/2}.
\]

\medskip
\noindent(ii) {\it Bayes-Hilbert space}. Recall the definition of $\mathfrak{B}^2(\mathcal{S},\mathcal{B}(\mathcal{S}),\nu)$ in Example \ref{Density valued}.
For this space, $\bzero$ is the density function $f$ defined by
$f(s)=\nu(\mathcal{S})^{-1}$ for all $s\in\mathcal{S}$.
For $c \in \mbR$ and $f, g\in \mathfrak{B}^2(\mathcal{S},\mathcal{B}(\mathcal{S}),\nu)$, the vector addition $f \oplus g$ and scalar multiplication $c \odot f$ are defined by
\[
f \oplus g = \frac{f\cdot g}{\int_{\mathcal{S}} f(s)\cdot g(s)d\nu(s)}, \quad
c \odot f = \frac{f^c}{\int_{\mathcal{S}} (f(s))^cd\nu(s)}.
\]
The inner product and norm are defined by
\begin{eqnarray*}\begin{split}
\lng f, g\rng &=\frac{1}{2\nu(\mathcal{S})}\int_{\mathcal{S}}\int_{\mathcal{S}} \log\left(\frac{f(s)}{f(s')}\right)\log\left(\frac{g(s)}
{g(s')}\right)d\nu(s)d\nu(s'),\\
\|f\|&=\left(\frac{1}{2\nu(\mathcal{S})}\int_{\mathcal{S}}\int_{\mathcal{S}} \left(\log\left(\frac{f(s)}{f(s')}\right)\right)^2d\nu(s) d\nu(s')\right)^{1/2}.
\end{split}
\end{eqnarray*}

\medskip
\noindent(iii) {\it Lebesgue-Bochner-Hilbert space}. Recall the definition of $L^2((\mathcal{S},\mathcal{B}(\mathcal{S}),\nu),\mathcal{H})$ in Example \ref{Fully observed fd}.
For this space, $\bzero$ is the function $\mf$ such that $\mf(s)=\bzero_{\mcH}$ for all $s\in\mathcal{S}$, where $\bzero_{\mcH}$ is a zero vector of $\mcH$. For $c\in \mbR$ and $\mf, \mg\in L^2((\mathcal{S},\mathcal{B}(\mathcal{S}),\nu),\mathcal{H})$, the vector addition $\mf \oplus \mg$ is the function defined by $(\mf \oplus \mg)(s)=\mf(s)\op_{\mcH}\mg(s)$ for all $s\in\mathcal{S}$ and the scalar multiplication $c \odot \mf$ is the function defined by $(c \odot \mf)(s)=c\od_{\mcH}\mf(s)$ for all $s\in\mathcal{S}$, where $\op_{\mcH}$ and $\od_{\mcH}$ are respectively a vector addition and a scalar multiplication on $\mcH$. For an inner product $\lng\cdot,\cdot\rng_{\mcH}$ on $\mcH$, an inner product and a norm on $L^2((\mathcal{S},\mathcal{B}(\mathcal{S}),\nu),\mathcal{H})$ are defined by
\[
\lng \mf,\mg\rng=\int_{\mathcal{S}}\lng \mf(s),\mg(s)\rng_{\mathcal{H}}d\nu(s), \quad \|\mf\|=\left(\int_{\mathcal{S}}\|\mf(s)\|_{\mcH}^2d\nu(s)\right)^{1/2}.
\]

\subsection{Definition of Bochner integral}\label{bochner definition}

We let $(\mathcal{S},\Sigma,\nu)$ be a measure space and $\mathbb{B}$ be a separable Banach space equipped with a vector addition $\op_{\mbB}$, a scalar multiplication $\od_{\mbB}$ and a norm $\|\cdot\|_{\mbB}$. Let $\mf$ be a $\mbB$-valued $\nu$-simple map defined on $(\mathcal{S},\Sigma,\nu)$,
that is,
\[
\mf(s)=\sideset{}{_\mbB}\bigoplus_{i=1}^N (\I(s\in A_i)\odot_{\mbB}\mb_i)
\]
for some $\mb_i\in\mbB$ and disjoint $A_i\in\Sigma$ with $\nu(A_i)<\infty$. Its Bochner integral is defined by
\[
\int_{\mathcal{S}} \mf(s)d\nu(s) = \sideset{}{_\mbB}\bigoplus_{i=1}^N (\nu(A_i)\odot_{\mbB}\mb_i).
\]
Now, let $\mf^*:\mathcal{S}\rightarrow\mbB$ be a measurable map such that $\|\mf^*\|_{\mbB}$ is Lebesgue integrable with respect to $\nu$. Then, there exists a sequence of $\nu$-simple maps
$\mf_\ell$ such that $\mf^*(s)=\underset{\ell\rightarrow\infty}{\lim}\mf_\ell(s)$ and $\|\mf_\ell(s)\|_{\mbB}\leq\|\mf^*(s)\|_{\mbB}$
for all $s$ and $\ell$. For such sequence of $\nu$-simple maps, the Bochner integral of $\mf^*$ is defined by
\ba
\int_{\mathcal{S}}\mf^*(s)d\nu(s) = \underset{\ell\rightarrow\infty}{\lim}\int_{\mathcal{S}}\mf_\ell(s)d\nu(s).
\ea
The definition is independent of the choice of a sequence of $\nu$-simple maps.

\subsection{Notions on manifolds}\label{Riemannian definition}

In this section, we give the definitions of several notions on manifolds. More details can be found in standard textbooks on manifolds; see, e.g., Tu (2017) and Lee (2018). We first introduce the notions used in Example \ref{Riemannian fd}. Throughout this section, we let $\mathcal{M}$ and $c:I\rightarrow\mathcal{M}$ denote a smooth manifold and a smooth curve on an interval $I\subset\mathbb{R}$ unless otherwise specified. We also let $C^\infty(\mathcal{M})$ denote the set of smooth functions from $\mathcal{M}$ to $\mathbb{R}$. A tangent vector at $p\in\mathcal{M}$ is a linear function $v:C^\infty(\mathcal{M})\rightarrow\mathbb{R}$ satisfying $v(f\cdot g)=v(f)\cdot g(p)+f(p)\cdot v(g)$ for all $f,g\in C^\infty(\mathcal{M})$. The set $T_p\mathcal{M}$ of all tangent vectors at $p$ is called the tangent space of $\mathcal{M}$ at $p$.
The set $\bigcup_{p\in\mathcal{M}}T_p\mathcal{M}$ is denoted by $T\mathcal{M}$ and called tangent bundle of $\mathcal{M}$. A smooth map $X:\mathcal{M}\rightarrow T\mathcal{M}$ is called a smooth vector field if $X(p)\in T_p\mathcal{M}$ for all $p\in\mathcal{M}$. Denote by $\mathfrak{X}(\mathcal{M})$ the set of all smooth vector fields on $\mathcal{M}$. A connection in $T\mathcal{M}$ is a map $\nabla:\mathfrak{X}(\mathcal{M})\times\mathfrak{X}(\mathcal{M})\rightarrow\mathfrak{X}(\mathcal{M})$ such that $\nabla_XY:=\nabla(X,Y)$ satisfies (i) $\nabla_{f_1\cdot X_1+f_2\cdot X_2}Y=f_1\cdot\nabla_{X_1}Y+f_2\cdot\nabla_{X_2}Y$ for all $f_1,f_2\in C^{\infty}(\mathcal{M})$ and $X_1,X_2,Y\in\mathfrak{X}(\mathcal{M})$; (ii) $\nabla_X(a_1\cdot Y_1+a_2\cdot Y_2)=a_1\cdot\nabla_XY_1+a_2\cdot\nabla_XY_2$ for all $a_1,a_2\in\mathbb{R}$ and $Y_1,Y_2\in\mathfrak{X}(\mathcal{M})$; (iii) $\nabla_X(f\cdot Y)=(Xf)\cdot Y+f\cdot\nabla_XY$ for all $f\in C^{\infty}(\mathcal{M})$ and $X,Y\in\mathfrak{X}(\mathcal{M})$. Here, for $f\in C^{\infty}(\mathcal{M})$ and $X\in \mathfrak{X}(\mathcal{M})$, $f\cdot X\in\mathfrak{X}(\mathcal{M})$ is defined by $(f\cdot X)(p)=f(p)\cdot X(p)$ and $Xf\in C^{\infty}(\mathcal{M})$ is defined by
\begin{align}\label{Xf map}
Xf(p)=X(p)(f).
\end{align}
A smooth vector field along $c$ is a smooth map $V:I\rightarrow T\mathcal{M}$ such that $V(t)\in T_{c(t)}\mathcal{M}$ for all $t\in I$. Denote by $\mathfrak{X}(c)$ the set of all smooth vector fields along $c$. Define the velocity $c'(t)\in T_{c(t)}\mathcal{M}$ of $c$ at $t\in I$ by $c'(t)(f)=(f\circ c)'(t)$.
For each connection $\nabla$ in $T\mathcal{M}$, there exists a unique map $D_{\nabla}:\mathfrak{X}(c)\rightarrow\mathfrak{X}(c)$, called the covariant derivative along $c$, satisfying (i) $D_{\nabla}(a_1\cdot V_1+a_2\cdot V_2)=a_1\cdot D_{\nabla}(V_1)+a_2\cdot D_{\nabla}(V_2)$ for all $a_1,a_2\in\mathbb{R}$ and $V_1,V_2\in\mathfrak{X}(c)$; (ii) $D_{\nabla}(f\cdot V)=f'\cdot V+f\cdot D_{\nabla}(V)$ for all $f\in C^{\infty}(I)$ and $V\in\mathfrak{X}(c)$; (iii) if $V=Y\circ c$ for some $Y\in\mathfrak{X}(\mathcal{M})$, then $D_{\nabla}(V)(t)=\nabla_XY(c(t))$ for some $X\in\mathfrak{X}(\mathcal{M})$ with $X(c(t))=c'(t)$; see e.g., Theorem 4.24 in Lee (2018). In condition (iii) above, $\nabla_XY(c(t))$ does not depend on $X$ itself but only depends on its value at $c(t)$ (e.g., Proposition 4.5 in Lee (2018)) and thus $\nabla_XY(c(t))$ is often written as $\nabla_{c'(t)}Y(c(t))$. A smooth vector field $V$ along $c$ is said to be parallel along $c$ if $D_\nabla(V)(t)=0_{c(t)}$ for all $t\in I$, where $0_{c(t)}$ is the zero vector of $T_{c(t)}\mathcal{M}$. For each $t\in I$ and $v\in T_{c(t)}\mathcal{M}$, there exists a unique smooth vector field $V_v$ parallel along $c$ such that $V_v(t)=v$ (e.g., Theorem 4.32 in Lee (2018)). Such $V_v$ is called the parallel transport of $v$ along $c$. For each $t_0,t_1\in I$, define a map $\mathcal{P}_{t_0t_1}^c:T_{c(t_0)}\mathcal{M}\rightarrow T_{c(t_1)}\mathcal{M}$ by $\mathcal{P}_{t_0t_1}^c(v)=V_v(t_1)$. Such map $\mathcal{P}_{t_0t_1}^c$ is called the parallel transport map along $c$.

For two measurable functions $g_1,g_2:\mathcal{T}\rightarrow\mathcal{M}$ on a compact set $\mathcal{T}\subset\mathbb{R}$, let $\{\gamma(t,\cdot):t\in\mathcal{T}\}$ be a family of smooth curves $\gamma(t,\cdot):[0,1]\rightarrow\mathcal{M}$ such that $\gamma(\cdot,u):\mathcal{T}\rightarrow\mathcal{M}$ is measurable for each $u\in[0,1]$ and that $\gamma(t,0)=g_1(t)$ and $\gamma(t,1)=g_2(t)$ for each $t\in\mathcal{T}$. Then, for a given connection $\nabla$ in $T\mathcal{M}$, there exists the parallel transport map $\mathcal{P}_{01}^{\gamma(t,\cdot)}:T_{g_1(t)}\mathcal{M}\rightarrow T_{g_2(t)}\mathcal{M}$ along $\gamma(t,\cdot)$ defined by $\mathcal{P}_{01}^{\gamma(t,\cdot)}(v)=V_v(1)$, where $V_v$ is the parallel transport of $v$ along $\gamma(t,\cdot)$. Define the parallel transport of vector fields $\Gamma_{g_1,g_2}:\mathfrak{T}(g_1)\rightarrow\mathfrak{T}(g_2)$ by $\Gamma_{g_1,g_2}(V)(t)=\mathcal{P}_{01}^{\gamma(t,\cdot)}(V(t))$. The latter map was introduced by Lin and Yao (2019). 

A smooth curve $c:I\rightarrow\mathcal{M}$ is said to be regular if $c'(t)\neq 0_{c(t)}$ for all $t\in I$. A continuous curve $c:[a,b]\rightarrow\mathcal{M}$ is said to be admissible if there exists a partition $a=t_0<\cdots<t_k=b$ such that $c$ restricted to each $[t_{j-1},t_j]$ is regular. Such partition is called an admissible partition for $c$. A Riemannian metric is a map $g$ on $\mathcal{M}$ that assigns each $p\in\mathcal{M}$ to an inner product $\lng\cdot,\cdot\rng_p$ on $T_p\mathcal{M}$ such that $p\mapsto \lng X(p),Y(p)\rng_p$ is smooth on $\mathcal{M}$ for any $X,Y\in\mathfrak{X}(\mathcal{M})$. A Riemannian manifold is a smooth manifold equipped with a Riemannian metric. If $(\mathcal{M},g)$ is a Riemannian manifold and $c:[a,b]\rightarrow\mathcal{M}$ is an admissible curve with an admissible partition $a=t_0<\cdots<t_k=b$ for $c$, then the length $L_g$ of $c$ is defined by $L_g(c)=\sum_{j=1}^k\int_{t_{j-1}}^{t_j}\lng c'(t),c'(t)\rng_{c(t)}dt$. An admissible curve $c:[a,b]\rightarrow\mathcal{M}$ is said to be minimizing if $L_g(c)\leq L_g(\tilde{c})$ for all admissible curves $\tilde{c}:[\tilde{a},\tilde{b}]\rightarrow\mathcal{M}$ with $c(a)=\tilde{c}(\tilde{a})$ and $c(b)=\tilde{c}(\tilde{b})$. If a Riemannian manifold $(\mathcal{M},g)$ is connected, then any two points of $\mathcal{M}$ can be joined by an admissible curve (e.g., Proposition 2.50 in Lee (2018)) and the Riemannian distance $d_{\mathcal{M}}$ is defined by $d_{\mathcal{M}}(p,q)=\inf\{L_g(c):c\text{~is an admissible curve joining~}p\text{~and~}q\}$. 

For a smooth manifold $\mathcal{M}$ and a given connection $\nabla$ in $T\mathcal{M}$, a smooth curve $c:I\rightarrow\mathcal{M}$ is called a geodesic if $D_\nabla(c')(t)=0_{c(t)}$ for all $t\in I$. A geodesic $c:I\rightarrow\mathcal{M}$ is said to be maximal if it cannot be extended to a geodesic on a larger interval. For each $p\in\mathcal{M}$ and $v\in T_p\mathcal{M}$, there exists a unique maximal geodesic $c_{p,v}$ with $c(0)=p$ and $c'(0)=v$ defined on some open interval $I_{p,v}$ containing $0$ (e.g., Corollary 4.28 in Lee (2018)). If $I_{p,v}$ contains $[0,1]$, define $\Exp_p(v)=c_{p,v}(1)$. A Riemannian manifold is said to be geodesically complete if every maximal geodesic is defined on $\mathbb{R}$. For a geodesically complete Riemannian manifold $(\mathcal{M},g)$, the (Riemannian) exponential map $\Exp_p$ at $p$ is defined on $T_p\mathcal{M}$. A connected Riemannian manifold is geodesically complete if and only if it is complete as a metric space by the Hopf-Rinow theorem. Suppose that $(\mathcal{M},g)$ is a complete and connected Riemannian manifold. For each $p\in\mathcal{M}$ and $v\in T_p\mathcal{M}$, the cut time $t_{\rm cut}(p,v)\in(0,\infty]$ is defined by $t_{\rm cut}(p,v)=\sup\{t>0:c_{p,v}\text{~restricted to~}[0,t]\text{~is minimizing}\}$. If $t_{\rm cut}(p,v)<\infty$, then the cut point of $p$ along $c_{p,v}$ is defined by $c_{p,v}(t_{\rm cut}(p,v))\in\mathcal{M}$. The cut locus of $p$ denoted by ${{\rm Cut}}(p)$ is the set of all $q\in\mathcal{M}$ such that $q$ is the cut point of $p$ along some maximal geodesic. In this case, the inverse image $\Exp_p^{-1}(q)$ of $q\in\mathcal{M}\setminus{{\rm Cut}}(p)$ under the exponential map $\Exp_p$ is uniquely determined by $\Exp_p^{-1}(q)=c'(1)$, where $c:[0,1]\rightarrow\mathcal{M}$ is the unique minimizing geodesic joining $p$ and $q$ (e.g., Kendall and Le (2011)). Write the inverse map $\Exp_p^{-1}:\mathcal{M}\setminus{{\rm Cut}}(p)\rightarrow T_p\mathcal{M}$ by $\Log_p$.

For each $X,Y\in\mathfrak{X}(\mathcal{M})$, define a map $[X,Y]$ that assigns each $f\in C^{\infty}(\mathcal{M})$ to another smooth function by $[X,Y](f)=X(Yf)-Y(Xf)$, where $X(Yf)$ and $Y(Xf)$ are smooth functions on $\mathcal{M}$ defined in the sense of (\ref{Xf map}). A connection $\nabla$ in $T\mathcal{M}$ is said to be symmetric if $(\nabla_XY)f-(\nabla_YX)f=[X,Y](f)$ for all $X,Y\in\mathfrak{X}(\mathcal{M})$ and $f\in C^{\infty}(\mathcal{M})$, where $(\nabla_XY)f$ and $(\nabla_YX)f$ are smooth functions on $\mathcal{M}$ defined in the sense of (\ref{Xf map}). For a Riemannian manifold $(\mathcal{M},g)$, a connection $\nabla$ in $T\mathcal{M}$ is said to be compatible with $g$ if every parallel transport map along a given smooth curve $c:I\rightarrow\mathcal{M}$ is a linear isometry. For a Riemannian manifold $(\mathcal{M},g)$, the Levi-Civita connection is the connection $\nabla$ in $T\mathcal{M}$ symmetric and compatible with $g$. 

We now introduce the notions used in Example \ref{B.3}. Let $(\mathcal{M},g)$ be a $m$-dimensional compact Riemannian manifold with $m\in\mathbb{N}$. A pair $(U,\varphi)$ of an open set
$U \subset \mcM$ and a homeomorphism $\varphi: U \to \varphi(U)\subset\mbR^m$
is called a chart of $\mcM$. Let $\{(U_\beta,\varphi_\beta)\}_\beta$ be a finite collection of charts
such that $\mcM=\bigcup_\beta U_\beta$. Such a finite collection exists since $\mcM$ is compact.
We write $\varphi_\beta=(\varphi_{\beta,1},\ldots,\varphi_{\beta, m})$, where $\varphi_{\beta,l}:U\rightarrow\mathbb{R}$.
Let $\{\rho_\beta\}_\beta$ be a partition of unity subordinate to $\{U_\beta\}_\beta$.
Also, let $\frac{\partial}{\partial \varphi_{\beta,l}}\big|_{\varphi_\beta^{-1}(t)}$ denote the function defined by
$\frac{\partial}{\partial \varphi_{\beta,l}}\big|_{\varphi_\beta^{-1}(t)}(h)=\mathfrak{D}^1_l(h\circ\varphi^{-1}_\beta)(t)$ for each smooth function
$h:U_\beta\rightarrow\mathbb{R}$, where $\mathfrak{D}^1_l$ for $1 \le l \le m$ are the usual partial differential operators acting on real-valued functions on
$\mbR^m$.
Define $G_\beta\big(\varphi_\beta^{-1}(t)\big)$
to be the $m\times m$ matrix whose $(l,l')$th entry is given by
\begin{align*}
\bigg<\frac{\partial}{\partial \varphi_{\beta,l}}\bigg|_{\varphi_\beta^{-1}(t)},\;
\frac{\partial}{\partial \varphi_{\beta,l'}}\bigg|_{\varphi_\beta^{-1}(t)}\bigg>_{\varphi_\beta^{-1}(t)}.
\end{align*}
Let $C(\mathcal{M},\mbR)$ denote the class of real-valued continuous functions defined on $\mcM$.
Define $F:C(\mathcal{M},\mbR)\rightarrow\mbR$ by
\ba
F(f)=\sum_{\beta}\int_{\varphi_\beta(U_\beta)} f\big(\varphi_\beta^{-1}(t)\big) \cdot \rho_\beta\big(\varphi_\beta^{-1}(t)\big)
\cdot\sqrt{{\rm det}(G_\beta(\varphi_\beta^{-1}(t)))}dt.
\ea
The value of $F(f)$ is independent of the choices of $\{(U_\beta,\varphi_\beta)\}_\beta$ and $\{\rho_\beta\}_\beta$.
Due to the Riesz representation theorem for measures,
there exists a unique regular Borel measure $\nu$ on $\mathcal{M}$ such that
\ba
F(f)=\int_{\mathcal{M}}f(p)d\nu(p), \quad f\in C(\mathcal{M},\mathbb{R}).
\ea
We call such a $\nu$ the Riemannian volume measure induced by $g$.

Let $(\mcM,g)$ be a $m$-dimensional Riemannian manifold. For $p\in\mcM$, let $0_p$ denote the zero vector in $T_p\mcM$ and $B_{T_p\mcM}(0_p,r)$ denote the open ball in $T_p\mcM$ centered at $0_p$ with radius $r>0$. A diffeomorphism between two smooth manifolds is a smooth bijective map whose inverse is also smooth. Then, the injectivity radius of $\mathcal{M}$ at $p$ is defined by
\begin{align*}
\inj_{p}(\mathcal{M})=\sup\{r>0:\Exp_{p}|_{B_{T_{p}(\mathcal{M})}(0_p,r)}~\text{is a diffeomorphism onto its image}\},
\end{align*}
where $\Exp_{p}|_{B_{T_{p}(\mathcal{M})}(0_p,r)}$ is the exponential map at $p$ restricted to $B_{T_{p}(\mathcal{M})}(0_p,r)$.
For an orthonormal basis $\{e_{p,1},\ldots,e_{p,m}\}$ of $T_p\mathcal{M}$, define an isometric isomorphism $\iota_p:\mathbb{R}^m\rightarrow T_p\mathcal{M}$ by $\iota_p(t_1,\ldots,t_m)=\sum_{k=1}^mt_k\,e_{p,k}$. For $r<\inj_p(\mathcal{M})$, define $\Exp^*_p=\Exp_p\circ \iota_p:B(\bzero_{\mathbb{R}^m},r)\rightarrow \Exp_{p}(B_{T_{p}(\mathcal{M})}(0_p,r))$, where $B(\bzero_{\mathbb{R}^m},r)$ is the open ball centered at the zero vector $\bzero_{\mathbb{R}^m}$ of $\mathbb{R}^m$ with radius $r$. Then, the volume density function $\theta(\cdot;p):\Exp_{p}(B_{T_{p}(\mathcal{M})}(0_p,r))\rightarrow(0,\infty)$ of $\mathcal{M}$ at $p$ is defined by $\theta(q;p)=\sqrt{\det(G_{p}(q))}$, where $G_{p}(q)$ is the matrix whose $(i,j)$th entry is given by
\begin{align*}
\bigg<\frac{\partial}{\partial((\Exp^*_{p})^{-1}_i)}\bigg|_{q},\frac{\partial}{\partial((\Exp^*_{p})^{-1}_j)}\bigg|_{q}\bigg>_q.
\end{align*}
Here, $(\Exp^*_{p})^{-1}_i$ and $(\Exp^*_{p})^{-1}_j$ are $i$th and $j$th coordinates of $(\Exp^*_{p})^{-1}$, respectively, and $\frac{\partial}{\partial((\Exp^*_{p})^{-1}_k)}\big|_{q}\in T_q\mathcal{M}$ for $1\leq k\leq m$ is defined by
\ba
\frac{\partial}{\partial((\Exp^*_{p})^{-1}_k)}\bigg|_{q}(f)=\frac{\partial(f\circ\Exp^*_p)(t)}{\partial t_k}\bigg|_{t=(\Exp^*_{p})^{-1}(q)}
\ea
for $f\in C^{\infty}(\mathcal{M})$ and $t=(t_1,\ldots,t_m)$.

\subsection{Condition (L)}\label{Lin condition}

In this section, we let $(\mathcal{M}^*,g)$ be a complete and connected Riemannian manifold and $Z^*$ be a $\mathcal{M}^*$-valued random function defined on a compact set $\mathcal{T}\subset\mathbb{R}$ equipped with a finite Borel measure. We let $\mathfrak{X}(\mathcal{M}^*)$ denote the set of all smooth vector fields on $\mathcal{M}^*$ and $\nabla$ denote the Levi-Civita connection. We also let $\Log_{(\cdot)}q$ for $q\in\mathcal{M}^*$ denote the map that assigns each $p\in\mathcal{M}^*\setminus{\rm Cut}(q)$ to $\Log_pq\in T_p(\mathcal{M}^*)$ and each $p\in{\rm Cut}(q)$ to an element in the preimage $\Exp_p^{-1}(q)\in T_p(\mathcal{M}^*)$.
For each $q\in\mathcal{M}^*$, it holds that $\nu({{\rm Cut}}(q))=0$ for the Riemannian volume measure $\nu$ induced by $g$; see e.g., Theorem 10.34 in Lee (2018).

\begin{customcon}{(L)}\label{L}\leavevmode
\begin{itemize}
\item[{\rm(L1)}] The exponential map $\Exp_p:T_p\mathcal{M}^*\rightarrow\mathcal{M}^*$ is surjective for every $p\in\mathcal{M}^*$.
\item[{\rm(L2)}] $F(p,t):=\E((d_{\mathcal{M}^*}(p,Z^*(t)))^2)<\infty$ for any $(p,t)\in\mathcal{M}^*\times\mathcal{T}$ and $\sup_{(p,t)\in\mathcal{K}\times\mathcal{T}}F(p,t)<\infty$ for any compact set $\mathcal{K}\subset\mathcal{M}^*$.
\item[{\rm(L3)}] The intrinsic mean function $\mu_{Z^*}$ of $Z^*$ defined by $\mu_{Z^*}(t)=\argmin_{p\in\mathcal{M}^*}F(p,t)$ exists and the image of $\mu_{Z^*}$ is bounded.
\item[{\rm(L4)}] For $({{\rm Cut}}(\mu_{Z^*}(t)))^\varepsilon=\{p\in\mathcal{M}:\inf_{q\in{{\rm Cut}}(\mu_{Z^*}(t))}d_{\mathcal{M}^*}(p,q)<\varepsilon\}$, $\Prob(Z^*(t)\in\mathcal{M}^*\setminus({{\rm Cut}}(\mu_{Z^*}(t))\cup({{\rm Cut}}(\mu_{Z^*}(t)))^\varepsilon)\text{~for all~}t\in\mathcal{T})=1$ for some $\varepsilon>0$ and the sample paths of $Z^*$ are continuous.
\item[{\rm(L5)}] $\inf_{t\in\mathcal{T}}\inf_{p\in\mathcal{M}^*:d_{\mathcal{M}^*}(p,\mu_{Z^*}(t))\geq\varepsilon}(F(p,t)-F(\mu_{Z^*}(t),t))>0$ for any $\varepsilon>0$.
\item[{\rm(L6)}] $\inf_{t\in\mathcal{T}}{\lambda_{\rm{min}}\E(H_t)}>0$, where $\lambda_{\rm{min}}\E(H_t)$ is the smallest eigenvalue of $\E(H_t)$ and $H_t$ is defined by $H_t(V)=-\nabla_V\Log_{(\cdot)}Z^*(t)$ for $V\in\mathfrak{X}(\mathcal{M}^*)$.
\item[{\rm(L7)}] $\sup_{t_1\neq t_2}d_{\mathcal{M}^*}(\mu_{Z^*}(t_1),\mu_{Z^*}(t_2))/|t_1-t_2|<\infty$ and $\E((\sup_{t_1\neq t_2}d_{\mathcal{M}^*}(Z^*(t_1),Z^*(t_2))/|t_1-t_2|)^2)<\infty$.
\item[{\rm(L8)}] There exist positive constants $c$ and $\varepsilon$ such that, for all $q\in\mathcal{M}^*$ and $p\in\mathcal{M}^*\setminus({{\rm Cut}}(q)\cup({{\rm Cut}}(q))^\varepsilon)$, the linear map $H_{p,q}:T_p\mathcal{M}^*\rightarrow T_p\mathcal{M}^*$ defined by $H_{p,q}(v)=(-\nabla_V\Log_{(\cdot)}q)(p)$ for $V\in\mathfrak{X}(\mathcal{M}^*)$ with $V(p)=v$ has the operator norm bounded by $c\cdot(1+d_{\mathcal{M}^*}(p,q))$.
\end{itemize}
\end{customcon}

Condition \ref{L} is a weak condition. In particular, conditions (L1) and (L8) on $\mathcal{M}^*$ are superfluous if $\mathcal{M}^*$ is compact. We note that $H_{p,q}(v)$ in (L8) can be written as $(-\nabla_v\Log_{(\cdot)}q)(p)$ and $H_{p,q}$ in (L8) can be written as $(-\nabla\Log_{(\cdot)}q)(p)$ in light of Proposition 4.5 in Lee et al. (2018). A full discussion of condition \ref{L} can be found in \cite{Lin and Yao (2019)} and \cite{Lin et al. (2022)}.

\subsection{Group additive models}\label{Lie case}

Suppose that $\mathcal{N}$ is a Lie group equipped with a group operation $\bullet$. Consider the group additive model
\begin{align}\label{group additive model}
W(t)=g_0(t)\bullet g_1(\xi_1)(t)\bullet\cdots\bullet g_d(\xi_d)(t)\bullet\zeta(t),\quad t\in\mathcal{T},
\end{align}
where the intrinsic mean function of the random function $g_j(\xi_j):\mathcal{T}\rightarrow\mathcal{N}$ exists and is given by the function that is identically equal to the identity element $e$ of $\mathcal{N}$ for all $1\leq j\leq d$, and the conditional intrinsic mean function of $\zeta$ given $\mxi=\mx$ defined by $\mu_{\zeta|\mxi=\mx}(t)=\argmin_{p\in\mathcal{N}}\E((d_{\mathcal{N}}(p,\zeta(t)))^2|\mxi=\mx)$ exists and is given by the function that is identically equal to $e$ for all $\mx$. Proposition 1 in Lin et al. (2022) implies that (\ref{group additive model}) is equivalent to
\begin{align}\label{group additive model2}
\Log_{\mu_{W}(t)}W(t)=\Log_{\mu_{W}(t)}g_0(t)+\sum_{j=1}^d\mathcal{P}^{\gamma(t,\cdot)}_{01}(\Log_e(g_j(\xi_j)(t)))+\mathcal{P}^{\gamma(t,\cdot)}_{01}(\Log_e(\zeta(t))),\quad t\in\mathcal{T},
\end{align}
provided that $\mathcal{N}$ is an abelian Lie group endowed with a bi-invariant Riemannian metric that turns $\mathcal{N}$ into a Hadamard manifold, where $\mathcal{P}^{\gamma(t,\cdot)}_{01}$ is the parallel transport map along the minimizing geodesic $\gamma(t,\cdot):[0,1]\rightarrow\mathcal{N}$ between $\gamma(t,0)=e$ and $\gamma(t,1)=\mu_{W}(t)$ for each $t\in\mathcal{T}$. We note that (\ref{group additive model2}) can be written as
\begin{align}\label{group additive model3}
\Log_{\mu_{W}}W=f_0+\sum_{j=1}^df_j(\xi_j)+\epsilon,
\end{align}
where the function $f_0$ is defined by $f_0(t)=\Log_{\mu_{W}(t)}g_0(t)\in T_{\mu_{W}(t)}\mathcal{N}$ and the random functions $f_j(\xi_j)$ and $\epsilon$ are defined by $f_j(\xi_j)(t)=\mathcal{P}^{\gamma(t,\cdot)}_{01}(\Log_e(g_j(\xi_j)(t)))\in T_{\mu_{W}(t)}\mathcal{N}$ and $\epsilon(t)=\mathcal{P}^{\gamma(t,\cdot)}_{01}(\Log_e(\zeta(t)))\in T_{\mu_{W}(t)}\mathcal{N}$. We note that $\E(f_j(\xi_j)(t))=0_{\mu_{W}(t)}$ and $\E(\epsilon(t)|\mxi)=0_{\mu_{W}(t)}$ for all $t\in\mathcal{T}$, where $0_{\mu_{W}(t)}$ is the zero vector of $T_{\mu_{W}(t)}\mathcal{N}$. Suppose that $\gamma(\cdot,u):\mathcal{T}\rightarrow\mathcal{N}$ is measurable for each $u\in[0,1]$, $g_j(\xi_j)$ and $\zeta$ have continuous sample paths, and
\ba
\E\bigg(\int_{\mathcal{T}}\lng\Log_e(g_j(\xi_j)(t)),\Log_e(g_j(\xi_j)(t))\rng_ed\nu(t)\bigg)&<\infty,\\
\E\bigg(\int_{\mathcal{T}}\lng\Log_e(\zeta(t)),\Log_e(\zeta(t))\rng_ed\nu(t)\bigg)&<\infty,
\ea
where $\lng\cdot,\cdot\rng_e$ is the inner product on $T_e\mathcal{N}$. Then, $f_j(\xi_j)$ and $\epsilon$ can be viewed as random elements taking values in $\mathfrak{T}(\mu_{W})$, and thus model (\ref{group additive model3}) becomes a Hilbertian additive model on $\mathfrak{T}(\mu_{W})$. This indicates that the additive model treated in Section \ref{extension} is equivalent to a group additive model for some Lie groups.

\subsection{Implementation details}\label{suppl-imple}

In this section, we give details on the implementation of our estimator. Suppose that $\hat{\mf}^{[0]}_j(x_j)$ takes the form $\bigoplus_{i=1}^n(w^{i,[0]}_j(x_j)\od\tilde{\mY}^i)$ for some weights $w^{i,[0]}_j(x_j)\in\mathbb{R}$. In this case, (\ref{algorithm}) is equivalent to taking
\begin{align}\label{simple algorithm}
\hat{\mf}^{[r]}_j(x_j)=\oplusin(w_j^{i,[r]}(x_j)\odot\tilde{\mY}^i), \quad 1 \le j \le d,
\end{align}
where $w_j^{i,[r]}(x_j)$ are obtained from the iterative algorithm
\begin{align*}
w_j^{i,[r]}(x_j)=&\frac{K_{h_j}(x_j,\tilde{\xi}_j^i)\I(\tilde{\mxi}^i\in D)}{\hat{p}_j^D(x_j)\cdot\hat{p}_0^D\cdot n}-\frac{\I(\tilde{\mxi}^i\in D)}{\hat{p}_0^D\cdot n}-\underset{k<j}{\sum}\int_{D_k}w_k^{i,[r]}(x_k)\frac{\hat{p}^D_{jk}(x_j,x_k)}{\hat{p}^D_j(x_j)}dx_k\\
&-\underset{k>j}{\sum}\int_{D_k} w_k^{i,[r-1]}(x_k)\frac{\hat{p}^D_{jk}(x_j,x_k)}{\hat{p}^D_j(x_j)}dx_k,\quad1\leq j\leq d,
\end{align*}
containing usual Lebesgue integrals. This follows from the property of Bochner integration that
\begin{align*}
\mbox{(Bochner)}\int_\mathcal{S} f(s)\odot\mh\, d\nu(s)=\mbox{(Lebesgue)}\int_\mathcal{S} f(s)d\nu(s)\odot\mh,
\end{align*}
where $(\mathcal{S},\Sigma,\nu)$ is any measure space, $f$ is any real-valued integrable function defined on $\mathcal{S}$ and $\mh$ is any constant in $\mbH$. In Section \ref{numerical study}, we specifically choose
\[
\hat{\mf}_j^{[0]}(x_j)=n^{-1}\odot \bigoplus_{i=1}^nw_j^{i,[0]}(x_j)\odot\mY^i=n^{-1}\odot \bigoplus_{i=1}^n \Big(\frac{K_{h_j}(x_j,\xi_j^i)}{\hat{p}_j^D(x_j)}-1\Big)\odot\mY^i,
\]
so that they satisfy (\ref{estimated constraint}). We use the biweight-type kernel defined at (\ref{biweight kernel}) for $K_j$. For the convergence criterion of the SBF algorithm, we set
\[
\max_{1 \le j \le d}\int_{D_j}\big\|\hat{\mf}^{[r]}_j(x_j)\ominus\hat{\mf}^{[r-1]}_j(x_j)\big\|^2dx_j<10^{-4}.
\]
For the bandwidth selection, we use the following CBS (Coordinate-wise Bandwidth Selection) scheme introduced by Jeon and Park (2020):

{\bf CBS algorithm}. Let $\CV(h_1,\ldots,h_d)$ denote a cross-validatory criterion for bandwidths $h_1,\ldots,h_d$. Take a bandwidth grid $\{g_{j1},\ldots,g_{jG_j}\}$ with $G_j\in\mathbb{N}$ for each $1\le j\le d$. Choose an initial bandwidth $h^{(0)}_j$ from $\{g_{j1},\ldots,g_{jG_j}\}$ for each $1\leq j\leq d$.
For $t=1,2,\ldots$, find
\begin{align*}
h^{(t)}_j=\underset{g_j\in\{g_{j1},\ldots,g_{jG_j}\}}{\argmin}\CV(h^{(t)}_1,\ldots,h^{(t)}_{j-1},g_j,h^{(t-1)}_{j+1},\ldots,h^{(t-1)}_d), \quad 1 \le j \le d.
\end{align*}
Repeat the procedure until $(h^{(t)}_1,\ldots,h^{(t)}_d)=(h^{(t-1)}_1,\ldots,h^{(t-1)}_d)$.\qed

We take a 5-fold cross-validation for the CV criterion. We choose $\{g_{j1},\ldots,g_{jG_j}\}=\{c_j+b_j\times k:k=0,\ldots,20\}$ for some small constant $c_j>0$ that satisfies condition \ref{A}. In the simulation studies, we choose $b_j=0.01$ and $b_j=0.05$ when $L_j=1$ and $L_j=2$, respectively. In the real data applications, we choose $b_j=0.025$ and $b_j=a_j/40$ when $L_j=1$ and $L_j=2$, respectively. Here, $a_j=\max\{\|\tilde{\xi}^i_j-\tilde{\xi}^{i'}_j\|_j:1\leq i,i'\leq n\}$.

For the kernel $K^*$ and the bandwidth $h_i^*$ used to construct $\tilde{Y}^i(\cdot)$ in (\ref{Pelletier density}), we take the Epanechnikov kernel and 10-fold cross-validatory bandwidth selected from $\{s_i+0.1\times k:k=0,\ldots,40\}$ for some small $s_i$ inducing that $\tilde{Y}^i(\cdot)>0$. The latter standard bandwidth selection scheme, which targets to minimize $\int_{\mathcal{S}}(\tilde{Y}^i(s)-Y^i(s))^2d\nu(s)$, is reasonable since $\|\tilde{Y}^i(\cdot)\om Y^i(\cdot)\|^2\leq\const\cdot\int_{\mathcal{S}}(\tilde{Y}^i(s)-Y^i(s))^2d\nu(s)$ with probability tending to one; see the proof of Proposition \ref{imperfect density prop 2}.

\subsection{Some general lemmas}\label{collection of lemmas}

We introduce two lemmas on the uniform rates of convergence for kernel-type estimators with general semi-metric space-valued predictors. The lemmas generalize Lemma S.1 in Jeon et al. (2021a) since the former lemmas involve an additional indicator function. Before we state them, we introduce the notion of upper Minkowski dimension. The upper Minkowski dimension of a semi-metric space $\mathcal{M}$ is defined as
$\limsup_{\delta\rightarrow0}(\log N(\mathcal{M},\delta)/\log\delta^{-1})$, where $N(\mathcal{M},\delta)$ is the $\delta$-covering number of $\mathcal{M}$. Examples of semi-metric spaces with finite upper Minkowski dimension include compact subsets of Euclidean spaces.

Let $(\mathcal{X}_j,\rho_j)$ for $1\leq j\leq d$ and $d\geq1$ be semi-metric spaces. Let also $\mu_j$ be the Borel measure on $\mathcal{X}_j$ for $1\leq j\leq d$. Define the product space $\mathcal{X}=\prod_{j=1}^d\mathcal{X}_j$ and the product measure $\mu=\bigotimes_{j=1}^d\mu_j$. Let $\mX=(X_1,\ldots,X_d)$ be a random element taking values in $\mathcal{X}$ and $p_\mX$ be the density of $\mX$ with respect to $\mu$. Let $\mathcal{M}(n)$ be a sequence of subsets of $\mathcal{X}$.

Let $J$ be any subset of $\{1,\ldots,d\}$ and define $\mathcal{X}_J=\prod_{j\in J}\mathcal{X}_j$. For $\mx_J=(x_j:j\in J)\in\mathcal{X}_J$, define $B_{\mathcal{X}_J}(\mx_J,R)=\{(u_j:j\in J)\in\mathcal{X}_J:\sum_{j\in J}\rho^2_j(x_j,u_j)<R^2\}$ for $R>0$. Let $\mathcal{M}_j$ be a subset of $\mathcal{X}_j$ for $j\in J$ and $\mathcal{M}_J:=\prod_{j\in J}\mathcal{M}_j$ have a finite upper Minkowski dimension. For $j\in J$, let $A_{nj}$ and $B_{nj}$ be positive sequences and $L_{nj}:\mathcal{X}_j\times\mathcal{X}_j\rightarrow\mathbb{R}$ be a sequence of measurable functions.

\begin{customlemma}{S.1}\label{uniform rates general}
Let $\mdelta_n$ be a sequence of $\mbH$-valued random elements such that $\E(\|\mdelta_n\|^\alpha)=O(U_{n\alpha})$ for some $2\leq\alpha<\infty$ and some positive sequence $U_{n\alpha}$. Assume that (i) for all $j\in J$, $\sup_{x_j\in \mathcal{M}_j, \bu\in\bigcup_n\mathcal{M}(n)}|L_{nj}(x_j,u_j)|=O(A_{nj})$, and that there exists a constant $C_j>0$ such that $|L_{nj}(x_j,u_j)-L_{nj}(x_j^*,u_j)|\leq C_j\cdot A_{nj}\cdot B_{nj}\cdot \rho_j(x_j,x_j^*)$ for all $x_j, x_j^*\in \mathcal{M}_j$ and $\bu\in\bigcup_n\mathcal{M}(n)$; (ii) $\sup_{\mx_J\in \mathcal{M}_J}\int_{\mathcal{M}(n)}\E(\|\mdelta_n\|^2|\mX=\bu)\cdot\prod_{j\in J}(L_{nj}(x_j,u_j))^2d\mu(\bu)=O(V_{nJ}\cdot\prod_{j\in J}A_{nj})$ for some positive sequence $V_{nJ}$; (iii) $n^{-1+2\beta+2/\alpha}\cdot\prod_{j\in J}A_{nj}=o(1)$ and $n^{-\beta}\cdot(\log{n})^{1/2}\cdot U_{n\alpha}^{1/\alpha}\cdot V_{nJ}^{-1/2}=O(1)$ for some constant $\beta>0$, and $n^c\cdot\max_{j\in J}B_{nj}\cdot U_{n\alpha}^{1/\alpha}\cdot V_{nJ}^{-1/2}=O(1)$ for some constant $c\in\mathbb{R}$; (iv) $p_\mX$ is bounded on $\bigcup_n\mathcal{M}(n)$. Let $\{(\mX^i,\mdelta_n^i):1\leq i\leq n\}$ be a set of i.i.d. copies of $(\mX,\mdelta_n)$. Then, it holds that
\begin{align*}
\sup_{\mx_J\in \mathcal{M}_J}\|\mS_n(\mx_J)\ominus \E(\mS_n(\mx_J))\|=O_p\left(\sqrt{\frac{\log{n}\cdot V_{nJ}\cdot\prod_{j\in J}A_{nj}}{n}}\right)
\end{align*}
for $\mS_n(\mx_J):=n^{-1}\odot\bigoplus_{i=1}^n(\I(\mX^i\in \mathcal{M}(n))\cdot\prod_{j\in J}L_{j,n}(x_j,X_j^i))\odot\mdelta_n^i$.
\end{customlemma}

\begin{proof}
Since the upper Minkowski dimension of $\mathcal{M}_J$ is finite, for sufficiently large $\kappa>0$, there exist $N(n^{-\kappa})\in\mathbb{N}$, $0<m<\infty$ and $\{\mx_J^{(1)},\cdots,\mx_J^{(N(n^{-\kappa}))}\}\subset \mathcal{M}_J$ such that $N(n^{-\kappa})=O(n^{\kappa m})$ and $\{B_{\mathcal{X}_J}(\mx_J^{(l)},n^{-\kappa}):1\leq l\leq N(n^{-\kappa})\}$ covers $\mathcal{M}_J$. For the constant $\beta$ in the condition (iii), it holds that
\begin{align*}
\Prob\bigg(\|\mdelta_n^i\|\leq n^{1/2-\beta}\cdot U_{n\alpha}^{1/\alpha}\cdot\prod_{j\in J}A_{nj}^{-1/2}\text{~for all~}1\leq i\leq n\bigg)\geq1-\frac{\E(\|\mdelta_n\|^\alpha)n}{U_{n\alpha}(n^{1/2-\beta}\cdot\prod_{j\in J}A_{nj}^{-1/2})^\alpha}\rightarrow1
\end{align*}
from the moment condition on $\mdelta_n$ and the first condition in (iii).
Hence,
\begin{align*}
\sup_{\mx_J\in\mathcal{M}_J}&\bigg\|\frac{1}{n}\od\bigoplus_{i=1}^n\bigg(\I(\mX^i\in \mathcal{M}(n))\prod_{j\in J}L_{nj}(x_j,X_j^i)\I\bigg(\|\mdelta_n^i\|>n^{1/2-\beta} U_{n\alpha}^{1/\alpha}\prod_{j\in J}A_{nj}^{-1/2}\bigg)\bigg)\od\mdelta_n^i\\
&\om \E\bigg(\bigg(\I(\mX\in \mathcal{M}(n))\prod_{j\in J}L_{nj}(x_j,X_j)\I\bigg(\|\mdelta_n\|>n^{1/2-\beta}U_{n\alpha}^{1/\alpha}\prod_{j\in J}A_{nj}^{-1/2}\bigg)\bigg)\od\mdelta_n\bigg)\bigg\|=0
\end{align*}
with probability tending to one. Thus,
\begin{align*}
\mS_n(\mx_J)\om \E(\mS_n(\mx_J))=\mU_n(\mx_J)\om \E(\mU_n(\mx_J))\op o_p\left(\sqrt{\frac{\log{n}\cdot V_{nJ}\cdot\prod_{j\in J}A_{nj}}{n}}\right),
\end{align*}
where
\begin{align*}
\mU_n(\mx_J)=\frac{1}{n}\od\bigoplus_{i=1}^n\bigg(\I(\mX^i\in \mathcal{M}(n))\prod_{j\in J}L_{nj}(x_j,X_j^i)\I\bigg(\|\mdelta_n^i\|\leq n^{1/2-\beta} U_{n\alpha}^{1/\alpha}\prod_{j\in J}A_{nj}^{-1/2}\bigg)\bigg)\od\mdelta_n^i.
\end{align*}
Note that
\begin{align*}
&\sup_{\mx_J\in\mathcal{M}_J}\|\mU_n(\mx_J)\om \E(\mU_n(\mx_J))\|\\
&\leq\max_{1\leq l\leq N(n^{-\kappa})}\|\mU_n(\mx_J^{(l)})\om \E(\mU_n(\mx_J^{(l)}))\|+\max_{1\leq l\leq N(n^{-\kappa})}\sup_{\mx_J\in \mathcal{M}_J\cap B_{\mathcal{X}_J}(\mx_J^{(l)},n^{-\kappa})}\|\mU_n(\mx_J)\om \mU_n(\mx_J^{(l)})\|\\
&\qquad\qquad+\max_{1\leq l\leq N(n^{-\kappa})}\sup_{\mx_J\in \mathcal{M}_J\cap B_{\mathcal{X}_J}(\mx_J^{(l)},n^{-\kappa})}\|\E(\mU_n(\mx_J))\om \E(\mU_n(\mx_J^{(l)}))\|.
\end{align*}
It follows that
\begin{align*}
&\max_{1\leq l\leq N(n^{-\kappa})}\sup_{\mx_J\in \mathcal{M}_J\cap B_{\mathcal{X}_J}(\mx_J^{(l)},n^{-\kappa})}\|\mU_n(\mx_J)\om \mU_n(\mx_J^{(l)})\|\\
&\qquad\qquad \leq \const\cdot n^{1/2-\beta-\kappa}\cdot U_{n\alpha}^{1/\alpha}\cdot\prod_{j\in J}A_{nj}^{1/2}\cdot\max_{j\in J}B_{nj},\textrm{ and }\\
&\max_{1\leq l\leq N(n^{-\kappa})}\sup_{\mx_J\in \mathcal{M}_J\cap B_{\mathcal{X}_J}(\mx_J^{(l)},n^{-\kappa})}\|\E(\mU_n(\mx_J))\om \E(\mU_n(\mx_J^{(l)}))\|\\
&\qquad\qquad \leq \const\cdot n^{1/2-\beta-\kappa}\cdot U_{n\alpha}^{1/\alpha}\cdot\prod_{j\in J}A_{nj}^{1/2}\cdot\max_{j\in J}B_{nj},
\end{align*}
from condition (i). Hence, by taking sufficiently large $\kappa$, we get
\begin{align*}
\sup_{\mx_J\in\mathcal{M}_J}\|\mU_n(\mx_J)\om \E(\mU_n(\mx_J))\|\leq&\max_{1\leq l\leq N(n^{-\kappa})}\|\mU_n(\mx_J^{(l)})\om \E(\mU_n(\mx_J^{(l)}))\|\\
&+o_p\left(\sqrt{\frac{\log{n}\cdot V_{nJ}\cdot\prod_{j\in J}A_{nj}}{n}}\right)
\end{align*}
from the second condition in (iii). Thus, it suffices to show that
\begin{align}\label{uniform rate claim}
\max_{1\leq l\leq N(n^{-\kappa})}\|\mU_n(\mx_J^{(l)})\om \E(\mU_n(\mx_J^{(l)}))\|=O_p\left(\sqrt{\frac{\log{n}\cdot V_{nJ}\cdot\prod_{j\in J}A_{nj}}{n}}\right).
\end{align}
For any $\mx_J\in\mathcal{M}_J$, we define
\begin{align*}
\mZ_n^i(\mx_J)=&\bigg(\frac{1}{n}\I(\mX^i\in \mathcal{M}(n))\prod_{j\in J}L_{nj}(x_j,X_j^i)\I\bigg(\|\mdelta_n^i\|\leq n^{1/2-\beta}U_{n\alpha}^{1/\alpha}\prod_{j\in J}A_{nj}^{-1/2}\bigg)\bigg)\od\mdelta_n^i\\
&\om \E\bigg(\bigg(\frac{1}{n}\I(\mX\in \mathcal{M}(n))\prod_{j\in J}L_{nj}(x_j,X_j)\I\bigg(\|\mdelta_n\|\leq n^{1/2-\beta}U_{n\alpha}^{1/\alpha}\prod_{j\in J}A_{nj}^{-1/2}\bigg)\bigg)\od\mdelta_n\bigg).
\end{align*}
Then, we have
\begin{align*}
\E(\mZ_n^i(\mx_J))&= \bzero, \quad  \|\mZ_n^i(\mx_J)\|\leq \const\cdot n^{-1/2-\beta}\cdot U_{n\alpha}^{1/\alpha}\cdot\prod_{j\in J}A_{nj}^{1/2},\quad\mbox{and}\\
\sum_{i=1}^n\E(\|\mZ_n^i(\mx_J)\|^2)&\leq \frac{1}{n}\cdot\E\bigg(\I(\mX\in \mathcal{M}(n))\cdot\|\mdelta_n\|^2\cdot\prod_{j\in J}(L_{nj}(x_j,X_j))^2\bigg)\\
&=\frac{1}{n}\cdot\int_{\mathcal{X}}\I(\bu\in \mathcal{M}(n))\cdot\E(\|\mdelta_n\|^2|\mX=\bu)\cdot p_\mX(\bu)\cdot\prod_{j\in J}(L_{nj}(x_j,u_j))^2d\mu(\bu)\\
&\leq \frac{\const}{n}\cdot\int_{\mathcal{M}(n)}\E(\|\mdelta_n\|^2|\mX=\bu)\cdot\prod_{j\in J}(L_{nj}(x_j,u_j))^2d\mu(\bu)\\
&\leq \frac{\const \cdot V_{nJ}\cdot\prod_{j\in J}A_{nj}}{n},
\end{align*}
where the first inequality follows from the first condition in (i), the second inequality follows from condition (iv) and the last inequality follows from condition (ii). Now, Theorem 2.6.2 in Bosq (2000) with the third condition in (iii) and the fact $N(n^{-\kappa})=O(n^{\kappa m})$ gives (\ref{uniform rate claim}).
\end{proof}

The next lemma considers the case where $\mdelta_n$ is a bounded random element possibly depending on $n$.

\begin{customlemma}{S.2}\label{uniform rates general density}
Let $\mdelta_n$ be a $\mbH$-valued random element such that $\|\mdelta_n\|\cdot\prod_{j\in J}L_{nj}(x_j,u_j)\leq C\cdot U_n\cdot\prod_{j\in J}L_{nj}(x_j,u_j)$ for some constant $C>0$ and positive sequence $U_n$. Assume conditions (i) and (iv) of Lemma \ref{uniform rates general} hold, that $\sup_{\mx_J\in \mathcal{M}_J}\int_{\mathcal{M}(n)}\prod_{j\in J}(L_{nj}(x_j,u_j))^2d\mu(\bu)=O(V_{nJ}\cdot\prod_{j\in J}A_{nj})$ for some positive sequence $V_{nJ}$, that $n^c\cdot\max_{j\in J}B_{nj}\cdot V_{nJ}^{-1/2}\cdot\prod_{j\in J}A_{nj}^{1/2}=O(1)$ for some constant $c\in\mathbb{R}$ and that $n^{-1/2}\cdot(\log n)^{1/2}\cdot V_{nJ}^{-1/2}\cdot\prod_{j\in J}A^{1/2}_{nj}=O(1)$. Let $\{(\mX^i,\mdelta_n^i):1\leq i\leq n\}$ be a set of i.i.d. copies of $(\mX,\mdelta_n)$. Then, it holds that
\begin{align*}
\sup_{\mx_J\in \mathcal{M}_J}\|\mS_n(\mx_J)\ominus \E(\mS_n(\mx_J))\|=O_p\left(U_n\cdot\sqrt{\frac{\log{n}\cdot V_{nJ}\cdot\prod_{j\in J}A_{nj}}{n}}\right)
\end{align*}
for $\mS_n(\mx_J):=n^{-1}\odot\bigoplus_{i=1}^n(\I(\mX^i\in \mathcal{M}(n))\cdot\prod_{j\in J}L_{j,n}(x_j,X_j^i))\odot\mdelta_n^i$.
\end{customlemma}

\begin{proof}
The proof of this lemma is similar to that of Lemma \ref{uniform rates general}. We therefore omit the reasoning for each argument. For sufficiently large $\kappa>0$, there exist $N(n^{-\kappa})\in\mathbb{N}$, $0<m<\infty$ and $\{\mx_J^{(1)},\cdots,\mx_J^{(N(n^{-\kappa}))}\}\subset \mathcal{M}_J$ such that $N(n^{-\kappa})=O(n^{\kappa m})$ and $\{B_{\mathcal{X}_J}(\mx_J^{(l)},n^{-\kappa}):1\leq l\leq N(n^{-\kappa})\}$ covers $\mathcal{M}_J$.
Note that
\begin{align*}
\sup_{\mx_J\in\mathcal{M}_J}\|\mS_n(\mx_J)\om \E(\mS_n(\mx_J))\|\leq&\max_{1\leq l\leq N(n^{-\kappa})}\|\mS_n(\mx_J^{(l)})\om \E(\mS_n(\mx_J^{(l)}))\|+\\
&\max_{1\leq l\leq N(n^{-\kappa})}\sup_{\mx_J\in \mathcal{M}_J\cap B_{\mathcal{X}_J}(\mx_J^{(l)},n^{-\kappa})}\|\mS_n(\mx_J)\om \mS_n(\mx_J^{(l)})\|+\\
&\max_{1\leq l\leq N(n^{-\kappa})}\sup_{\mx_J\in \mathcal{M}_J\cap B_{\mathcal{X}_J}(\mx_J^{(l)},n^{-\kappa})}\|\E(\mS_n(\mx_J))\om \E(\mS_n(\mx_J^{(l)}))\|.
\end{align*}
It follows that
\begin{align*}
&\max_{1\leq l\leq N(n^{-\kappa})}\sup_{\mx_J\in \mathcal{M}_J\cap B_{\mathcal{X}_J}(\mx_J^{(l)},n^{-\kappa})}\|\mS_n(\mx_J)\om \mS_n(\mx_J^{(l)})\|\\
&\qquad\leq\const\cdot n^{-\kappa}\cdot U_n\cdot\prod_{j\in J}A_{nj}\cdot\max_{j\in J}B_{nj},\textrm{ and }\\
&\max_{1\leq l\leq N(n^{-\kappa})}\sup_{\mx_J\in \mathcal{M}_J\cap B_{\mathcal{X}_J}(\mx_J^{(l)},n^{-\kappa})}\|\E(\mS_n(\mx_J))\om \E(\mS_n(\mx_J^{(l)}))\|\\
&\qquad\leq\const\cdot n^{-\kappa}\cdot U_n\cdot\prod_{j\in J}A_{nj}\cdot\max_{j\in J}B_{nj}.
\end{align*}
Hence, by taking sufficiently large $\kappa$, we get
\begin{align*}
\sup_{\mx_J\in\mathcal{M}_J}\|\mS_n(\mx_J)\om \E(\mS_n(\mx_J))\|\leq&\max_{1\leq l\leq N(n^{-\kappa})}\|\mS_n(\mx_J^{(l)})\om \E(\mS_n(\mx_J^{(l)}))\|\\
&+o_p\left(U_n\cdot\sqrt{\frac{\log{n}\cdot V_{nJ}\cdot\prod_{j\in J}A_{nj}}{n}}\right).
\end{align*}
Thus, it suffices to show that
\begin{align}\label{uniform rate claim density}
\max_{1\leq l\leq N(n^{-\kappa})}\|\mS_n(\mx_J^{(l)})\om \E(\mS_n(\mx_J^{(l)}))\|=O_p\left(U_n\cdot\sqrt{\frac{\log{n}\cdot V_{nJ}\cdot\prod_{j\in J}A_{nj}}{n}}\right).
\end{align}
For any $\mx_J\in\mathcal{M}_J$, we define
\begin{align*}
\mZ_n^i(\mx_J)=&\bigg(\frac{1}{n}\I(\mX^i\in \mathcal{M}(n))\prod_{j\in J}L_{nj}(x_j,X_j^i)\bigg)\od\mdelta_n^i\\
&\om \E\bigg(\bigg(\frac{1}{n}\I(\mX\in \mathcal{M}(n))\prod_{j\in J}L_{nj}(x_j,X_j)\bigg)\od\mdelta_n\bigg).
\end{align*}
Then, we have
\begin{align*}
\E(\mZ_n^i(\mx_J))&=\bzero,\quad \|\mZ_n^i(\mx_J)\|\leq \const\cdot U_n\cdot\frac{1}{n}\cdot\prod_{j\in J}A_{nj},\quad\mbox{and}\\
\sum_{i=1}^n\E(\|\mZ_n^i(\mx_J)\|^2)&\leq\frac{\const}{n}\cdot U_n^2\cdot\E\bigg(\I(\mX\in \mathcal{M}(n))\prod_{j\in J}(L_{nj}(x_j,X_j))^2\bigg)\\
&=\frac{\const}{n}\cdot U_n^2\cdot\int_{\mathcal{X}}\I(\bu\in \mathcal{M}(n))\prod_{j\in J}(L_{nj}(x_j,u_j))^2p_\mX(\bu)d\mu(\bu)\\
&\leq\frac{\const}{n}\cdot U_n^2\cdot \int_{\mathcal{M}(n)}\prod_{j\in J}(L_{nj}(x_j,u_j))^2d\mu(\bu)\\
&\leq\frac{\const \cdot U_n^2\cdot V_{nJ}\cdot\prod_{j\in J}A_{nj}}{n}.
\end{align*}
Now, Theorem 2.6.2 in Bosq (2000) gives (\ref{uniform rate claim density}).
\end{proof}

\begin{customremark}{S.1}\label{uniform rate remark}\leavevmode
\begin{itemize}
\item[1.] In Lemma \ref{uniform rates general}, one can take $\mathcal{M}(n)$ or $\mdelta_n$ independent of $n$. Suppose that $\mathcal{M}(n)\equiv\mathcal{M}:=\prod_{j=1}^d\mathcal{M}_j$ and $\mdelta_n\equiv\mdelta$, where $\mathcal{M}_j$ for each $j\notin J$ is a subset of $\mathcal{X}_j$ such that $\mu_j(\mathcal{M}_j)<\infty$. In that case, one obtains the uniform rate $O_p\left(\sqrt{n^{-1}\cdot\log{n}\cdot\prod_{j\in J}A_{nj}}\right)$ under the following conditions: $\E(\|\mdelta\|^\alpha)<\infty$ for some constant $2\leq\alpha<\infty$, conditions (i) and (iv) of Lemma \ref{uniform rates general} hold, $\max_{j\in J}\sup_{x_j\in \mathcal{M}_j}\int_{\mathcal{M}_j}|L_{nj}(x_j,u_j)|d\mu_j(u_j)=O(1)$, $n^{-1+2\beta+2/\alpha}\cdot\prod_{j\in J}A_{nj}=o(1)$ for some constant $\beta>0$, $n^c\cdot\max_{j\in J}B_{nj}=O(1)$ for some constant $c\in\mathbb{R}$ and $\E(\|\mdelta\|^2|\mX=\cdot)$ is bounded on $\mathcal{M}$.
\item[2.] In Lemma \ref{uniform rates general density}, one can replace the conditions that $n^c\cdot\max_{j\in J}B_{nj}\cdot V_{nJ}^{-1/2}\cdot\prod_{j\in J}A^{1/2}_{nj}=O(1)$  for some constant $c\in\mathbb{R}$ and that $n^{-1/2}\cdot(\log n)^{1/2}\cdot V_{nJ}^{-1/2}\cdot\prod_{j\in J}A^{1/2}_{nj}=O(1)$ by the following alternative sufficient condition using Lemma \ref{uniform rates general} with $\mdelta_n=\I(\mX\in\mathcal{M}(n))$: There exist positive constants $\alpha$ and $\beta$ and a constant $c\in\mathbb{R}$ such that $n^{-1+2\beta+2/\alpha}\cdot\prod_{j\in J}A_{nj}=o(1)$, $n^{-\beta}\cdot(\log{n})^{1/2}\cdot \Prob(\mX\in\mathcal{M}(n))^{1/\alpha}\cdot V_{nJ}^{-1/2}=O(1)$, and $n^c\cdot\max_{j\in J}B_{nj}\cdot \Prob(\mX\in\mathcal{M}(n))^{1/\alpha}\cdot V_{nJ}^{-1/2}=O(1)$.
\item[3.] In Lemma \ref{uniform rates general density} and in Remark \ref{uniform rate remark}-2, one can take $\mathcal{M}(n)\equiv\mathcal{M}$. In that case, one obtains the uniform rate $O_p\left(U_n\cdot\sqrt{n^{-1}\cdot\log{n}\cdot\prod_{j\in J}A_{nj}}\right)$ under the respective conditions with $V_{nJ}\equiv1$.
\item[4.] In Lemmas \ref{uniform rates general}-\ref{uniform rates general density} and in the above remarks, the second condition in (i) of Lemma \ref{uniform rates general} can be replaced by the following alternative sufficient condition: There exist positive constants $C_j$ and $R_j$ such that $|L_{nj}(x_j,u_j)-L_{nj}(x_j^*,u_j)|\leq C_j\cdot A_{nj}\cdot B_{nj}\cdot \rho_j(x_j,x_j^*)$ for all $\bu\in\bigcup_n\mathcal{M}(n)$ and $x_j, x_j^*\in \mathcal{M}_j\cap B_{\mathcal{X}_j}(u_j,R_j)$, and there exists a positive sequence $C_{nj}$ eventually less than $R_j$ such that $L_{nj}(x_j,u_j)=0$ for all $x_j\in \mathcal{M}_j$ and $\bu\in\bigcup_n\mathcal{M}(n)$ with $\rho_j(x_j,u_j)>C_{nj}$.
\end{itemize}
\end{customremark}

\subsection{Specific lemmas}

Define $\check{p}_0^D$, $\check{p}^D$, $\check{p}_j^D$ and $\check{p}_{jk}^D$, the ``oracle'' equivalent of $p_0^D,p^D,p_j^D$ and $p_{jk}^D$, as $\hat{p}_0^D$, $\hat{p}^D$, $\hat{p}_j^D$ and $\hat{p}_{jk}^D$ with $\tilde{\mxi}^i$, $\tilde{\xi}^i_j$ and $\tilde{\xi}^i_{jk}$ in their respective definitions replaced with $\mxi^i$, $\xi^i_j$ and $\xi^i_{jk}$, respectively. We also define $\check{\mf}_0$ and $\check{\mm}_j$ as $\hat{\mf}_0$ and $\hat{\mm}_j$ with $\hat{p}_0^D$, $\hat{p}_j^D$, $\tilde{\mxi}^i$ and $\tilde{\xi}^i_j$ in the definitions of $\hat{\mf}_0$ and $\hat{\mm}_j$ being replaced by $\check{p}_0^D$, $\check{p}_j^D$, $\mxi^i$ and $\xi^i_j$, respectively.

Recall that $\max_{1\leq i\leq n}\|\tilde{\xi}_j^i-\xi_j^i\|_j=O_p(a_{nj})$ for all $1\leq j\leq d$. Throughout the proofs in the Supplementary Material, we assume that
\begin{align}\label{Op event}
\max_{1\leq i\leq n}\|\tilde{\xi}_j^i-\xi_j^i\|_j\leq C\cdot a_{nj}\text{~for all~}1\leq j\leq d\text{~for some constant~}C>0
\end{align}
since if a rate of convergence holds under the event (\ref{Op event}), then the rate of convergence holds under the whole probability space. We note that the event (\ref{Op event}) guarantees (\ref{indicator inequality}) presented in the following lemma.

\begin{customlemma}{S.3}\label{p0 rate}
Assume that condition (B3) holds and that $p$ is bounded on $D^+(\varepsilon)$ for some $\varepsilon>0$. Then, it holds that
\ba
|\hat{p}_0^D-\check{p}_0^D|&=O_p\bigg(\sum_{j=1}^da_{nj}+n^{-1/2}\cdot\bigg(\sum_{j=1}^da_{nj}\bigg)^{1/2}\bigg),\\
|\hat{p}_0^D-p_0^D|&=O_p\bigg(\sum_{j=1}^da_{nj}+n^{-1/2}\bigg).
\ea
\end{customlemma}

\begin{proof}
We define
\begin{align}\label{Dn definition}
\begin{split}
D(j,n)&=\bigg\{\mx\in\prod_{j=1}^dD_j^+(C\cdot a_{nj}):x_j\in \bigcup_{u_j\in \partial D_j}\bar{B}_j(u_j,C\cdot a_{nj})\bigg\},\\
D(n)&=\bigcup_{j=1}^dD(j,n).
\end{split}
\end{align}
Then, one can check that
\begin{align}\label{indicator inequality}
|\I(\tilde{\mxi}^i\in D)-\I(\mxi^i\in D)|\leq\I(\mxi^i\in D(n)).
\end{align}
Hence, $|\hat{p}_0^D-\check{p}_0^D|\leq n^{-1}\sum_{i=1}^n\I(\mxi^i\in D(n))$. Also, it holds that
\begin{align}\label{indicator rate}
\begin{split}
\E\bigg(n^{-1}\sum_{i=1}^n\I(\mxi\in D(n))\bigg)&=\Prob(\mxi\in D(n))\leq\const\Leb(D(n))\\
&\leq\const\sum_{j=1}^d\Leb(D(j,n))=O\bigg(\sum_{j=1}^da_{nj}\bigg),
\end{split}
\end{align}
where the first inequality follows from the fact that $D(n)\subset D^+(\varepsilon)$ for sufficiently large $n$ with the boundedness of $p$ on $D^+(\varepsilon)$ and the last equality follows from Lemma A4.3 in Kallenberg (2017) with (B3). Also, we have $\Var\left(n^{-1}\sum_{i=1}^n\I(\mxi^i\in D(n))\right)\leq n^{-1}\E(\I(\mxi\in D(n)))=O(n^{-1}\cdot\sum_{j=1}^da_{nj})$. Thus, the first assertion follows. For the second assertion, we note that $|\check{p}_0^D-p_0^D|=O_p(n^{-1/2})$. This with the first assertion gives the second assertion.
\end{proof}

We now provide six new lemmas (Lemmas \ref{uniform rates}-\ref{uniform consistency general}) regarding the uniform rates of convergence. Each lemma is based on a minimal assumption. They have different applications to different stochastic terms in our asymptotic analysis and they are of interest in their own right. In the lemmas, $J$ is any subset of $\{1,\ldots,d\}$.

\begin{customlemma}{S.4}\label{uniform rates}
Let $\{\mdelta^i:1\leq i\leq n\}$ be a set of i.i.d. copies of a $\mbH$-valued random element $\mdelta$ such that $\E(\|\mdelta\|^\alpha)<\infty$ for some $2<\alpha<\infty$.
Assume that condition (B3) holds, that $K_j(\|\cdot\|_j):\mathbb{R}^{L_j}\rightarrow[0,\infty)$ are Lipschitz continuous for all $j\in J$, that $\max_{j\in J}h_j=O(1)$, that $n^{-1+2\beta+2/\alpha}\cdot\prod_{j\in J}h_j^{-L_j}=o(1)$ for some constant $\beta>0$ and that $p$ and $\E(\|\mdelta\|^2|\mxi=\cdot)$ are bounded on $D$. Then, it holds that
\begin{align*}
\sup_{\mx_J\in\prod_{j\in J}D_j}\|\mS_n(\mx_J)\ominus \E(\mS_n(\mx_J))\|=O_p\left(\sqrt{\frac{\log{n}}{n\cdot\prod_{j\in J}h_j^{L_j}}}\right)
\end{align*}
for $\mS_n(\mx_J):=n^{-1}\odot\bigoplus_{i=1}^n(\I(\mxi^i\in D)\cdot\prod_{j\in J}K_{h_j}(x_j,\xi_j^i))\odot\mdelta^i$. The same uniform rate of convergence holds for $\mS^*_n(\mx_J):=n^{-1}\odot\bigoplus_{i=1}^n(\I(\mxi^i\in D)\cdot\prod_{j\in J}h_j^{-L_j}K_j(\|x_j-\xi_j^i\|_j/h_j))\odot\mdelta^i$ without condition (B3).
\end{customlemma}

\begin{proof}
We apply Remark \ref{uniform rate remark}-1 for the proof. We first prove the first part of the lemma. For this, it suffices to verify that the following three assertions hold with $\mathcal{M}(n)\equiv D$, $\mathcal{M}_j=D_j$, $L_{nj}=K_{h_j}$, $A_{nj}=h_j^{-L_j}$ and $B_{nj}=h_j^{-1}$: (1) the condition (i) of Lemma \ref{uniform rates general}; (2) $\max_{j\in J}\sup_{x_j\in \mathcal{M}_j}\int_{\mathcal{M}_j}|L_{nj}(x_j,u_j)|d\mu_j(u_j)=O(1)$; (3) $n^c\cdot\max_{j\in J}B_{nj}=O(1)$ for some constant $c\in\mathbb{R}$. For the first assertion, we note that
\begin{align}\label{lower bound}
\begin{split}
&\inf_n\inf_{u_j\in D_j}\int_{D_j}h_j^{-L_j}K_j\bigg(\frac{\|t_j-u_j\|_j}{h_j}\bigg)dt_j\\
&=\inf_n\inf_{u_j\in D_j}\int_{(D_j-u_j)/h_j}K_j(\|s_j\|_j)ds_j\\
&\geq \inf_n\inf_{u_j\in D_j}\int_{(D_j-u_j)/h_j\cap B_j(\bzero_j,\tau_j)}K_j(\|s_j\|_j)ds_j\\
&\geq \inf_{t\in[0,\tau_j]}K_j(t)\cdot\inf_n\inf_{u_j\in D_j}\Leb_j((D_j-u_j)/h_j\cap B_j(\bzero_j,\tau_j))\\
&\geq \const,
\end{split}
\end{align}
where $(D_j-u_j)/h_j=\{(x_j-u_j)/h_j:x_j\in D_j\}$ and the last inequality follows from (B3) and that $\max_{j\in J}h_j=O(1)$. In fact,
\ba
\inf_n\inf_{u_j\in D_j}\Leb_j((D_j-u_j)/h_j\cap B_j(\bzero_j,\tau_j))\geq\const
\ea
under (B3) and that $\max_{j\in J}h_j=O(1)$, as demonstrated in Jeon et al. (2021a). Hence,
\begin{align}\label{denominator bound}
K_{h_j}(x_j,u_j)\leq \const\cdot h_j^{-L_j}K_j\bigg(\frac{\|x_j-u_j\|_j}{h_j}\bigg),
\end{align}
where the $\const$ is independent of $x_j, u_j\in D_j$ and $n$. We note that the Lipschitz continuity of $K_j(\|\cdot\|_j)$ implies that $K_j(\|\cdot\|_j)$ is bounded on $\mathbb{R}^{L_j}$. Hence, $\sup_{x_j,u_j\in D_j}K_{h_j}(x_j,u_j)=O(h_j^{-L_j})$ from (\ref{denominator bound}) and the boundedness of $K_j(\|\cdot\|_j)$.
Also, (\ref{denominator bound}) with the Lipschitz continuity of $K_j(\|\cdot\|_j)$ implies
\ba
|K_{h_j}(x_j,u_j)-K_{h_j}(x^*_j,u_j)|&\leq\const\cdot h_j^{-L_j}\cdot\bigg|K_j\bigg(\frac{\|x_j-u_j\|_j}{h_j}\bigg)-K_j\bigg(\frac{\|x_j^*-u_j\|_j}{h_j}\bigg)\bigg|\\
&\leq \const\cdot h_j^{-L_j}\cdot h_j^{-1}\cdot\|x_j-x^*_j\|_j
\ea
for any $x_j,x_j^*\in D_j$ and $u_j\in\mathbb{R}^{L_j}$. Hence, the first assertion follows. For the second assertion, we note that (\ref{denominator bound}) implies
\begin{align}\label{nominator bound}
\begin{split}
\sup_{x_j\in D_j}\int_{D_j}K_{h_j}(x_j,u_j)du_j\leq&\const\cdot\sup_{x_j\in D_j}\int_{D_j}h_j^{-L_j}K_j\bigg(\frac{\|x_j-u_j\|_j}{h_j}\bigg)du_j\\
=&\const\cdot\sup_{x_j\in D_j}\int_{(D_j-x_j)/h_j\cap B_j(\bzero_j,1)}K_j(\|s_j\|_j)ds_j\\
\leq&\const\cdot\int_{B_j(\bzero_j,1)}K_j(\|s_j\|_j)ds_j.
\end{split}
\end{align}
This with the boundedness of $K_j(\|\cdot\|_j)$ gives the second assertion. For the third assertion, we note that $n^{-1+2\beta+2/\alpha}\cdot\max_{j\in J}h_j^{-1}\leq\const\cdot n^{-1+2\beta+2/\alpha}\cdot\prod_{j\in J}h_j^{-L_j}=o(1)$. The inequality follows from the condition that $\max_{j\in J}h_j=O(1)$. This verifies the third assertion. Hence, the first part of the lemma follows. The second part of the lemma similarly follows.
\end{proof}

For the next lemma, we recall the definition of $D(n)$ given in (\ref{Dn definition}).

\begin{customlemma}{S.5}\label{uniform rates Dn}
Let $\{\mdelta^i:1\leq i\leq n\}$ be a set of i.i.d. copies of a $\mbH$-valued random element $\mdelta$ such that $\E(\|\mdelta\|^\alpha)<\infty$ for some $2<\alpha<\infty$.
Assume that condition (B3) holds, that $K_j(\|\cdot\|_j):\mathbb{R}^{L_j}\rightarrow[0,\infty)$ are Lipschitz continuous for all $j\in J$, that $\max_{j\in J}h_j=O(1)$, that $n^{-1+2\beta+2/\alpha}\cdot\prod_{j\in J}h_j^{-L_j}=o(1)$ for some constant $\beta>0$, that $n^{-\beta}\cdot(\log{n})^{1/2}\cdot(\sum_{j=1}^da_{nj})^{1/\alpha}\cdot(\sum_{j\in J}h_j^{-L_j}\cdot a_{nj}+\sum_{j\notin J}a_{nj})^{-1/2}=O(1)$ and that $p$ and $\E(\|\mdelta\|^{\alpha}|\mxi=\cdot)$ are bounded on $D(\varepsilon)$ for some $\varepsilon>0$. Then, it holds that
\begin{align*}
\sup_{\mx_J\in\prod_{j\in J}D_j}\|\mS_n(\mx_J)\ominus \E(\mS_n(\mx_J))\|=O_p\left(\sqrt{\frac{(\sum_{j\in J}h_j^{-L_j}\cdot a_{nj}+\sum_{j\notin J}a_{nj})\cdot\log{n}}{n\cdot\prod_{j\in J}h_j^{L_j}}}\right)
\end{align*}
for $\mS_n(\mx_J):=n^{-1}\odot\bigoplus_{i=1}^n(\I(\mxi^i\in D(n))\cdot\prod_{j\in J}h_j^{-L_j}K_j(\|x_j-\xi_j^i\|_j/h_j))\odot\mdelta^i$.
\end{customlemma}

\begin{proof}
We apply Lemma \ref{uniform rates general} with $\mathcal{M}(n)\equiv D(n)$, $L_{nj}(x_j,u_j)=h_j^{-L_j}K_j(\|x_j-u_j\|_j/h_j)$ and $\mdelta_n=\I(\mxi\in D(n))\od\mdelta$ for the proof.
We note that the condition (i) of Lemma \ref{uniform rates general} holds with $\mathcal{M}_j=D_j$, $A_{nj}=h_j^{-L_j}$ and $B_{nj}=h_j^{-1}$ by arguing as in the proof of Lemma \ref{uniform rates}. Hence, it suffices to verify that the following three assertions hold with $U_{n\alpha}=\sum_{j=1}^da_{nj}$ and $V_{nJ}=\sum_{j\in J}h_j^{-L_j}\cdot a_{nj}+\sum_{j\notin J}a_{nj}$: (1) $\E(\|\mdelta_n\|^\alpha)=O(U_{n\alpha})$; (2) the condition (ii) of Lemma \ref{uniform rates general}; (3) $n^c\cdot\max_{j\in J}B_{nj}\cdot U_{n\alpha}^{1/\alpha}\cdot V_{nJ}^{-1/2}=O(1)$ for some constant $c\in\mathbb{R}$.
For the first assertion, we note that
\begin{align*}
\E(\I(\mxi\in D(n))\cdot\|\mdelta\|^\alpha)&=\E(\I(\mxi\in D(n))\cdot\E(\|\mdelta\|^\alpha|\mxi))\\
&=\int_{D(n)}\E(\|\mdelta\|^\alpha|\mxi=\mx)\cdot p(\mx)d\mx\\
&\leq\sup_{\mx\in D(n)}\E(\|\mdelta\|^\alpha|\mxi=\mx)\cdot\sup_{\mx\in D(n)}p(\mx)\cdot\Leb(D(n))\\
&=O\bigg(\sum_{j=1}^da_{nj}\bigg),
\end{align*}
where the last equality follows from (\ref{indicator rate}). This verifies the first assertion. For the second assertion, we note that
\begin{align*}
&\sup_{\mx_J\in \prod_{j\in J}D_j}\int_{D(n)}\E(\I(\mxi\in D(n))\cdot\|\mdelta\|^2|\mxi=\bu)\cdot\prod_{j\in J}h_j^{-2L_j}K^2_j\bigg(\frac{\|x_j-u_j\|_j}{h_j}\bigg)d\bu\\
&=\prod_{j\in J}h_j^{-L_j}\sup_{\mx_J\in \prod_{j\in J}D_j}\int_{D(n)}\E(\|\mdelta\|^2|\mxi=\bu)\cdot\prod_{j\in J}h_j^{-L_j}K^2_j\bigg(\frac{\|x_j-u_j\|_j}{h_j}\bigg)d\bu\\
&\leq\const\cdot \prod_{j\in J}h_j^{-L_j}\cdot\sup_{\mx_J\in \prod_{j\in J}D_j}\int_{D(n)}\prod_{j\in J}h_j^{-L_j}K^2_j\bigg(\frac{\|x_j-u_j\|_j}{h_j}\bigg)d\bu,
\end{align*}
where the inequality follows from the boundedness of $\E(\|\mdelta\|^{\alpha}|\mxi=\cdot)$ on $D(n)$ with H\"{o}lder's inequality. We note that
\begin{align*}
\int_{D(n)}\prod_{j\in J}h_j^{-L_j}K^2_j\bigg(\frac{\|x_j-u_j\|_j}{h_j}\bigg)d\bu\leq&\sum_{k\in J}\int_{D(k,n)}\prod_{j\in J}h_j^{-L_j}K^2_j\bigg(\frac{\|x_j-u_j\|_j}{h_j}\bigg)d\bu\\
&+\sum_{k\notin J}\int_{D(k,n)}\prod_{j\in J}h_j^{-L_j}K^2_j\bigg(\frac{\|x_j-u_j\|_j}{h_j}\bigg)d\bu,
\end{align*}
where $D(k,n)=\big\{\mx\in\prod_{j=1}^dD_j^+(C\cdot a_{nj}):x_k\in \bigcup_{u_k\in \partial D_k}\bar{B}_k(u_k,C\cdot a_{nk})\big\}$.
Recall the definition $D_k^-(C\cdot a_{nk})=D_k\setminus\bigcup_{u_k\in \partial D_k}B_k(u_k,C\cdot a_{nk})$ given in Section \ref{asymptotic properties sbf}.
Then,
\begin{align*}
&\int_{D(k,n)}\prod_{j\in J}h_j^{-L_j}K^2_j\bigg(\frac{\|x_j-u_j\|_j}{h_j}\bigg)d\bu\\
&=\int_{D_k^+(C\cdot a_{nk})\setminus D_k^-(C\cdot a_{nk})}h_k^{-L_k}K^2_k\bigg(\frac{\|x_k-u_k\|_k}{h_k}\bigg)du_k\\
&\quad\cdot\bigg(\prod_{j\in J:j\neq k}\int_{D_j^+(C\cdot a_{nj})}h_j^{-L_j}K^2_j\bigg(\frac{\|x_j-u_j\|_j}{h_j}\bigg)du_j\bigg)\cdot\prod_{j\notin J}\Leb_j(D_j^+(C\cdot a_{nj}))\\
&\leq\const\cdot\int_{(D_k^+(C\cdot a_{nk})\setminus D_k^-(C\cdot a_{nk})-x_k)/h_k\cap B_k(\bzero_k,1)}K^2_k(\|s_k\|_k)ds_k\cdot\prod_{j\in J:j\neq k}\int_{B_j(\bzero_j,1)}K^2_j(\|s_j\|_j)ds_j\\
&\leq\const\cdot\Leb_k((D_k^+(C\cdot a_{nk})\setminus D_k^-(C\cdot a_{nk})-x_k)/h_k)\\
&=\const\cdot\Leb_k((D_k^+(C\cdot a_{nk})\setminus D_k^-(C\cdot a_{nk}))/h_k)
\end{align*}
for $k\in J$, where the first inequality follows by arguing as in (\ref{nominator bound}), the last equality follows from the translation invariance property of the Lebesgue measure and the $\const$ are independent of $\mx_J$. We note that $\Leb_k((D_k^+(C\cdot a_{nk})\setminus D_k^-(C\cdot a_{nk}))/h_k)=O(h_k^{-L_k}\cdot a_{nk})$ by Lemma A4.3 in Kallenberg (2017) and (B3). Similarly,
\begin{align*}
&\int_{D(k,n)}\prod_{j\in J}h_j^{-L_j}K^2_j\bigg(\frac{\|x_j-u_j\|_j}{h_j}\bigg)d\bu\\
&=\Leb_k(D_k^+(C\cdot a_{nk})\setminus D_k^-(C\cdot a_{nk}))\cdot\bigg(\prod_{j\in J}\int_{D_j^+(C\cdot a_{nj})}h_j^{-L_j}K^2_j\bigg(\frac{\|x_j-u_j\|_j}{h_j}\bigg)du_j\bigg)\\
&\quad\cdot\prod_{j\notin J:j\neq k}\Leb_j(D_j^+(C\cdot a_{nj}))\\
&\leq\const\cdot\Leb_k(D_k^+(C\cdot a_{nk})\setminus D_k^-(C\cdot a_{nk}))
\end{align*}
for $k\notin J$, where the $\const$ is independent of $\mx_J$. We note that $\Leb_k(D_k^+(C\cdot a_{nk})\setminus D_k^-(C\cdot a_{nk}))=O(a_{nk})$ by Lemma A4.3 in Kallenberg (2017) and (B3). Hence, we have
\begin{align*}
\sup_{\mx_J\in \prod_{j\in J}D_j}\int_{D(n)}\prod_{j\in J}h_j^{-L_j}K^2_j\bigg(\frac{\|x_j-u_j\|_j}{h_j}\bigg)d\bu=O\bigg(\sum_{j\in J}h_j^{-L_j}\cdot a_{nj}+\sum_{j\notin J}a_{nj}\bigg).
\end{align*}
This gives the second assertion. For the third assertion, we get
\begin{align*}
n^{-1+\beta+2/\alpha}\cdot\bigg(\sum_{j=1}^da_{nj}\bigg)^{1/\alpha}\cdot\bigg(\sum_{j\in J}h_j^{-L_j}\cdot a_{nj}+\sum_{j\notin J}a_{nj}\bigg)^{-1/2}\cdot\prod_{j\in J}h_j^{-L_j}=o(1)
\end{align*}
by multiplying the two $o(1)$ terms $n^{-1+2\beta+2/\alpha}\cdot\prod_{j\in J}h_j^{-L_j}$ and $n^{-\beta}\cdot(\sum_{j=1}^da_{nj})^{1/\alpha}\cdot(\sum_{j\in J}h_j^{-L_j}\cdot a_{nj}+\sum_{j\notin J}a_{nj})^{-1/2}$. Since $\max_{j\in J}h_j^{-1}\leq\const\cdot\prod_{j\in J}h_j^{-L_j}$, the third assertion follows. This completes the proof.
\end{proof}

\begin{customlemma}{S.6}\label{uniform rates Dn2}
Let $\{\mdelta^i:1\leq i\leq n\}$ be identically distributed copies of a $\mbH$-valued random element $\mdelta$ such that $\E(\|\mdelta\|^\alpha)<\infty$ for some $2<\alpha<\infty$.
Assume that condition (B3) holds, that $K_j(\|\cdot\|_j):\mathbb{R}^{L_j}\rightarrow[0,\infty)$ are Lipschitz continuous for all $j\in J$, that $\max_{j\in J}h_j=O(1)$, that $n^{-1+2\beta+2/\alpha}\cdot\prod_{j\in J}h_j^{-L_j}=o(1)$ for some constant $\beta>0$, that $n^{-\beta}\cdot(\log{n})^{1/2}\cdot(\sum_{j=1}^da_{nj})^{1/\alpha}\cdot(\sum_{j\in J}h_j^{-L_j}\cdot a_{nj}+\sum_{j\notin J}a_{nj})^{-1/2}=O(1)$ that $p$ and $\E(\|\mdelta\|^{\alpha}|\mxi=\cdot)$ are bounded on $D$, that $\mdelta^i$ are conditionally independent given $(\mxi^1,\ldots,\mxi^n,\tilde{\mxi}^1,\ldots,\tilde{\mxi}^n)$   and that $\max_{1\leq i\leq n}\E(\|\mdelta^i\|^2|\mxi^1,\ldots,\mxi^n,\tilde{\mxi}^1,\ldots,\tilde{\mxi}^n)<C$ almost surely for a constant $C>0$. Then, it holds that
\begin{align*}
&\sup_{\mx_J\in\prod_{j\in J}D_j}\|\mS_n(\mx_J)\ominus \E(\mS_n(\mx_J)|\mxi^1,\ldots,\mxi^n,\tilde{\mxi}^1,\ldots,\tilde{\mxi}^n)\|\\
&=O_p\left(\sqrt{\frac{(\sum_{j\in J}h_j^{-L_j}\cdot a_{nj}+\sum_{j\notin J}a_{nj})\cdot\log{n}}{n\cdot\prod_{j\in J}h_j^{L_j}}}\right)
\end{align*}
for $\mS_n(\mx_J):=n^{-1}\odot\bigoplus_{i=1}^n(\I(\mxi^i\in D)\cdot(\I(\tilde{\mxi}^i\in D)-\I(\mxi^i\in D))\cdot\prod_{j\in J}K_{h_j}(x_j,\xi_j^i))\odot\mdelta^i$.
\end{customlemma}

\begin{proof}
By arguing as in the proofs of Lemmas \ref{uniform rates general} and \ref{uniform rates Dn} with $\mathcal{M}(n)\equiv D$, $L_{nj}=K_{h_j}$ and $\mdelta_n=(\I(\mxi\in D)\cdot(\I(\tilde{\mxi}\in D)-\I(\mxi\in D)))\od\mdelta$, we get
\begin{align*}
&\mS_n(\mx_J)\om \E(\mS_n(\mx_J)|\mxi^1,\ldots,\mxi^n,\tilde{\mxi}^1,\ldots,\tilde{\mxi}^n)\\
&=\max_{1\leq l\leq N(n^{-\kappa})}\|\mU_n(\mx_J^{(l)})\om \E(\mU_n(\mx_J^{(l)})|\mxi^1,\ldots,\mxi^n,\tilde{\mxi}^1,\ldots,\tilde{\mxi}^n)\|\\
&\quad\op o_p\left(\sqrt{\frac{\log{n}\cdot \big(\sum_{j\in J}h_j^{-L_j}\cdot a_{nj}+\sum_{j\notin J}a_{nj}\big)}{n\cdot\prod_{j\in J}h_j^{L_j}}}\right),
\end{align*}
where $N(n^{-\kappa})$ and $\{\mx_J^{(1)},\cdots,\mx_J^{(N(n^{-\kappa}))}\}\subset \prod_{j\in J}D_j$ are defined as in the proof of Lemma \ref{uniform rates general}, and
\begin{align*}
&\mU_n(\mx_J)\\
&=\frac{1}{n}\od\bigoplus_{i=1}^n\bigg(\I(\mxi^i\in D)\cdot\prod_{j\in J}K_{h_j}(x_j,\xi_j^i)\\
&\qquad\qquad\cdot\I\bigg(\|(\I(\mxi^i\in D)\cdot(\I(\tilde{\mxi}^i\in D)-\I(\mxi^i\in D)))\od\mdelta^i\|\leq n^{1/2-\beta}\bigg(\sum_{j=1}^da_{nj}\bigg)^{1/\alpha}\prod_{j\in J}h_j^{L_j/2}\bigg)\\
&\qquad\qquad\cdot\I(\mxi^i\in D)\cdot(\I(\tilde{\mxi}^i\in D)-\I(\mxi^i\in D))\bigg)\od\mdelta^i.
\end{align*}
We define $\mZ_n^i(\mx_J)$ so that $\mU_n(\mx_J)\om \E(\mU_n(\mx_J)|\mxi^1,\ldots,\mxi^n,\tilde{\mxi}^1,\ldots,\tilde{\mxi}^n)=\bigoplus_{i=1}^n\mZ_n^i(\mx_J)$. Since $\E(\mZ_n^i(\mx_J)|\mxi^1,\ldots,\mxi^n,\tilde{\mxi}^1,\ldots,\tilde{\mxi}^n)=\bzero$ and $\mdelta^i$ are conditionally independent, the conditional versions of (2.18) and (2.74) in Bosq (2000) give
\begin{align}\label{conditional concentration}
\begin{split}
&\Prob\left(\bigg\|\bigoplus_{i=1}^n\mZ_n^i(\mx_J)\bigg\|>T_n\bigg|\mxi^1,\ldots,\mxi^n,\tilde{\mxi}^1,\ldots,\tilde{\mxi}^n\right)\\
&\leq 2\exp(-\lambda_n\cdot T_n)\cdot\prod_{i=1}^n\E\big(\exp(\lambda_n\cdot\|\mZ_n^i(\mx_J)\|)-\lambda_n\cdot\|\mZ_n^i(\mx_J)\||\mxi^1,\ldots,\mxi^n,\tilde{\mxi}^1,\ldots,\tilde{\mxi}^n\big)\\
&=2n^{-C_0}\cdot\prod_{i=1}^n\E\big(1+\lambda_n^2\cdot\|\mZ_n^i(\mx_J)\|^2\cdot\exp(\tilde{\mZ}_n^i(\mx_J))|\mxi^1,\ldots,\mxi^n,\tilde{\mxi}^1,\ldots,\tilde{\mxi}^n\big),
\end{split}
\end{align}
where
\ba
T_n=C_0\cdot\sqrt{\frac{\log{n}\cdot \big(\sum_{j\in J}h_j^{-L_j}\cdot a_{nj}+\sum_{j\notin J}a_{nj}\big)}{n\cdot\prod_{j\in J}h_j^{L_j}}}
\ea
for a constant $C_0>0$,
\ba
\lambda_n=\sqrt{\frac{n\cdot\prod_{j\in J}h_j^{L_j}\cdot \log{n}}{\sum_{j\in J}h_j^{-L_j}\cdot a_{nj}+\sum_{j\notin J}a_{nj}}}
\ea
and $\tilde{\mZ}^i_n(\mx_J)$ is a $\mbH$-valued random element satisfying $\|\tilde{\mZ}^i_n(\mx_J)\|\leq\lambda_n\cdot \|\mZ_n^i(\mx_J)\|$. The equality in (\ref{conditional concentration}) follows from the fact that $\exp(y)=1+y+y^2\exp(y^*)/2$ for some $y^*$ with $|y^*|\leq |y|$. We note that
\ba
\|\mZ_n^i(\mx_J)\|\leq&\const\cdot n^{-1/2-\beta}\cdot \bigg(\sum_{j=1}^da_{nj}\bigg)^{1/\alpha}\cdot\prod_{j\in J}h_j^{-L_j/2}.
\ea
Hence,
\ba
&\E\big(1+\lambda_n^2\cdot\|\mZ_n^i(\mx_J)\|^2\cdot\exp(\tilde{\mZ}_n^i(\mx_J))|\mxi^1,\ldots,\mxi^n,\tilde{\mxi}^1,\ldots,\tilde{\mxi}^n\big)\\
&\leq1+\lambda_n^2\cdot\E(\|\mZ_n^i(\mx_J)\|^2|\mxi^1,\ldots,\mxi^n,\tilde{\mxi}^1,\ldots,\tilde{\mxi}^n)\\
&\qquad\qquad\cdot\exp\bigg(\const\cdot \lambda_n\cdot n^{-1/2-\beta}\cdot \bigg(\sum_{j=1}^da_{nj}\bigg)^{1/\alpha}\cdot\prod_{j\in J}h_j^{-L_j/2}\bigg)\\
&\leq1+O(1)\cdot \lambda_n^2\cdot n^{-2}\cdot\I(\mxi^i\in D\cap D(n))\cdot\E(\|\mdelta^i\|^2|\mxi^1,\ldots,\mxi^n,\tilde{\mxi}^1,\ldots,\tilde{\mxi}^n)\cdot\prod_{j\in J}(K_{h_j}(x_j,\xi_j^i))^2\\
&\leq1+O(1)\cdot \lambda_n^2\cdot n^{-2}\cdot\I(\mxi^i\in D\cap D(n))\cdot\prod_{j\in J}h_j^{-2L_j}K^2_j\bigg(\frac{\|x_j-\xi_j^i\|_j}{h_j}\bigg)
\ea
almost surely, where the second inequality follows from that $n^{-\beta}\cdot(\log{n})^{1/2}\cdot(\sum_{j=1}^da_{nj})^{1/\alpha}\cdot(\sum_{j\in J}h_j^{-L_j}\cdot a_{nj}+\sum_{j\notin J}a_{nj})^{-1/2}=O(1)$ and the last inequality follows from the boundedness of $\max_{1\leq i\leq n}\E(\|\mdelta^i\|^2|\mxi^1,\ldots,\mxi^n,\tilde{\mxi}^1,\ldots,\tilde{\mxi}^n)$ and (\ref{denominator bound}). Thus,
\ba
&\Prob\left(\bigg\|\bigoplus_{i=1}^n\mZ_n^i(\mx_J)\bigg\|>T_n\right)\\
&\leq2n^{-C_0}\cdot \E\bigg(\prod_{i=1}^n\bigg(1+O(1)\cdot \lambda_n^2\cdot n^{-2}\cdot\I(\mxi^i\in D\cap D(n))\cdot\prod_{j\in J}h_j^{-2L_j}K^2_j\bigg(\frac{\|x_j-\xi_j^i\|_j}{h_j}\bigg)\bigg)\bigg)\\
&\leq2n^{-C_0}\cdot \bigg(\E\bigg(1+O(1)\cdot \lambda_n^2\cdot n^{-2}\cdot\I(\mxi\in D\cap D(n))\cdot\prod_{j\in J}h_j^{-2L_j}K^2_j\bigg(\frac{\|x_j-\xi_j\|_j}{h_j}\bigg)\bigg)\bigg)^n.
\ea
It holds that
\ba
&\E\bigg(\I(\mxi\in D\cap D(n))\cdot\prod_{j\in J}h_j^{-2L_j}K^2_j\bigg(\frac{\|x_j-\xi_j\|_j}{h_j}\bigg)\bigg)\\
&=O\bigg(\prod_{j\in J}h_j^{-L_j}\bigg(\sum_{j\in J}h_j^{-L_j}\cdot a_{nj}+\sum_{j\notin J}a_{nj}\bigg)\bigg)
\ea
by arguing as in the proof of Lemma \ref{uniform rates Dn}. Therefore, we have
\ba
\Prob\left(\bigg\|\bigoplus_{i=1}^n\mZ_n^i(\mx_J)\bigg\|>T_n\right)&\leq2n^{-C_0}\cdot(1+O(1)\cdot n^{-1}\cdot\log{n})^n\\
&\leq2n^{-C_0}\cdot\exp(O(1)\cdot\log{n})\\
&=2n^{-C_0}\cdot n^{O(1)}.
\ea
This with the fact that $N(n^{-\kappa})=O(n^{\kappa m})$ for some $m>0$ completes the proof.
\end{proof}

In the next two lemmas, we consider the case where there is no $\mdelta$ or $\delta$ term.

\begin{customlemma}{S.7}\label{uniform rates density}
Assume that condition (B3) holds, that $K_j(\|\cdot\|_j):\mathbb{R}^{L_j}\rightarrow[0,\infty)$ are Lipschitz continuous for all $j\in J$, that $\max_{j\in J}h_j=O(1)$, that $n^{-1/2}\cdot(\log n)^{1/2}\cdot\prod_{j\in J}h_j^{-L_j/2}=O(1)$ and that $p$ is bounded on $D$. Then, it holds that
\begin{align*}
\sup_{\mx_J\in\prod_{j\in J}D_j}|S_n(\mx_J)-\E(S_n(\mx_J))|=O_p\left(\sqrt{\frac{\log{n}}{n\cdot\prod_{j\in J}h_j^{L_j}}}\right)
\end{align*}
for $S_n(\mx_J):=n^{-1}\sum_{i=1}^n\I(\mxi^i\in D)\cdot\prod_{j\in J}K_{h_j}(x_j,\xi_j^i)$. The same uniform rate of convergence holds for $S^*_n(\mx_J):=n^{-1}\sum_{i=1}^n\I(\mxi^i\in D)\cdot\prod_{j\in J}h_j^{-L_j}K_j(\|x_j-\xi_j^i\|_j/h_j)$ without condition (B3).
\end{customlemma}

\begin{proof}
We prove that $n^c\cdot\max_{j\in J}h_j^{-1}\cdot\prod_{j\in J}h_j^{-L_j/2}=O(1)$ for some constant $c\in\mathbb{R}$. Then, the lemma follows by applying Remark \ref{uniform rate remark}-3 and arguing as in the proof of Lemma \ref{uniform rates}. It holds that $\prod_{j\in J}h_j^{-L_j/2}=o(n^{1/2})$ from the condition that $n^{-1/2}\cdot(\log n)^{1/2}\cdot\prod_{j\in J}h_j^{-L_j/2}=O(1)$. Since $\max_{j\in J}h_j^{-1}\leq\const\cdot\prod_{j\in J}h_j^{-L_j}$, we have $\max_{j\in J}h_j^{-1}\cdot\prod_{j\in J}h_j^{-L_j/2}=o(n^{3/2})$. This completes the proof.
\end{proof}

\begin{customlemma}{S.8}\label{uniform rates density Dn}
Assume that condition (B3) holds, that $K_j(\|\cdot\|_j):\mathbb{R}^{L_j}\rightarrow[0,\infty)$ are Lipschitz continuous for all $j\in J$, that $\max_{j\in J}h_j=O(1)$, that $n^{-1/2}\cdot(\log n)^{1/2}\cdot\prod_{j\in J}h_j^{-L_j/2}=O(1)$, that $\sum_{j\in J}h_j^{-L_j}\cdot a_{nj}=O(1)$, that $n^{-1/2}\cdot(\log n)^{1/2}\cdot(\sum_{j\in J}h_j^{-L_j}\cdot a_{nj}+\sum_{j\notin J}a_{nj})^{-1/2}\cdot\prod_{j\in J}h_j^{-L_j/2}=O(1)$ and that $p$ is bounded on $D$. Then, it holds that
\begin{align*}
\sup_{\mx_J\in\prod_{j\in J}D_j}|S_n(\mx_J)-\E(S_n(\mx_J))|=O_p\left(\sqrt{\frac{(\sum_{j\in J}h_j^{-L_j}\cdot a_{nj}+\sum_{j\notin J}a_{nj})\cdot\log{n}}{n\cdot\prod_{j\in J}h_j^{L_j}}}\right)
\end{align*}
for $S_n(\mx_J):=n^{-1}\sum_{i=1}^n\I(\mxi^i\in D\cap D(n))\cdot\prod_{j\in J}K_{h_j}(x_j,\xi_j^i)$. If we further assume that $p$ is bounded on $D(\varepsilon)$ for some $\varepsilon>0$, then the same uniform rate of convergence holds for $S^*_n(\mx_J):=n^{-1}\sum_{i=1}^n\I(\mxi^i\in D(n))\cdot\prod_{j\in J}h_j^{-L_j}K_j(\|x_j-\xi_j^i\|_j/h_j)$.
\end{customlemma}

\begin{proof}
We prove that $n^c\cdot\max_{j\in J}h_j^{-1}\cdot(\sum_{j\in J}h_j^{-L_j}\cdot a_{nj}+\sum_{j\notin J}a_{nj})^{-1/2}\cdot\prod_{j\in J}h_j^{-L_j/2}=O(1)$ for some constant $c\in\mathbb{R}$. Then, the lemma follows by applying Lemma \ref{uniform rates general density} and arguing as in the proof of Lemma \ref{uniform rates Dn}. It holds that $n^{-1/2}\cdot\max_{j\in J}h_j^{-L_j/2}=o(1)$ from the conditions that $n^{-1/2}\cdot(\log n)^{1/2}\cdot(\sum_{j\in J}h_j^{-L_j}\cdot a_{nj}+\sum_{j\notin J}a_{nj})^{-1/2}\cdot\prod_{j\in J}h_j^{-L_j/2}=O(1)$, that $\max_{j\in J}h_j=O(1)$ and that $\sum_{j\in J}h_j^{-L_j}\cdot a_{nj}=O(1)$. Since $n^{-1/2}\cdot\max_{j\in J}h_j^{-1/2}\leq\const\cdot n^{-1/2}\cdot\max_{j\in J}h_j^{-L_j/2}$, we have $n^{-1}\cdot\max_{j\in J}h_j^{-1}=o(1)$. By multiplying this term and $n^{-1/2}\cdot(\sum_{j\in J}h_j^{-L_j}\cdot a_{nj}+\sum_{j\notin J}a_{nj})^{-1/2}\cdot\prod_{j\in J}h_j^{-L_j/2}=o(1)$, we get $n^{-3/2}\cdot\max_{j\in J}h_j^{-1}\cdot(\sum_{j\in J}h_j^{-L_j}\cdot a_{nj}+\sum_{j\notin J}a_{nj})^{-1/2}\cdot\prod_{j\in J}h_j^{-L_j/2}=o(1)$. This completes the proof.
\end{proof}


\begin{customremark}{S.2}
Lemmas \ref{uniform rates density} and \ref{uniform rates density Dn} are based on the conditions in Lemma \ref{uniform rates general density} not on the alternative condition in Remark \ref{uniform rate remark}-2 since the latter condition is stronger than the former conditions as long as $(n\cdot \sum_{j=1}^da_{nj})^{-1}=O(1)$.
\end{customremark}

The next lemma is useful to get the uniform consistency of our marginal density estimators. To state the lemma, we define $p^D_J(\mx_J)=\int_{\prod_{j\notin J}D_j}p^D(\mx)d\mx_{-J}$ for $\mx_{-J}=(x_j:j\notin J)$.

\begin{customlemma}{S.9}\label{uniform consistency general}
Assume that condition (B3) holds, that $\max_{j\in J}h_j=o(1)$ and that $p^D_J$ is continuous on $\prod_{j\in J}D_j$. Then, it holds that
\ba
\sup_{\mx_J\in\prod_{j\in J}D_j}\bigg|\E\bigg(n^{-1}\sum_{i=1}^n\prod_{j\in J}K_{h_j}(x_j,\xi_j^i)\I(\mxi^i\in D)\bigg)-p^D_0p^D_J(\mx_J)\prod_{j\in J}\int_{D_j}K_{h_j}(x_j,u_j)du_j\bigg|&=o(1),\\
\sup_{\mx_J\in\prod_{j\in J}D^-_j(2h_j)}\bigg|\E\bigg(n^{-1}\sum_{i=1}^n\prod_{j\in J}K_{h_j}(x_j,\xi_j^i)\I(\mxi^i\in D)\bigg)-p^D_0p^D_J(\mx_J)\bigg|&=o(1).
\ea
\end{customlemma}

\begin{proof}
We first prove the first assertion. We note that
\ba
\E\bigg(n^{-1}\sum_{i=1}^n\prod_{j\in J}K_{h_j}(x_j,\xi_j^i)\I(\mxi^i\in D)\bigg)&=\int_Dp(\bu)\prod_{j\in J}K_{h_j}(x_j,u_j)d\bu\\
&=p_0^D\cdot\int_{\prod_{j\in J}D_j}p^D_J(\bu_J)\prod_{j\in J}K_{h_j}(x_j,u_j)d\bu_J.
\ea
We also note that the continuity of $p$ on $D$ implies the continuity of $p^D_J$ on $\prod_{j\in J}D_j$. Hence,
\ba
&\sup_{\mx_J\in\prod_{j\in J}D_j}\bigg|\E\bigg(n^{-1}\sum_{i=1}^n\prod_{j\in J}K_{h_j}(x_j,\xi_j^i)\I(\mxi^i\in D)\bigg)-p_0^D\cdot p^D_J(\mx_J)\cdot\prod_{j\in J}\int_{D_j}K_{h_j}(x_j,u_j)du_j\bigg|\\
&=p_0^D\cdot\sup_{\mx_J\in\prod_{j\in J}D_j}\bigg|\int_{\prod_{j\in J}D_j}(p^D_J(\bu_J)-p^D_J(\mx_J))\prod_{j\in J}K_{h_j}(x_j,u_j)d\bu_J\bigg|\\
&=p_0^D\cdot\sup_{\mx_J\in\prod_{j\in J}D_j}\bigg|\int_{\prod_{j\in J}D_j\cap B_j(x_j,h_j)}(p^D_J(\bu_J)-p^D_J(\mx_J))\prod_{j\in J}K_{h_j}(x_j,u_j)d\bu_J\bigg|\\
&\leq p_0^D\cdot\sup_{\bu_J,\bu^*_J\in\prod_{j\in J}D_j:\sum_{j\in J}\|u_j-u^*_j\|_j^2<\sum_{j\in J}h_j^2}|p^D_J(\bu_J)-p^D_J(\bu^*_J)|\cdot\prod_{j\in J}\sup_{x_j\in D_j}\int_{D_j}K_{h_j}(x_j,u_j)du_j\\
&=o(1),
\ea
where the last equality follows from the continuity of $p^D_J$ on the compact set $\prod_{j\in J}D_j$ and (\ref{nominator bound}). For the second assertion, it suffices to prove that $\int_{D_j}K_{h_j}(x_j,u_j)du_j=1$ for $x_j\in D^-_j(2h_j)$ and sufficiently large $n$. We note that, if $x_j\in D^-_j(2h_j)$ and $u_j\in B_j(x_j,h_j)$, then $B(\bzero_j,1)\subset(D_j-x_j)/h_j$ and $B(\bzero_j,1)\subset(D_j-u_j)/h_j$ for sufficiently large $n$. Hence,
\begin{align}\label{integral one in interior}
\begin{split}
\int_{D_j}\frac{1}{h_j^{L_j}}K_j\bigg(\frac{\|t_j-u_j\|_j}{h_j}\bigg)dt_j=&\int_{(D_j-u_j)/h_j\cap B(\bzero_j,1)}K_j(\|s_j\|_j)ds_j\\
=&\int_{B(\bzero_j,1)}K_j(\|s_j\|_j)ds_j\in(0,\infty)\\
\int_{D_j}\frac{1}{h_j^{L_j}}K_j\bigg(\frac{\|x_j-u_j\|_j}{h_j}\bigg)du_j=&\int_{(D_j-x_j)/h_j\cap B(\bzero_j,1)}K_j(\|s_j\|_j)ds_j\\
=&\int_{B(\bzero_j,1)}K_j(\|s_j\|_j)ds_j\in(0,\infty)
\end{split}
\end{align}
for $x_j\in D^-_j(2h_j)$ and sufficiently large $n$. This completes the proof.
\end{proof}

The next lemma provides some properties on the partial derivatives of $K_j(\|\cdot\|_j)$. The lemma is frequently used in the proofs of Lemmas \ref{marginal density approximation 1}, \ref{marginal density approximation 2} and \ref{marginal regression approximation} to be presented and the theorems in this paper without referring it. Below, $u_{jl}$ denotes the $l$th element of $u_j$.

\begin{customlemma}{S.10}\label{K derivative}
Assume that conditions (B3) and (B4) hold. Then, for all $1\leq j\leq d$, it holds that
\begin{align*}
\sup_{(x_j,u_j)\in\mathbb{R}^{L_j}\times\mathbb{R}^{L_j}}\bigg|\frac{\partial}{\partial u_{jl}}K_j\left(\frac{\|x_j-u_j\|_j}{h_j}\right)\bigg|=O(h_j^{-1}).
\end{align*}
Also,
\begin{align*}
\sup_{u_j\in\mathbb{R}^{L_j}}\bigg|\bigg|\frac{\partial}{\partial u_{jl}}K_j\left(\frac{\|x_j-u_j\|_j}{h_j}\right)\bigg|-\bigg|\frac{\partial}{\partial u_{jl}}K_j\left(\frac{\|x^*_j-u_j\|_j}{h_j}\right)\bigg|\bigg|\leq O(h_j^{-2})\cdot\|x_j-x_j^*\|_j
\end{align*}
for all $x_j, x_j^*\in\mathbb{R}^{L_j}$. In addition,
\ba
\sup_{x_j\in D_j}\int_{D_j^+(C\cdot a_{nj})}\bigg|\frac{\partial}{\partial u_{jl}}K_j\bigg(\frac{\|x_j-u_j\|_j}{h_j}\bigg)\bigg|^mdu_j&=O(h_j^{L_j-m}),\\
\sup_{x_j\in D_j}\int_{\bigcup_{z_j\in \partial D_j}\bar{B}_j(z_j,C\cdot a_{nj})}\bigg|\frac{\partial}{\partial u_{jl}}K_j\bigg(\frac{\|x_j-u_j\|_j}{h_j}\bigg)\bigg|^mdu_j&=O(h_j^{-m}\cdot a_{nj}).
\ea
for any $m\geq 1$, where $C>0$ is the constant in (\ref{Op event}). Moreover,
\ba
\sup_{u_j\in\mathbb{R}^{L_j}}\int_{D_j}\bigg|\frac{\partial}{\partial u_{jl}}K_j\bigg(\frac{\|x_j-u_j\|_j}{h_j}\bigg)\bigg|dx_j=O(h_j^{L_j-1}).
\ea
\end{customlemma}

\begin{proof}
We define $K_j^{1,l}(x_j)=\partial K_j(\|u_j\|_j)/\partial u_{jl}|_{u_j=x_j}$ and $g_{x_j,h_j}(u_j)=(x_j-u_j)/h_j$. Then, $\partial K_j(\|x_j-u_j\|_j/h_j)/\partial u_{jl}=-h_j^{-1}\cdot K_j^{1,l}(g_{x_j,h_j}(u_j))$. Hence,
\begin{align*}
\sup_{(x_j,u_j)\in\mathbb{R}^{L_j}\times\mathbb{R}^{L_j}}\bigg|\frac{\partial}{\partial u_{jl}}K_j\left(\frac{\|x_j-u_j\|_j}{h_j}\right)\bigg|\leq\const\cdot h_j^{-1}
\end{align*}
by the boundedness of $K_j^{1,l}$. Also, the Lipschitz continuity of $K_j^{1,l}$ implies the Lipschitz continuity of $|K_j^{1,l}|$ and the latter implies that
\begin{align*}
&\sup_{u_j\in\mathbb{R}^{L_j}}\bigg|\bigg|\frac{\partial}{\partial u_{jl}}K_j\left(\frac{\|x_j-u_j\|_j}{h_j}\right)\bigg|-\bigg|\frac{\partial}{\partial u_{jl}}K_j\left(\frac{\|x^*_j-u_j\|_j}{h_j}\right)\bigg|\bigg|\\
&=h_j^{-1}\cdot\sup_{u_j\in\mathbb{R}^{L_j}}||K_j^{1,l}(g_{x_j,h_j}(u_j))|-|K_j^{1,l}(g_{x^*_j,h_j}(u_j))||\\
&\leq\const\cdot h_j^{-1}\cdot\sup_{u_j\in\mathbb{R}^{L_j}}\|g_{x_j,h_j}(u_j)-g_{x^*_j,h_j}(u_j)\|_j\\
&=\const\cdot h_j^{-2}\cdot\|x_j-x_j^*\|_j
\end{align*}
for all $x_j, x_j^*\in\mathbb{R}^{L_j}$. Also,
\begin{align*}
&\sup_{x_j\in D_j}\int_{D_j^+(C\cdot a_{nj})}\bigg|\frac{\partial}{\partial u_{jl}}K_j\bigg(\frac{\|x_j-u_j\|_j}{h_j}\bigg)\bigg|^mdu_j\\
&=h_j^{L_j-m}\cdot\sup_{x_j\in D_j}\int_{D_j^+(C\cdot a_{nj})}h_j^{-L_j}|K_j^{1,l}(g_{x_j,h_j}(u_j))|^mdu_j\\
&=h_j^{L_j-m}\cdot\sup_{x_j\in D_j}\int_{(D_j^+(C\cdot a_{nj})-x_j)/h_j\cap B_j(\bzero_j,1)}|K^{1,l}_j(s_j)|^mds_j\\
&\leq h_j^{L_j-m}\cdot\int_{B_j(\bzero_j,1)}|K^{1,l}_j(s_j)|^mds_j\\
&\leq\const\cdot h_j^{L_j-m}
\end{align*}
for any $m\geq 1$. Also,
\begin{align*}
&\sup_{x_j\in D_j}\int_{\bigcup_{z_j\in \partial D_j}\bar{B}_j(z_j,C\cdot a_{nj})}\bigg|\frac{\partial}{\partial u_{jl}}K_j\bigg(\frac{\|x_j-u_j\|_j}{h_j}\bigg)\bigg|^mdu_j\\
&=h_j^{L_j-m}\cdot\sup_{x_j\in D_j}\int_{\bigcup_{z_j\in \partial D_j}\bar{B}_j(z_j,C\cdot a_{nj})}h_j^{-L_j}|K_j^{1,l}(g_{x_j,h_j}(u_j))|^mdu_j\\
&=h_j^{L_j-m}\cdot\sup_{x_j\in D_j}\int_{(\bigcup_{z_j\in \partial D_j}\bar{B}_j(z_j,C\cdot a_{nj})-x_j)/h_j\cap B_j(\bzero_j,1)}|K^{1,l}_j(s_j)|^mds_j\\
&\leq\const\cdot h_j^{L_j-m}\cdot\Leb_j\bigg(\bigg(\bigcup_{z_j\in \partial D_j}\bar{B}_j(z_j,C\cdot a_{nj})-x_j\bigg)/h_j\bigg)\\
&=\const\cdot h_j^{L_j-m}\cdot\Leb_j\bigg(\bigg(\bigcup_{z_j\in \partial D_j}\bar{B}_j(z_j,C\cdot a_{nj})\bigg)/h_j\bigg)\\
&\leq\const\cdot h_j^{-m}\cdot a_{nj}
\end{align*}
for any $m\geq 1$, where the last equality follows from the translation invariance property of the Lebesgue measure. Also,
\ba
&\sup_{u_j\in\mathbb{R}^{L_j}}\int_{D_j}\bigg|\frac{\partial}{\partial u_{jl}}K_j\bigg(\frac{\|x_j-u_j\|_j}{h_j}\bigg)\bigg|dx_j\\
&=h_j^{L_j-1}\cdot\sup_{u_j\in\mathbb{R}^{L_j}}\int_{D_j}h_j^{-L_j}|K_j^{1,l}(g_{x_j,h_j}(u_j))|dx_j\\
&=h_j^{L_j-1}\cdot\sup_{u_j\in\mathbb{R}^{L_j}}\int_{(D_j-u_j)/h_j\cap B_j(\bzero_j,1)}|K_j^{1,l}(s_j)|ds_j\\
&\leq\const\cdot h_j^{L_j-1}.
\ea
This completes the proof.
\end{proof}

In the next two lemmas, we derive the uniform rates of convergence between the marginal density estimators based on $\tilde{\mxi}^i$ and those based on $\mxi^i$. We only use the first parts of the lemmas throughout the paper but the second parts are also of interest in their own right.

\begin{customlemma}{S.11}\label{marginal density approximation 1}
Assume that conditions (B3) and (B4) hold, that $h_j=O(1)$, that $n^{-1/2}\cdot(\log{n})^{1/2}\cdot h_j^{-L_j/2}=O(1)$, that $h_j^{-1}\cdot a_{nj}=O(1)$, that $h_j^{-L_j-2}\cdot a^2_{nj}=O(1)$ and that $p$ is bounded on $D^+(\varepsilon)$ for some $\varepsilon>0$. Then, it holds that
\ba
\sup_{x_j\in D_j}|\hat{p}^D_j(x_j)-\check{p}^D_j(x_j)|=O_p\left(\sqrt{\frac{\log{n}}{n\cdot h_j^{L_j}}}+h_j^{-L_j}\cdot a_{nj}+\sum_{k\neq j}a_{nk}+h_j^{-L_j-2}\cdot a_{nj}^2\right).
\ea
If we further assume that $n^{-1/2}\cdot(\log{n})^{1/2}\cdot(h_j^{-L_j}\cdot a_{nj}+\sum_{k\neq j}a_{nk})^{-1/2}\cdot h_j^{-L_j/2}=O(1)$, then it holds that
\ba
\sup_{x_j\in D_j}|\hat{p}^D_j(x_j)-\check{p}^D_j(x_j)|=O_p\bigg(h_j^{-L_j}\cdot a_{nj}+\sum_{k\neq j}a_{nk}+h_j^{-L_j-2}\cdot a_{nj}^2\bigg).
\ea
\end{customlemma}

\begin{proof}
We first prove the second part of the lemma. We note that
\begin{align}\label{pj approximation}
\begin{split}
\sup_{x_j\in D_j}|\hat{p}^D_j(x_j)-\check{p}^D_j(x_j)|&\leq\frac{|\hat{p}_0^D-\check{p}_0^D|}{\hat{p}_0^D\check{p}_0^D}\sup_{x_j\in D_j}\bigg|\frac{1}{n}\sum_{i=1}^nK_{h_j}(x_j,\xi_j^i)\I(\mxi^i\in D)\bigg|\\
&\quad+\frac{1}{\hat{p}_0^D}\sup_{x_j\in D_j}\bigg|\frac{1}{n}\sum_{i=1}^n(K_{h_j}(x_j,\tilde{\xi}_j^i)\I(\tilde{\mxi}^i\in D)-K_{h_j}(x_j,\xi_j^i)\I(\mxi^i\in D))\bigg|.
\end{split}
\end{align}
We first approximate the first term on the right hand side of (\ref{pj approximation}). For this, we approximate $\sup_{x_j\in D_j}|n^{-1}\sum_{i=1}^nK_{h_j}(x_j,\xi_j^i)\I(\mxi^i\in D)|$.
By Lemma \ref{uniform rates density}, it holds that
\ba
\sup_{x_j\in D_j}\bigg|\frac{1}{n}\sum_{i=1}^nK_{h_j}(x_j,\xi_j^i)\I(\mxi^i\in D)-\E(K_{h_j}(x_j,\xi_j)\I(\mxi\in D))\bigg|=O_p(n^{-1/2}\cdot(\log{n})^{1/2}\cdot h_j^{-L_j/2})
\ea
provided that $n^{-1/2}\cdot(\log{n})^{1/2}\cdot h_j^{-L_j/2}=O(1)$.
Also,
\ba
\sup_{x_j\in D_j}\E(K_{h_j}(x_j,\xi_j)\I(\mxi\in D))&=\sup_{x_j\in D_j}\int_DK_{h_j}(x_j,u_j)p(\bu)d\bu\\
&=p_0^D\cdot\sup_{x_j\in D_j}\int_{D_j}K_{h_j}(x_j,u_j)p_j^D(u_j)du_j\\
&\leq\const\cdot\sup_{x_j\in D_j}\int_{D_j}h_j^{-L_j}K_j\bigg(\frac{\|x_j-u_j\|_j}{h_j}\bigg)du_j\\
&\leq\const\cdot\int_{B_j(\bzero_j,1)}K_j(\|s_j\|_j)ds_j\\
&=O(1),
\ea
where the first inequality follows from (\ref{denominator bound}) and the boundedness of $p$. Hence, we have
\ba
\sup_{x_j\in D_j}\bigg|\frac{1}{n}\sum_{i=1}^nK_{h_j}(x_j,\xi_j^i)\I(\mxi^i\in D)\bigg|=O_p(1).
\ea
Since $\hat{p}_0^D$ and $\check{p}_0^D$ are bounded away from zero with probability tending to one by Lemma \ref{p0 rate} and its proof, the first term on the right hand side of (\ref{pj approximation}) has the same rate with $|\hat{p}_0^D-\check{p}_0^D|$, which is
\ba
O_p\bigg(\sum_{j=1}^da_{nj}+n^{-1/2}\cdot\bigg(\sum_{j=1}^da_{nj}\bigg)^{1/2}\bigg)
\ea
by Lemma \ref{p0 rate}.

We now approximate the second term on the right hand side of (\ref{pj approximation}). We note that
\begin{align}\label{decomposition j}
\begin{split}
&\sup_{x_j\in D_j}\bigg|\frac{1}{n}\sum_{i=1}^n(K_{h_j}(x_j,\tilde{\xi}_j^i)\I(\tilde{\mxi}^i\in D)-K_{h_j}(x_j,\xi_j^i)\I(\mxi^i\in D))\bigg|\\
&\leq\sup_{x_j\in D_j}\bigg|\frac{1}{n}\sum_{i=1}^nK_{h_j}(x_j,\xi_j^i)\I(\mxi^i\in D)(\I(\tilde{\mxi}^i\in D)-\I(\mxi^i\in D))\bigg|\\
&\quad+\sup_{x_j\in D_j}\bigg|\frac{1}{n}\sum_{i=1}^n(K_{h_j}(x_j,\tilde{\xi}_j^i)\I(\tilde{\mxi}^i\in D)-K_{h_j}(x_j,\xi_j^i)\I(\mxi^i\in D))\I(\tilde{\mxi}^i\in D)\bigg|\\
&\leq\sup_{x_j\in D_j}\bigg|\frac{1}{n}\sum_{i=1}^nK_{h_j}(x_j,\xi_j^i)\I(\mxi^i\in D)(\I(\tilde{\mxi}^i\in D)-\I(\mxi^i\in D))\bigg|\\
&\quad+\sup_{x_j\in D_j}\bigg|\frac{1}{n}\sum_{i=1}^nK_{h_j}(x_j,\tilde{\xi}_j^i)\I(\tilde{\mxi}^i\in D)(\I(\tilde{\mxi}^i\in D)-\I(\mxi^i\in D))\bigg|\\
&\quad+\sup_{x_j\in D_j}\bigg|\frac{1}{n}\sum_{i=1}^n(K_{h_j}(x_j,\tilde{\xi}_j^i)-K_{h_j}(x_j,\xi_j^i))\I(\mxi^i\in D)\I(\tilde{\mxi}^i\in D)\bigg|\\
&\leq\sup_{x_j\in D_j}\frac{1}{n}\sum_{i=1}^nK_{h_j}(x_j,\xi_j^i)\I(\mxi^i\in D\cap D(n))\\
&\quad+\sup_{x_j\in D_j}\frac{1}{n}\sum_{i=1}^nK_{h_j}(x_j,\tilde{\xi}_j^i)\I(\tilde{\mxi}^i\in D)\I(\mxi^i\in D(n))\\
&\quad+\sup_{x_j\in D_j}\frac{1}{n}\sum_{i=1}^n|K_{h_j}(x_j,\tilde{\xi}_j^i)-K_{h_j}(x_j,\xi_j^i)|\I(\mxi^i\in D)\I(\tilde{\mxi}^i\in D),
\end{split}
\end{align}
where the last inequality follows from (\ref{indicator inequality}). We first approximate the first term on the right hand side of (\ref{decomposition j}). We note that
\ba
&\sup_{x_j\in D_j}\frac{1}{n}\sum_{i=1}^nK_{h_j}(x_j,\xi_j^i)\I(\mxi^i\in D\cap D(n))\\
&\leq\sup_{x_j\in D_j}\bigg|\frac{1}{n}\sum_{i=1}^nK_{h_j}(x_j,\xi_j^i)\I(\mxi^i\in D\cap D(n))-\E(K_{h_j}(x_j,\xi_j)\I(\mxi\in D\cap D(n)))\bigg|\\
&\quad+\sup_{x_j\in D_j}\E(K_{h_j}(x_j,\xi_j)\I(\mxi\in D\cap D(n))).
\ea
We define $V_{nj}=h_j^{-L_j}\cdot a_{nj}+\sum_{k\neq j}a_{nk}$. The first term on the right hand side of the above inequality has the rate $O_p(n^{-1/2}\cdot(\log{n})^{1/2}\cdot V_{nj}^{1/2}\cdot h_j^{-L_j/2})$ provided that $n^{-1/2}\cdot(\log{n})^{1/2}\cdot V_{nj}^{-1/2}\cdot h_j^{-L_j/2}=O(1)$ by Lemma \ref{uniform rates density Dn}. Then,
\begin{align}\label{sup E}
\begin{split}
&\sup_{x_j\in D_j}\E(K_{h_j}(x_j,\xi_j)\I(\mxi\in D\cap D(n)))\\
&=\sup_{x_j\in D_j}\int_{D\cap D(n)}K_{h_j}(x_j,u_j)p(\bu)d\bu\\
&\leq\sup_{x_j\in D_j}\sum_{k=1}^d\int_{D\cap D(k,n)}K_{h_j}(x_j,u_j)p(\bu)d\bu\\
&=p_0^D\cdot\sup_{x_j\in D_j}\int_{D_j\setminus D_j^-(C\cdot a_{nj})}K_{h_j}(x_j,u_j)p_j^D(u_j)du_j\\
&\quad+p_0^D\cdot\sup_{x_j\in D_j}\sum_{k\neq j}\int_{D_j}\int_{D_k\setminus D_k^-(C\cdot a_{nk})}K_{h_j}(x_j,u_j)p_{jk}^D(u_j,u_k)du_kdu_j\\
&\leq\const\cdot\sup_{x_j\in D_j}\int_{D_j\setminus D_j^-(C\cdot a_{nj})}K_{h_j}(x_j,u_j)du_j\\
&\quad+\const\cdot\sum_{k\neq j}\Leb_k(D_k\setminus D_k^-(C\cdot a_{nk}))\cdot\sup_{x_j\in D_j}\int_{D_j}K_{h_j}(x_j,u_j)du_j\\
&\leq\const\cdot\sup_{x_j\in D_j}\int_{(D_j\setminus D_j^-(C\cdot a_{nj})-x_j)/h_j}K_j(\|s_j\|_j)ds_j\\
&\quad+\const\cdot\sum_{k\neq j}\Leb_k(D_k\setminus D_k^-(C\cdot a_{nk}))\cdot\int_{B_j(\bzero_j,1)}K_j(\|s_j\|_j)ds_j\\
&\leq\const\cdot\bigg(\Leb_j((D_j\setminus D_j^-(C\cdot a_{nj}))/h_j)+\sum_{k\neq j}\Leb_k(D_k\setminus D_k^-(C\cdot a_{nk}))\bigg)\\
&=O(V_{nj}).
\end{split}
\end{align}
Hence,
\ba
\sup_{x_j\in D_j}\frac{1}{n}\sum_{i=1}^nK_{h_j}(x_j,\xi_j^i)\I(\mxi^i\in D\cap D(n))=O_p(n^{-1/2}\cdot(\log{n})^{1/2}\cdot V_{nj}^{1/2}\cdot h_j^{-L_j/2}+V_{nj}).
\ea
We now approximate the second term on the right hand side of (\ref{decomposition j}). We note that
\ba
&\sup_{x_j\in D_j}\frac{1}{n}\sum_{i=1}^nK_{h_j}(x_j,\tilde{\xi}_j^i)\I(\tilde{\mxi}^i\in D)\I(\mxi^i\in D(n))\\
&\leq\const\sup_{x_j\in D_j}\frac{1}{n}\sum_{i=1}^n\frac{1}{h_j^{L_j}}K_j\bigg(\frac{\|x_j-\tilde{\xi}_j^i\|_j}{h_j}\bigg)\I(\mxi^i\in D(n)).
\ea
Here,
\begin{align}\label{K taylor}
K_j\bigg(\frac{\|x_j-\tilde{\xi}_j^i\|_j}{h_j}\bigg)&=K_j\bigg(\frac{\|x_j-\xi_j^i\|_j}{h_j}\bigg)+\sum_{l=1}^{L_j}\frac{\partial}{\partial u_{jl}}K_j\bigg(\frac{\|x_j-u_j\|_j}{h_j}\bigg)\bigg|_{u_j=\bar{\xi}_j^i}(\tilde{\xi}_{jl}^i-\xi_{jl}^i)
\end{align}
for some random vector $\bar{\xi}_j^i$ lying on the line connecting $\tilde{\xi}_j^i$ and $\xi_j^i$. We note that
\ba
&\sup_{x_j\in D_j}\frac{1}{n}\sum_{i=1}^n\frac{1}{h_j^{L_j}}K_j\bigg(\frac{\|x_j-\xi_j^i\|_j}{h_j}\bigg)\I(\mxi^i\in D(n))\\
&\leq\sup_{x_j\in D_j}\bigg|\frac{1}{n}\sum_{i=1}^n\frac{1}{h_j^{L_j}}K_j\bigg(\frac{\|x_j-\xi_j^i\|_j}{h_j}\bigg)\I(\mxi^i\in D(n))-\E\bigg(\frac{1}{h_j^{L_j}}K_j\bigg(\frac{\|x_j-\xi_j\|_j}{h_j}\bigg)\I(\mxi\in D(n))\bigg)\bigg|\\
&\quad+\sup_{x_j\in D_j}\E\bigg(\frac{1}{h_j^{L_j}}K_j\bigg(\frac{\|x_j-\xi_j\|_j}{h_j}\bigg)\I(\mxi\in D(n))\bigg)\\
&=O_p(n^{-1/2}\cdot(\log{n})^{1/2}\cdot V_{nj}^{1/2}\cdot h_j^{-L_j/2})+O_p(V_{nj}),
\ea
provided that $n^{-1/2}\cdot(\log{n})^{1/2}\cdot V_{nj}^{-1/2}\cdot h_j^{-L_j/2}=O(1)$ by Lemma \ref{uniform rates density Dn} and arguing as in (\ref{sup E}). Also, it holds that
\begin{align}\label{pj approximation pf}
\begin{split}
&\sup_{x_j\in D_j}\bigg|\frac{1}{n}\sum_{i=1}^n\frac{1}{h_j^{L_j}}\frac{\partial}{\partial u_{jl}}K_j\bigg(\frac{\|x_j-u_j\|_j}{h_j}\bigg)\bigg|_{u_j=\bar{\xi}_j^i}(\tilde{\xi}_{jl}^i-\xi_{jl}^i)\I(\mxi^i\in D(n))\bigg|\\
&\leq\sup_{x_j\in D_j}\frac{1}{n}\sum_{i=1}^n\frac{1}{h_j^{L_j}}\bigg|\frac{\partial}{\partial u_{jl}}K_j\bigg(\frac{\|x_j-u_j\|_j}{h_j}\bigg)\bigg|_{u_j=\bar{\xi}_j^i}(\tilde{\xi}_{jl}^i-\xi_{jl}^i)\I(\mxi^i\in D(n))\bigg|\\
&\leq\max_{1\leq i\leq n}|\tilde{\xi}_{jl}^i-\xi_{jl}^i|\sup_{x_j\in D_j}\frac{1}{n}\sum_{i=1}^n\frac{1}{h_j^{L_j}}\bigg|\frac{\partial}{\partial u_{jl}}K_j\bigg(\frac{\|x_j-u_j\|_j}{h_j}\bigg)\bigg|_{u_j=\bar{\xi}_j^i}\bigg|\I(\mxi^i\in D(n))\\
&\leq\max_{1\leq i\leq n}|\tilde{\xi}_{jl}^i-\xi_{jl}^i|\sup_{x_j\in D_j}\frac{1}{n}\sum_{i=1}^n\frac{1}{h_j^{L_j}}\bigg|\frac{\partial}{\partial u_{jl}}K_j\bigg(\frac{\|x_j-u_j\|_j}{h_j}\bigg)\bigg|_{u_j=\bar{\xi}_j^i}-\frac{\partial}{\partial u_{jl}}K_j\bigg(\frac{\|x_j-u_j\|_j}{h_j}\bigg)\bigg|_{u_j=\xi_j^i}\bigg|\\
&\qquad\qquad\qquad\qquad\qquad\qquad\qquad\cdot\I(\mxi^i\in D(n))\\
&\quad+\max_{1\leq i\leq n}|\tilde{\xi}_{jl}^i-\xi_{jl}^i|\sup_{x_j\in D_j}\frac{1}{n}\sum_{i=1}^n\frac{1}{h_j^{L_j}}\bigg|\frac{\partial}{\partial u_{jl}}K_j\bigg(\frac{\|x_j-u_j\|_j}{h_j}\bigg)\bigg|_{u_j=\xi_j^i}\bigg|\I(\mxi^i\in D(n))\\
&\leq h_j^{-L_j-2}\max_{1\leq i\leq n}|\tilde{\xi}_{jl}^i-\xi_{jl}^i|^2\frac{1}{n}\sum_{i=1}^n\I(\mxi^i\in D(n))\\
&\quad+\max_{1\leq i\leq n}|\tilde{\xi}_{jl}^i-\xi_{jl}^i|\sup_{x_j\in D_j}\bigg|\frac{1}{n}\sum_{i=1}^n\frac{1}{h_j^{L_j}}\bigg|\frac{\partial}{\partial u_{jl}}K_j\bigg(\frac{\|x_j-u_j\|_j}{h_j}\bigg)\bigg|_{u_j=\xi_j^i}\bigg|\I(\mxi^i\in D(n))\\
&\qquad\qquad\qquad\qquad\qquad\qquad-\E\bigg(\frac{1}{h_j^{L_j}}\bigg|\frac{\partial}{\partial u_{jl}}K_j\bigg(\frac{\|x_j-u_j\|_j}{h_j}\bigg)\bigg|_{u_j=\xi_j}\bigg|\I(\mxi\in D(n))\bigg)\bigg|\\
&\quad+\max_{1\leq i\leq n}|\tilde{\xi}_{jl}^i-\xi_{jl}^i|\sup_{x_j\in D_j}\E\bigg(\frac{1}{h_j^{L_j}}\bigg|\frac{\partial}{\partial u_{jl}}K_j\bigg(\frac{\|x_j-u_j\|_j}{h_j}\bigg)\bigg|_{u_j=\xi_j}\bigg|\I(\mxi\in D(n))\bigg).
\end{split}
\end{align}
The first term has the rate
\ba
O_p(h_j^{-L_j-2}\cdot a_{nj}^2)\cdot O_p\bigg(\sum_{j=1}^da_{nj}+n^{-1/2}\cdot\bigg(\sum_{j=1}^da_{nj}\bigg)^{1/2}\bigg)
\ea
by the proof of Lemma \ref{p0 rate}. We now approximate the second term on the right hand side of (\ref{pj approximation pf}). For $\kappa>0$, we let $\{x_j^{(l)}:1\leq l\leq N(n^{-\kappa})\}\subset D_j$ be a set of points such that $N(n^{-\kappa})=O(n^{\kappa L_j})$ and $\{B_j(x_j^{(l)},n^{-\kappa}):1\leq l\leq N(n^{-\kappa})\}$ covers $D_j$.
Then, it holds that 
\begin{align*}
&\sup_{x_j\in D_j}\bigg|\frac{1}{n}\sum_{i=1}^n\frac{1}{h_j^{L_j}}\bigg|\frac{\partial}{\partial u_{jl}}K_j\bigg(\frac{\|x_j-u_j\|_j}{h_j}\bigg)\bigg|_{u_j=\xi^i_j}\bigg|\I(\mxi^i\in D(n))\\
&\qquad\quad-\E\bigg(\frac{1}{h_j^{L_j}}\bigg|\frac{\partial}{\partial u_{jl}}K_j\bigg(\frac{\|x_j-u_j\|_j}{h_j}\bigg)\bigg|_{u_j=\xi_j}\bigg|\I(\mxi\in D(n))\bigg)\bigg|\\
&\leq\max_{1\leq l\leq N(n^{-\kappa})}\bigg|\frac{1}{n}\sum_{i=1}^n\frac{1}{h_j^{L_j}}\bigg|\frac{\partial}{\partial u_{jl}}K_j\bigg(\frac{\|x_j^{(l)}-u_j\|_j}{h_j}\bigg)\bigg|_{u_j=\xi^i_j}\bigg|\I(\mxi^i\in D(n))\\
&\qquad\qquad\qquad-\E\bigg(\frac{1}{h_j^{L_j}}\bigg|\frac{\partial}{\partial u_{jl}}K_j\bigg(\frac{\|x^{(l)}_j-u_j\|_j}{h_j}\bigg)\bigg|_{u_j=\xi_j}\bigg|\I(\mxi\in D(n))\bigg)\bigg|\\
&+\max_{1\leq l\leq N(n^{-\kappa})}\sup_{x_j\in D_j\cap B_j(x_j^{(l)},n^{-\kappa})}\bigg|\frac{1}{n}\sum_{i=1}^n\frac{1}{h_j^{L_j}}\bigg|\bigg|\frac{\partial}{\partial u_{jl}}K_j\bigg(\frac{\|x_j-u_j\|_j}{h_j}\bigg)\bigg|_{u_j=\xi^i_j}\bigg|\\
&\qquad\qquad\qquad\qquad\qquad\qquad\qquad\qquad\qquad-\bigg|\frac{\partial}{\partial u_{jl}}K_j\bigg(\frac{\|x_j^{(l)}-u_j\|_j}{h_j}\bigg)\bigg|_{u_j=\xi^i_j}\bigg|\bigg|\I(\mxi^i\in D(n))\\
&+\max_{1\leq l\leq N(n^{-\kappa})}\sup_{x_j\in D_j\cap B_j(x_j^{(l)},n^{-\kappa})}\bigg|\E\bigg(\frac{1}{h_j^{L_j}}\bigg|\frac{\partial}{\partial u_{jl}}K_j\bigg(\frac{\|x_j-u_j\|_j}{h_j}\bigg)\bigg|_{u_j=\xi_j}\bigg|\I(\mxi\in D(n))\bigg)\\
&\qquad\qquad\qquad\qquad\qquad\qquad\qquad-\E\bigg(\frac{1}{h_j^{L_j}}\bigg|\frac{\partial}{\partial u_{jl}}K_j\bigg(\frac{\|x^{(l)}_j-u_j\|_j}{h_j}\bigg)\bigg|_{u_j=\xi_j}\bigg|\I(\mxi\in D(n))\bigg).
\end{align*}
The second and third terms on the right hand side of the above inequality are bounded by $\const\cdot n^{-\kappa}\cdot h_j^{-L_j-2}$.
Hence, the second and third terms are negligible by taking sufficiently large $\kappa$ and thus it suffices to approximate the first term of the above inequality. We claim that the first term has the rate
\begin{align}\label{claim for pj}
O_p(n^{-1/2}\cdot(\log{n})^{1/2}\cdot V_{nj}^{1/2}\cdot h_j^{-L_j/2-1}).
\end{align}
For this, we define
\begin{align*}
Z_n^i(x_j)=&\frac{1}{n}\bigg(\frac{1}{h_j^{L_j}}\bigg|\frac{\partial}{\partial u_{jl}}K_j\bigg(\frac{\|x_j-u_j\|_j}{h_j}\bigg)\bigg|_{u_j=\xi_j^i}\bigg|\I(\mxi^i\in D(n))\\
&\qquad-\E\bigg(\frac{1}{h_j^{L_j}}\bigg|\frac{\partial}{\partial u_{jl}}K_j\bigg(\frac{\|x_j-u_j\|_j}{h_j}\bigg)\bigg|_{u_j=\xi_j}\bigg|\I(\mxi\in D(n))\bigg)\bigg).
\end{align*}
It holds that
\begin{align*}
\E(Z_n^i(x_j))=0,&\quad|Z_n^i(x_j)|\leq\const\cdot n^{-1}\cdot h_j^{-L_j-1},\quad\mbox{and}\\
\sum_{i=1}^n\E(|Z_n^i(x_j)|^2)&\leq n^{-1}\cdot \E\bigg(\bigg(\frac{1}{h_j^{L_j}}\frac{\partial}{\partial u_{jl}}K_j\bigg(\frac{\|x_j-u_j\|_j}{h_j}\bigg)\bigg|_{u_j=\xi_j}\bigg)^2\I(\mxi\in D(n))\bigg)\\
&=O(n^{-1}\cdot h_j^{-L_j-2}\cdot V_{nj}).
\end{align*}
Now, Theorem 2.6.2 in Bosq (2000) gives (\ref{claim for pj}) if $n^{-1/2}\cdot(\log{n})^{1/2}\cdot V_{nj}^{-1/2}\cdot h_j^{-L_j/2}=O(1)$.
Hence, the second term on the right hand side of (\ref{pj approximation pf}) has the rate
\begin{align}\label{concentration of 2nd pj}
O_p(n^{-1/2}\cdot(\log{n})^{1/2}\cdot V_{nj}^{1/2}\cdot h_j^{-L_j/2-1}\cdot a_{nj}).
\end{align}
We now approximate the third term on the right hand side of (\ref{pj approximation pf}). We note that
\ba
\sup_{x_j\in D_j}\E\bigg(\frac{1}{h_j^{L_j}}\bigg|\frac{\partial}{\partial u_{jl}}K_j\bigg(\frac{\|x_j-u_j\|_j}{h_j}\bigg)\bigg|_{u_j=\xi_j}\bigg|\I(\mxi\in D(n))\bigg)=O(h_j^{-1}\cdot V_{nj}),
\ea
Hence, the third term on the right hand side of (\ref{pj approximation pf}) has the rate
\begin{align}\label{pj 2nd 2nd rate}
O_p(h_j^{-1}\cdot a_{nj}\cdot V_{nj}).
\end{align}
Thus, the second term on the right hand side of (\ref{decomposition j}) has the rate
\begin{align}\label{pj 2nd rate}
\begin{split}
&O_p(n^{-1/2}\cdot(\log{n})^{1/2}\cdot V_{nj}^{1/2}\cdot h_j^{-L_j/2}+V_{nj})\\
&+O_p(h_j^{-L_j-2}\cdot a_{nj}^2)\cdot O_p\bigg(\sum_{j=1}^da_{nj}+n^{-1/2}\cdot\bigg(\sum_{j=1}^da_{nj}\bigg)^{1/2}\bigg)\\
&+O_p(n^{-1/2}\cdot(\log{n})^{1/2}\cdot V_{nj}^{1/2}\cdot h_j^{-L_j/2-1}\cdot a_{nj})+O_p(h_j^{-1}\cdot a_{nj}\cdot V_{nj})\\
&=O_p(n^{-1/2}\cdot(\log{n})^{1/2}\cdot V_{nj}^{1/2}\cdot h_j^{-L_j/2}+V_{nj}),
\end{split}
\end{align}
where the equality follows from that $h_j^{-1}\cdot a_{nj}=O(1)$ and $h_j^{-L_j-2}\cdot a^2_{nj}=O(1)$.
We now approximate the third term on the right hand side of (\ref{decomposition j}). We note that
\begin{align}\label{pj 3rd part}
\begin{split}
&\frac{1}{n}\sum_{i=1}^n|K_{h_j}(x_j,\tilde{\xi}_j^i)-K_{h_j}(x_j,\xi_j^i)|\I(\mxi^i\in D)\I(\tilde{\mxi}^i\in D)\\
&\leq\const\cdot\frac{1}{n}\sum_{i=1}^n\frac{1}{h_j^{L_j}}\bigg|K_j\bigg(\frac{\|x_j-\tilde{\xi}_j^i\|_j}{h_j}\bigg)-K_j\bigg(\frac{\|x_j-\xi_j^i\|_j}{h_j}\bigg)\bigg|\I(\mxi^i\in D)\\
&\quad+\const\cdot\frac{1}{n}\sum_{i=1}^n\int_{D_j}\frac{1}{h_j^{L_j}}\bigg|K_j\bigg(\frac{\|t_j-\tilde{\xi}_j^i\|_j}{h_j}\bigg)-K_j\bigg(\frac{\|t_j-\xi_j^i\|_j}{h_j}\bigg)\bigg|dt_j\\
&\qquad\qquad\qquad\qquad\qquad\cdot \frac{1}{h_j^{L_j}}K_j\bigg(\frac{\|x_j-\xi_j^i\|_j}{h_j}\bigg)\I(\mxi^i\in D)\\
&\leq\const\cdot\frac{1}{n}\sum_{i=1}^n\frac{1}{h_j^{L_j}}\bigg|\sum_{l=1}^{L_j}\frac{\partial}{\partial u_{jl}}K_j\bigg(\frac{\|x_j-u_j\|_j}{h_j}\bigg)\bigg|_{u_j=\bar{\xi}_j^i}(\tilde{\xi}_{jl}^i-\xi_{jl}^i)\bigg|\I(\mxi^i\in D)\\
&\quad+\const\cdot\frac{1}{n}\sum_{i=1}^n\int_{D_j}\frac{1}{h_j^{L_j}}\bigg|\sum_{l=1}^{L_j}\frac{\partial}{\partial u_{jl}}K_j\bigg(\frac{\|t_j-u_j\|_j}{h_j}\bigg)\bigg|_{u_j=\bar{\xi}_j^i}(\tilde{\xi}_{jl}^i-\xi_{jl}^i)\bigg|dt_j\\
&\qquad\qquad\qquad\qquad\qquad\cdot \frac{1}{h_j^{L_j}}K_j\bigg(\frac{\|x_j-\xi_j^i\|_j}{h_j}\bigg)\I(\mxi^i\in D)\\
&\leq\const\cdot\sum_{l=1}^{L_j}\max_{1\leq i\leq n}|\tilde{\xi}_{jl}^i-\xi_{jl}^i|\frac{1}{n}\sum_{i=1}^n\frac{1}{h_j^{L_j}}\bigg|\frac{\partial}{\partial u_{jl}}K_j\bigg(\frac{\|x_j-u_j\|_j}{h_j}\bigg)\bigg|_{u_j=\bar{\xi}_j^i}\bigg|\I(\mxi^i\in D)\\
&\quad+\const\cdot h_j^{-1}\max_{1\leq i\leq n,1\leq l\leq L_j}|\tilde{\xi}_{jl}^i-\xi_{jl}^i|\frac{1}{n}\sum_{i=1}^n\frac{1}{h_j^{L_j}}K_j\bigg(\frac{\|x_j-\xi_j^i\|_j}{h_j}\bigg)\I(\mxi^i\in D),
\end{split}
\end{align}
where the $\const$ is independent of $x_j$. We may show that
the first term on the right hand side of (\ref{pj 3rd part}) has the uniform rate
\ba
O_p(h_j^{-L_j-2}\cdot a^2_{nj}+n^{-1/2}\cdot(\log{n})^{1/2}\cdot h_j^{-L_j/2-1}\cdot a_{nj}+h_j^{-1}\cdot a_{nj})
\ea
over $x_j\in D_j$ provided that $n^{-1/2}\cdot(\log{n})^{1/2}\cdot h_j^{-L_j/2}=O(1)$ by arguing as the proof of (\ref{pj 2nd rate}). We may also show that the second term on the right hand side of (\ref{pj 3rd part}) has the uniform rate
\ba
O_p(n^{-1/2}\cdot(\log{n})^{1/2}\cdot h_j^{-L_j/2-1}\cdot a_{nj}+h_j^{-1}\cdot a_{nj})
\ea
over $x_j\in D_j$ provided that $n^{-1/2}\cdot(\log{n})^{1/2}\cdot h_j^{-L_j/2}=O(1)$ by Lemma \ref{uniform rates density}.
Hence, the third term on the right hand side of (\ref{decomposition j}) has the rate
\begin{align*}
O_p(h_j^{-L_j-2}\cdot a^2_{nj}+h_j^{-1}\cdot a_{nj})
\end{align*}
since $n^{-1/2}\cdot(\log{n})^{1/2}\cdot h_j^{-L_j/2}=O(1)$.
Therefore,
\ba
\sup_{x_j\in D_j}|\hat{p}^D_j(x_j)-\check{p}^D_j(x_j)|&=O_p\bigg(n^{-1/2}\cdot(\log{n})^{1/2}\cdot V_{nj}^{1/2}\cdot h_j^{-L_j/2}+V_{nj}+h_j^{-L_j-2}\cdot a_{nj}^2\bigg)\\
&=O_p(V_{nj}+h_j^{-L_j-2}\cdot a_{nj}^2)
\ea
since $n^{-1/2}\cdot(\log{n})^{1/2}\cdot V_{nj}^{-1/2}\cdot h_j^{-L_j/2}=O(1)$. This completes the proof for the second part of the lemma.

The first part of the lemma follows similarly as in the above proof by obtaining the rates with $V_{nj}^{1/2}$ being replaced by 1 using the bound $\I(\mxi^i\in D\cap D(n))\leq\I(\mxi^i\in D)$ or $\I(\mxi^i\in D(n))\leq\I(\mxi^i\in\prod_{j=1}^dD_j^+(C\cdot a_{nj}))$ whenever $V_{nj}^{1/2}$ appears. In this case, the condition that $n^{-1/2}\cdot(\log{n})^{1/2}\cdot h_j^{-L_j/2}=O(1)$ is required instead of the condition that $n^{-1/2}\cdot(\log{n})^{1/2}\cdot V_{nj}^{-1/2}\cdot h_j^{-L_j/2}=O(1)$.
\end{proof}

\begin{customlemma}{S.12}\label{marginal density approximation 2}
Assume that conditions (B3) and (B4) hold, that $h_j\vee h_k=O(1)$, that $n^{-1/2}\cdot(\log{n})^{1/2}\cdot h_j^{-L_j/2}\cdot h_k^{-L_k/2}=O(1)$, that $(h_j^{-1}\cdot a_{nj})\vee(h_k^{-1}\cdot a_{nk})=O(1)$, that $(h_j^{-L_j-2}\cdot a^2_{nj})\vee(h_k^{-L_k-2}\cdot a^2_{nk})=O(1)$ and that $p$ is bounded on $D^+(\varepsilon)$ for some $\varepsilon>0$. Then, it holds that
\ba
&\sup_{(x_j,x_k)\in D_j\times D_k}|\hat{p}^D_{jk}(x_j,x_k)-\check{p}^D_{jk}(x_j,x_k)|\\
&=O_p\left(\sqrt{\frac{\log{n}}{n\cdot h_j^{L_j}\cdot h_k^{L_k}}}+h_j^{-L_j}\cdot a_{nj}+h_k^{-L_k}\cdot a_{nk}+\sum_{m\neq j,k}a_{nm}+h_j^{-L_j-2}\cdot a_{nj}^2+h_k^{-L_k-2}\cdot a_{nk}^2\right).
\ea
If we further assume that $n^{-1/2}\cdot(\log{n})^{1/2}\cdot(h_j^{-L_j}\cdot a_{nj}+h_k^{-L_k}\cdot a_{nk}+\sum_{m\neq j,k}a_{nm})^{-1/2}\cdot h_j^{-L_j/2}\cdot h_k^{-L_k/2}=O(1)$, then it holds that
\ba
&\sup_{(x_j,x_k)\in D_j\times D_k}|\hat{p}^D_{jk}(x_j,x_k)-\check{p}^D_{jk}(x_j,x_k)|\\
&=O_p\bigg(h_j^{-L_j}\cdot a_{nj}+h_k^{-L_k}\cdot a_{nk}+\sum_{m\neq j,k}a_{nm}+h_j^{-L_j-2}\cdot a_{nj}^2+h_k^{-L_k-2}\cdot a_{nk}^2\bigg).
\ea
\end{customlemma}

\begin{proof}
Since the outline of the proof is similar to the proof of Lemma \ref{marginal density approximation 1}, we omit some details for simplicity. We first prove the second part of the lemma.
We note that
\begin{align}\label{pjk approximation}
\begin{split}
&\sup_{(x_j,x_k)\in D_j\times D_k}|\hat{p}^D_{jk}(x_j,x_k)-\check{p}^D_{jk}(x_j,x_k)|\\
&\leq\frac{|\hat{p}_0^D-\check{p}_0^D|}{\hat{p}_0^D\check{p}_0^D}\sup_{(x_j,x_k)\in D_j\times D_k}\bigg|\frac{1}{n}\sum_{i=1}^nK_{h_j}(x_j,\xi_j^i)K_{h_k}(x_k,\xi_k^i)\I(\mxi^i\in D)\bigg|\\
&\quad+\frac{1}{\hat{p}_0^D}\sup_{(x_j,x_k)\in D_j\times D_k}\bigg|\frac{1}{n}\sum_{i=1}^n(K_{h_j}(x_j,\tilde{\xi}_j^i)K_{h_k}(x_k,\tilde{\xi}_k^i)\I(\tilde{\mxi}^i\in D)\\
&\qquad\qquad\qquad\qquad\qquad\qquad-K_{h_j}(x_j,\xi_j^i)K_{h_k}(x_k,\xi_k^i)\I(\mxi^i\in D))\bigg|.
\end{split}
\end{align}
We first approximate the first term on the right hand side of (\ref{pjk approximation}). For this, we approximate $\sup_{(x_j,x_k)\in D_j\times D_k}|n^{-1}\sum_{i=1}^nK_{h_j}(x_j,\xi_j^i)K_{h_k}(x_k,\xi_k^i)\I(\mxi^i\in D)|$.
By Lemma \ref{uniform rates density}, it holds that
\ba
&\sup_{(x_j,x_k)\in D_j\times D_k}\bigg|\frac{1}{n}\sum_{i=1}^nK_{h_j}(x_j,\xi_j^i)K_{h_k}(x_k,\xi_k^i)\I(\mxi^i\in D)-\E(K_{h_j}(x_j,\xi_j)K_{h_k}(x_k,\xi_k)\I(\mxi\in D))\bigg|\\
&=O_p(n^{-1/2}\cdot(\log{n})^{1/2}\cdot h_j^{-L_j/2}\cdot h_k^{-L_k/2}),
\ea
provided that $n^{-1/2}\cdot(\log{n})^{1/2}\cdot h_j^{-L_j/2}\cdot h_k^{-L_k/2}=O(1)$.
Also,
\ba
&\sup_{(x_j,x_k)\in D_j\times D_k}\E(K_{h_j}(x_j,\xi_j)K_{h_k}(x_k,\xi_k)\I(\mxi\in D))\\
&=p_0^D\cdot\sup_{(x_j,x_k)\in D_j\times D_k}\int_{D_k}\int_{D_j}K_{h_j}(x_j,u_j)K_{h_k}(x_k,u_k)p_{jk}^D(u_j,u_k)du_jdu_k\\
&\leq\const\cdot\sup_{x_j\in D_j}\int_{D_j}h_j^{-L_j}K_j\bigg(\frac{\|x_j-u_j\|_j}{h_j}\bigg)du_j\cdot\sup_{x_k\in D_k}\int_{D_k}h_k^{-L_k}K_k\bigg(\frac{\|x_k-u_k\|_k}{h_k}\bigg)du_k\\
&=O(1).
\ea
Hence, we have
\begin{align*}
\sup_{(x_j,x_k)\in D_j\times D_k}\bigg|\frac{1}{n}\sum_{i=1}^nK_{h_j}(x_j,\xi_j^i)K_{h_k}(x_k,\xi_k^i)\I(\mxi^i\in D)\bigg|=O_p(1).
\end{align*}
Thus, the first term on the right hand side of (\ref{pjk approximation}) has the same rate with $|\hat{p}_0^D-\check{p}_0^D|$, which is
\ba
O_p\bigg(\sum_{j=1}^da_{nj}+n^{-1/2}\cdot\bigg(\sum_{j=1}^da_{nj}\bigg)^{1/2}\bigg).
\ea

We now approximate the second term on the right hand side of (\ref{pjk approximation}). We note that
\begin{align}\label{decomposition jk}
\begin{split}
&\sup_{(x_j,x_k)\in D_j\times D_k}\bigg|\frac{1}{n}\sum_{i=1}^n(K_{h_j}(x_j,\tilde{\xi}_j^i)K_{h_k}(x_k,\tilde{\xi}_k^i)\I(\tilde{\mxi}^i\in D)-K_{h_j}(x_j,\xi_j^i)K_{h_k}(x_k,\xi_k^i)\I(\mxi^i\in D))\bigg|\\
&\leq\sup_{(x_j,x_k)\in D_j\times D_k}\frac{1}{n}\sum_{i=1}^nK_{h_j}(x_j,\xi_j^i)K_{h_k}(x_k,\xi_k^i)\I(\mxi^i\in D\cap D(n))\\
&\quad+\sup_{(x_j,x_k)\in D_j\times D_k}\frac{1}{n}\sum_{i=1}^nK_{h_j}(x_j,\tilde{\xi}_j^i)K_{h_k}(x_k,\tilde{\xi}_k^i)\I(\tilde{\mxi}^i\in D)\I(\mxi^i\in D(n))\\
&\quad+\sup_{(x_j,x_k)\in D_j\times D_k}\frac{1}{n}\sum_{i=1}^n|K_{h_j}(x_j,\tilde{\xi}_j^i)K_{h_k}(x_k,\tilde{\xi}_k^i)-K_{h_j}(x_j,\xi_j^i)K_{h_k}(x_k,\xi_k^i)|\I(\mxi^i\in D)\I(\tilde{\mxi}^i\in D).
\end{split}
\end{align}
We first approximate the first term on the right hand side of (\ref{decomposition jk}). We note that
\ba
&\sup_{(x_j,x_k)\in D_j\times D_k}\frac{1}{n}\sum_{i=1}^nK_{h_j}(x_j,\xi_j^i)K_{h_k}(x_k,\xi_k^i)\I(\mxi^i\in D\cap D(n))\\
&\leq\sup_{(x_j,x_k)\in D_j\times D_k}\bigg|\frac{1}{n}\sum_{i=1}^nK_{h_j}(x_j,\xi_j^i)K_{h_k}(x_k,\xi_k^i)\I(\mxi^i\in D\cap D(n))\\
&\qquad\qquad\qquad\quad-\E(K_{h_j}(x_j,\xi_j)K_{h_k}(x_k,\xi_k)\I(\mxi\in D\cap D(n)))\bigg|\\
&\quad+\sup_{(x_j,x_k)\in D_j\times D_k}\E(K_{h_j}(x_j,\xi_j)K_{h_k}(x_k,\xi_k)\I(\mxi\in D\cap D(n))).
\ea
We define $V_{nj}=h_j^{-L_j}\cdot a_{nj}+\sum_{k\neq j}a_{nk}$ and $V_{njk}=h_j^{-L_j}\cdot a_{nj}+h_k^{-L_k}\cdot a_{nk}+\sum_{m\neq j,k}a_{nm}$. The first term on the right hand side of the above inequality has the rate
\ba
O_p(n^{-1/2}\cdot(\log{n})^{1/2}\cdot V_{njk}^{1/2}\cdot h_j^{-L_j/2}\cdot h_k^{-L_k/2})
\ea
provided that $n^{-1/2}\cdot(\log{n})^{1/2}\cdot V_{njk}^{-1/2}\cdot h_j^{-L_j/2}\cdot h_j^{-L_k/2}=O(1)$ by Lemma \ref{uniform rates density Dn}. Also,
\begin{align*}
&\sup_{(x_j,x_k)\in D_j\times D_k}\E(K_{h_j}(x_j,\xi_j)K_{h_k}(x_k,\xi_k)\I(\mxi\in D\cap D(n)))\\
&\leq\sup_{(x_j,x_k)\in D_j\times D_k}\sum_{k=1}^d\int_{D\cap D(k,n)}K_{h_j}(x_j,u_j)K_{h_k}(x_k,u_k)p(\bu)d\bu\\
&=p_0^D\cdot\sup_{(x_j,x_k)\in D_j\times D_k}\int_{D_k}\int_{D_j\setminus D_j^-(C\cdot a_{nj})}K_{h_j}(x_j,u_j)K_{h_k}(x_k,u_k)p_{jk}^D(u_j,u_k)du_jdu_k\\
&\quad+p_0^D\cdot\sup_{(x_j,x_k)\in D_j\times D_k}\int_{D_j}\int_{D_k\setminus D_k^-(C\cdot a_{nk})}K_{h_j}(x_j,u_j)K_{h_k}(x_k,u_k)p_{jk}^D(u_j,u_k)du_kdu_j\\
&\quad+p_0^D\cdot\sup_{(x_j,x_k)\in D_j\times D_k}\sum_{m\neq j,k}\int_{D_j}\int_{D_k}\int_{D_m\setminus D_m^-(Ca_{nm})}K_{h_j}(x_j,u_j)K_{h_k}(x_k,u_k)\\
&\qquad\qquad\qquad\qquad\qquad\qquad\qquad\qquad\qquad\qquad\qquad\cdot p_{jkm}^D(u_j,u_k,u_m)du_mdu_kdu_j\\
&\leq\const\cdot\sup_{(x_j,x_k)\in D_j\times D_k}\int_{D_k}\int_{D_j\setminus D_j^-(C\cdot a_{nj})}K_{h_j}(x_j,u_j)K_{h_k}(x_k,u_k)du_jdu_k\\
&\quad+\const\cdot\sup_{(x_j,x_k)\in D_j\times D_k}\int_{D_j}\int_{D_k\setminus D_k^-(C\cdot a_{nk})}K_{h_j}(x_j,u_j)K_{h_k}(x_k,u_k)du_kdu_j\\
&\quad+\const\cdot\sum_{m\neq j,k}\Leb_m(D_m\setminus D_m^-(Ca_{nm}))\\
&\leq\const\cdot\sup_{x_j\in D_j}\int_{(D_j\setminus D_j^-(C\cdot a_{nj})-x_j)/h_j}K_j(\|s_j\|_j)ds_j\\
&\quad+\const\cdot\sup_{x_k\in D_k}\int_{(D_k\setminus D_k^-(C\cdot a_{nk})-x_k)/h_k}K_k(\|s_k\|_k)ds_k\\
&\quad+\const\cdot\sum_{m\neq j,k}\Leb_m(D_m\setminus D_m^-(Ca_{nm}))\\
&\leq\const\cdot\bigg(\Leb_j((D_j\setminus D_j^-(C\cdot a_{nj}))/h_j)+\Leb_k((D_k\setminus D_k^-(C\cdot a_{nk}))/h_k)\\
&\qquad\qquad\qquad+\sum_{m\neq j,k}\Leb_m(D_m\setminus D_m^-(Ca_{nm}))\bigg)\\
&=O(V_{njk}),
\end{align*}
where $p_{jkm}^D$ is similarly defined as $p_{jk}^D$.
Hence,
\ba
&\sup_{(x_j,x_k)\in D_j\times D_k}\frac{1}{n}\sum_{i=1}^nK_{h_j}(x_j,\xi_j^i)K_{h_k}(x_k,\xi_k^i)\I(\mxi^i\in D\cap D(n))\\
&=O_p(n^{-1/2}\cdot(\log{n})^{1/2}\cdot V_{njk}^{1/2}\cdot h_j^{-L_j/2}\cdot h_k^{-L_k/2}+V_{njk}).
\ea
We now approximate the second term on the right hand side of (\ref{decomposition jk}). We note that
\ba
&\sup_{(x_j,x_k)\in D_j\times D_k}\frac{1}{n}\sum_{i=1}^nK_{h_j}(x_j,\tilde{\xi}_j^i)K_{h_k}(x_k,\tilde{\xi}_k^i)\I(\tilde{\mxi}^i\in D)\I(\mxi^i\in D(n))\\
&\leq\const\cdot\sup_{(x_j,x_k)\in D_j\times D_k}\frac{1}{n}\sum_{i=1}^n\frac{1}{h_j^{L_j}}K_j\bigg(\frac{\|x_j-\tilde{\xi}_j^i\|_j}{h_j}\bigg)\frac{1}{h_k^{L_k}}K_k\bigg(\frac{\|x_k-\tilde{\xi}_k^i\|_k}{h_k}\bigg)\I(\mxi^i\in D(n)).
\ea
Here,
\ba
&K_j\bigg(\frac{\|x_j-\tilde{\xi}_j^i\|_j}{h_j}\bigg)K_k\bigg(\frac{\|x_k-\tilde{\xi}_k^i\|_k}{h_k}\bigg)\\
&=\bigg(K_j\bigg(\frac{\|x_j-\xi_j^i\|_j}{h_j}\bigg)+\sum_{l=1}^{L_j}\frac{\partial}{\partial u_{jl}}K_j\bigg(\frac{\|x_j-u_j\|_j}{h_j}\bigg)\bigg|_{u_j=\bar{\xi}_j^i}(\tilde{\xi}_{jl}^i-\xi_{jl}^i)\bigg)\\
&\quad\cdot\bigg(K_k\bigg(\frac{\|x_k-\xi_k^i\|_k}{h_k}\bigg)+\sum_{l=1}^{L_k}\frac{\partial}{\partial u_{kl}}K_k\bigg(\frac{\|x_k-u_k\|_k}{h_k}\bigg)\bigg|_{u_k=\bar{\xi}_k^i}(\tilde{\xi}_{kl}^i-\xi_{kl}^i)\bigg)
\ea
for some random vector $\bar{\xi}_j^i$ lying on the line connecting $\tilde{\xi}_j^i$ and $\xi_j^i$, and some random vector $\bar{\xi}_k^i$ lying on the line connecting $\tilde{\xi}_k^i$ and $\xi_k^i$. We note that
\ba
&\sup_{(x_j,x_k)\in D_j\times D_k}\frac{1}{n}\sum_{i=1}^n\frac{1}{h_j^{L_j}}K_j\bigg(\frac{\|x_j-\xi_j^i\|_j}{h_j}\bigg)\frac{1}{h_k^{L_k}}K_k\bigg(\frac{\|x_k-\xi_k^i\|_k}{h_k}\bigg)\I(\mxi^i\in D(n))\\
&\leq\sup_{(x_j,x_k)\in D_j\times D_k}\bigg|\frac{1}{n}\sum_{i=1}^n\frac{1}{h_j^{L_j}}K_j\bigg(\frac{\|x_j-\xi_j^i\|_j}{h_j}\bigg)\frac{1}{h_k^{L_k}}K_k\bigg(\frac{\|x_k-\xi_k^i\|_k}{h_k}\bigg)\I(\mxi^i\in D(n))\\
&\qquad\qquad\qquad\qquad-\E\bigg(\frac{1}{h_j^{L_j}}K_j\bigg(\frac{\|x_j-\xi_j\|_j}{h_j}\bigg)\frac{1}{h_k^{L_k}}K_k\bigg(\frac{\|x_k-\xi_k\|_k}{h_k}\bigg)\I(\mxi\in D(n))\bigg)\bigg|\\
&\quad+\sup_{(x_j,x_k)\in D_j\times D_k}\E\bigg(\frac{1}{h_j^{L_j}}K_j\bigg(\frac{\|x_j-\xi_j\|_j}{h_j}\bigg)\frac{1}{h_k^{L_k}}K_k\bigg(\frac{\|x_k-\xi_k\|_k}{h_k}\bigg)\I(\mxi\in D(n))\bigg)\\
&=O_p(n^{-1/2}\cdot(\log{n})^{1/2}\cdot V_{njk}^{1/2}\cdot h_j^{-L_j/2}\cdot h_k^{-L_k/2}+V_{njk}),
\ea
provided that $n^{-1/2}\cdot(\log{n})^{1/2}\cdot V_{njk}^{-1/2}\cdot h_j^{-L_j/2}\cdot h_k^{-L_k/2}=O(1)$ by Lemma \ref{uniform rates density Dn}. Also, it holds that
\begin{align}\label{pjk 2nd term}
\begin{split}
&\bigg|\frac{1}{n}\sum_{i=1}^n\frac{1}{h_j^{L_j}}\frac{\partial}{\partial u_{jl}}K_j\bigg(\frac{\|x_j-u_j\|_j}{h_j}\bigg)\bigg|_{u_j=\bar{\xi}_j^i}(\tilde{\xi}_{jl}^i-\xi_{jl}^i)\frac{1}{h_k^{L_k}}K_k\bigg(\frac{\|x_k-\xi_k^i\|_k}{h_k}\bigg)\I(\mxi^i\in D(n))\bigg|\\
&\leq\max_{1\leq i\leq n}|\tilde{\xi}_{jl}^i-\xi_{jl}^i|\frac{1}{n}\sum_{i=1}^n\frac{1}{h_j^{L_j}}\bigg|\frac{\partial}{\partial u_{jl}}K_j\bigg(\frac{\|x_j-u_j\|_j}{h_j}\bigg)\bigg|_{u_j=\bar{\xi}_j^i}\bigg|\frac{1}{h_k^{L_k}}K_k\bigg(\frac{\|x_k-\xi_k^i\|_k}{h_k}\bigg)\I(\mxi^i\in D(n))\\
&\leq\max_{1\leq i\leq n}|\tilde{\xi}_{jl}^i-\xi_{jl}^i|\frac{1}{n}\sum_{i=1}^n\frac{1}{h_j^{L_j}}\bigg|\frac{\partial}{\partial u_{jl}}K_j\bigg(\frac{\|x_j-u_j\|_j}{h_j}\bigg)\bigg|_{u_j=\bar{\xi}_j^i}-\frac{\partial}{\partial u_{jl}}K_j\bigg(\frac{\|x_j-u_j\|_j}{h_j}\bigg)\bigg|_{u_j=\xi_j^i}\bigg|\\
&\qquad\qquad\qquad\qquad\qquad\cdot\frac{1}{h_k^{L_k}}K_k\bigg(\frac{\|x_k-\xi_k^i\|_k}{h_k}\bigg)\I(\mxi^i\in D(n))\\
&~~+\max_{1\leq i\leq n}|\tilde{\xi}_{jl}^i-\xi_{jl}^i|\bigg|\frac{1}{n}\sum_{i=1}^n\frac{1}{h_j^{L_j}}\bigg|\frac{\partial}{\partial u_{jl}}K_j\bigg(\frac{\|x_j-u_j\|_j}{h_j}\bigg)\bigg|_{u_j=\xi_j^i}\bigg|\frac{1}{h_k^{L_k}}K_k\bigg(\frac{\|x_k-\xi_k^i\|_k}{h_k}\bigg)\I(\mxi^i\in D(n))\\
&\qquad\qquad\qquad-\E\bigg(\frac{1}{h_j^{L_j}}\bigg|\frac{\partial}{\partial u_{jl}}K_j\bigg(\frac{\|x_j-u_j\|_j}{h_j}\bigg)\bigg|_{u_j=\xi_j}\bigg|\frac{1}{h_k^{L_k}}K_k\bigg(\frac{\|x_k-\xi_k\|_k}{h_k}\bigg)\I(\mxi\in D(n))\bigg)\bigg|\\
&~~+\max_{1\leq i\leq n}|\tilde{\xi}_{jl}^i-\xi_{jl}^i|\E\bigg(\frac{1}{h_j^{L_j}}\bigg|\frac{\partial}{\partial u_{jl}}K_j\bigg(\frac{\|x_j-u_j\|_j}{h_j}\bigg)\bigg|_{u_j=\xi_j}\bigg|\frac{1}{h_k^{L_k}}K_k\bigg(\frac{\|x_k-\xi_k\|_k}{h_k}\bigg)\I(\mxi\in D(n))\bigg).
\end{split}
\end{align}
The first term on the right hand side of (\ref{pjk 2nd term}) has the uniform rate
\ba
O_p(h_j^{-L_j-2}\cdot a_{nj}^2\cdot (n^{-1/2}\cdot(\log{n})^{1/2}\cdot V_{nk}^{1/2}\cdot h_k^{-L_k/2}+V_{nk}))
\ea
over $(x_j,x_k)\in D_j\times D_k$ provided that $n^{-1/2}\cdot(\log{n})^{1/2}\cdot V_{nk}^{-1/2}\cdot h_k^{-L_k/2}=O(1)$ by Lemma \ref{uniform rates density Dn}. We now approximate the second term on the right hand side of (\ref{pjk 2nd term}) uniformly for $(x_j,x_k)\in D_j\times D_k$.
For this, we define
\begin{align*}
Z_n^i(x_j,x_k)=&\frac{1}{n}\bigg(\frac{1}{h_j^{L_j}}\bigg|\frac{\partial}{\partial u_{jl}}K_j\bigg(\frac{\|x_j-u_j\|_j}{h_j}\bigg)\bigg|_{u_j=\xi_j^i}\bigg|\frac{1}{h_k^{L_k}}K_k\bigg(\frac{\|x_k-\xi_k^i\|_k}{h_k}\bigg)\I(\mxi^i\in D(n))\\
&-\E\bigg(\frac{1}{h_j^{L_j}}\bigg|\frac{\partial}{\partial u_{jl}}K_j\bigg(\frac{\|x_j-u_j\|_j}{h_j}\bigg)\bigg|_{u_j=\xi_j}\bigg|\frac{1}{h_k^{L_k}}K_k\bigg(\frac{\|x_k-\xi_k\|_k}{h_k}\bigg)\I(\mxi\in D(n))\bigg)\bigg).
\end{align*}
It holds that
\begin{align*}
&\E(Z_n^i(x_j,x_k))=0,\quad|Z_n^i(x_j,x_k)|\leq\const\cdot n^{-1}\cdot h_j^{-L_j-1}\cdot h_k^{-L_k},\quad\mbox{and}\\
&\sum_{i=1}^n\E(|Z_n^i(x_j,x_k)|^2)\\
&\leq n^{-1}\cdot \E\bigg(\bigg(\frac{1}{h_j^{L_j}}\frac{\partial}{\partial u_{jl}}K_j\bigg(\frac{\|x_j-u_j\|_j}{h_j}\bigg)\bigg|_{u_j=\xi_j}\bigg)^2\bigg(\frac{1}{h_k^{L_k}}K_k\bigg(\frac{\|x_k-\xi_k\|_k}{h_k}\bigg)\bigg)^2\I(\mxi\in D(n))\bigg)\\
&=O(n^{-1}\cdot h_j^{-L_j-2}\cdot h_k^{-L_k}\cdot V_{njk}).
\end{align*}
Now, by arguing as in the proof of (\ref{concentration of 2nd pj}) and using Theorem 2.6.2 in Bosq (2000), we may prove that the second term on the right hand side of (\ref{pjk 2nd term}) has the uniform rate
\ba
O_p(n^{-1/2}\cdot(\log{n})^{1/2}\cdot V_{njk}^{1/2}\cdot h_j^{-L_j/2-1}\cdot h_k^{-L_k/2}\cdot a_{nj})
\ea
over $(x_j,x_k)\in D_j\times D_k$ provided that $n^{-1/2}\cdot(\log{n})^{1/2}\cdot V_{njk}^{-1/2}\cdot h_j^{-L_j/2}\cdot h_k^{-L_k/2}=O(1)$. For the third term on the right hand side of (\ref{pjk 2nd term}), we note that
\ba
&\sup_{(x_j,x_k)\in D_j\times D_k}\E\bigg(\frac{1}{h_j^{L_j}}\bigg|\frac{\partial}{\partial u_{jl}}K_j\bigg(\frac{\|x_j-u_j\|_j}{h_j}\bigg)\bigg|_{u_j=\xi_j}\bigg|\frac{1}{h_k^{L_k}}K_k\bigg(\frac{\|x_k-\xi_k\|_k}{h_k}\bigg)\I(\mxi\in D(n))\bigg)\\
&=O(h_j^{-1}\cdot V_{njk}).
\ea
Hence, the third term on the right hand side of (\ref{pjk 2nd term}) has the uniform rate
\ba
O(h_j^{-1}\cdot a_{nj}\cdot V_{njk})
\ea
over $(x_j,x_k)\in D_j\times D_k$. Thus,
\ba
&\bigg|\frac{1}{n}\sum_{i=1}^n\frac{1}{h_j^{L_j}}\frac{\partial}{\partial u_{jl}}K_j\bigg(\frac{\|x_j-u_j\|_j}{h_j}\bigg)\bigg|_{u_j=\bar{\xi}_j^i}(\tilde{\xi}_{jl}^i-\xi_{jl}^i)\frac{1}{h_k^{L_k}}K_k\bigg(\frac{\|x_k-\xi_k^i\|_k}{h_k}\bigg)\I(\mxi^i\in D(n))\bigg|\\
&=O_p(h_j^{-L_j-2}\cdot a^2_{nj}\cdot V_{nk}+n^{-1/2}\cdot(\log{n})^{1/2}\cdot V_{njk}^{1/2}\cdot h_j^{-L_j/2-1}\cdot h_k^{-L_k/2}\cdot a_{nj}+h_j^{-1}\cdot a_{nj}\cdot V_{njk})
\ea
uniformly for $(x_j,x_k)\in D_j\times D_k$ since $V_{nk}^{1/2}\cdot h_j^{-L_j-2}\cdot a_{nj}^2=O(V_{njk}^{1/2}\cdot h_j^{-L_j/2-1}\cdot a_{nj})$. Similarly, we may show that
\ba
&\bigg|\frac{1}{n}\sum_{i=1}^n\frac{1}{h_j^{L_j}}K_j\bigg(\frac{\|x_j-\xi_j^i\|_j}{h_j}\bigg)\frac{1}{h_k^{L_k}}\frac{\partial}{\partial u_{kl}}K_k\bigg(\frac{\|x_k-u_k\|_k}{h_k}\bigg)\bigg|_{u_k=\bar{\xi}_k^i}(\tilde{\xi}_{kl}^i-\xi_{kl}^i)\I(\mxi^i\in D(n))\bigg|\\
&=O_p(h_k^{-L_k-2}\cdot a^2_{nk}\cdot V_{nj}+n^{-1/2}\cdot(\log{n})^{1/2}\cdot V_{njk}^{1/2}\cdot h_j^{-L_j/2}\cdot h_k^{-L_k/2-1}\cdot a_{nk}+h_k^{-1}\cdot a_{nk}\cdot V_{njk})
\ea
uniformly for $(x_j,x_k)\in D_j\times D_k$. We now approximate
\ba
&\bigg|\frac{1}{n}\sum_{i=1}^n\frac{1}{h_j^{L_j}}\frac{\partial}{\partial u_{jl}}K_j\bigg(\frac{\|x_j-u_j\|_j}{h_j}\bigg)\bigg|_{u_j=\bar{\xi}_j^i}(\tilde{\xi}_{jl}^i-\xi_{jl}^i)\frac{1}{h_k^{L_k}}\frac{\partial}{\partial u_{kl'}}K_k\bigg(\frac{\|x_k-u_k\|_k}{h_k}\bigg)\bigg|_{u_k=\bar{\xi}_k^i}(\tilde{\xi}_{kl'}^i-\xi_{kl'}^i)\\
&\qquad\qquad\cdot\I(\mxi^i\in D(n))\bigg|\\
&\leq\frac{1}{n}\sum_{i=1}^n\frac{1}{h_j^{L_j}}\bigg|\frac{\partial}{\partial u_{jl}}K_j\bigg(\frac{\|x_j-u_j\|_j}{h_j}\bigg)\bigg|_{u_j=\bar{\xi}_j^i}-\frac{\partial}{\partial u_{jl}}K_j\bigg(\frac{\|x_j-u_j\|_j}{h_j}\bigg)\bigg|_{u_j=\xi_j^i}\\
&\qquad\qquad\qquad\qquad+\frac{\partial}{\partial u_{jl}}K_j\bigg(\frac{\|x_j-u_j\|_j}{h_j}\bigg)\bigg|_{u_j=\xi_j^i}\bigg||\tilde{\xi}_{jl}^i-\xi_{jl}^i|\\
&\qquad\qquad\cdot\frac{1}{h_k^{L_k}}\bigg|\frac{\partial}{\partial u_{kl'}}K_k\bigg(\frac{\|x_k-u_k\|_k}{h_k}\bigg)\bigg|_{u_k=\bar{\xi}_k^i}-\frac{\partial}{\partial u_{kl'}}K_k\bigg(\frac{\|x_k-u_k\|_k}{h_k}\bigg)\bigg|_{u_k=\xi_k^i}\\
&\qquad\qquad\qquad\qquad+\frac{\partial}{\partial u_{kl'}}K_k\bigg(\frac{\|x_k-u_k\|_k}{h_k}\bigg)\bigg|_{u_k=\xi_k^i}\bigg||\tilde{\xi}_{kl'}^i-\xi_{kl'}^i|\I(\mxi^i\in D(n))
\ea
uniformly for $(x_j,x_k)\in D_j\times D_k$. We note that
\ba
&\frac{1}{n}\sum_{i=1}^n\frac{1}{h_j^{L_j}}\bigg|\frac{\partial}{\partial u_{jl}}K_j\bigg(\frac{\|x_j-u_j\|_j}{h_j}\bigg)\bigg|_{u_j=\bar{\xi}_j^i}-\frac{\partial}{\partial u_{jl}}K_j\bigg(\frac{\|x_j-u_j\|_j}{h_j}\bigg)\bigg|_{u_j=\xi_j^i}\bigg||\tilde{\xi}_{jl}^i-\xi_{jl}^i|\\
&\cdot\frac{1}{h_k^{L_k}}\bigg|\frac{\partial}{\partial u_{kl'}}K_k\bigg(\frac{\|x_k-u_k\|_k}{h_k}\bigg)\bigg|_{u_k=\bar{\xi}_k^i}-\frac{\partial}{\partial u_{kl'}}K_k\bigg(\frac{\|x_k-u_k\|_k}{h_k}\bigg)\bigg|_{u_k=\xi_k^i}\bigg||\tilde{\xi}_{kl'}^i-\xi_{kl'}^i|\I(\mxi^i\in D(n))\\
&=O_p(h_j^{-L_j-2}\cdot a_{nj}^2\cdot h_k^{-L_k-2}\cdot a_{nk}^2)\cdot O_p\bigg(\sum_{j=1}^da_{nj}+n^{-1/2}\cdot\bigg(\sum_{j=1}^da_{nj}\bigg)^{1/2}\bigg)
\ea
uniformly for $(x_j,x_k)\in D_j\times D_k$. We also note that
\ba
&\frac{1}{n}\sum_{i=1}^n\frac{1}{h_j^{L_j}}\bigg|\frac{\partial}{\partial u_{jl}}K_j\bigg(\frac{\|x_j-u_j\|_j}{h_j}\bigg)\bigg|_{u_j=\bar{\xi}_j^i}-\frac{\partial}{\partial u_{jl}}K_j\bigg(\frac{\|x_j-u_j\|_j}{h_j}\bigg)\bigg|_{u_j=\xi_j^i}\bigg||\tilde{\xi}_{jl}^i-\xi_{jl}^i|\\
&\qquad\qquad\cdot\frac{1}{h_k^{L_k}}\bigg|\frac{\partial}{\partial u_{kl'}}K_k\bigg(\frac{\|x_k-u_k\|_k}{h_k}\bigg)\bigg|_{u_k=\xi_k^i}\bigg||\tilde{\xi}_{kl'}^i-\xi_{kl'}^i|\I(\mxi^i\in D(n))\\
&=O_p(h_j^{-L_j-2}\cdot a_{nj}^2)\cdot O_p(n^{-1/2}\cdot(\log{n})^{1/2}\cdot V_{nk}^{1/2}\cdot h_k^{-L_k/2-1}\cdot a_{nk}+h_k^{-1}\cdot a_{nk}\cdot V_{nk})
\ea
and
\ba
&\frac{1}{n}\sum_{i=1}^n\frac{1}{h_j^{L_j}}\bigg|\frac{\partial}{\partial u_{jl}}K_j\bigg(\frac{\|x_j-u_j\|_j}{h_j}\bigg)\bigg|_{u_j=\xi_j^i}\bigg||\tilde{\xi}_{jl}^i-\xi_{jl}^i|\\
&\cdot\frac{1}{h_k^{L_k}}\bigg|\frac{\partial}{\partial u_{kl'}}K_k\bigg(\frac{\|x_k-u_k\|_k}{h_k}\bigg)\bigg|_{u_k=\bar{\xi}_k^i}-\frac{\partial}{\partial u_{kl'}}K_k\bigg(\frac{\|x_k-u_k\|_k}{h_k}\bigg)\bigg|_{u_k=\xi_k^i}\bigg||\tilde{\xi}_{kl'}^i-\xi_{kl'}^i|\I(\mxi^i\in D(n))\\
&=O_p(h_k^{-L_k-2}\cdot a_{nk}^2)\cdot O_p(n^{-1/2}\cdot(\log{n})^{1/2}\cdot V_{nj}^{1/2}\cdot h_j^{-L_j/2-1}\cdot a_{nj}+h_j^{-1}\cdot a_{nj}\cdot V_{nj})
\ea
uniformly for $(x_j,x_k)\in D_j\times D_k$ by arguing as in the proof of (\ref{pj 2nd rate}).
We note that
\begin{align}\label{4th term pjk}
\begin{split}
&\sup_{(x_j,x_k)\in D_j\times D_k}|F_{nll'}(x_j,x_k)|\\
&\leq\sup_{(x_j,x_k)\in D_j\times D_k}|F_{nll'}(x_j,x_k)-\E(F_{nll'}(x_j,x_k))|+\sup_{(x_j,x_k)\in D_j\times D_k}|\E(F_{nll'}(x_j,x_k))|,
\end{split}
\end{align}
where
\ba
F_{nll'}(x_j,x_k)&=\frac{1}{n}\sum_{i=1}^n\frac{1}{h_j^{L_j}}\bigg|\frac{\partial}{\partial u_{jl}}K_j\bigg(\frac{\|x_j-u_j\|_j}{h_j}\bigg)\bigg|_{u_j=\xi_j^i}\bigg||\tilde{\xi}_{jl}^i-\xi_{jl}^i|\\
&\qquad\qquad\cdot\frac{1}{h_k^{L_k}}\bigg|\frac{\partial}{\partial u_{kl'}}K_k\bigg(\frac{\|x_k-u_k\|_k}{h_k}\bigg)\bigg|_{u_k=\xi_k^i}\bigg||\tilde{\xi}_{kl'}^i-\xi_{kl'}^i|\I(\mxi^i\in D(n)).
\ea
We first approximate the first term on the right hand side of (\ref{4th term pjk}). For this, we define
\begin{align*}
&Z_n^i(x_j,x_k)\\
&=\frac{1}{n}\bigg(\frac{1}{h_j^{L_j}}\bigg|\frac{\partial}{\partial u_{jl}}K_j\bigg(\frac{\|x_j-u_j\|_j}{h_j}\bigg)\bigg|_{u_j=\xi_j^i}\bigg|\frac{1}{h_k^{L_k}}\bigg|\frac{\partial}{\partial u_{kl'}}K_k\bigg(\frac{\|x_k-u_k\|_k}{h_k}\bigg)\bigg|_{u_k=\xi_k^i}\bigg|\I(\mxi^i\in D(n))\\
&\quad-\E\bigg(\frac{1}{h_j^{L_j}}\bigg|\frac{\partial}{\partial u_{jl}}K_j\bigg(\frac{\|x_j-u_j\|_j}{h_j}\bigg)\bigg|_{u_j=\xi_j}\bigg|\frac{1}{h_k^{L_k}}\bigg|\frac{\partial}{\partial u_{kl'}}K_k\bigg(\frac{\|x_k-u_k\|_k}{h_k}\bigg)\bigg|_{u_k=\xi_k}\bigg|\I(\mxi\in D(n))\bigg)\bigg).
\end{align*}
It holds that
\begin{align*}
&\E(Z_n^i(x_j,x_k))=0,\quad|Z_n^i(x_j,x_k)|\leq\const\cdot n^{-1}\cdot h_j^{-L_j-1}\cdot h_k^{-L_k-1},\quad\mbox{and}\\
&\sum_{i=1}^n\E(|Z_n^i(x_j,x_k)|^2)\\
&\leq n^{-1}\cdot\E\bigg(\bigg(\frac{1}{h_j^{L_j}}\frac{\partial}{\partial u_{jl}}K_j\bigg(\frac{\|x_j-u_j\|_j}{h_j}\bigg)\bigg|_{u_j=\xi^i_j}\frac{1}{h_k^{L_k}}\frac{\partial}{\partial u_{kl'}}K_k\bigg(\frac{\|x_k-u_k\|_k}{h_k}\bigg)\bigg|_{u_k=\xi^i_k}\bigg)^2\\
&\qquad\qquad\qquad\cdot\I(\mxi^i\in D(n))\bigg)\\
&=O(n^{-1}\cdot h_j^{-L_j-2}\cdot h_k^{-L_k-2}\cdot V_{njk}).
\end{align*}
Now, by arguing as in the proof of (\ref{concentration of 2nd pj}) and using Theorem 2.6.2 in Bosq (2000), we may prove that the first term on the right hand side of (\ref{4th term pjk}) has the rate
\ba
O_p(n^{-1/2}\cdot(\log{n})^{1/2}\cdot V_{njk}^{1/2}\cdot h_j^{-L_j/2-1}\cdot h_k^{-L_k/2-1}\cdot a_{nj}\cdot a_{nk}),
\ea
provided that $n^{-1/2}\cdot(\log{n})^{1/2}\cdot V_{njk}^{-1/2}\cdot h_j^{-L_j/2}\cdot h_k^{-L_k/2}=O(1)$. For the second term on the right hand side of (\ref{4th term pjk}), we note that
\ba
&\sup_{(x_j,x_k)\in D_j\times D_k}\E\bigg(\frac{1}{h_j^{L_j}}\bigg|\frac{\partial}{\partial u_{jl}}K_j\bigg(\frac{\|x_j-u_j\|_j}{h_j}\bigg)\bigg|_{u_j=\xi^i_j}\bigg|\frac{1}{h_k^{L_k}}\bigg|\frac{\partial}{\partial u_{kl'}}K_j\bigg(\frac{\|x_k-u_k\|_k}{h_k}\bigg)\bigg|_{u_k=\xi^i_k}\bigg|\\
&\qquad\qquad\qquad\qquad\cdot\I(\mxi^i\in D(n))\bigg)\\
&=O(h_j^{-1}\cdot h_k^{-1}\cdot V_{njk}).
\ea
Hence, the second term on the right hand side of (\ref{4th term pjk}) has the rate
\ba
&O(h_j^{-1}\cdot a_{nj}\cdot h_k^{-1}\cdot a_{nk}\cdot V_{njk}).
\ea
Thus,
\ba
&\sup_{(x_j,x_k)\in D_j\times D_k}|F_{nll'}(x_j,x_k)|\\
&=O_p(n^{-1/2}\cdot(\log{n})^{1/2}\cdot V_{njk}^{1/2}\cdot h_j^{-L_j/2-1}\cdot h_k^{-L_k/2-1}\cdot a_{nj}\cdot a_{nk})+O(h_j^{-1}\cdot a_{nj}\cdot h_k^{-1}\cdot a_{nk}\cdot V_{njk}).
\ea
Combining the obtained rates gives that the second term on the right hand side of (\ref{decomposition jk}) has the rate
\begin{align}\label{2nd term in pjk}
O_p(n^{-1/2}\cdot(\log{n})^{1/2}\cdot V_{njk}^{1/2}\cdot h_j^{-L_j/2}\cdot h_k^{-L_k/2}+V_{njk}).
\end{align}
We now approximate the third term on the right hand side of (\ref{decomposition jk}). We note that
\begin{align}\label{pjk 3rd part}
\begin{split}
&\frac{1}{n}\sum_{i=1}^n|K_{h_j}(x_j,\tilde{\xi}_j^i)K_{h_k}(x_k,\tilde{\xi}_k^i)-K_{h_j}(x_j,\xi_j^i)K_{h_k}(x_k,\xi_k^i)|\I(\mxi^i\in D)\I(\tilde{\mxi}^i\in D)\\
&\leq\const\cdot\frac{1}{n}\sum_{i=1}^n\frac{1}{h_j^{L_j}}\frac{1}{h_k^{L_k}}\bigg|K_j\bigg(\frac{\|x_j-\tilde{\xi}_j^i\|_j}{h_j}\bigg)K_k\bigg(\frac{\|x_k-\tilde{\xi}_k^i\|_k}{h_k}\bigg)\\
&\qquad\qquad\qquad\qquad\qquad-K_j\bigg(\frac{\|x_j-\xi_j^i\|_j}{h_j}\bigg)K_k\bigg(\frac{\|x_k-\xi_k^i\|_k}{h_k}\bigg)\bigg|\I(\mxi^i\in D)\\
&\quad+\const\cdot\frac{1}{n}\sum_{i=1}^n\bigg|\int_{D_j}\frac{1}{h_j^{L_j}}K_j\bigg(\frac{\|t_j-\tilde{\xi}_j^i\|_j}{h_j}\bigg)dt_j\int_{D_k}\frac{1}{h_k^{L_k}}K_k\bigg(\frac{\|t_k-\tilde{\xi}_k^i\|_k}{h_k}\bigg)dt_k\\
&\qquad\qquad\qquad\qquad\qquad-\int_{D_j}\frac{1}{h_j^{L_j}}K_j\bigg(\frac{\|t_j-\xi_j^i\|_j}{h_j}\bigg)dt_j\int_{D_k}\frac{1}{h_k^{L_k}}K_k\bigg(\frac{\|t_k-\xi_k^i\|_k}{h_k}\bigg)dt_k\bigg|\\
&\qquad\qquad\qquad\qquad\qquad\cdot \frac{1}{h_j^{L_j}}K_j\bigg(\frac{\|x_j-\xi_j^i\|_j}{h_j}\bigg)\frac{1}{h_k^{L_k}}K_k\bigg(\frac{\|x_k-\xi_k^i\|_k}{h_k}\bigg)\I(\mxi^i\in D).
\end{split}
\end{align}
The first term on the right hand side of (\ref{pjk 3rd part}) has the uniform rate
\ba
O_p(h_j^{-L_j-2}\cdot a_{nj}^2+h_k^{-L_k-2}\cdot a_{nk}^2+h_j^{-1}\cdot a_{nj}+h_k^{-1}\cdot a_{nk})
\ea
over $(x_j,x_k)\in D_j\times D_k$ provided that $n^{-1/2}\cdot(\log{n})^{1/2}\cdot h_j^{-L_j/2}\cdot h_k^{-L_k/2}=O(1)$. This follows by arguing as in the proof of (\ref{2nd term in pjk}). The second term on the right hand side of (\ref{pjk 3rd part}) is bounded by
\ba
&\const\cdot\frac{1}{n}\sum_{i=1}^n\int_{D_j}\frac{1}{h_j^{L_j}}\bigg|K_j\bigg(\frac{\|t_j-\tilde{\xi}_j^i\|_j}{h_j}\bigg)-K_j\bigg(\frac{\|t_j-\xi_j^i\|_j}{h_j}\bigg)\bigg|dt_j\\
&\qquad\qquad\qquad\qquad\qquad\cdot \frac{1}{h_j^{L_j}}K_j\bigg(\frac{\|x_j-\xi_j^i\|_j}{h_j}\bigg)\frac{1}{h_k^{L_k}}K_k\bigg(\frac{\|x_k-\xi_k^i\|_k}{h_k}\bigg)\I(\mxi^i\in D)\\
+&\const\cdot\frac{1}{n}\sum_{i=1}^n\int_{D_k}\frac{1}{h_k^{L_k}}\bigg|K_k\bigg(\frac{\|t_k-\tilde{\xi}_k^i\|_k}{h_k}\bigg)-K_k\bigg(\frac{\|t_k-\xi_k^i\|_k}{h_k}\bigg)\bigg|dt_k\\
&\qquad\qquad\qquad\qquad\qquad\cdot \frac{1}{h_j^{L_j}}K_j\bigg(\frac{\|x_j-\xi_j^i\|_j}{h_j}\bigg)\frac{1}{h_k^{L_k}}K_k\bigg(\frac{\|x_k-\xi_k^i\|_k}{h_k}\bigg)\I(\mxi^i\in D)
\ea
and it has the uniform rate
\ba
O_p(h_j^{-1}\cdot a_{nj}+h_k^{-1}\cdot a_{nk})
\ea
over $(x_j,x_k)\in D_j\times D_k$ provided that $n^{-1/2}\cdot(\log{n})^{1/2}\cdot h_j^{-L_j/2}\cdot h_k^{-L_k/2}=O(1)$.
Hence, the third term on the right hand side of (\ref{decomposition jk}) has the rate
\ba
O_p(h_j^{-L_j-2}\cdot a_{nj}^2+h_k^{-L_k-2}\cdot a_{nk}^2+h_j^{-1}\cdot a_{nj}+h_k^{-1}\cdot a_{nk}).
\ea
This completes the proof for the second part of the lemma.

The first part of the lemma follows similarly as in the above proof by obtaining the rates with $V_{nj}^{1/2}$, $V_{nk}^{1/2}$ and $V_{njk}^{1/2}$ being replaced by 1 using the bound $\I(\mxi^i\in D\cap D(n))\leq\I(\mxi^i\in D)$ or $\I(\mxi^i\in D(n))\leq\I(\mxi^i\in\prod_{j=1}^dD_j^+(C\cdot a_{nj}))$ whenever $V_{nj}^{1/2}$, $V_{nk}^{1/2}$ and $V_{njk}^{1/2}$ appear. In this case, the condition that $n^{-1/2}\cdot(\log{n})^{1/2}\cdot h_j^{-L_j/2}\cdot h_k^{-L_k/2}=O(1)$ is required instead of the condition that $n^{-1/2}\cdot(\log{n})^{1/2}\cdot V_{njk}^{-1/2}\cdot h_j^{-L_j/2}\cdot h_k^{-L_k/2}=O(1)$.
\end{proof}

Now, we are ready to state the uniform consistency of our marginal density estimators.

\begin{customlemma}{S.13}\label{marginal density convergence 1}
Assume that conditions (B3)--(B4) and (B6) hold, that $h_j=o(1)$, that $n^{-1/2}\cdot(\log{n})^{1/2}\cdot h_j^{-L_j/2}=o(1)$, that $p^D_j$ is continuous on $D_j$ and that $p$ is bounded on $D^+(\varepsilon)$ for some $\varepsilon>0$. Then, it holds that 
\ba
\sup_{x_j\in D_j}\bigg|\hat{p}^D_j(x_j)-p^D_j(x_j)\int_{D_j}K_{h_j}(x_j,u_j)du_j\bigg|&=o_p(1),\quad\mbox{and}\\
\sup_{x_j\in D_j^-(2h_j)}|\hat{p}^D_j(x_j)-p^D_j(x_j)|&=o_p(1).
\ea
\end{customlemma}

\begin{proof}
Note that
\ba
&\bigg|\hat{p}^D_j(x_j)-p^D_j(x_j)\int_{D_j}K_{h_j}(x_j,u_j)du_j\bigg|\\
&\leq|\hat{p}^D_j(x_j)-\check{p}^D_j(x_j)|\\
&\quad+\bigg|(\check{p}_0^D)^{-1}\bigg(n^{-1}\sum_{i=1}^nK_{h_j}(x_j,\xi_j^i)\I(\mxi^i\in D)-\E\bigg(n^{-1}\sum_{i=1}^nK_{h_j}(x_j,\xi_j^i)\I(\mxi^i\in D)\bigg)\bigg)\bigg|\\
&\quad+\bigg|(\check{p}_0^D)^{-1}\E\bigg(n^{-1}\sum_{i=1}^nK_{h_j}(x_j,\xi_j^i)\I(\mxi^i\in D)\bigg)-(p_0^D)^{-1}p_0^Dp^D_j(x_j)\int_{D_j}K_{h_j}(x_j,u_j)du_j\bigg|\\
&\leq|\hat{p}^D_j(x_j)-\check{p}^D_j(x_j)|\\
&\quad+\bigg|(\check{p}_0^D)^{-1}\bigg(n^{-1}\sum_{i=1}^nK_{h_j}(x_j,\xi_j^i)\I(\mxi^i\in D)-\E\bigg(n^{-1}\sum_{i=1}^nK_{h_j}(x_j,\xi_j^i)\I(\mxi^i\in D)\bigg)\bigg)\bigg|\\
&\quad+\bigg|(\check{p}_0^D)^{-1}\bigg(\E\bigg(n^{-1}\sum_{i=1}^nK_{h_j}(x_j,\xi_j^i)\I(\mxi^i\in D)\bigg)-p_0^Dp^D_j(x_j)\int_{D_j}K_{h_j}(x_j,u_j)du_j\bigg)\bigg|\\
&\quad+\bigg|(\check{p}_0^D)^{-1}p^D_j(x_j)\int_{D_j}K_{h_j}(x_j,u_j)du_j(\check{p}_0^D-p_0^D)\bigg|.
\ea
The first term on the right hand side of the above inequality has the uniform rate
\ba
O_p\left(\sqrt{\frac{\log{n}}{n\cdot h_j^{L_j}}}+h_j^{-L_j}\cdot a_{nj}+\sum_{k\neq j}a_{nk}+h_j^{-L_j-2}\cdot a_{nj}^2\right)=o_p(1)
\ea
over $x_j\in D_j$ by the first part of Lemma \ref{marginal density approximation 1}.
We note that $(\check{p}_0^D)^{-1}$ is bounded away from zero with probability tending to one by Lemma \ref{p0 rate}. Hence, the second term on the right hand side has the uniform rate
\ba
O_p(n^{-1/2}\cdot(\log{n})^{1/2}\cdot h_j^{-L_j})=o_p(1)
\ea
over $x_j\in D_j$ by Lemma \ref{uniform rates density}. Also, the third term on the right hand side has the uniform rate $o_p(1)$ over $x_j\in D_j$ by Lemma \ref{uniform consistency general}. The fourth term has the uniform rate $O_p(n^{-1/2})=o_p(1)$ over $x_j\in D_j$ by (\ref{nominator bound}) and the proof of Lemma \ref{p0 rate}. Hence, the first assertion follows. The second assertion follows from (\ref{integral one in interior}) and the first assertion.
\end{proof}

\begin{customlemma}{S.14}\label{marginal density convergence 2}
Assume that conditions (B3)--(B4) and (B6) hold, that $h_j\vee h_k=o(1)$, that $n^{-1/2}\cdot(\log{n})^{1/2}\cdot h_j^{-L_j/2}\cdot h_k^{-L_k/2}=o(1)$, that $p^D_{jk}$ is continuous on $D_j\times D_k$ and that $p$ is bounded on $D^+(\varepsilon)$ for some $\varepsilon>0$. Then, it holds that 
\ba
\sup_{(x_j,x_k)\in D_j\times D_k}\bigg|\hat{p}^D_{jk}(x_j,x_k)-p^D_{jk}(x_j,x_k)\int_{D_j}K_{h_j}(x_j,u_j)du_j\int_{D_k}K_{h_k}(x_k,u_k)du_k\bigg|&=o_p(1),\quad\mbox{and}\\
\sup_{(x_j,x_k)\in D_j^-(2h_j)\times D_k^-(2h_k)}|\hat{p}^D_{jk}(x_j,x_k)-p^D_{j,k}(x_j,x_k)|&=o_p(1).
\ea
\end{customlemma}

\begin{proof}
Since the proof is similar to that of Lemma \ref{marginal density convergence 1}, we only it. Note that
\ba
&\bigg|\hat{p}^D_{jk}(x_j,x_k)-p^D_{jk}(x_j,x_k)\int_{D_j}K_{h_j}(x_j,u_j)du_j\int_{D_k}K_{h_k}(x_k,u_k)du_k\bigg|\\
&\leq|\hat{p}^D_{jk}(x_j,x_k)-\check{p}^D_{jk}(x_j,x_k)|\\
&\quad+\bigg|(\check{p}_0^D)^{-1}\bigg(n^{-1}\sum_{i=1}^nK_{h_j}(x_j,\xi_j^i)K_{h_k}(x_k,\xi_k^i)\I(\mxi^i\in D)\\
&\qquad\qquad\qquad-\E\bigg(n^{-1}\sum_{i=1}^nK_{h_j}(x_j,\xi_j^i)K_{h_k}(x_k,\xi_k^i)\I(\mxi^i\in D)\bigg)\bigg)\bigg|\\
&\quad+\bigg|(\check{p}_0^D)^{-1}\bigg(\E\bigg(n^{-1}\sum_{i=1}^nK_{h_j}(x_j,\xi_j^i)K_{h_k}(x_k,\xi_k^i)\I(\mxi^i\in D)\bigg)\\
&\qquad\qquad\qquad-p_0^D p^D_{jk}(x_j,x_k)\int_{D_j}K_{h_j}(x_j,u_j)du_j\int_{D_k}K_{h_k}(x_k,u_k)du_k\bigg)\bigg|\\
&\quad+\bigg|(\check{p}_0^D)^{-1}p^D_{jk}(x_j,x_k)\int_{D_j}K_{h_j}(x_j,u_j)du_j\int_{D_k}K_{h_k}(x_k,u_k)du_k(\check{p}_0^D-p_0^D)\bigg|.
\ea
The first term on the right hand side of the above inequality has the uniform rate $o_p(1)$ over $(x_j,x_k)\in D_j\times D_k$ by the first part of Lemma \ref{marginal density approximation 2}.
The second term on the right hand side has the uniform rate $o_p(1)$ over $(x_j,x_k)\in D_j\times D_k$ by Lemma \ref{uniform rates density}. The third term on the right hand side has the uniform rate $o_p(1)$ over $(x_j,x_k)\in D_j\times D_k$ by Lemma \ref{uniform consistency general}. The fourth term has the uniform rate $o_p(1)$ over $(x_j,x_k)\in D_j\times D_k$ by (\ref{nominator bound}) and the proof of Lemma \ref{p0 rate}. Hence, the first assertion follows. The second assertion follows from (\ref{integral one in interior}) and the first assertion.
\end{proof}

The next lemma is useful to construct the uniform consistency of our marginal regression estimators.

\begin{customlemma}{S.15}\label{marginal regression approximation}
Assume that conditions (B1) and (B3)--(B4) hold, that $h_j=O(1)$, that $n^{-1/2}\cdot(\log{n})^{1/2}\cdot h_j^{-L_j/2}=O(1)$, that $n^{-1+2\beta+2/\alpha}\cdot h_j^{-L_j}=o(1)$ for some constant $\beta>0$, that $h_j^{-1}\cdot a_{nj}=O(1)$, that $h_j^{-L_j-2}\cdot a^2_{nj}=O(1)$ and that $p$ is bounded on $D^+(\varepsilon)$ for some $\varepsilon>0$. Then, it holds that
\ba
&\sup_{x_j\in D_j}\bigg\|n^{-1}\od\bigoplus_{i=1}^n((K_{h_j}(x_j,\tilde{\xi}_j^i)\I(\tilde{\mxi}^i\in D))\od\tilde{\mY}^i)\om n^{-1}\od\bigoplus_{i=1}^n((K_{h_j}(x_j,\xi_j^i)\I(\mxi^i\in D))\od\mY^i)\bigg\|\\
&=O_p\left(\sqrt{\frac{\log{n}}{n\cdot h_j^{L_j}}}+h_j^{-L_j}\cdot a_{nj}+\sum_{k\neq j}a_{nk}+h_j^{-L_j-2}\cdot a_{nj}^2+b_n\right).
\ea
If we further assume that $n^{-\beta}\cdot(\log{n})^{1/2}\cdot(\sum_{j=1}^da_{nj})^{1/\alpha}\cdot(h_j^{-L_j}\cdot a_{nj}+\sum_{k\neq j}a_{nk})^{-1/2}=O(1)$, then it holds that
\ba
&\sup_{x_j\in D_j}\bigg\|n^{-1}\od\bigoplus_{i=1}^n((K_{h_j}(x_j,\tilde{\xi}_j^i)\I(\tilde{\mxi}^i\in D))\od\tilde{\mY}^i)\om n^{-1}\od\bigoplus_{i=1}^n((K_{h_j}(x_j,\xi_j^i)\I(\mxi^i\in D))\od\mY^i)\bigg\|\\
&=O_p\bigg(h_j^{-L_j}\cdot a_{nj}+\sum_{k\neq j}a_{nk}+h_j^{-L_j-2}\cdot a_{nj}^2+b_n\bigg).
\ea
\end{customlemma}

\begin{proof}
We first prove the second part of the lemma. We note that
\begin{align}\label{decomposition in marginal regression}
\begin{split}
&\bigg\|n^{-1}\od\bigoplus_{i=1}^n((K_{h_j}(x_j,\tilde{\xi}_j^i)\I(\tilde{\mxi}^i\in D))\od\tilde{\mY}^i)\om n^{-1}\od\bigoplus_{i=1}^n((K_{h_j}(x_j,\xi_j^i)\I(\mxi^i\in D))\od\mY^i)\bigg\|\\
&\leq\bigg\|n^{-1}\od\bigoplus_{i=1}^n(K_{h_j}(x_j,\xi_j^i)\I(\mxi^i\in D)(\I(\tilde{\mxi}^i\in D)-\I(\mxi^i\in D)))\od\mY^i\bigg\|\\
&\quad+\bigg\|n^{-1}\od\bigoplus_{i=1}^n(K_{h_j}(x_j,\tilde{\xi}_j^i)\I(\tilde{\mxi}^i\in D)(\I(\tilde{\mxi}^i\in D)-\I(\mxi^i\in D)))\od\tilde{\mY}^i\bigg\|\\
&\quad+\bigg\|n^{-1}\od\bigoplus_{i=1}^n(K_{h_j}(x_j,\tilde{\xi}_j^i)\od\tilde{\mY}^i\om K_{h_j}(x_j,\xi_j^i)\od\mY^i)\od(\I(\mxi^i\in D)\I(\tilde{\mxi}^i\in D))\bigg\|.
\end{split}
\end{align}
We first approximate the first term on the right hand side of (\ref{decomposition in marginal regression}) uniformly for $x_j\in D_j$. Define $V_{nj}=h_j^{-L_j}\cdot a_{nj}+\sum_{k\neq j}a_{nk}$. Note that
\begin{align}\label{some inequality}
\begin{split}
&\sup_{x_j\in D_j}\bigg\|n^{-1}\od\bigoplus_{i=1}^n(K_{h_j}(x_j,\xi_j^i)\I(\mxi^i\in D)(\I(\tilde{\mxi}^i\in D)-\I(\mxi^i\in D)))\od\mY^i\bigg\|\\
&\leq\const\cdot\sup_{x_j\in D_j}n^{-1}\sum_{i=1}^n\frac{1}{h_j^{L_j}}K_j\bigg(\frac{\|x_j-\xi_j^i\|_j}{h_j}\bigg)\I(\mxi^i\in D(n))\|\mY^i\|\\
&\leq\const\cdot\sup_{x_j\in D_j}\bigg|n^{-1}\sum_{i=1}^n\frac{1}{h_j^{L_j}}K_j\bigg(\frac{\|x_j-\xi_j^i\|_j}{h_j}\bigg)\I(\mxi^i\in D(n))\|\mY^i\|\\
&\qquad\qquad\qquad\qquad-\E\bigg(\frac{1}{h_j^{L_j}}K_j\bigg(\frac{\|x_j-\xi_j\|_j}{h_j}\bigg)\I(\mxi\in D(n))\|\mY\|\bigg)\bigg|\\
&\quad+\const\cdot\sup_{x_j\in D_j}\E\bigg(\frac{1}{h_j^{L_j}}K_j\bigg(\frac{\|x_j-\xi_j\|_j}{h_j}\bigg)\I(\mxi\in D(n))\E(\|\mY\||\mxi)\bigg)\\
&=O_p(n^{-1/2}\cdot(\log{n})^{1/2}\cdot V_{nj}^{-1/2}\cdot h_j^{-L_j/2}+V_{nj})
\end{split}
\end{align}
provided that $n^{-1+2\beta+2/\alpha}\cdot h_j^{-L_j}=o(1)$ and $n^{-\beta}\cdot(\log{n})^{1/2}\cdot(\sum_{j=1}^da_{nj})^{1/\alpha}\cdot V_{nj}^{-1/2}=O(1)$, where the equality in (\ref{some inequality}) follows from Lemma \ref{uniform rates Dn}. The second term on the right hand side of (\ref{decomposition in marginal regression}) satisfies
\begin{align}\label{second term in marginal regression}
\begin{split}
&\sup_{x_j\in D_j}\bigg\|n^{-1}\od\bigoplus_{i=1}^n(K_{h_j}(x_j,\tilde{\xi}_j^i)\I(\tilde{\mxi}^i\in D)(\I(\tilde{\mxi}^i\in D)-\I(\mxi^i\in D)))\od\tilde{\mY}^i\bigg\|\\
&\leq\sup_{x_j\in D_j}n^{-1}\sum_{i=1}^nK_{h_j}(x_j,\tilde{\xi}_j^i)\I(\tilde{\mxi}^i\in D)\I(\mxi^i\in D(n))\|\tilde{\mY}^i\|\\
&\leq\const\cdot\sup_{x_j\in D_j}n^{-1}\sum_{i=1}^n\frac{1}{h_j^{L_j}}K_j\bigg(\frac{\|x_j-\tilde{\xi}_j^i\|_j}{h_j}\bigg)\I(\mxi^i\in D(n))\|\tilde{\mY}^i\|\\
&\leq\const\cdot\max_{1\leq i\leq n}\|\tilde{\mY}^i\om\mY^i\|\cdot\sup_{x_j\in D_j}n^{-1}\sum_{i=1}^n\frac{1}{h_j^{L_j}}K_j\bigg(\frac{\|x_j-\tilde{\xi}_j^i\|_j}{h_j}\bigg)\I(\mxi^i\in D(n))\\
&\quad+\const\cdot\sup_{x_j\in D_j}n^{-1}\sum_{i=1}^n\frac{1}{h_j^{L_j}}K_j\bigg(\frac{\|x_j-\tilde{\xi}_j^i\|_j}{h_j}\bigg)\I(\mxi^i\in D(n))\|\mY^i\|.
\end{split}
\end{align}
Note that the first term on the right hand side of (\ref{second term in marginal regression}) has the rate
\ba
&O_p(b_n)\cdot O_p(n^{-1/2}\cdot(\log{n})^{1/2}\cdot h_j^{-L_j/2}+V_{nj})
\ea
provided that $n^{-1/2}\cdot(\log{n})^{1/2}\cdot h_j^{-L_j/2}=O(1)$ by arguing as in the proof for the first part of Lemma \ref{marginal density approximation 1}. We now approximate the second term on the right hand side of (\ref{second term in marginal regression}). Recall (\ref{K taylor}). It holds that
\ba
&\sup_{x_j\in D_j}\frac{1}{n}\sum_{i=1}^n\frac{1}{h_j^{L_j}}K_j\bigg(\frac{\|x_j-\xi_j^i\|_j}{h_j}\bigg)\I(\mxi^i\in D(n))\|\mY^i\|\\
&=O_p(n^{-1/2}\cdot(\log{n})^{1/2}\cdot V_{nj}^{1/2}\cdot h_j^{-L_j/2}+V_{nj}),
\ea
provided that $n^{-1+2\beta+2/\alpha}\cdot h_j^{-L_j}=o(1)$ and $n^{-\beta}\cdot(\log{n})^{1/2}\cdot(\sum_{j=1}^da_{nj})^{1/\alpha}\cdot V_{nj}^{-1/2}=O(1)$ by (\ref{some inequality}). Also, it holds that
\begin{align}\label{Y approximation pf}
\begin{split}
&\sup_{x_j\in D_j}\bigg|\frac{1}{n}\sum_{i=1}^n\frac{1}{h_j^{L_j}}\frac{\partial}{\partial u_{jl}}K_j\bigg(\frac{\|x_j-u_j\|_j}{h_j}\bigg)\bigg|_{u_j=\bar{\xi}_j^i}(\tilde{\xi}_{jl}^i-\xi_{jl}^i)\I(\mxi^i\in D(n))\|\mY^i\|\bigg|\\
&\leq\max_{1\leq i\leq n}|\tilde{\xi}_{jl}^i-\xi_{jl}^i|\sup_{x_j\in D_j}\frac{1}{n}\sum_{i=1}^n\frac{1}{h_j^{L_j}}\bigg|\frac{\partial}{\partial u_{jl}}K_j\bigg(\frac{\|x_j-u_j\|_j}{h_j}\bigg)\bigg|_{u_j=\bar{\xi}_j^i}\bigg|\I(\mxi^i\in D(n))\|\mY^i\|\\
&\leq\max_{1\leq i\leq n}|\tilde{\xi}_{jl}^i-\xi_{jl}^i|\sup_{x_j\in D_j}\frac{1}{n}\sum_{i=1}^n\frac{1}{h_j^{L_j}}\bigg|\frac{\partial}{\partial u_{jl}}K_j\bigg(\frac{\|x_j-u_j\|_j}{h_j}\bigg)\bigg|_{u_j=\bar{\xi}_j^i}-\frac{\partial}{\partial u_{jl}}K_j\bigg(\frac{\|x_j-u_j\|_j}{h_j}\bigg)\bigg|_{u_j=\xi_j^i}\bigg|\\
&\qquad\qquad\qquad\qquad\qquad\qquad\qquad\cdot\I(\mxi^i\in D(n))\|\mY^i\|\\
&\quad+\max_{1\leq i\leq n}|\tilde{\xi}_{jl}^i-\xi_{jl}^i|\sup_{x_j\in D_j}\bigg|\frac{1}{n}\sum_{i=1}^n\frac{1}{h_j^{L_j}}\bigg|\frac{\partial}{\partial u_{jl}}K_j\bigg(\frac{\|x_j-u_j\|_j}{h_j}\bigg)\bigg|_{u_j=\xi_j^i}\bigg|\I(\mxi^i\in D(n))\|\mY^i\|\\
&\qquad\qquad\qquad\qquad\qquad\qquad-\E\bigg(\frac{1}{h_j^{L_j}}\bigg|\frac{\partial}{\partial u_{jl}}K_j\bigg(\frac{\|x_j-u_j\|_j}{h_j}\bigg)\bigg|_{u_j=\xi_j}\bigg|\I(\mxi\in D(n))\|\mY\|\bigg)\bigg|\\
&\quad+\max_{1\leq i\leq n}|\tilde{\xi}_{jl}^i-\xi_{jl}^i|\sup_{x_j\in D_j}\E\bigg(\frac{1}{h_j^{L_j}}\bigg|\frac{\partial}{\partial u_{jl}}K_j\bigg(\frac{\|x_j-u_j\|_j}{h_j}\bigg)\bigg|_{u_j=\xi_j}\bigg|\I(\mxi\in D(n))\|\mY\|\bigg).
\end{split}
\end{align}
The first term on the right hand side of (\ref{Y approximation pf}) has the rate
\ba
O_p(h_j^{-L_j-2}\cdot a_{nj}^2)\cdot O_p\bigg(\sum_{j=1}^da_{nj}+n^{-1/2}\cdot\bigg(\sum_{j=1}^da_{nj}\bigg)^{1/2}\bigg).
\ea
We now approximate the second term on the right hand side of (\ref{Y approximation pf}). For $\kappa>0$, we let $\{x_j^{(l)}:1\leq l\leq N(n^{-\kappa})\}\subset D_j$ be a set of points such that $N(n^{-\kappa})=O(n^{\kappa L_j})$ and $\{B_j(x_j^{(l)},n^{-\kappa}):1\leq l\leq N(n^{-\kappa})\}$ covers $D_j$.
Define
\ba
U_n(x_j)&=\frac{1}{n}\sum_{i=1}^n\frac{1}{h_j^{L_j}}\bigg|\frac{\partial}{\partial u_{jl}}K_j\bigg(\frac{\|x_j-u_j\|_j}{h_j}\bigg)\bigg|_{u_j=\xi_j^i}\bigg|\I(\mxi^i\in D(n))\|\mY^i\|\\
&\qquad\qquad\cdot\I\bigg(\I(\mxi^i\in D(n))\|\mY^i\|\leq n^{1/2-\beta}\bigg(\sum_{k=1}^da_{nk}\bigg)^{1/\alpha}h_j^{L_j/2}\bigg).
\ea
It holds that
\begin{align*}
&\sup_{x_j\in D_j}\bigg|\frac{1}{n}\sum_{i=1}^n\frac{1}{h_j^{L_j}}\bigg|\frac{\partial}{\partial u_{jl}}K_j\bigg(\frac{\|x_j-u_j\|_j}{h_j}\bigg)\bigg|_{u_j=\xi_j^i}\bigg|\I(\mxi^i\in D(n))\|\mY^i\|\\
&\qquad\quad-\E\bigg(\frac{1}{h_j^{L_j}}\bigg|\frac{\partial}{\partial u_{jl}}K_j\bigg(\frac{\|x_j-u_j\|_j}{h_j}\bigg)\bigg|_{u_j=\xi_j}\bigg|\I(\mxi\in D(n))\|\mY\|\bigg)\bigg|\\
&\leq\max_{1\leq l\leq N(n^{-\kappa})}|U_n(x^{(l)}_j)-\E(U_n(x^{(l)}_j))|+o_p(n^{-1/2}\cdot(\log{n})^{1/2}\cdot V_{nj}^{1/2}\cdot h_j^{-L_j/2-1})
\end{align*}
provided that $n^{-1+2\beta+2/\alpha}\cdot h_j^{-L_j}=o(1)$ by taking sufficiently large $\kappa$. We claim that the first term on the right hand side of the above inequality has the rate
\begin{align}\label{claim Y}
O_p(n^{-1/2}\cdot(\log{n})^{1/2}\cdot V_{nj}^{1/2}\cdot h_j^{-L_j/2-1}).
\end{align}
To see this, define
\begin{align*}
Z_n^i(x_j)=&\frac{1}{n}\bigg(\frac{1}{h_j^{L_j}}\bigg|\frac{\partial}{\partial u_{jl}}K_j\bigg(\frac{\|x_j-u_j\|_j}{h_j}\bigg)\bigg|_{u_j=\xi_j^i}\bigg|\I(\mxi^i\in D(n))\|\mY^i\|\\
&\qquad\cdot\I\bigg(\I(\mxi^i\in D(n))\|\mY^i\|\leq n^{1/2-\beta}\bigg(\sum_{k=1}^da_{nk}\bigg)^{1/\alpha}h_j^{L_j/2}\bigg)\\
&-\E\bigg(\frac{1}{h_j^{L_j}}\bigg|\frac{\partial}{\partial u_{jl}}K_j\bigg(\frac{\|x_j-u_j\|_j}{h_j}\bigg)\bigg|_{u_j=\xi_j}\bigg|\I(\mxi\in D(n))\|\mY\|\\
&\qquad\cdot\I\bigg(\I(\mxi\in D(n))\|\mY\|\leq n^{1/2-\beta}\bigg(\sum_{k=1}^da_{nk}\bigg)^{1/\alpha}h_j^{L_j/2}\bigg)\bigg)\bigg).
\end{align*}
It holds that
\begin{align*}
&\E(Z_n^i(x_j))=0,\quad|Z_n^i(x_j)|\leq\const\cdot n^{-1/2-\beta}\cdot\bigg(\sum_{k=1}^da_{nk}\bigg)^{1/\alpha}\cdot h_j^{-L_j/2-1},\quad\mbox{and}\\
&\sum_{i=1}^n\E(|Z_n^i(x_j)|^2)\\
&\leq n^{-1}\cdot \E\bigg(\bigg(\frac{1}{h_j^{L_j}}\frac{\partial}{\partial u_{jl}}K_j\bigg(\frac{\|x_j-u_j\|_j}{h_j}\bigg)\bigg|_{u_j=\xi_j}\bigg)^2\I(\mxi\in D(n))\|\mY\|\bigg)\\
&=O(n^{-1}\cdot h_j^{-L_j-2}\cdot V_{nj}).
\end{align*}
Now, Theorem 2.6.2 in Bosq (2000) gives (\ref{claim Y}) provided that $n^{-\beta}\cdot(\log{n})^{1/2}\cdot(\sum_{j=1}^da_{nj})^{1/\alpha}\cdot V_{nj}^{-1/2}=O(1)$.
Hence, the second term on the right hand side of (\ref{Y approximation pf}) has the rate
\ba
O_p(n^{-1/2}\cdot(\log{n})^{1/2}\cdot V_{nj}^{1/2}\cdot h_j^{-L_j/2-1}\cdot a_{nj}).
\ea
The third term on the right hand side of (\ref{Y approximation pf}) has the rate
\ba
O_p(h_j^{-1}\cdot a_{nj}\cdot V_{nj})
\ea
by arguing as in the proof of (\ref{pj 2nd 2nd rate}). Hence, the second term on the right hand side of (\ref{decomposition in marginal regression}) has the uniform rate
\ba
O_p(n^{-1/2}\cdot(\log{n})^{1/2}\cdot V_{nj}^{1/2}\cdot h_j^{-L_j/2}+V_{nj})
\ea
over $x_j\in D_j$. The third term on the right hand side of (\ref{decomposition in marginal regression}) satisfies
\begin{align}\label{consistency final}
\begin{split}
&\sup_{x_j\in D_j}\bigg\|n^{-1}\od\bigoplus_{i=1}^n(K_{h_j}(x_j,\tilde{\xi}_j^i)\od\tilde{\mY}^i\om K_{h_j}(x_j,\xi_j^i)\od\mY^i)\od(\I(\mxi^i\in D)\I(\tilde{\mxi}^i\in D))\bigg\|\\
&\leq \sup_{x_j\in D_j}\bigg\|n^{-1}\od\bigoplus_{i=1}^n(K_{h_j}(x_j,\tilde{\xi}_j^i)\od(\tilde{\mY}^i\om\mY^i))\od(\I(\mxi^i\in D)\I(\tilde{\mxi}^i\in D))\bigg\|\\
&\quad+\sup_{x_j\in D_j}\bigg\|n^{-1}\od\bigoplus_{i=1}^n((K_{h_j}(x_j,\tilde{\xi}_j^i)-K_{h_j}(x_j,\xi_j^i))\od\mY^i)\od(\I(\mxi^i\in D)\I(\tilde{\mxi}^i\in D))\bigg\|\\
&\leq\max_{1\leq i\leq n}\|\tilde{\mY}^i\om\mY^i\|\cdot\sup_{x_j\in D_j}n^{-1}\sum_{i=1}^nK_{h_j}(x_j,\tilde{\xi}_j^i)\I(\mxi^i\in D)\I(\tilde{\mxi}^i\in D)\\
&\quad+\const\cdot\sup_{x_j\in D_j}n^{-1}\sum_{i=1}^n\frac{1}{h_j^{L_j}}\bigg|K_j\bigg(\frac{\|x_j-\tilde{\xi}_j^i\|_j}{h_j}\bigg)-K_j\bigg(\frac{\|x_j-\xi_j^i\|_j}{h_j}\bigg)\bigg|\|\mY^i\|\I(\mxi^i\in D)\\
&\quad+\const\cdot\sup_{x_j\in D_j}n^{-1}\sum_{i=1}^n\int_{D_j}\frac{1}{h_j^{L_j}}\bigg|K_j\bigg(\frac{\|t_j-\tilde{\xi}_j^i\|_j}{h_j}\bigg)-K_j\bigg(\frac{\|t_j-\xi_j^i\|_j}{h_j}\bigg)\bigg|dt_j\\
&\qquad\qquad\qquad\qquad\qquad\qquad\quad\cdot \frac{1}{h_j^{L_j}}K_j\bigg(\frac{\|x_j-\xi_j^i\|_j}{h_j}\bigg)\|\mY^i\|\I(\mxi^i\in D)\\
&=O_p(b_n)+O_p(h_j^{-L_j-2}\cdot a^2_{nj}+h_j^{-1}\cdot a_{nj})
\end{split}
\end{align}
provided that $n^{-1+2\beta+2/\alpha}\cdot h_j^{-L_j}=o(1)$ and $n^{-1/2}\cdot(\log{n})^{1/2}\cdot h_j^{-L_j/2}=O(1)$, where the equality in (\ref{consistency final}) follows by arguing as in the proof of Lemma \ref{marginal density approximation 1} using Lemma \ref{uniform rates} instead of Lemma \ref{uniform rates density}. This completes the proof for the second part of the lemma.

The first part of the lemma follows similarly as in the above proof by obtaining the rates with $V_{nj}^{1/2}$ being replaced by 1 using the bound $\I(\mxi^i\in D\cap D(n))\leq\I(\mxi^i\in D)$ or $\I(\mxi^i\in D(n))\leq\I(\mxi^i\in\prod_{j=1}^dD_j^+(C\cdot a_{nj}))$ whenever $V_{nj}^{1/2}$ appears. In this case, the condition that $n^{-\beta}\cdot(\log{n})^{1/2}=O(1)$ is required instead of the condition that $n^{-\beta}\cdot(\log{n})^{1/2}\cdot(\sum_{j=1}^da_{nj})^{1/\alpha}\cdot(h_j^{-L_j}\cdot a_{nj}+\sum_{k\neq j}a_{nk})^{-1/2}=O(1)$.
\end{proof}

We now present a high-level sufficient condition for Propositions \ref{existence} and \ref{convergence1}.

\begin{customcon}{(F)}\label{F}\leavevmode
For all $1 \le j \neq k \le d$ and $x_j\in D_j$, it holds that $\hat{p}^D_j(x_j)>0$,
\[
\int_{D_k}\frac{(\hat{p}^D_{jk}(x_j,x_k))^2}{\hat{p}^D_k(x_k)}dx_k<\infty\text{~and~}\int_{D_k}\int_{D_j}\frac{(\hat{p}^D_{jk}(x_j,x_k))^2}{\hat{p}^D_j(x_j)\hat{p}^D_k(x_k)}dx_jdx_k<\infty.
\]
\end{customcon}

The next two lemmas are variations of Propositions 2.1 and 2.2 in Jeon et al. (2021a) and Theorems 1 and 2 in Jeon et al. (2021b) for flexible estimation domains.

\begin{customlemma}{S.16}\label{general existence}
Assume that condition \ref{F} holds. Then, the existence of the SBF estimator and the uniqueness of the sum of the individual SBF estimators in Proposition \ref{existence} hold. If we further assume that $\hat{p}^D>0$ on $D$, then the uniqueness of the individual SBF estimators in Proposition \ref{existence} holds.
\end{customlemma}

\begin{proof}
We define
\begin{align*}
L_2^\mbH(\hat{p}^D)=\bigg\{\mg:D\rightarrow\mbH\bigg|\|\mg\|_{2,D,n}:=\bigg(\int_{D}\|\mg(\mx)\|^2d\hat{P}_{\mxi}^D(\mx)\bigg)^{1/2}<\infty\bigg\}.
\end{align*}
We also define its subspaces
\begin{align*}
L_2^\mbH(\hat{p}_j^D)=\{\mg_j\in L_2^\mbH(\hat{p}^D):\exists~\text{a univariate map}~\mg^*_j:D_j\rightarrow\mbH~\text{such that}~\mg_j(\mx)=\mg^*_j(x_j)\}
\end{align*}
for $1\leq j\leq d$ and their sum space $S^\mbH(\hat{p}^D)=\{\bigoplus_{j=1}^d\mg_j:\mg_j\in L_2^\mbH(\hat{p}_j^D)\}$. By arguing as in the proof of Theorem 3.3 in Jeon and Park (2020), we may show that $S^\mbH(\hat{p}^D)$ is a closed subspace of $L_2^\mbH(\hat{p}^D)$ under condition \ref{F}.

Define functional $\hat{F}:S^\mbH(\hat{p}^D)\rightarrow\mathbb{R}$ by
\ba
\hat{F}(\mg)=\int_Dn^{-1}\sum_{i=1}^n\|\tilde{\mY}^i\om\mg(\mx)\|^2\I(\tilde{\mxi}^i\in D)\prod_{j=1}^dK_{h_j}(x_j,\tilde{\xi}_j^i)d\mx.
\ea
One can check that $\hat{F}$ is a well-defined, strictly convex, continuous and G\^{a}teaux differentiable functional satisfying $\hat{F}(\mg)\rightarrow\infty$
as $\|\mg\|_{2,D,n}\rightarrow\infty$.
Hence, Lemma~4 in Beltrami (1967) entails
that there exists a minimizer of $\hat F$ in $S^\mbH(\hp^D)$. By Theorem 5.3.19 in Atkinson and Han (2009), $\hat{\mg}$ being a minimizer of
$\hat F$ is equivalent to $\mathfrak{D}^1\hat F(\hat\mg)(\hat\mk)=0$ for all $\hat\mk \in S^\mbH(\hp^D)$,
where
\[
\mathfrak{D}^1\hat F(\hat\mg)(\hat\mk)=-2\int_{D}n^{-1}\sumin \lng\tilde{\mY}^i \ominus \hat\mg(\mx),\hat\mk(\mx)\rng\I(\tilde{\mxi}^i\in D)\prod_{j=1}^dK_{h_j}(x_j,\tilde{\xi}_j^i)d\mx
\]
is the G\^{a}teaux derivative of $\hat F$ at $\hat\mg$ in the direction of $\hat\mk$. With specification of $\hat\mk\in S^\mbH(\hp^D)$
to $\hat\mk_j \in L_2^\mbH(\hp_j^D)$ for each $1 \le j \le d$, we find that
\begin{align}\label{equivalence2 in existence}
\int_{D_{-j}}n^{-1}\odot \bigg(\oplusin(\tilde{\mY}^i \ominus \hat{\mg}(\mx))
\odot \bigg(\I(\tilde{\mxi}^i\in D)\prod_{j=1}^dK_{h_j}(x_j,\tilde{\xi}_j^i)\bigg)\bigg)d\mx_{-j}=\bzero
\end{align}
almost everywhere $x_j \in D_j$ with respect to $\Leb_j$ for all $1 \le j \le d$.
Let $\hat{\mg}=\hat{\mg}_0 \oplus \bigoplus_{j=1}^d \hat{\mg}_j$ be
a decomposition of $\hat\mg$ such that $(\hat{\mg}_j:1\leq j\leq d)$ satisfies (\ref{estimated constraint}).
Plugging the decomposition into the left hand side of
(\ref{equivalence2 in existence}) and using (\ref{normalization}),
we may show that $\hat{\mg}_0=\hat{\mf}_0$ and $(\hat\mg_j: 1 \le j \le d)$ satisfies
\begin{align}\label{SBF almost everywhere equation}
\hat{\mg}_j(x_j)=\hat\mm_j(x_j) \ominus \hat{\mf}_0 \ominus \opluskj\int_{D_k}
\hat{\mg}_k(x_k) \odot \frac{\hat{p}^D_{jk}(x_j,x_k)}{\hat{p}^D_j(x_j)}dx_k
\end{align}
almost everywhere $x_j \in D_j$ with respect to $\Leb_j$ for all $1 \le j \le d$. By taking the right hand side of (\ref{SBF almost everywhere equation})
as $\hat\mf_j(x_j)$, we see that $(\hat\mf_j:1\leq j\leq d)$ satisfies both (\ref{estimated Bochner equation}) and (\ref{estimated constraint}).


The second and third conclusions of the lemma follow by arguing as in the proof of Proposition 2.1 in Jeon et al. (2021a).
\end{proof}

\begin{customlemma}{S.17}\label{general algorithm convergence 1}
Assume that condition \ref{F} holds. Then, the conclusion of Proposition \ref{convergence1} holds.
\end{customlemma}

\begin{proof}
The first conclusion of the lemma follows by arguing as in the proof of Proposition 2.2 in Jeon et al. (2021a). The first conclusion implies that
\ba
\sum_{r=1}^\infty\int_D\|\hat\mf^{[r]}(\mx)\om\hat\mf(\mx)\|^2d\hat{P}^D_{\mxi}(\mx)=\int_D\sum_{r=1}^\infty\|\hat\mf^{[r]}(\mx)\om\hat\mf(\mx)\|^2d\hat{P}^D_{\mxi}(\mx)<\infty.
\ea
This implies that $\sum_{r=1}^\infty\|\hat\mf^{[r]}(\mx)\om\hat\mf(\mx)\|^2<\infty$ almost everywhere with respect to $\hat{P}^D_{\mxi}$. The latter implies the second conclusion of the lemma. The third conclusion of the lemma follows from the second conclusion by Egorov's theorem.
\end{proof}

The next lemma provides a new result that does not exist in the literature.

\begin{customlemma}{S.18}\label{general convergence1 individual}
Assume that $\hat{p}^D$ is bounded away from zero and infinity on $D$. Then, the conclusion of Proposition \ref{convergence1 individual} holds.
\end{customlemma}

\begin{proof}
Since $S^\mbH(\hat{p}^D)$ is a closed subspace of $L_2^\mbH(\hat{p}^D)$ as demonstrated in the proof of Lemma \ref{general existence}, Lemma S.7 in Jeon and Park (2020) implies that there exist a map $\bigoplus_{j=1}^d\hat{\mg}_j^{[r]}\in S^\mbH(\hat{p}^D)$ and a constant $\hat{c}>0$ that depends only on $\hat{p}^D$ such that $\bigoplus_{j=1}^d\hat{\mg}_j^{[r]}(x_j)=\hat\mf^{[r]}(\mx)\om\hat\mf(\mx)$ almost everywhere with respect to $\hat{P}^D_{\mxi}$ and
\begin{align}\label{individual inequality 1}
\max_{1\leq j\leq d}\int_{D_j}\|\hat{\mg}_j^{[r]}(x_j)\|^2\hat{p}^D_j(x_j)dx_j\leq \hat{c}\cdot\int_D\|\hat\mf^{[r]}(\mx)\om\hat\mf(\mx)\|^2d\hat{P}^D_{\mxi}(\mx).
\end{align}
Define $\hat{\mc}_j^{[r]}=\int_{D_j}\hat{\mg}_j^{[r]}(x_j)\od\hat{p}^D_j(x_j)dx_j$. Then,
\begin{align}\label{individual inequality 2}
\begin{split}
&\int_{D_j}\|\hat{\mg}_j^{[r]}(x_j)\|^2\hat{p}^D_j(x_j)dx_j\\
&=\int_{D_j}\|\hat{\mg}_j^{[r]}(x_j)\om\hat{\mc}^{[r]}_j\|^2\hat{p}^D_j(x_j)dx_j+2\int_{D_j}\lng\hat{\mc}^{[r]}_j,\hat{\mg}_j^{[r]}(x_j)\om\hat{\mc}^{[r]}_j\rng\hat{p}^D_j(x_j)dx_j+\|\hat{\mc}^{[r]}_j\|^2\\
&=\int_{D_j}\|\hat{\mg}_j^{[r]}(x_j)\om\hat{\mc}^{[r]}_j\|^2\hat{p}^D_j(x_j)dx_j+\|\hat{\mc}^{[r]}_j\|^2\\
&\geq\int_{D_j}\|\hat{\mg}_j^{[r]}(x_j)\om\hat{\mc}^{[r]}_j\|^2\hat{p}^D_j(x_j)dx_j.
\end{split}
\end{align}
Note that
\ba
\bigoplus_{j=1}^d\hat{\mc}_j^{[r]}=\int_D(\hat\mf^{[r]}(\mx)\om\hat\mf(\mx))d\hat{P}^D_{\mxi}(\mx)=\bigoplus_{j=1}^d\int_{D_j}(\hat\mf^{[r]}_j(x_j)\om\hat\mf_j(x_j))\od\hat{p}^D_j(x_j)dx_j=\bzero.
\ea
By arguing as in the proof of Lemma S.8 in Jeon and Park (2020) with the condition that $\hat{p}^D$ is bounded away from zero and infinity on $D$, we get that $\hat{\mg}_j^{[r]}(x_j)\om\hat{\mc}^{[r]}_j=\hat{\mf}_j^{[r]}(x_j)\om\hat{\mf}_j(x_j)\op\hat{\mc}_j$ almost everywhere with respect to $\Leb_j$, where $\hat{\mc}_j\in\mbH$ are constants such that $\bigoplus_{j=1}^d\hat{\mc}_j=\bzero$. Since
\begin{align*}
\hat{\mc}_j=\int_{D_j}(\hat{\mg}_j^{[r]}(\mx)\om\hat{\mc}_j^{[r]})\od\hat{p}^D_j(x_j)dx_j\om\int_{D_j}(\hat{\mf}_j^{[r]}(x_j)\om\hat{\mf}_j(x_j))\od\hat{p}^D_j(x_j)dx_j=\bzero,
\end{align*}
we get
\begin{align}\label{individual inequality 3}
\begin{split}
\int_{D_j}\|\hat{\mg}_j^{[r]}(x_j)\om\hat{\mc}^{[r]}_j\|^2\hat{p}^D_j(x_j)dx_j&=\int_{D_j}\|\hat{\mf}_j^{[r]}(x_j)\om\hat{\mf}_j(x_j)\|^2\hat{p}^D_j(x_j)dx_j\\
&\geq\int_{D_j}\|\hat{\mf}_j^{[r]}(x_j)\om\hat{\mf}_j(x_j)\|^2dx_j.
\end{split}
\end{align}
Now, (\ref{individual inequality 1}), (\ref{individual inequality 2}), (\ref{individual inequality 3}), Lemma \ref{general algorithm convergence 1} and the proof of Lemma \ref{general algorithm convergence 1} give the desired result.
\end{proof}


Below, we provide a high-level sufficient condition for Theorem \ref{convergence2}.

\begin{customcon}{(G)}\label{G}\leavevmode
For all $1\leq j\neq k\leq d$, the following holds.
\begin{itemize}
\item[{\rm(G1)}] $\lim_{n\rightarrow\infty}\Prob\big(\sup_{x_j\in D_j}\hat{p}^D_j(x_j)^{-1}<C\big)=1$ for some $C>0$.
\item[{\rm(G2)}] $\lim_{n\rightarrow\infty}\Prob\big(\sup_{(x_j,x_k)\in D_j\times D_k}\hat{p}^D_{jk}(x_j,x_k)<C\big)=1$ for some $C>0$.
\item[{\rm(G3)}] $\lim_{n\rightarrow\infty}\Prob\big(\sup_{x_j\in D_j}\|\hat{\mm}_j(x_j)\|<C\big)=1$ for some $C>0$.
\item[{\rm(G4)}] $\int_{D_j}(\hat{p}_j^D(x_j)-p_j^D(x_j))^2dx_j=o_p(1)$.
\item[{\rm(G5)}] $\int_{D_k}\int_{D_j}(\hat{p}_{jk}^D(x_j,x_k)-p_{jk}^D(x_j,x_k))^2dx_jdx_k=o_p(1)$.
\end{itemize}
\end{customcon}

The following lemma is a variation of Proposition 2.3 in Jeon et al. (2021a) and Theorem 3 in Jeon et al. (2021b) for flexible estimation domains. We note that the two results in Jeon et al. (2021a) and Jeon et al. (2021b) do not contain the uniqueness of their individual SBF estimators.

\begin{customlemma}{S.19}\label{general algorithm convergence 2}
Assume that conditions \ref{G} and (B7) hold and that $p$ is bounded away from zero and infinity on $D$. Then, the conclusions of Theorem \ref{convergence2} hold.
\end{customlemma}

\begin{proof}
Define a probability measure $P_{\mxi}^D$ on the product sigma-field $\bigotimes_{j=1}^{d}\mathcal{B}(D_j)$ by $P_{\mxi}^D(A)=\int_Ap^D(\mx)d\mx$ and define
\begin{align*}
L_2^\mbH(p^D)=\bigg\{\mg:D\rightarrow\mbH\bigg|\|\mg\|_{2,D}:=\bigg(\int_{D}\|\mg(\mx)\|^2dP_{\mxi}^D(\mx)\bigg)^{1/2}<\infty\bigg\}.
\end{align*}
We also define its subspaces
\begin{align*}
L_2^\mbH(p_j^D)=\{\mg_j\in L_2^\mbH(p^D):\exists~\text{a univariate map}~\mg^*_j:D_j\rightarrow\mbH~\text{such that}~\mg_j(\mx)=\mg^*_j(x_j)\}
\end{align*}
for $1\leq j\leq d$ and their sum space $S^\mbH(p^D)=\{\bigoplus_{j=1}^d\mg_j:\mg_j\in L_2^\mbH(p_j^D)\}$. By arguing as in the proof of Theorem 3.3 in Jeon and Park (2020), we may show that $S^\mbH(p^D)$ is a closed subspace of $L_2^\mbH(p^D)$ under the population version of condition \ref{F}. For all $1\leq j\leq d$, we define $\hat{\pi}_j:S^\mbH(p^D)\rightarrow L_2^\mbH(p^D_j)$ by
\ba
\hat{\pi}_j(\mg)(\mx)=\int_{D_{-j}}\mg(\mx)\od\frac{\hp^D(\mx)}{\hp_j^D(x_j)}d\mx_{-j}.
\ea
We note that the operators are well defined with probability tending to one.
We also define $\hat{T}:S^\mbH(p^D)\rightarrow S^\mbH(p^D)$ by $\hat{T}=(I-\hat{\pi}_d)\circ\cdots\circ(I-\hat{\pi}_1)$, where $I$ is the identity operator on $S^\mbH(p^D)$. The operator $\hat{T}$ is also well defined with probability tending to one. The following arguments hold with probability tending to one.

By arguing as in the proof of Theorem 3.5 in Jeon and Park (2020), we can show that $\|\hat{T}\|_{\rm op}\leq\tau$ for some $\tau<1$, where $\|\cdot\|_{\rm op}$ is the operator norm on $S^\mbH(p^D)$. Hence, the contraction mapping theorem (see e.g., Chapter 4.4 in Sacks (2017)) implies that, up to measure 0 with respect to $P_{\mxi}^D$, there exists a unique solution of the functional equation $\mg=\hat{T}(\mg)\oplus\hat\mm$ over $\mg\in S^\mbH(p^D)$, where
\ba
\hat\mm=\hat\mm_d \ominus \hat{\mf}_0 \oplus (I-\hat{\pi}_d)(\hat\mm_{d-1}) \oplus \cdots \oplus (I-\hat{\pi}_d)\circ \cdots \circ (I-\hat{\pi}_2)
(\hat\mm_1).
\ea
We note that $\hat\mm$ belongs to $S^\mbH(p^D)$. We may show that a solution $(\hat{\mf}_j:1\leq j\leq d)$ of (\ref{estimated Bochner equation}) subject to (\ref{estimated constraint}) satisfies that $\hat{\mf}_{\op}\in S^\mbH(p^D)$ and $\hat\mf_\op=\hat{T}(\hat\mf_\op)\oplus\hat\mm$, where $\hat\mf_\op=\bigoplus_{j=1}^d\hat{\mf}_j$. Hence, if $(\hat{\mf}^\star_j:1\leq j\leq d)$ is another solution, then $\bigoplus_{j=1}^d\hat{\mf}_j=\bigoplus_{j=1}^d\hat{\mf}^\star_j$ almost everywhere with respect to $P_{\mxi}^D$. By arguing as in the proof of Lemma S.8 in Jeon and Park (2020) with the condition that $p$ is bounded away from zero and infinity on $D$, we get that $\hat{\mf}_j(x_j)=\hat{\mf}^\star_j(x_j)\op\hat{\mc}_j$ almost everywhere with respect to $\Leb_j$, where $\hat{\mc}_j\in\mbH$ are stochastic constants such that $\bigoplus_{j=1}^d\hat{\mc}_j=\bzero$. Since
\begin{align*}
\hat{\mc}_j=\int_{D_j}(\hat{\mf}_j(x_j)\om\hat{\mf}^\star_j(x_j))\od\hat{p}^D_j(x_j)dx_j=\bzero
\end{align*}
by (\ref{estimated constraint}), we get the desired result.
The remaining parts of the lemma follow by arguing as in the proof of Proposition 2.3 in Jeon et al. (2021a) with Egorov's theorem.
\end{proof}

We now state a lemma for Theorems \ref{holder rate} and \ref{differentiable rate}. It is a variation of Theorem 5.1 in Jeon et al. (2021a) for flexible estimation domains. For the statement, we define
\begin{align*}
\hat{\mdelta}^{A}_j(x_j)&=n^{-1}\odot\bigoplus_{i=1}^n\big((K_{h_j}(x_j,\tilde{\xi}^i_j)\I(\tilde{\mxi^i}\in D))\odot\mepsilon_i\big),\\
\hat{\mdelta}^{B}_j(x_j)&=n^{-1}\odot\bigoplus_{i=1}^n\big((K_{h_j}(x_j,\tilde{\xi}^i_j)\I(\tilde{\mxi^i}\in D))\odot(\mf_j(\xi^i_j)\om\mf_j(x_j))\big),\quad\mbox{and}\\
\hat{\mdelta}^{C}_j(x_j)&=\bigoplus_{k\neq j}\bigg(n^{-1}\odot\bigoplus_{i=1}^n
\int_{D_k}(\mf_k(\xi^i_k)\om\mf_k(x_k))\odot(K_{h_k}(x_k,\tilde{\xi}^i_k)K_{h_j}(x_j,\tilde{\xi}^i_j)\I(\tilde{\mxi^i}\in D))dx_k\bigg).
\end{align*}
Note that we have $\mf_j(\xi^i_j)$ instead of $\mf_j(\tilde{\xi}^i_j)$ in the definition of $\hat{\mdelta}^{B}_j(x_j)$. We also have $\mf_k(\xi^i_k)$ instead of $\mf_k(\tilde{\xi}^i_k)$ in the definition of $\hat{\mdelta}^{C}_j(x_j)$.
For $1 \le j \le d$, we let $S_j\in\mathcal{B}(D_j)$, and $A^{(1)}_{nj}(x_j)$, $A^{(2)}_{nj}$,
$A^{(3)}_{nj}(S_j)$, $B^{(1)}_{nj}(x_j)$, $B^{(2)}_{nj}(S_j)$, $B^{(3)}_{nj}(S_j)$,
$B^{(4)}_{nj}$, $C^{(1)}_{nj}(x_j)$, $C^{(2)}_{nj}$, $C^{(3)}_{nj}(S_j)$, $D^{(1)}_n$ and $D^{(2)}_{nj}$
be nonnegative sequences converging to zero such that
\begin{itemize}
\item[(i)] rates for $\hat{\mdelta}^{A}_j$:
\begin{align*}
\begin{split}
&\hat{\mdelta}^{A}_j(x_j)=O_p(A^{(1)}_{nj}(x_j)),\quad \bigg(\int_{D_j}\|\hat{\mdelta}^{A}_j(x_j)\|^2dx_j\bigg)^{1/2}=O_p(A^{(2)}_{nj}),\quad\mbox{and}\\
&\sup_{x_j\in S_j}\|\hat{\mdelta}^{A}_j(x_j)\|=O_p(A^{(3)}_{nj}( S_j));
\end{split}
\end{align*}
\item[(ii)] rates for $\hat{\mdelta}^{B}_j$:
\begin{align*}
\begin{split}
&\hat{\mdelta}^{B}_j(x_j)=O_p(B^{(1)}_{nj}(x_j)),\quad \bigg(\int_{ S_j}\|\hat{\mdelta}^{B}_j(x_j)\|^2dx_j\bigg)^{1/2}=O_p(B^{(2)}_{nj}( S_j)),\\
&\sup_{x_j\in S_j}\|\hat{\mdelta}^{B}_j(x_j)\|=O_p(B^{(3)}_{nj}( S_j)),\quad\mbox{and}\\
&\sum_{k\neq j}\sup_{x_j\in D_j}\bigg\|\int_{D_k}\hat{\mdelta}^{B}_k(x_k)\odot
\frac{\hat{p}^D_{jk}(x_j,x_k)}{\hat{p}^D_k(x_k)}dx_k\bigg\|=O_p(B^{(4)}_{nj});
\end{split}
\end{align*}
\item[(iii)] rates for $\hat{\mdelta}^{C}_j$:
\begin{align*}
\begin{split}
&\hat{\mdelta}^{C}_j(x_j)=O_p(C^{(1)}_{nj}(x_j)),\quad \bigg(\int_{D_j}\|\hat{\mdelta}^{C}_j(x_j)\|^2dx_j\bigg)^{1/2}=O_p(C^{(2)}_{nj}),\quad\mbox{and}\\
&\sup_{x_j\in S_j}\|\hat{\mdelta}^{C}_j(x_j)\|=O_p(C^{(3)}_{nj}( S_j));
\end{split}
\end{align*}
\item[(iv)] other rates:
\ba
&\hat{\mf}_0\ominus\mf_0=O_p(D^{(1)}_n),\quad n^{-1}\od\bigoplus_{i=1}^n\big(\mf_j(\xi^i_j)\od\I(\tilde{\mxi}^i\in D)\big)=O_p(D^{(2)}_{nj}).
\ea
\end{itemize}
Put $R_n=\sum_{j=1}^d(A^{(2)}_{nj}+B^{(4)}_{nj}+C^{(2)}_{nj})+D^{(1)}_n+b_n$, $E^{(1)}_{nj}(x_j)=A^{(1)}_{nj}(x_j)+B^{(1)}_{nj}(x_j)+C^{(1)}_{nj}(x_j)$ and $E^{(3)}_{nj}( S_j)=A^{(3)}_{nj}( S_j)+B^{(3)}_{nj}( S_j)+C^{(3)}_{nj}( S_j)$.

\begin{customlemma}{S.20}\label{general rate}
Assume that conditions (G1), (G2) and (G5) hold, that $p_j^D$ is bounded away from zero and infinity on $D_j$ for all $1\leq j\leq d$ and that $p_{jk}^D$ is bounded on $D_j\times D_k$ for all $1\leq j\neq k\leq d$. Then, it holds that
\begin{align*}
\bigg(\int_{\prod_{j=1}^d S_j}\|\hat{\mf}(\mx)\om\mf(\mx)\|^2p(\mx)d\mx\bigg)^{1/2}
=O_p\bigg(R_n+\sum_{j=1}^dB^{(2)}_{nj}( S_j)\bigg).
\end{align*}
If we further assume that $p$ is bounded away from zero and infinity on $D$, then it holds that, for all $1\leq j\leq d$,
\begin{align*}
&\hat{\mf}_j(x_j)\om\mf_j(x_j)=O_p\bigg(R_n+E^{(1)}_{nj}(x_j)+\sum_{k\neq j}D^{(2)}_{nk}\bigg)\text{~for all~}x_j\in D_j,\\
&\bigg(\int_{ S_j}\|\hat{\mf}_j(x_j)\om\mf_j(x_j)\|^2dx_j\bigg)^{1/2}
=O_p(R_n+B^{(2)}_{nj}( S_j)+D^{(2)}_{nj}),\\
&\sup_{x_j\in S_j}\|\hat{\mf}_j(x_j)\om\mf_j(x_j)\|=O_p\bigg(R_n+E^{(3)}_{nj}( S_j)+\sum_{k\neq j}D^{(2)}_{nk}\bigg).
\end{align*}
\end{customlemma}

\begin{proof}
We define
\ba
\hat{\mf}_j^A(x_j)&=(\hp^D_j(x_j)\cdot\hat{p}^D_0)^{-1}\od\hat{\mdelta}_j^A(x_j),\quad\hat{\mf}_j^B(x_j)=(\hp^D_j(x_j)\cdot\hat{p}^D_0)^{-1}\od\hat{\mdelta}_j^B(x_j),\\
\hat{\mf}_j^C(x_j)&=(\hp^D_j(x_j)\cdot\hat{p}^D_0)^{-1}\od\hat{\mdelta}_j^C(x_j),\quad\hat{\mDelta}_j(x_j)=\hat{\mf}_j(x_j)\om\mf_j(x_j)\om\hat{\mf}_j^B(x_j).
\ea
We also define $\check{\mm}_j(x_j)$ as $\hat{\mm}_j(x_j)$ with $\tilde{\mY}^i$ in the definition of $\hat{\mm}_j(x_j)$ being replaced by $\mY^i$.
Then, equation (\ref{estimated Bochner equation}) can be written as
\begin{align}\label{equation 3 new}
\begin{split}
\hat{\mDelta}_j(x_j)=&\hat{\mf}_j^A(x_j)\op\hat{\mf}_j^C(x_j)\op \mf_0\om\hat{\mf}_0\op\hat{\mm}_j(x_j)\om\check{\mm}_j(x_j)\om\bigoplus_{k\neq j}
\int_{D_k}\hat{\mf}_k^B(x_k)\od\frac{\hp^D_{jk}(x_j,x_k)}{\hp^D_j(x_j)}dx_k\\
&\om\bigoplus_{k\neq j}\int_{D_k}\hat{\mDelta}_k(x_k)\od\frac{\hp^D_{jk}(x_j,x_k)}{\hp^D_j(x_j)}dx_k,\quad1\leq j\leq d.
\end{split}
\end{align}
Note that
\begin{align*}
\bigg(\int_{D_j}\|\hat{\mf}_j^W(x_j)\|^2p^D_j(x_j)dx_j\bigg)^{1/2}&\leq(\hat{p}^D_0)^{-1}\bigg(\sup_{x_j\in D_j}\frac{p^D_j(x_j)}{(\hp^D_j(x_j))^2}\bigg)^{1/2}\bigg(\int_{D_j}\|\hat{\mf}_j^W(x_j)\|^2dx_j\bigg)^{1/2}\\
&=O_p(W_{nj}^{(2)})
\end{align*}
for $W=A,B$ and $C$. Also, note that
\begin{align*}
&\bigg(\int_{D_j}\bigg\|\bigoplus_{k\neq j}\int_{D_k}\hat{\mf}_k^B(x_k)\od\frac{\hp^D_{jk}(x_j,x_k)}{\hp^D_j(x_j)}dx_k
\bigg\|^2p^D_j(x_j)dx_j\bigg)^{1/2}\\
&\le (\hat{p}^D_0)^{-1}\bigg(\int_{D_j}\frac{p^D_j(x_j)}{(\hp^D_j(x_j))^2}dx_j\bigg)^{1/2}\sum_{k\neq j}\sup_{x_j\in D_j}\bigg\|\int_{D_k}\hat{\mdelta}_k^B(x_k)\od\frac{\hp^D_{jk}(x_j,x_k)}{\hp^D_k(x_k)}dx_k\bigg\|\\
&=O_p(B_{nj}^{(4)})
\end{align*}
and
\ba
\sup_{x_j\in D_j}\|\hat{\mm}_j(x_j)\om\check{\mm}_j(x_j)\|\leq\max_{1\leq i\leq n}\|\tilde{\mY}^i\om\mY^i\|=O_p(b_n).
\ea
Then, the map $\hat{\mLambda}_j$ defined by
\begin{align*}
\hat{\mLambda}_j(x_j)=\hat{\mf}_j^A(x_j)\op\hat{\mf}_j^C(x_j)\op \mf_0\om\hat{\mf}_0\op\hat{\mm}_j(x_j)\om\check{\mm}_j(x_j)\om\bigoplus_{k\neq j}\int_{D_k}
\hat{\mf}_k^B(x_k)\od\frac{\hp^D_{jk}(x_j,x_k)}{\hp^D_j(x_j)}dx_k
\end{align*}
satisfies
\begin{align}\label{Lambda order}
\bigg(\int_{D_j}\|\hat{\mLambda}_j(x_j)\|^2p^D_j(x_j)dx_j\bigg)^{1/2}=O_p(A_{nj}^{(2)}+B_{nj}^{(4)}+C_{nj}^{(2)}+D_n^{(1)}+b_n).
\end{align}
For a map $\mg:D\rightarrow\mbH$, we define $\|\mg\|_{2, S}=\left(\int_{ S}\|\mg(\mx)\|^2p(\mx)d\mx\right)^{1/2}$ for $ S=\prod_{j=1}^d S_j$. Using (\ref{Lambda order}) and by arguing as in the proof of Theorem 5.1 in Jeon et al. (2021a), we may show that $\|\bigoplus_{j=1}^d\hat{\mDelta}_j\|_{2,D}=O_p(R_n)$.
Since $\|\hat{\mf}\om\mf\|_{2, S}\leq\|\bigoplus_{j=1}^d(\hat{\mf}_j\om\mf_j)\|_{2, S}+\|\hat{\mf}_0\om\mf_0\|_{2, S}\leq\|\bigoplus_{j=1}^d\hat{\mf}^B_j\|_{2, S}+\|\bigoplus_{j=1}^d\hat{\mDelta}_j\|_{2,D}+O_p(D_n^{(1)})$, the first conclusion of the lemma follows.

Now, we further assume that $p$ is bounded away from zero and infinity on $D$. By Lemma S.7 in Jeon and Park (2020), there exist a constant $c>0$ and maps $(\hat{\mDelta}^*_1,\ldots,\hat{\mDelta}^*_d)$
such that $\bigoplus_{j=1}^d\hat{\mDelta}^*_j=\bigoplus_{j=1}^d\hat{\mDelta}_j$ almost everywhere with respect to $\Leb$ and
\begin{align*}
\max_{1\leq j\leq d}\bigg(\int_{D_j}\|\hat{\mDelta}^*_j(x_j)\|^2p^D_j(x_j)dx_j\bigg)^{1/2}\leq c\cdot\bigg\|\bigoplus_{j=1}^d\hat{\mDelta}_j\bigg\|_{2,D}.
\end{align*}
By arguing as in the proof of Lemma S.8 in Jeon and Park (2020), one may show that there exist stochastic constants $(\hat{\mc}_1,\ldots,\hat{\mc}_d)$
satisfying $\hat{\mDelta}_j=\hat{\mDelta}_j^*\op\hat{\mc}_j$ almost everywhere with respect to $\Leb_j$ for $1\leq j\leq d$ and $\bigoplus_{j=1}^d\hat{\mc}_j=\bzero$.
Note that
\begin{align*}
\hat{\mc}_j&=\int_{D_j}(\hat{\mDelta}_j(x_j)\om\hat{\mDelta}^*_j(x_j))\od\hp^D_j(x_j)dx_j\\
&=-1\od\int_{D_j}(\mf_j(x_j)\op\hat{\mf}_j^B(x_j))\od\hp^D_j(x_j)dx_j\om\int_{D_j}\hat{\mDelta}^*_j(x_j)\od\hp^D_j(x_j)dx_j\\
&=-(\hat{p}^D_0\cdot n)^{-1}\od\bigoplus_{i=1}^n(\mf_j(\xi^i_j)\od\I(\tilde{\mxi}^i\in D))\op O_p(R_n)\\
&=O_p(D_{nj}^{(2)}+R_n).
\end{align*}
Hence, $\left(\int_{D_j}\|\hat{\mDelta}_j(x_j)\|^2p^D_j(x_j)dx_j\right)^{1/2}=O_p(D_{nj}^{(2)}+R_n)$. Thus,
\begin{align*}
&\bigg(\int_{ S_j}\|\hat{\mf}_j(x_j)\om\mf_j(x_j)\|^2dx_j\bigg)^{1/2}\\
&\leq\Big(\inf_{x_j\in D_j}p^D_j(x_j)\Big)^{-1}\bigg(\int_{ S_j}\|\hat{\mf}_j(x_j)\om\mf_j(x_j)\|^2p^D_j(x_j)dx_j\bigg)^{1/2}\\
&\leq\Big(\inf_{x_j\in D_j}p^D_j(x_j)\Big)^{-1}\bigg(\bigg(\int_{D_j}
\|\hat{\mDelta}_j(x_j)\|^2p^D_j(x_j)dx_j\bigg)^{1/2}+\bigg(\int_{ S_j}\|\hat{\mf}_j^B(x_j)\|^2p^D_j(x_j)dx_j\bigg)^{1/2}\bigg)\\
&=O_p(D_{nj}^{(2)}+R_n+B_{nj}^{(2)}( S_j)).
\end{align*}
As for the pointwise and uniform rates of $\hat{\mf}_j\om\mf_j$, Note that
\begin{align*}
\|\hat{\mf}_j(x_j)\om\mf_j(x_j)\|\leq&\|\hat{\mf}_j^A(x_j)\|+\|\hat{\mf}_j^B(x_j)\|+\|\hat{\mf}_j^C(x_j)\|+\|\hat{\mf}_0\om\mf_0\|+\|\hat{\mm}_j(x_j)\om\check{\mm}_j(x_j)\|\\
&+\sum_{k\neq j}\bigg\|\int_{D_k}\hat{\mf}_k^B(x_k)\od\frac{\hp^D_{jk}(x_j,x_k)}{\hp^D_j(x_j)}dx_k\bigg\|+\sum_{k\neq j}\int_{D_k}\|\hat{\mDelta}_k(x_k)\|\frac{\hp^D_{jk}(x_j,x_k)}{\hp^D_j(x_j)}dx_k.
\end{align*}
This follows from (\ref{equation 3 new}). This together with the $L_2$ rates of $\hat{\mDelta}_k$ and the conditions in the lemma gives the desired pointwise and uniform rates.
\end{proof}

\begin{customlemma}{S.21}\label{A.I.}
Let $\{(\mdelta^i,\mxi^i):1\leq i\leq n\}$ be i.i.d. copies of $(\mdelta,\mxi)$, where $\mdelta$ is a $\mbH$-valued random element such that $\E(\|\mdelta\|^\alpha)<\infty$ for some $2<\alpha<\infty$. For a fixed $\mx\in\prod_{j=1}^d( D_j\setminus\partial D_j)$, assume that (i) $\E(\|\mdelta\|^\alpha|\xi_j=\cdot,\mxi\in D)$ is bounded on a neighborhood of $x_j$, and $\E(\lng\mdelta,\me_m\rng \lng\mdelta,\me_{m'}\rng |\xi_j=\cdot,\xi_k=\cdot,\mxi\in D)$ and $p^D_{jk}$ are bounded on a neighborhood of $(x_j,x_k)$ for all $1\leq j\neq k\leq d$; (ii) there exists a neighborhood $U_{x_j}$ of $x_j$ such that, for all $m$ and $m'$, $\E(\lng\mdelta,\me_m\rng \lng\mdelta,\me_{m'}\rng |\xi_j=\cdot,\mxi\in D)$ are continuous on $U_{x_j}$ for all $1\leq j\leq d$; (iii) $p^D_j$ is continuous on a neighborhood of $x_j$ and $p^D_j(x_j)>0$ for all $1\leq j\leq d$; (iv) $h_j=o(1)$ and $n^{-L_{\rm max}/(L_{\rm max}+4)}\cdot h_j^{-L_j}\rightarrow\kappa_j\in[0,\infty)$ for all $1\leq j\leq d$. Then,
\begin{align*}
\bigg(\big(n^{2/(L_{\rm max}+4)}\cdot (\hat{p}^D_j(x_j)\cdot \hat{p}_0^D\cdot n)^{-1}\big)\odot\bigg(\bigoplus_{i=1}^n\big((K_{h_j}(x_j,\xi^i_j)\I(\mxi^i\in D))\odot\mdelta^i\big)\bigg):1\leq j\leq d\bigg)\\
\overset{d}{\longrightarrow}\big(\mG(\bzero,C^*_{j,x_j}):1\leq j\leq d\big),
\end{align*}
where the covariance operator $C^*_{j,x_j}:\mbH\rightarrow\mbH$ is characterized by
\begin{align*}
\lng C^*_{j,x_j}(\mh),\me_m\rng =\sum_{m'\geq 1}\lng\mh,\me_{m'}\rng\cdot\mathcal{V}^*_{j,mm'}(x_j),\quad\mh\in\mbH,\;m\geq 1,
\end{align*}
where $\mathcal{V}^*_{j,mm'}(x_j)=\kappa_j\cdot (p_j^D(x_j)\cdot p_0^D)^{-1}\cdot \E(\lng\mdelta,\me_m\rng \lng\mdelta,\me_{m'}\rng |\xi_j=x_j,\mxi\in D)\cdot\int_{\mathbb{R}^{L_j}}K_j^2(\|\mathbf{t}\|)dt_j$.
In addition, $\mG(\bzero,C^*_{1,x_1}),\cdots,\mG(\bzero,C^*_{d,x_d})$ are independent.
\end{customlemma}

\begin{proof}
Since $\mx$ is a fixed interior point of $ D$, we may assume that $x_j\in  D^-_j(2h_j)$ for $1\leq j\leq d$. Let $\mbH^d$ denote the space of tuples $(\mh_j: 1 \le j \le d)$
with $\mh_j \in \mbH$.
Let $\|\cdot\|_{\mbH^d}$ and $\lng \cdot,\cdot \rng _{\mbH^d}$ denote the norm
and inner product on $\mbH^d$, respectively, defined in the standard way.
Let $\me_{jm} \in \mbH^d$ denote
$(\mathbf{0},\ldots,\mathbf{0}, \me_m,\mathbf{0},\ldots,\mathbf{0})$
where $\me_m$ is placed at the $j$th entry.
Then, $(\me_{jm}: 1 \le j \le d, \,m\ge1)$
forms an orthonormal basis of $\mbH^d$. Define
\begin{align*}
\meta_n^i(\mx)=\left(\frac{n^{2/(L_{\rm max}+4)}K_{h_j}(x_j,\xi^i_j)\I(\mxi^i\in D)}{np^D_j(x_j)p_0^D}\odot\mdelta^i:1\leq j\leq d\right)\in\mbH^d.
\end{align*}
Note that $\E(\lng \meta_n^i(\mx),\me_{jm} \rng _{\mbH^d})=0$ and
$\E(\|\meta_n^i(\mx)\|^2_{\mbH^d})<\infty$. Put $\mS_n(\mx)=\bigoplus_{i=1}^n
\meta_n^i(\mx)$ and let
\begin{align*}
\lim_{n\rightarrow\infty}\E\left(\lng\mS_n(\mx),\me_{jm}\rng _{\mbH^d} \cdot \lng\mS_n(\mx),\me_{km'}\rng _{\mbH^d}\right)=\mathcal{V}^*_{jmkm'}(\mx).
\end{align*}
Then, $\mathcal{V}^*_{jmkm'}(\mx)=\mathcal{V}^*_{j,mm'}(x_j)$ if $j=k$, and $\mathcal{V}^*_{jmkm'}(\mx)=0$ otherwise. We also get
\begin{align*}
&\limn\,\sumjd\sum_{m\geq l}\E\left(\lng\mS_n(\mx),\me_{jm}\rng _{\mbH^d}\cdot \lng\mS_n(\mx),\me_{jm}\rng _{\mbH^d}\right)\\
&=\sumjd\frac{\kappa_j}{p_j^D(x_j)\cdot p_0^D}\cdot \E(\|\mdelta\|^2 |\xi_j=x_j,\mxi\in D)\cdot\int_{B_j(\bzero_j,1)}K_j^2(\|t_j\|_j)dt_j.
\end{align*}
In addition, for $\tau=\alpha-2$,
\begin{align*}
\sumin \E\big(\big\|\meta_n^i(\mx)\big\|_{\mbH^d}^{2+\tau}\big)\leq&n^{-1-\tau+(4+2\tau)/(L_{\rm max}+4)}\cdot d^{1+\tau/2}\\
&\cdot\E\bigg(\I(\mxi\in D)\|\mdelta\|^{2+\tau}\sumjd\left(p_j^D(x_j)\cdot p_0^D\cdot h_j^{L_j}\right)^{-2-\tau}K^{2+\tau}_j\bigg(\frac{\|x_j-\xi_j\|_j}{h_j}\bigg)\bigg)\\
=&o(1).
\end{align*}
Therefore, by applying Theorem 1.1 in Kundu et al. (2000) for
infinite-dimensional $\mbH$ and Proposition S.2 in Jeon and Park (2020)
for finite-dimensional $\mbH$, we obtain $\mS_n(\mx)\overset{d}{\rightarrow}
\mG(\mathbf{0},C^*_\mx)$,
where $C^*_\mx:\mbH^d\rightarrow\mbH^d$ is a covariance operator such that, for all
$\mh=(\mh_1,\ldots,\mh_d)\in\mbH^d$,
\begin{align}\label{C and Cj}
\begin{split}
\lng C^*_\mx(\mh),\me_{jm} \rng _{\mbH^d}
&=\sumkd\sum_{m'\geq1} \lng \mh,\me_{km'} \rng _{\mbH^d}\cdot \mathcal{V}^*_{km'jm}(\mx)
=\sumkd\sum_{m'\geq1} \lng \mh_k,\me_{m'} \rng  \cdot \mathcal{V}^*_{km'jm}(\mx)\\
&=\sum_{m'\geq1} \lng \mh_j,\me_{m'} \rng  \cdot \mathcal{V}^*_{j,mm'}(x_j), \quad
1\leq j\leq d,\; m\geq 1.
\end{split}
\end{align}
Since $(\hp_j^D(x_j)\cdot\hp_0^D)^{-1}-(p_j^D(x_j)\cdot p_0^D)^{-1}=o_p(1)$ by Lemmas \ref{p0 rate} and \ref{marginal density convergence 1},
we get
\begin{align*}
\bigg(\big(n^{2/(L_{\rm max}+4)}\cdot (\hat{p}^D_j(x_j)\cdot \hat{p}_0^D\cdot n)^{-1}\big)\odot\bigg(\bigoplus_{i=1}^n\big((K_{h_j}(x_j,\xi^i_j)\I(\mxi^i\in D))\odot\mdelta^i\big)\bigg):1\leq j\leq d\bigg)\\
\overset{d}{\longrightarrow}\mG(\mathbf{0},C^*_\mx).
\end{align*}
The remaining proof follows by arguing as in the proof of Lemma S.5 in Jeon et al. (2021a).
\end{proof}

\begin{customlemma}{S.22}\label{markov}
Let $\{\mdelta^i:1\leq i\leq n\}$ be i.i.d. copies of a $\mbH$-valued random element $\mdelta$ such that $\E(\|\mdelta\|^\alpha)<\infty$ for some $0<\alpha<\infty$. Then, it holds that $\max_{1\leq i\leq n}\|\mdelta^i\|=O_p(n^{1/\alpha})$. If there exists a constant $c>0$ such that $\E\big(\exp(c\cdot\|\mdelta\|)\big)<\infty$, then it holds that $\max_{1\leq i\leq n}\|\mdelta^i\|=O_p(\log{n})$.
\end{customlemma}

\begin{proof}
We first prove the first result. We note that
\ba
\Prob(\max_{1\leq i\leq n}\|\mdelta^i\|>T\cdot n^{1/\alpha})&\leq n\cdot\Prob(\|\mdelta\|>T\cdot n^{1/\alpha})= n\cdot\Prob(\|\mdelta\|^{\alpha}>T^{\alpha}\cdot n)\\
&\leq n\cdot\E(\|\mdelta\|^{\alpha})/(T^{\alpha}\cdot n)=\E(\|\mdelta\|^{\alpha})/T^{\alpha},
\ea
where the third inequality follows from Markov's inequality. Hence,
\ba
\lim_{T\rightarrow\infty}\limsup_n\Prob(\max_{1\leq i\leq n}\|\mdelta^i\|>T\cdot n^{1/\alpha})\leq \lim_{T\rightarrow\infty}\E(\|\mdelta\|^{\alpha})/T^{\alpha}=0.
\ea
We now prove the second result. Markov's inequality again implies that
\ba
\Prob(\max_{1\leq i\leq n}\|\mdelta^i\|>T\cdot \log{n})&\leq n\cdot\Prob(\|\mdelta\|>T\cdot \log{n})\\
&= n\cdot\Prob(\exp(c\cdot\|\mdelta\|)>\exp(c\cdot T\cdot \log{n}))\\
&\leq \E\big(\exp(c\cdot\|\mdelta\|)\big)/(n^{c\cdot T-1}).
\ea
This gives the lemma.
\end{proof}

\subsection{Proof of Proposition \ref{HPC}}\label{HPC proof}

We first find the rate of $\|\hat{\mpsi}_r\om_*\mpsi_r\|_*$. We note that $C_{\mX}$ is a nonnegative-definite compact self-adjoint operator. Hence, if the eigenvalues of $C_{\mX}$ satisfy that $\lambda_1>\lambda_2>\cdots>\lambda_{r+1}$, then Lemma 2.3 in Horv\'{a}th and Kokoszka (2012) gives that
\begin{align}\label{Koko book}
\|\hat{\mpsi}_r\om_*\mpsi_r\|_*=\|\hat{\mpsi}_r\om_*{\rm{sgn}}(\lng\hat{\mpsi}_r,\mpsi_r\rng_*)\od_*\mpsi_r\|_*\leq\frac{2\sqrt{2}}{\delta_r}\cdot\|\hat{C}_{\mX}-C_{\mX}\|_{\rm op},
\end{align}
where $\delta_r=\lambda_1-\lambda_2$ if $r=1$ and $\delta_r=\min\{\lambda_{r-1}-\lambda_r,\lambda_r-\lambda_{r+1}\}$ otherwise, and $\|\cdot\|_{\rm op}$ denotes the operator norm. Hereafter, we let $\|\cdot\|_{\HS}$ denote the Hilbert-Schmidt norm of an operator. It is induced by an inner product and it holds that $\|C\|_{\rm op}\leq\|C\|_{\rm{HS}}$ for any Hilbert-Schmidt operator $C$. We note that $C_{\mX}$ and $\hat{C}_{\mX}$ are Hilbert-Schmidt operators so is $\hat{C}_{\mX}-C_{\mX}$. We also note that $\hat{C}_{\mX}$ is an unbiased estimator of $C_{\mX}$. These facts with some algebra entail that 
\begin{align}\label{cov bound}
\E(\|\hat{C}_{\mX}-C_{\mX}\|^2_{\rm op})\leq\E(\|\hat{C}_{\mX}-C_{\mX}\|^2_{\rm{HS}})=\E(\|\hat{C}_{\mX}\|^2_{\rm{HS}})-\|C\|^2_{\rm{HS}}\leq\E(\|\hat{C}_{\mX}\|^2_{\rm{HS}})=O(n^{-1})
\end{align}
provided that $\E(\|\mX\|_*^4)<\infty$. A version of (\ref{cov bound}) for mean zero $\mX$ can be found in Theorem 2.5 of Horv\'{a}th and Kokoszka (2012). Hence, we obtain $\|\hat{C}_{\mX}-C_{\mX}\|_{\rm op}\leq \|\hat{C}_{\mX}-C_{\mX}\|_{\rm HS}=O_p(n^{-1/2})$. Now, (\ref{Koko book}) and (\ref{cov bound}) imply that $\|\hat{\mpsi}_r\om_*\mpsi_r\|_*=O_p(n^{-1/2})$. Hence, we have
\ba
&\max_{1\leq i\leq n}|\lng\mX^i\om_*\bar{\mX},\hat{\mpsi}_r\rng_*-\lng\mX^i\om_*\E(\mX),\mpsi_r\rng_*|\\
&=\max_{1\leq i\leq n}|\lng \mX^i,\hat{\mpsi}_r\om_*\mpsi_r\rng_*-\lng\bar{\mX},\hat{\mpsi}_r\rng_*+\lng\E(\mX),\mpsi_r\rng_*|\\
&=\max_{1\leq i\leq n}|\lng \mX^i\om_*\E(\mX),\hat{\mpsi}_r\om_*\mpsi_r\rng_*+\lng\bar{\mX}\om_*\E(\mX),\hat{\mpsi}_r\rng_*|\\
&\leq\|\hat{\mpsi}_r\om_*\mpsi_r\|_*\cdot\max_{1\leq i\leq n}\|\mX^i\om_*\E(\mX)\|_*+\|\bar{\mX}\om_*\E(\mX)\|_*\cdot(\|\hat{\mpsi}_r\om_*\mpsi_r\|_*+\|\mpsi_r\|_*)\\
&=O_p(n^{-1/2})\cdot \max_{1\leq i\leq n}\|\mX^i\om_*\E(\mX)\|_*+O_p(n^{-1/2})\cdot O_p(1).
\ea
In case $\E(\|\mX\|_*^{\tau})<\infty$ for some $\tau\geq4$, we have $\max_{1\leq i\leq n}\|\mX^i\om_*\E(\mX)\|_*=O_p(n^{1/\tau})$ by Lemma \ref{markov}. In case $\mX$ has the exponential moment $\E\big(\exp(c\cdot\|\mX\|_*)\big)<\infty$, we have $\max_{1\leq i\leq n}\|\mX^i\om_*\E(\mX)\|_*=O_p(\log{n})$ by Lemma \ref{markov}. In case $\mX$ is a bounded random element, it holds that $\max_{1\leq i\leq n}\|\mX^i\om_*\E(\mX)\|_*=O_p(1)$. This proves the proposition. \qed

\subsection{Proof of Proposition \ref{HSC}}\label{HSC proof}

Since $C_{\mX\mY\mX}$ is a nonnegative-definite compact self-adjoint operator, if the eigenvalues of $C_{\mX\mY\mX}$ satisfy $\sigma_1>\sigma_2>\cdots>\sigma_{r+1}$, then Lemma 2.3 in Horv\'{a}th and Kokoszka (2012) implies that
\begin{align}\label{Koko book 2}
\|\hat{\mphi}_r\om_*\mphi_r\|_*\leq\frac{2\sqrt{2}}{\delta_r}\cdot\|\hat{C}_{\mX\mY\mX}-C_{\mX\mY\mX}\|_{\rm op},
\end{align}
where $\delta_r=\sigma_1-\sigma_2$ if $r=1$ and $\delta_r=\min\{\sigma_{r-1}-\sigma_r,\sigma_r-\sigma_{r+1}\}$ otherwise. Define the operators $\check{C}_{\mX\mY}:\mbH\rightarrow\mbH_*$ and $\check{C}_{\mY\mX}:\mbH_*\rightarrow\mbH$ by
\ba
\check{C}_{\mX\mY}(\mh)&=\frac{1}{n-1}\od_*\sideset{}{_*}\bigoplus_{i=1}^n(\lng\mY^i\om\bar{\mY},\mh\rng\od_*(\mX^i\om_*\bar{\mX})),\quad\mbox{and}\\
\check{C}_{\mY\mX}(\mh_*)&=\frac{1}{n-1}\od\bigoplus_{i=1}^n(\lng\mX^i\om_*\bar{\mX},\mh_*\rng_*\od(\mY^i\om\bar{\mY})),
\ea
where $\bar{\mY}=n^{-1}\od\bigoplus_{i=1}^n\mY^i$.
Then, it holds that
\begin{align}\label{cross covariance decomposition}
\begin{split}
&\|\hat{C}_{\mX\mY\mX}-C_{\mX\mY\mX}\|_{\HS}\\
&=\|(\hat{C}_{\mX\mY}-C_{\mX\mY})\circ(\hat{C}_{\mY\mX}-C_{\mY\mX}+C_{\mY\mX})+C_{\mX\mY}\circ(\hat{C}_{\mY\mX}-C_{\mY\mX})\|_{\HS}\\
&\leq\|\hat{C}_{\mX\mY}-C_{\mX\mY}\|_{\HS}\|\hat{C}_{\mY\mX}-C_{\mY\mX}\|_{\HS}+\|\hat{C}_{\mX\mY}-C_{\mX\mY}\|_{\HS}\|C_{\mY\mX}\|_{\HS}\\
&\quad+\|C_{\mX\mY}\|_{\HS}\|\hat{C}_{\mY\mX}-C_{\mY\mX}\|_{\HS}\\
&\leq(\|\hat{C}_{\mX\mY}-\check{C}_{\mX\mY}\|_{\HS}+\|\check{C}_{\mX\mY}-C_{\mX\mY}\|_{\HS})\cdot(\|\hat{C}_{\mY\mX}-\check{C}_{\mY\mX}\|_{\HS}+\|\check{C}_{\mY\mX}-C_{\mY\mX}\|_{\HS})\\
&\quad+(\|\hat{C}_{\mX\mY}-\check{C}_{\mX\mY}\|_{\HS}+\|\check{C}_{\mX\mY}-C_{\mX\mY}\|_{\HS})\cdot\|C_{\mY\mX}\|_{\HS}\\
&\quad+\|C_{\mX\mY}\|_{\HS}\cdot(\|\hat{C}_{\mY\mX}-\check{C}_{\mY\mX}\|_{\HS}+\|\check{C}_{\mY\mX}-C_{\mY\mX}\|_{\HS}).
\end{split}
\end{align}
We note that $\check{C}_{\mX\mY}$ and $\check{C}_{\mY\mX}$ are unbiased estimators of $C_{\mX\mY}$ and $C_{\mY\mX}$, respectively and that $\check{C}_{\mX\mY}$, $\check{C}_{\mY\mX}$, $C_{\mX\mY}$ and $C_{\mY\mX}$ are Hilbert-Schmidt operators. These facts with argument similar to (\ref{cov bound}) entail that $\E(\|\check{C}_{\mX\mY}-C_{\mX\mY}\|^2_{\rm HS})=O(n^{-1})$ and $\E(\|\check{C}_{\mY\mX}-C_{\mY\mX}\|^2_{\rm HS})=O(n^{-1})$ provided that $\E(\|\mX\|_*^2\|\mY\|^2)<\infty$. We can also prove that $\|\hat{C}_{\mX\mY}-\check{C}_{\mX\mY}\|_{\rm HS}=O_p(b_n)$ and $\|\hat{C}_{\mY\mX}-\check{C}_{\mY\mX}\|_{\rm HS}=O_p(b_n)$. Combining these results, we have $\|\hat{C}_{\mX\mY\mX}-C_{\mX\mY\mX}\|_{\rm HS}=O_p(n^{-1/2}+b_n)$. Hence, we get $\|\hat{\mphi}_r\om_*\mphi_r\|_*=O_p(n^{-1/2}+b_n)$ from (\ref{Koko book 2}).
Now, by arguing as in the proof of Proposition \ref{HPC}, we have $\max_{1\leq i\leq n}|\lng\mX^i\om_*\bar{\mX},\hat{\mphi}_r\rng_*-\lng\mX^i\om_*\E(\mX),\mphi_r\rng_*|=(n^{-1/2}+b_n)\cdot n^{1/\tau}$ in case $\E(\|\mX\|_*^{\tau})<\infty$ for some $\tau\geq2$, $\max_{1\leq i\leq n}|\lng\mX^i\om_*\bar{\mX},\hat{\mphi}_r\rng_*-\lng\mX^i\om_*\E(\mX),\mphi_r\rng_*|=(n^{-1/2}+b_n)\cdot \log{n}$ in case $\E\big(\exp(c\cdot\|\mX\|_*^{\tau})\big)<\infty$, and $\max_{1\leq i\leq n}|\lng\mX^i\om_*\bar{\mX},\hat{\mphi}_r\rng_*-\lng\mX^i\om_*\E(\mX),\mphi_r\rng_*|=n^{-1/2}+b_n$ in case $\mX$ is a bounded random element. \qed

\subsection{Proof of Proposition \ref{iRFPC}}\label{iRFPC proof}

By Proposition 2 in Lin and Yao (2019), $\Gamma_{\hat{\mu}_{Z},\mu_{Z}}$ is an unitary linear operator, that is, $\Gamma_{\hat{\mu}_{Z},\mu_{Z}}$ is a bijective bounded linear operator satisfying
\begin{align}\label{unitary property}
\lng\hat{V}_1,\hat{V}_2\rng_{\mathfrak{T}(\hat{\mu}_{Z})}=\lng\Gamma_{\hat{\mu}_{Z},\mu_{Z}}(\hat{V}_1),\Gamma_{\hat{\mu}_{Z},\mu_{Z}}(\hat{V}_2)\rng_{\mathfrak{T}(\mu_{Z})}
\end{align}
for all $\hat{V}_1,\hat{V}_2\in\mathfrak{T}(\hat{\mu}_{Z})$. By Proposition 2 and Theorem 7 in Lin and Yao (2019) and Lemma 2.3 in Horv\'{a}th and Kokoszka (2012), we have $\|\Gamma_{\hat{\mu}_{Z},\mu_{Z}}(\hat{\mpsi}_r)-\mpsi_r\|_{\mathfrak{T}(\mu_{Z})}=O(n^{-1/2})$. Also,
\begin{align}\label{parellel inequality all}
\begin{split}
&\max_{1\leq i\leq n}|\lng \Log_{\hat{\mu}_{Z}}Z^i,\hat{\mpsi}_r\rng_{\mathfrak{T}(\hat{\mu}_{Z})}-\lng\Log_{\mu_{Z}}Z^i,\mpsi_r\rng_{\mathfrak{T}(\mu_{Z})}|\\
&=\max_{1\leq i\leq n}|\lng \Gamma_{\hat{\mu}_{Z},\mu_{Z}}(\Log_{\hat{\mu}_{Z}}Z^i),\Gamma_{\hat{\mu}_{Z},\mu_{Z}}(\hat{\mpsi}_r)\rng_{\mathfrak{T}(\mu_{Z})}-\lng\Log_{\mu_{Z}}Z^i,\mpsi_r\rng_{\mathfrak{T}(\mu_{Z})}|\\
&=\max_{1\leq i\leq n}|\lng \Gamma_{\hat{\mu}_{Z},\mu_{Z}}(\Log_{\hat{\mu}_{Z}}Z^i)-\Log_{\mu_{Z}}Z^i,\Gamma_{\hat{\mu}_{Z},\mu_{Z}}(\hat{\mpsi}_r)-\mpsi_r+\mpsi_r\rng_{\mathfrak{T}(\mu_{Z})}\\
&\qquad\qquad+\lng\Log_{\mu_{Z}}Z^i,\Gamma_{\hat{\mu}_{Z},\mu_{Z}}(\hat{\mpsi}_r)-\mpsi_r\rng_{\mathfrak{T}(\mu_{Z})}|\\
&\leq(\|\Gamma_{\hat{\mu}_{Z},\mu_{Z}}(\hat{\mpsi}_r)-\mpsi_r\|_{\mathfrak{T}(\mu_{Z})}+\|\mpsi_r\|_{\mathfrak{T}(\mu_{Z})})\cdot\max_{1\leq i\leq n}\|\Gamma_{\hat{\mu}_{Z},\mu_{Z}}(\Log_{\hat{\mu}_{Z}}Z^i)-\Log_{\mu_{Z}}Z^i\|_{\mathfrak{T}(\mu_{Z})}\\
&\quad+\|\Gamma_{\hat{\mu}_{Z},\mu_{Z}}(\hat{\mpsi}_r)-\mpsi_r\|_{\mathfrak{T}(\mu_{Z})}\cdot\max_{1\leq i\leq n}\|\Log_{\mu_{Z}}Z^i\|_{\mathfrak{T}(\mu_{Z})}\\
&=O_p(1)\cdot\max_{1\leq i\leq n}\|\Gamma_{\hat{\mu}_{Z},\mu_{Z}}(\Log_{\hat{\mu}_{Z}}Z^i)-\Log_{\mu_{Z}}Z^i\|_{\mathfrak{T}(\mu_{Z})}+O_p(n^{-1/2})\cdot \max_{1\leq i\leq n}\|\Log_{\mu_{Z}}Z^i\|_{\mathfrak{T}(\mu_{Z})},
\end{split}
\end{align}
where the first equality follows from (\ref{unitary property}).
Using the first inequality at (5.7) in Kendall and Le (2011) and condition (L8), it holds that
\ba
&\|\mathcal{P}_{\hat{\mu}_{Z}(t),\mu_{Z}(t)}(\Log_{\hat{\mu}_{Z}(t)}Z^i(t))-\Log_{\mu_{Z}(t)}Z^i(t)\|_{\mu_{Z}(t)}\\
&\leq d_{\mathcal{M}}(\hat{\mu}_{Z}(t),\mu_{Z}(t))\cdot\sup_{\tilde{\mu}_{Z}(t)\in B_{\mathcal{M}}(\mu_{Z}(t),\rho)}\|(\nabla\Log_{(\cdot)}Z^i(t))(\tilde{\mu}_{Z}(t))\|_{\rm op}\\
&\leq d_{\mathcal{M}}(\hat{\mu}_{Z}(t),\mu_{Z}(t))\cdot\sup_{\tilde{\mu}_{Z}(t)\in B_{\mathcal{M}}(\mu_{Z}(t),\rho)}(c\cdot(1+d_{\mathcal{M}}(\tilde{\mu}_{Z}(t),Z^i(t))))\\
&\leq d_{\mathcal{M}}(\hat{\mu}_{Z}(t),\mu_{Z}(t))\cdot\sup_{\tilde{\mu}_{Z}(t)\in B_{\mathcal{M}}(\mu_{Z}(t),\rho)}(c\cdot(1+d_{\mathcal{M}}(\tilde{\mu}_{Z}(t),\mu_{Z}(t))+d_{\mathcal{M}}(\mu_{Z}(t),Z^i(t))))\\
&\leq d_{\mathcal{M}}(\hat{\mu}_{Z}(t),\mu_{Z}(t))\cdot(c\cdot(1+\rho+d_{\mathcal{M}}(\mu_{Z}(t),Z^i(t))))\\
&=d_{\mathcal{M}}(\hat{\mu}_{Z}(t),\mu_{Z}(t))\cdot(c\cdot(1+\rho+\|\Log_{\mu_{Z}(t)}Z^i(t)\|_{\mu_{Z}(t)}))\\
&\leq \const\cdot d_{\mathcal{M}}(\hat{\mu}_{Z}(t),\mu_{Z}(t))+\const\cdot d_{\mathcal{M}}(\hat{\mu}_{Z}(t),\mu_{Z}(t))\cdot\|\Log_{\mu_{Z}(t)}Z^i(t)\|_{\mu_{Z}(t)}
\ea
for sufficiently small $\rho>0$, where $B_{\mathcal{M}}(\mu_{Z}(t),\rho)=\{p\in\mathcal{M}:d_{\mathcal{M}}(p,\mu_{Z}(t))<\rho\}$.
Hence,
\begin{align}\label{parellel inequality2}
\begin{split}
&\|\Gamma_{\hat{\mu}_{Z},\mu_{Z}}(\Log_{\hat{\mu}_{Z}}Z^i)-\Log_{\mu_{Z}}Z^i\|^2_{\mathfrak{T}(\mu_{Z})}\\
&=\int_{\mathcal{T}}\|\mathcal{P}_{\hat{\mu}_{Z}(t),\mu_{Z}(t)}(\Log_{\hat{\mu}_{Z}(t)}Z^i(t))-\Log_{\mu_{Z}(t)}Z^i(t)\|^2_{\mu_{Z}(t)}d\nu(t)\\
&\leq \const\cdot\int_{\mathcal{T}}(d_{\mathcal{M}}(\hat{\mu}_{Z}(t),\mu_{Z}(t)))^2d\nu(t)\\
&\quad+\const\cdot\int_{\mathcal{T}}(d_{\mathcal{M}}(\hat{\mu}_{Z}(t),\mu_{Z}(t)))^2\cdot\|\Log_{\mu_{Z}(t)}Z^i(t)\|^2_{\mu_{Z}(t)}d\nu(t)\\
&\leq O_p(n^{-1})+\const\cdot\sup_{t\in\mathcal{T}}(d_{\mathcal{M}}(\hat{\mu}_{Z}(t),\mu_{Z}(t)))^2\cdot\|\Log_{\mu_{Z}}Z^i\|^2_{\mathfrak{T}(\mu_{Z})},
\end{split}
\end{align}
where $\mathcal{P}_{\hat{\mu}_{Z}(t),\mu_{Z}(t)}:T_{\hat{\mu}(t)}\mathcal{M}\rightarrow T_{\mu(t)}\mathcal{M}$ is the parallel transport map along $\gamma(t,\cdot)$ and the last inequality follows from Theorem 6 in Lin and Yao (2019). See also Lemma S3 in Lin et al. (2022) for the case where $\mathcal{T}$ is a singleton and $\nu$ is the counting measure. By Theorem 6 in Lin and Yao (2019), we have $\sup_{t\in\mathcal{T}}(d_{\mathcal{M}}(\hat{\mu}_{Z}(t),\mu_{Z}(t)))^2=O_p(n^{-1})$. In case $\E(\|\Log_{\mu_{Z}}Z\|_{\mathfrak{T}(\mu_{Z})}^{\tau})<\infty$ for some $\tau\geq4$, Lemma \ref{markov} gives that
\begin{align}\label{log max order}
\max_{1\leq i\leq n}\|\Log_{\mu_{Z}}Z^i\|_{\mathfrak{T}(\mu_{Z})}=O_p(n^{1/\tau}).
\end{align}
This with (\ref{parellel inequality2}) gives that $\max_{1\leq i\leq n}\|\Gamma_{\hat{\mu}_{Z},\mu_{Z}}(\Log_{\hat{\mu}_{Z}}Z^i)-\Log_{\mu_{Z}}Z^i\|_{\mathfrak{T}(\mu_{Z})}=O_p(n^{-1/2+1/\tau})$.
The latter with (\ref{parellel inequality all}) and (\ref{log max order}) gives the first assertion. The second and third assertions follow similarly. \qed

\subsection{Proof of Proposition \ref{iRHSC}}\label{iRHSC proof}

Define $\hat{C}_{X\mY}^{\Gamma}:\mbH\rightarrow\mathfrak{T}(\mu_{Z})$ and $\hat{C}_{\mY X}^{\Gamma}:\mathfrak{T}(\mu_{Z})\rightarrow\mbH$ by
\ba
\hat{C}_{X\mY}^{\Gamma}(\mh)&=\frac{1}{n-1}\sum_{i=1}^n(\lng\tilde{\mY}^i\om\bar{\tilde{\mY}},\mh\rng\cdot\Gamma_{\hat{\mu}_{Z},\mu_{Z}}(\Log_{\hat{\mu}_{Z}}Z^i)),\quad\mbox{and}\\
\hat{C}_{\mY X}^{\Gamma}(V)&=\frac{1}{n-1}\od\bigoplus_{i=1}^n(\lng\Gamma_{\hat{\mu}_{Z},\mu_{Z}}(\Log_{\hat{\mu}_{Z}}Z^i),V\rng_{\mathfrak{T}(\mu_{Z})}\od(\tilde{\mY}^i\om\bar{\tilde{\mY}})).
\ea
and define $\hat{C}_{X\mY X}^{\Gamma}=\hat{C}_{X\mY}^{\Gamma}\circ\hat{C}_{\mY X}^{\Gamma}:\mathfrak{T}(\mu_{Z})\rightarrow\mathfrak{T}(\mu_{Z})$. A direct computation with Proposition 2 in Lin and Yao (2019) shows that $\Phi_{\hat{\mu}_{Z},\mu_{Z}}(\hat{C}_{X\mY X})=\hat{C}_{X\mY X}^{\Gamma}$.
We also define $\check{C}_{X\mY}^{\Gamma}:\mbH\rightarrow\mathfrak{T}(\mu_{Z})$ and $\check{C}_{\mY X}^{\Gamma}:\mathfrak{T}(\mu_{Z})\rightarrow\mbH$ by
\ba
\check{C}^{\Gamma}_{X\mY}(\mh)&=\frac{1}{n-1}\sum_{i=1}^n(\lng\mY^i\om\bar{\mY},\mh\rng\cdot\Gamma_{\hat{\mu}_{Z},\mu_{Z}}(\Log_{\hat{\mu}_{Z}}Z^i)),\quad\mbox{and}\\
\check{C}^{\Gamma}_{\mY X}(V)&=\frac{1}{n-1}\od\bigoplus_{i=1}^n(\lng\Gamma_{\hat{\mu}_{Z},\mu_{Z}}(\Log_{\hat{\mu}_{Z}}Z^i),V\rng_{\mathfrak{T}(\mu_{Z})}\od(\mY^i\om\bar{\mY})).
\ea
By arguing as in (\ref{cross covariance decomposition}), we get
\begin{align}\label{cross covariance decomp2}
\begin{split}
&\|\Phi_{\hat{\mu}_{Z},\mu_{Z}}(\hat{C}_{X\mY X})-C_{X\mY X}\|_{\HS}\\
&=\|\hat{C}_{X\mY X}^{\Gamma}-C_{X\mY X}\|_{\HS}\\
&\leq(\|\hat{C}^{\Gamma}_{X\mY}-\check{C}^{\Gamma}_{X\mY}\|_{\HS}+\|\check{C}^{\Gamma}_{X\mY}-C_{X\mY}\|_{\HS})\cdot(\|\hat{C}^{\Gamma}_{\mY X}-\check{C}^{\Gamma}_{\mY X}\|_{\HS}+\|\check{C}^{\Gamma}_{\mY X}-C_{\mY X}\|_{\HS})\\
&\quad+(\|\hat{C}^{\Gamma}_{X\mY}-\check{C}^{\Gamma}_{X\mY}\|_{\HS}+\|\check{C}^{\Gamma}_{X\mY}-C_{X\mY}\|_{\HS})\cdot\|C_{\mY X}\|_{\HS}\\
&\quad+\|C_{X\mY}\|_{\HS}\cdot(\|\hat{C}^{\Gamma}_{\mY X}-\check{C}^{\Gamma}_{\mY X}\|_{\HS}+\|\check{C}^{\Gamma}_{\mY X}-C_{\mY X}\|_{\HS}).
\end{split}
\end{align}
We first claim that $\|\check{C}^{\Gamma}_{X\mY}-C_{X\mY}\|^2_{\HS}=O_p(n^{-1})$. Note that
\begin{align}\label{check c gamma decomp}
\begin{split}
&\check{C}^{\Gamma}_{X\mY}(\cdot)\\
&=\frac{1}{n-1}\sum_{i=1}^n\lng\mY^i\om\bar{\mY},\cdot\rng\cdot(\Gamma_{\hat{\mu}_{Z},\mu_{Z}}(\Log_{\hat{\mu}_{Z}}Z^i)-\Log_{\mu_{Z}}Z^i)+\frac{1}{n-1}\sum_{i=1}^n\lng\mY^i\om\bar{\mY},\cdot\rng\cdot\Log_{\mu_{Z}}Z^i
\end{split}
\end{align}
A direct computation shows that the second operator on the right hand side of (\ref{check c gamma decomp}) is an unbiased estimator of the cross-covariance operator $C_{X\mY}(\cdot)=\E(\lng\mY\om\E(\mY),\cdot\rng\cdot\Log_{\mu_{Z}}Z)$. Hence,
\ba
&\E\bigg(\bigg\|\frac{1}{n-1}\sum_{i=1}^n\lng\mY^i\om\bar{\mY},\cdot\rng\cdot\Log_{\mu_{Z}}Z^i-C_{X\mY}(\cdot)\bigg\|^2_{\HS}\bigg)\\
&\leq\E\bigg(\bigg\|\frac{1}{n-1}\sum_{i=1}^n\lng\mY^i\om\bar{\mY},\cdot\rng\cdot\Log_{\mu_{Z}}Z^i\bigg\|^2_{\HS}\bigg)\\
&=O(n^{-1}).
\ea
Also, the first operator on the right hand side of (\ref{check c gamma decomp}) satisfies that
\ba
&\bigg\|\frac{1}{n-1}\sum_{i=1}^n\lng\mY^i\om\bar{\mY},\cdot\rng\cdot(\Gamma_{\hat{\mu}_{Z},\mu_{Z}}(\Log_{\hat{\mu}_{Z}}Z^i)-\Log_{\mu_{Z}}Z^i)\bigg\|^2_{\HS}\\
&\leq\frac{\const}{(n-1)^2}\sum_{i=1}^n\sum_{j=1}^n(\|\mY^i\om\bar{\mY}\|^2+\|\mY^j\om\bar{\mY}\|^2)(\|\Gamma_{\hat{\mu}_{Z},\mu_{Z}}(\Log_{\hat{\mu}_{Z}}Z^i)-\Log_{\mu_{Z}}Z^i\|^2_{\mathfrak{T}(\mu_{Z})}\\
&\qquad\qquad\qquad\qquad\qquad\qquad\qquad+\|\Gamma_{\hat{\mu}_{Z},\mu_{Z}}(\Log_{\hat{\mu}_{Z}}Z^j)-\Log_{\mu_{Z}}Z^j\|^2_{\mathfrak{T}(\mu_{Z})})\\
&\leq\frac{\const}{(n-1)^2}\sum_{i=1}^n\sum_{j=1}^n(\|\mY^i\om\bar{\mY}\|^2+\|\mY^j\om\bar{\mY}\|^2)(O_p(n^{-1})+O_p(n^{-1})\cdot\|\Log_{\mu_{Z}}Z^i\|^2_{\mathfrak{T}(\mu_{Z})}\\
&\qquad\qquad\qquad\qquad\qquad\qquad\qquad+O_p(n^{-1})+O_p(n^{-1})\cdot\|\Log_{\mu_{Z}}Z^j\|^2_{\mathfrak{T}(\mu_{Z})})\\
&=O_p(n^{-1}),
\ea
where the second inequality follows from (\ref{parellel inequality2}). Thus, we have $\|\check{C}^{\Gamma}_{X\mY}-C_{X\mY}\|^2_{\rm HS}=O_p(n^{-1})$. We may similarly prove that $\|\check{C}^{\Gamma}_{\mY X}-C_{\mY X}\|^2_{\rm HS}=O_p(n^{-1})$.
We can also prove that
\ba
&\|\hat{C}^{\Gamma}_{X\mY}-\check{C}^{\Gamma}_{X\mY}\|_{\HS}\leq \max_{1\leq i\leq n}\|\tilde{\mY}^i\om\mY^i\|\cdot O_p(1)=O_p(b_n),\quad\mbox{and}\\
&\|\hat{C}^{\Gamma}_{\mY X}-\check{C}^{\Gamma}_{\mY X}\|_{\HS}\leq \max_{1\leq i\leq n}\|\tilde{\mY}^i\om\mY^i\|\cdot O_p(1)=O_p(b_n).
\ea
Plugging those rates into (\ref{cross covariance decomp2}) gives that $\|\hat{C}^{\Gamma}_{X\mY X}-C_{X\mY X}\|_{\HS}=O_p(n^{-1/2}+b_n)$. By Proposition 2 in Lin and Yao (2019), $\Gamma_{\hat{\mu}_{Z},\mu_{Z}}(\hat{\mphi}_r)$ is indeed an eigenfunction induced from the spectral theorem on $\hat{C}^{\Gamma}_{X\mY X}$. Hence, Lemma 2.3 in Horv\'{a}th and Kokoszka (2012) implies that $\|\Gamma_{\hat{\mu}_{Z},\mu_{Z}}(\hat{\mphi}_r)-\mphi_r\|_{\mathfrak{T}(\mu_{Z})}=O_p(n^{-1/2}+b_n)$.
Now, by arguing as in the proof of Proposition \ref{iRFPC}, we get the desired rates for
$\max_{1\leq i\leq n}|\lng \Log_{\hat{\mu}_{Z}}Z^i,\hat{\mphi}_r\rng_{\mathfrak{T}(\hat{\mu}_{Z})}-\lng\Log_{\mu_{Z}}Z^i,\mphi_r\rng_{\mathfrak{T}(\mu_{Z})}|$. \qed

\subsection{Proof of Proposition \ref{imperfect density prop}}

We first approximate $\max_{1\leq i\leq n}\sup_{s\in\mathcal{S}}|\tilde{Y}^i(s)-Y^i(s)|$. We write $Y^i(\cdot)$ by $g_i(\cdot)$ for simplicity. We first find the rate of $\max_{1\leq i\leq n}\sup_{s\in\mathcal{S}}|\hat{g}_i(s)-g_i(s)|$. Note that
\ba
\max_{1\leq i\leq n}\sup_{s\in\mathcal{S}}|\hat{g}_i(s)-g_i(s)|\leq \max_{1\leq i\leq n}\sup_{s\in\mathcal{S}}|\hat{g}_i(s)-\E(\hat{g}_i(s))|+\max_{1\leq i\leq n}\sup_{s\in\mathcal{S}}|\E(\hat{g}_i(s))-g_i(s)|.
\ea
Here,
\ba
\E(\hat{g}_i(s))-g_i(s)&=w(s,h_i^*)\int_{\mathcal{S}}(h_i^*)^{-q}K^*\bigg(\frac{\|y-s\|_{\mathbb{R}^q}}{h_i^*}\bigg)g_i(y)dy-g_i(s)\\
&=w(s,h_i^*)\int_{\mathcal{S}}(h_i^*)^{-q}K^*\bigg(\frac{\|y-s\|_{\mathbb{R}^q}}{h_i^*}\bigg)(g_i(y)-g_i(s))dy\\
&=w(s,h_i^*)\int_{(\mathcal{S}-s)/h_i^*\cap B(\bzero_{\mathbb{R}^q},1)}K^*(\|t\|_{\mathbb{R}^q})(g_i(s+th_i^*)-g_i(s))dt,
\ea
where $(\mathcal{S}-s)/h_i^*=\{(y-s)/h_i^*:y\in\mathcal{S}\}$, $B(\bzero_{\mathbb{R}^q},1)$ is the open ball centered at the zero vector $\bzero_{\mathbb{R}^q}$ of $\mathbb{R}^q$ with radius 1 and $s^*$ lies between $s+th_i^*$ and $s$. Since $w(s,h_i^*)$ is bounded uniformly over $n$ and $s\in\mathcal{S}$ and $g_i$ are Lipschitz continuous with a common Lipschitz constant, we get
\ba
\max_{1\leq i\leq n}\sup_{s\in\mathcal{S}}|\E(\hat{g}_i(s))-g_i(s)|&\leq \const\cdot h^*_{\rm max}\cdot\int_{B(\bzero_{\mathbb{R}^q},1)}K^*(\|t\|_{\mathbb{R}^q})\|t\|_{\mathbb{R}^q}dt\\
&=O(h^*_{\rm max}),
\ea
where $h^*_{\rm max}=\max_{1\leq i\leq n}h_i^*$. We now approximate $\max_{1\leq i\leq n}\sup_{s\in\mathcal{S}}|\hat{g}_i(s)-\E(\hat{g}_i(s))|$. Note that $\hat{g}_i(s)-\E(\hat{g}_i(s))=w(s,h_i^*)(\hat{f}_i(s)-\E(\hat{f}_i(s)))$, where
\ba
\hat{f}_i(s)=\frac{1}{n_i^*(h_i^*)^{q}}\sum_{j=1}^{n_i^*}K^*\bigg(\frac{\|s-Y^*_{ij}\|_{\mathbb{R}^{q}}}{h_i^*}\bigg).
\ea
Hence, it suffices to find the rate of $\max_{1\leq i\leq n}\sup_{s\in\mathcal{S}}|\hat{f}_i(s)-\E(\hat{f}_i(s))|$. For $\varsigma>0$, there exist $N(n^{-\varsigma})\in\mathbb{N}$ and $\{s^{(1)},\cdots,s^{(N(n^{-\varsigma}))}\}\subset \mathcal{S}$ such that $N(n^{-\varsigma})=O(n^{\varsigma q})$ and $\{B(s^{(l)},n^{-\varsigma}):1\leq l\leq N(n^{-\varsigma})\}$ covers $\mathcal{S}$, where $B(s^{(l)},n^{-\varsigma})$ is the open ball centered at $s^{(l)}$ with radius $n^{-\varsigma}$. Note that
\begin{align}\label{imperfect density proof part1}
\begin{split}
\max_{1\leq i\leq n}\sup_{s\in\mathcal{S}}|\hat{f}_i(s)-\E(\hat{f}_i(s))|&\leq\max_{1\leq i\leq n}\max_{1\leq l\leq N(n^{-\gamma})}|\hat{f}_i(s^{(l)})-\E(\hat{f}_i(s^{(l)}))|\\
&\quad+\max_{1\leq i\leq n}\max_{1\leq l\leq N(n^{-\gamma})}\sup_{s\in \mathcal{S}\cap B(s^{(l)},n^{-\varsigma})}|\hat{f}_i(s)-\hat{f}_i(s^{(l)})|\\
&\quad+\max_{1\leq i\leq n}\max_{1\leq l\leq N(n^{-\gamma})}\sup_{s\in \mathcal{S}\cap B(s^{(l)},n^{-\varsigma})}|\E(\hat{f}_i(s))-\E(\hat{f}_i(s^{(l)}))|.
\end{split}
\end{align}
Here,
\ba
|\hat{f}_i(s)-\hat{f}_i(s^{(l)})|\leq\const\frac{\|s-s^{(l)}\|_{\mathbb{R}^q}}{(h_i^*)^{q+1}},
\ea
where $\const$ depends only on $K^*$. Hence, the second and third terms on the right hand side of (\ref{imperfect density proof part1}) are negligible by taking sufficiently large $\varsigma>0$. We also note that
\ba
&\Prob\bigg(\max_{1\leq i\leq n}\max_{1\leq l\leq N(n^{-\gamma})}|\hat{f}_i(s^{(l)})-\E(\hat{f}_i(s^{(l)}))|>T\cdot H_n\bigg)\\
&\leq O(n^{\varsigma q+1})\max_{1\leq i\leq n}\max_{1\leq l\leq N(n^{-\gamma})}\Prob\big(|\hat{f}_i(s^{(l)})-\E(\hat{f}_i(s^{(l)}))|>T\cdot H_n\big),
\ea
where
\ba
H_n=\bigg(\frac{\log{n}}{\min_{1\leq i\leq n}(n_i^*(h_i^*)^{q})}\bigg)^{1/2}.
\ea
Define
\ba
Z_{ijl}=\frac{1}{n_i^*(h_i^*)^q}\bigg(K^*\bigg(\frac{\|s^{(l)}-Y_{ij}^*\|_{\mathbb{R}^{q}}}{h_i^*}\bigg)-\E\bigg(K^*\bigg(\frac{\|s^{(l)}-Y_{ij}^*\|_{\mathbb{R}^{q}}}{h_i^*}\bigg)\bigg)\bigg).
\ea
It holds that
\ba
\E(Z_{ijl})=0,\quad\|Z_{ijl}\|\leq\const(n_i^*(h_i^*)^q)^{-1},\quad\mbox{and}\quad\sum_{j=1}^{n_i^*}\E(Z_{ijl}^2)\leq\const(n_i^*(h_i^*)^q)^{-1},
\ea
where the $\const$ are independent of $i,j,l$ and $s$. Now, Theorem 2.6.2 in Bosq (2000) gives that
\ba
\Prob\big(|\hat{f}_i(s^{(l)})-\E(\hat{f}_i(s^{(l)}))|>T\cdot H_n\big)&\leq 2\exp(-T\cdot H_n^2 \cdot n_i^*(h_i^*)^q)\\
&\leq 2\exp(-T\cdot H_n^2\cdot\min_{1\leq i\leq n}(n_i^*(h_i^*)^q))\\
&\leq 2n^{-T}
\ea
for sufficiently large $T$ and $n$, which implies that
\ba
\max_{1\leq i\leq n}\max_{1\leq l\leq N(n^{-\gamma})}|\hat{f}_i(s^{(l)})-\E(\hat{f}_i(s^{(l)}))|=O_p(H_n).
\ea
Hence, we have $\max_{1\leq i\leq n}\sup_{s\in\mathcal{S}}|\hat{g}_i(s)-g_i(s)|=O_p(h^*_{\rm max}+H_n)$. Since
\ba
|\tilde{Y}^i(s)-Y^i(s)|&\leq\frac{|\hat{g}_i(s)-g_i(s)|}{\int_{\mathcal{S}}\hat{g}_i(s)ds}+\frac{g_i(s)\int_{\mathcal{S}}|\hat{g}_i(s)-g_i(s)|ds}{\int_{\mathcal{S}}\hat{g}_i(s)ds}\\
&\leq\frac{|\hat{g}_i(s)-g_i(s)|}{1-\Leb_{\mathbb{R}^q}(\mathcal{S})\sup_{s\in\mathcal{S}}|\hat{g}_i(s)-g_i(s)|}+\frac{g_i(s)\Leb_{\mathbb{R}^q}(\mathcal{S})\sup_{s\in\mathcal{S}}|\hat{g}_i(s)-g_i(s)|}{1-\Leb_{\mathbb{R}^q}(\mathcal{S})\sup_{s\in\mathcal{S}}|\hat{g}_i(s)-g_i(s)|},
\ea
we get
\begin{align}\label{density first result}
\max_{1\leq i\leq n}\sup_{s\in\mathcal{S}}|\tilde{Y}^i(s)-Y^i(s)|=O_p(h^*_{\rm max}+H_n).
\end{align}
We now approximate $\max_{1\leq i\leq n}\|\tilde{\mY}^i\om\mY^i\|$. We first claim that $\max_{1\leq i\leq n}\int_{\mathcal{S}}(\log(\tilde{Y}^i(s)))^2ds<\infty$ with probability tending to one. Since $\max_{1\leq i\leq n}\int_{\mathcal{S}}(\log(Y^i(s)))^2ds<\infty$ almost surely, it suffices to prove that
\ba
\max_{1\leq i\leq n}\int_{\mathcal{S}}(\log(\tilde{Y}^i(s))-\log(Y^i(s)))^2ds<\infty
\ea
with probability tending to one. Since
\ba
\log(\tilde{Y}^i(s))-\log(Y^i(s))=\frac{1}{\breve{Y}^i(s)}(\tilde{Y}^i(s)-Y^i(s))
\ea
for some $\breve{Y}^i(s)$ lying between $\tilde{Y}^i(s)$ and $Y^i(s)$, we get
\ba
\int_{\mathcal{S}}(\log(\tilde{Y}^i(s))-\log(Y^i(s)))^2ds&\leq\frac{1}{(\inf_{s\in\mathcal{S}}\breve{Y}^i(s))^2}\int_{\mathcal{S}}(\tilde{Y}^i(s)-Y^i(s))^2ds\\
&\leq\frac{\Leb_{\mathbb{R}^q}(\mathcal{S})\sup_{s\in\mathcal{S}}|\tilde{Y}^i(s)-Y^i(s)|^2}{(\inf_{s\in\mathcal{S}}Y^i(s)-\sup_{s\in\mathcal{S}}|\tilde{Y}^i(s)-Y^i(s)|)^2},
\ea
where the second inequality holds with probability tending to one by (\ref{density first result}). Hence,
\ba
&\max_{1\leq i\leq n}\int_{\mathcal{S}}(\log(\tilde{Y}^i(s))-\log(Y^i(s)))^2ds\\
&\leq\frac{\Leb_{\mathbb{R}^q}(\mathcal{S})\max_{1\leq i\leq n}\sup_{s\in\mathcal{S}}|\tilde{Y}^i(s)-Y^i(s)|^2}{(\min_{1\leq i\leq n}\inf_{s\in\mathcal{S}}Y^i(s)-\max_{1\leq i\leq n}\sup_{s\in\mathcal{S}}|\tilde{Y}^i(s)-Y^i(s)|)^2}
\ea
with probability tending to one by (\ref{density first result}). This with (\ref{density first result}) proves the claim. To find the rate of $\max_{1\leq i\leq n}\|\tilde{\mY}^i\om\mY^i\|$,
we let $\Leb_{\mathbb{R}^q}$ denote the Lebesgue measure on $\mathbb{R}^q$. A direct computation shows that
\ba \|\tilde{\mY}^i\om\mY^i\|^2&=\int_{\mathcal{S}}(\log(\tilde{Y}^i(s))-\log(Y^i(s)))^2ds-\frac{1}{\Leb_{\mathbb{R}^q}(\mathcal{S})}\bigg(\int_{\mathcal{S}}\log(\tilde{Y}^i(s))-\log(Y^i(s))ds\bigg)^2\\
&\leq\int_{\mathcal{S}}(\log(\tilde{Y}^i(s))-\log(Y^i(s)))^2ds.
\ea
Hence,
\ba
\max_{1\leq i\leq n}\|\tilde{\mY}^i\om\mY^i\|\leq\frac{\Leb_{\mathbb{R}^q}(\mathcal{S})^{1/2}\max_{1\leq i\leq n}\sup_{s\in\mathcal{S}}|\tilde{Y}^i(s)-Y^i(s)|}{\min_{1\leq i\leq n}\inf_{s\in\mathcal{S}}Y^i(s)-\max_{1\leq i\leq n}\sup_{s\in\mathcal{S}}|\tilde{Y}^i(s)-Y^i(s)|}
\ea
with probability tending to one. Thus, the desired rate follows from (\ref{density first result}). \qed

\subsection{Proof of Proposition \ref{imperfect density prop 2}}

We first find the rate of $\max_{1\leq i\leq n}\sup_{s\in\mathcal{S}}|\tilde{Y}^i(s)-Y^i(s)|$. We write $Y^i(\cdot)$ by $g_i(\cdot)$ and $\tilde{Y}^i(\cdot)$ by $\tilde{g}_i$ for simplicity. Note that
\ba
\max_{1\leq i\leq n}\sup_{s\in\mathcal{S}}|\tilde{g}_i(s)-g_i(s)|\leq \max_{1\leq i\leq n}\sup_{s\in\mathcal{S}}|\tilde{g}_i(s)-\E(\tilde{g}_i(s))|+\max_{1\leq i\leq n}\sup_{s\in\mathcal{S}}|\E(\tilde{g}_i(s))-g_i(s)|.
\ea
We first find the rate of the second term on the right hand side. Define $\Exp^*_s=\Exp_s\circ \iota_s:\mathbb{R}^q\rightarrow\mathcal{S}$, where $\Exp_s:T_s\mathcal{S}\rightarrow\mathcal{S}$ is the exponential map and $\iota_s:\mathbb{R}^q\rightarrow T_s\mathcal{S}$ is an isometric isomorphism for each $s\in\mathbb{S}$. Note that
\ba
&=\int_{\mathcal{S}}(h_i^*)^{-q}\frac{1}{\theta_{\mathcal{S}}(s;y)}K^*\bigg(\frac{d_{\mathcal{S}}(s,y)}{h_i^*}\bigg)d\nu(y)\\
&=\int_{B(\bzero_{\mathbb{R}^q},h_i^*)}(h_i^*)^{-q}\frac{1}{\theta_{\mathcal{S}}(s;\Exp_s^*(t))}K^*\bigg(\frac{d_{\mathcal{S}}(s,\Exp_s^*(t))}{h_i^*}\bigg)\sqrt{\det G(\Exp_s^*(t))}dt\\
&=\int_{B(\bzero_{\mathbb{R}^q},h_i^*)}(h_i^*)^{-q}K^*\bigg(\frac{\|t\|_{\mathbb{R}^q}}{h_i^*}\bigg)dt\\
&=\int_{B(\bzero_{\mathbb{R}^q},1)}K^*(\|u\|_{\mathbb{R}^q})du\\
&=1.
\ea
Hence,
\ba
\max_{1\leq i\leq n}\sup_{s\in\mathcal{S}}|\E(\tilde{g}_i(s))-g_i(s)|&\leq\max_{1\leq i\leq n}\sup_{s\in\mathcal{S}}\int_{\mathcal{S}}(h_i^*)^{-q}\frac{1}{\theta_{\mathcal{S}}(s;y)}K^*\bigg(\frac{d_{\mathcal{S}}(s,y)}{h_i^*}\bigg)|g_i(y)-g_i(s)|d\nu(y)\\
&\leq \const\max_{1\leq i\leq n}\sup_{s\in\mathcal{S}}\int_{\mathcal{S}}(h_i^*)^{-q}\frac{1}{\theta_{\mathcal{S}}(s;y)}K^*\bigg(\frac{d_{\mathcal{S}}(s,y)}{h_i^*}\bigg)d_{\mathcal{S}}(s,y)d\nu(y)\\
&\leq \const\max_{1\leq i\leq n} h_i^*\sup_{s\in\mathcal{S}}\int_{\mathcal{S}}(h_i^*)^{-q}\frac{1}{\theta_{\mathcal{S}}(s;y)}K^*\bigg(\frac{d_{\mathcal{S}}(s,y)}{h_i^*}\bigg)d\nu(y)\\
&=\const\max_{1\leq i\leq n}h_i^*.
\ea
We now find the rate of $\max_{1\leq i\leq n}\sup_{s\in\mathcal{S}}|\tilde{g}_i(s)-\E(\tilde{g}_i(s))|$. For $\varsigma>0$, there exist $N(n^{-\varsigma})\in\mathbb{N}$ and $\{s^{(1)},\cdots,s^{(N(n^{-\varsigma}))}\}\subset \mathcal{S}$ such that $N(n^{-\varsigma})=O(n^{\varsigma q})$ and $\{B_{\mathcal{S}}(s^{(l)},n^{-\varsigma}):1\leq l\leq N(n^{-\varsigma})\}$ covers $\mathcal{S}$, where $B_{\mathcal{S}}(s,r)=\{y\in\mathcal{S}:d_{\mathcal{S}}(s,y)<r\}$ for $s\in\mathcal{S}$ and $r>0$. Note that
\begin{align}\label{imperfect density proof part2}
\begin{split}
\max_{1\leq i\leq n}\sup_{s\in\mathcal{S}}|\tilde{g}_i(s)-\E(\tilde{g}_i(s))|&\leq\max_{1\leq i\leq n}\max_{1\leq l\leq N(n^{-\gamma})}|\tilde{g}_i(s^{(l)})-\E(\tilde{g}_i(s^{(l)}))|\\
&\quad+\max_{1\leq i\leq n}\max_{1\leq l\leq N(n^{-\gamma})}\sup_{s\in B_{\mathcal{S}}(s^{(l)},n^{-\varsigma})}|\tilde{g}_i(s)-\tilde{g}_i(s^{(l)})|\\
&\quad+\max_{1\leq i\leq n}\max_{1\leq l\leq N(n^{-\gamma})}\sup_{s\in B_{\mathcal{S}}(s^{(l)},n^{-\varsigma})}|\E(\tilde{g}_i(s))-\E(\tilde{g}_i(s^{(l)}))|.
\end{split}
\end{align}
Here,
\ba
|\tilde{g}_i(s)-\tilde{g}_i(s^{(l)})|\leq\frac{1}{n_i^*(h_i^*)^{q}}\sum_{j=1}^{n_i^*}|U_{ijl}|,
\ea
where
\ba
U_{ijl}=K^*\bigg(\frac{d_{\mathcal{S}}(s,Y^*_{ij})}{h_i^*}\bigg)\frac{1}{\theta_{\mathcal{S}}(s;Y^*_{ij})}-K^*\bigg(\frac{d_{\mathcal{S}}(s^{(l)},Y^*_{ij})}{h_i^*}\bigg)\frac{1}{\theta_{\mathcal{S}}(s^{(l)};Y^*_{ij})}.
\ea
If $d_{\mathcal{S}}(s,Y^*_{ij})>h_i^*$ and $d_{\mathcal{S}}(s^{(l)},Y^*_{ij})>h_i^*$, then $U_{ijl}=0$. If $d_{\mathcal{S}}(s,Y^*_{ij})\leq h_i^*$ or $d_{\mathcal{S}}(s^{(l)},Y^*_{ij})\leq h_i^*$, then $\max\{d_{\mathcal{S}}(s,Y^*_{ij}),d_{\mathcal{S}}(s^{(l)},Y^*_{ij})\}<h_i^*+n^{-\varsigma}\leq\max_{1\leq i\leq n}h_i^*+n^{-\varsigma}<R$ for $n>m$, where $m$ is some number independent of $i, j, l$ and $s$. In the latter case, $|(\theta_{\mathcal{S}}(s;Y^*_{ij}))^{-1}-(\theta_{\mathcal{S}}(s^{(l)};Y^*_{ij}))^{-1}|\leq L\cdot d_{\mathcal{S}}(s,s^{(l)})<L\cdot n^{-\varsigma}$ for $n>m$. Also,
\ba
(\theta_{\mathcal{S}}(s^{(l)};Y^*_{ij}))^{-1}&\leq|(\theta_{\mathcal{S}}(s^{(l)};Y^*_{ij}))^{-1}-(\theta_{\mathcal{S}}(Y^*_{ij};Y^*_{ij}))^{-1}|+|\theta_{\mathcal{S}}(Y^*_{ij};Y^*_{ij}))^{-1}|\\
&=|(\theta_{\mathcal{S}}(s^{(l)};Y^*_{ij}))^{-1}-(\theta_{\mathcal{S}}(Y^*_{ij};Y^*_{ij}))^{-1}|+1\\
&\leq L\cdot R+1
\ea
for $n>m$. Hence,
\ba
|U_{ijl}|\leq&K^*\bigg(\frac{d_{\mathcal{S}}(s,Y^*_{ij})}{h_i^*}\bigg)\bigg|\frac{1}{\theta_{\mathcal{S}}(s;Y^*_{ij})}-\frac{1}{\theta_{\mathcal{S}}(s^{(l)};Y^*_{ij})}\bigg|\\
&+\frac{1}{\theta_{\mathcal{S}}(s^{(l)};Y^*_{ij})}\bigg|K^*\bigg(\frac{d_{\mathcal{S}}(s,Y^*_{ij})}{h_i^*}\bigg)-K^*\bigg(\frac{d_{\mathcal{S}}(s^{(l)},Y^*_{ij})}{h_i^*}\bigg)\bigg|\\
\leq &\const\cdot L\cdot n^{-\varsigma}+(L\cdot R+1)\cdot\const\cdot n^{-\varsigma}\cdot (h_i^*)^{-1}
\ea
for $n>m$, where the above $\const$ depends only on $K^*$. Thus,
\ba
|\tilde{g}_i(s)-\tilde{g}_i(s^{(l)})|\leq\const\cdot n^{-\varsigma}\cdot (h_i^*)^{-q-1}
\ea
for $n>m$, where the above $\const$ does not depend on $i, j, l$ and $s$.
Hence, the second and third terms on the right hand side of (\ref{imperfect density proof part2}) are negligible by taking sufficiently large $\varsigma>0$.
Define
\ba
Z_{ijl}=\frac{1}{n_i^*(h_i^*)^q}\bigg(K^*\bigg(\frac{d_{\mathcal{S}}(s,Y^*_{ij})}{h_i^*}\bigg)\frac{1}{\theta_{\mathcal{S}}(s;Y^*_{ij})}-\E\bigg(K^*\bigg(\frac{d_{\mathcal{S}}(s,Y^*_{ij})}{h_i^*}\bigg)\frac{1}{\theta_{\mathcal{S}}(s;Y^*_{ij})}\bigg)\bigg).
\ea
It holds that
\ba
\E(Z_{ijl})=0,\quad\|Z_{ijl}\|\leq\const(n_i^*(h_i^*)^q)^{-1},\quad\mbox{and}\quad\sum_{j=1}^{n_i^*}\E(Z_{ijl}^2)\leq\const(n_i^*(h_i^*)^q)^{-1},
\ea
where the $\const$ are independent of $i,j,l$ and $s$. By arguing as in the proof of Proposition \ref{imperfect density prop}, we get
\ba
\max_{1\leq i\leq n}\max_{1\leq l\leq N(n^{-\gamma})}|\tilde{g}_i(s^{(l)})-\E(\tilde{g}_i(s^{(l)}))|=O_p(H_n).
\ea
Hence, we obtain
\begin{align}\label{density second result}
\max_{1\leq i\leq n}\sup_{s\in\mathcal{S}}|\tilde{Y}^i(s)-Y^i(s)|=O_p(h^*_{\rm max}+H_n).
\end{align}
By arguing as in the proof of Proposition \ref{imperfect density prop},
\ba
\max_{1\leq i\leq n}\|\tilde{\mY}^i\om\mY^i\|\leq\frac{\nu(\mathcal{S})^{1/2}\max_{1\leq i\leq n}\sup_{s\in\mathcal{S}}|\tilde{Y}^i(s)-Y^i(s)|}{\min_{1\leq i\leq n}\inf_{s\in\mathcal{S}}Y^i(s)-\max_{1\leq i\leq n}\sup_{s\in\mathcal{S}}|\tilde{Y}^i(s)-Y^i(s)|}
\ea
with probability tending to one by (\ref{density second result}). This with (\ref{density second result}) gives the desired rate for $\max_{1\leq i\leq n}\|\tilde{\mY}^i\om\mY^i\|$. \qed

\subsection{Proof of Proposition \ref{existence}}\label{proof existence}

We apply Lemma \ref{general existence} for the proof. For this, we first verify condition \ref{F}. We note that $\sup_{u_j\in D_j}\int_{D_j}h_j^{-L_j}K_j(\|t_j-u_j\|_j/h_j)dt_j\leq\const$ by the proof of (\ref{nominator bound}). Hence, $\hat{p}^D_j(x_j)>0$ for each $x_j\in D_j$ under condition \ref{A} since $K_{h_j}(x_j,\tilde{\xi}_j^i)\I(\tilde{\mxi}^i\in D)>0$ if and only if $\tilde{\mxi}^i\in D$ and $\tilde{\xi}^i_j\in B_j(x_j,h_j)$. Also, for each $\tilde{\mxi}^i$, the function $K_{h_j}(\cdot,\tilde{\xi}_j^i)\I(\tilde{\mxi}^i\in D):D_j\rightarrow[0,\infty)$ is continuous since $K_j$ is continuous. Hence, $\hat{p}^D_j$ is continuous on the compact set $D_j$ and thus $\inf_{x_j\in D_j}\hat{p}^D_j(x_j)>0$. Also, $\hat{p}^D_{jk}(x_j,x_k)<\infty$ for each $(x_j,x_k)\in D_j\times D_k$ and $K_{h_j}(\cdot,\tilde{\xi}_j^i)K_{h_k}(\cdot,\tilde{\xi}_k^i)\I(\tilde{\mxi}^i\in D):D_j\times D_k\rightarrow[0,\infty)$ is continuous for each $\tilde{\mxi}^i$. Hence, $\sup_{(x_j,x_k)\in D_j\times D_k}\hat{p}^D_{jk}(x_j,x_k)<\infty$. This verifies condition \ref{F}. For the second part, we note that $\prod_{j=1}^dK_{h_j}(x_j,\tilde{\xi}_j^i)\I(\tilde{\mxi}^i\in D)>0$ if and only if $\tilde{\mxi}^i\in D\cap\prod_{j=1}^dB_j(x_j,h_j)$. Hence, $\hat{p}^D(\mx)>0$ for each $\mx\in D$. This completes the proof. \qed

\subsection{Proof of Proposition \ref{convergence1}}\label{proof convergence1}

The proposition follows from Lemma \ref{general algorithm convergence 1} and the proof of Proposition \ref{existence}. \qed

\subsection{Proof of Proposition \ref{convergence1 individual}}\label{proof convergence1 individual}

We note that $0<\hat{p}^D<\infty$ on $D$ by the proof of Proposition \ref{existence}. Also, the function $\prod_{j=1}^dK_{h_j}(\cdot,\tilde{\xi}_j^i)\I(\tilde{\mxi}^i\in D):D\rightarrow[0,\infty)$ is continuous for each $\tilde{\mxi}^i$. Hence, $\hat{p}^D$ is continuous on the compact set $D$ and thus $\hat{p}^D$ is bounded away from zero and infinity on $D$. This with Lemma \ref{general convergence1 individual} proves the proposition. \qed

\subsection{Proof of Theorem \ref{convergence2}}

We apply Lemma \ref{general algorithm convergence 2} for the proof. Hence, it suffices to verify condition \ref{G}.

We note that $\sup_{u_j\in D_j}\int_{D_j}h_j^{-L_j}K_j(\|t_j-u_j\|_j/h_j)dt_j\leq\const$ by the proof of (\ref{nominator bound}). This with (\ref{lower bound}) implies that $\inf_{x_j\in D_j}\int_{D_j}K_{h_j}(x_j,u_j)du_j\geq\const$. The latter with the condition that $p_j^D$ is bounded away from zero on $D_j$ and Lemma \ref{marginal density convergence 1} implies that $\hat{p}_j^D$ is bounded away from zero on $D_j$ with probability tending to one. This verifies condition (G1). Also, (\ref{nominator bound}) with the condition that $p_{jk}^D$ is bounded on $D_j\times D_k$ and Lemma \ref{marginal density convergence 2} implies that $\hat{p}_{jk}^D$ is bounded on $D_j\times D_k$ with probability tending to one. This verifies condition (G2).

We now verify condition (G3).
Since $\hat{p}_j^D$ is bounded away from zero on $D_j$ with probability tending to one and $\hat{p}_0^D>0$ with probability tending to one, it suffices to prove that
\ba
\lim_{n\rightarrow\infty}\Prob\left(\sup_{x_j\in D_j}\bigg\|n^{-1}\od\bigoplus_{i=1}^n(K_{h_j}(x_j,\tilde{\xi}_j^i)\I(\tilde{\mxi}^i\in D))\od\tilde{\mY}^i\bigg\|<\const\right)=1.
\ea
Note that
\ba
&\sup_{x_j\in D_j}\bigg\|n^{-1}\od\bigoplus_{i=1}^n(K_{h_j}(x_j,\tilde{\xi}_j^i)\I(\tilde{\mxi}^i\in D))\od\tilde{\mY}^i\bigg\|\\
&\leq\sup_{x_j\in D_j}\bigg\|n^{-1}\od\bigoplus_{i=1}^n((K_{h_j}(x_j,\tilde{\xi}_j^i)\I(\tilde{\mxi}^i\in D))\od\tilde{\mY}^i)\\
&\qquad\qquad\qquad\om n^{-1}\od\bigoplus_{i=1}^n((K_{h_j}(x_j,\xi_j^i)\I(\mxi^i\in D))\od\mY^i)\bigg\|\\
&\quad+\sup_{x_j\in D_j}\bigg\|n^{-1}\od\bigoplus_{i=1}^n(K_{h_j}(x_j,\xi_j^i)\I(\mxi^i\in D))\od\mY^i\bigg\|.
\ea
By Lemma \ref{marginal regression approximation}, the first term on the right hand side has the rate $o_p(1)$. Hence, it suffices to show that
\ba
\lim_{n\rightarrow\infty}\Prob\left(\sup_{x_j\in D_j}\bigg\|n^{-1}\od\bigoplus_{i=1}^n(K_{h_j}(x_j,\xi_j^i)\I(\mxi^i\in D))\od\mY^i\bigg\|<\const\right)=1.
\ea
Note that
\ba
&\sup_{x_j\in D_j}\bigg\|n^{-1}\od\bigoplus_{i=1}^n(K_{h_j}(x_j,\xi_j^i)\I(\mxi^i\in D))\od\mY^i\bigg\|\\
&\leq\sup_{x_j\in D_j}\bigg\|n^{-1}\od\bigoplus_{i=1}^n(K_{h_j}(x_j,\xi_j^i)\I(\mxi^i\in D))\od\mY^i-\E((K_{h_j}(x_j,\xi_j)\I(\mxi^i\in D))\od\mY)\bigg\|\\
&\quad+\sup_{x_j\in D_j}\|\E((K_{h_j}(x_j,\xi_j)\I(\mxi\in D))\od\mY)\|\\
&=O_p(n^{-1/2}\cdot(\log{n})^{1/2}\cdot h_j^{-L_j/2})+\sup_{x_j\in D_j}\E(K_{h_j}(x_j,\xi_j)\I(\mxi\in D)\E(\|\mY\||\mxi))\\
&\leq o_p(1)+\const,
\ea
where the first equality follows from Lemma \ref{uniform rates}. This verifies condition (G3).

We now verify condition (G4). Note that
\begin{align*}
&\int_{D_j}(\hat{p}_j^D(x_j)-p_j^D(x_j))^2dx_j\\
&\le 2\int_{D_j}\bigg(\hat{p}_j^D(x_j)-p_j^D(x_j)\int_{D_j}K_{h_j}(x_j,u_j)du_j\bigg)^2dx_j\\
&\quad + 2\int_{D_j}(p_j^D(x_j))^2\bigg(\int_{D_j}K_{h_j}(x_j,u_j)du_j-1\bigg)^2dx_j\\
&\le 2\Leb_j(D_j)\sup_{x_j\in D_j}\bigg(\hat{p}_j^D(x_j)-p_j^D(x_j)\int_{D_j}K_{h_j}(x_j,u_j)du_j\bigg)^2\\
&\quad + 2\Leb_j(D_j\setminus D^-_j(2h_j))\sup_{x_j\in D_j\setminus D^-_j(2h_j)}(p_j^D(x_j))^2\sup_{x_j\in D_j\setminus D^-_j(2h_j)}\bigg(\int_{D_j}K_{h_j}(x_j,u_j)du_j-1\bigg)^2\\
&=o_p(1),
\end{align*}
where the last equality follows from Lemma \ref{marginal density convergence 1} and the fact that $\Leb_j(D_j\setminus D^-_j(2h_j))\rightarrow0$. The latter fact is guaranteed by (B3).

Condition (G5) is similarly verified as in the verification of (G4) using Lemma \ref{marginal density convergence 2} instead of Lemma \ref{marginal density convergence 1}. \qed

\subsection{Proof of Theorem \ref{holder rate}}

We apply Lemma \ref{general rate}. Conditions (G1), (G2) and (G5) are verified in the proof of Theorem \ref{convergence2}. We approximate the terms that appear in Lemma \ref{general rate} to complete the proof.

We first approximate $\hat{\mf}_0\om\mf_0$. Note that
\begin{align*}
\hat{\mf}_0\om\mf_0&=(\hat{p}_0^D\cdot n)^{-1}\od\bigoplus_{i=1}^n(\tilde{\mY}^i\od\I(\tilde{\mxi}^i\in D))\om(p^D_0)^{-1}\od\E(\mY\od\I(\mxi\in D))\\
&=(\hat{p}_0^D\cdot p_0^D)^{-1}\bigg(\bigg(n^{-1}\od\bigoplus_{i=1}^n(\tilde{\mY}^i\od\I(\tilde{\mxi}^i\in D))\om\E(\mY\od\I(\mxi\in D))\bigg)\od p^D_0\\
&\qquad\qquad\qquad\quad\op(p^D_0-\hat{p}^D_0)\od\E(\mY\od\I(\mxi\in D))\bigg).
\end{align*}
Hence, it suffices to approximate $n^{-1}\od\bigoplus_{i=1}^n(\tilde{\mY}^i\od\I(\tilde{\mxi}^i\in D))\om\E(\mY\od\I(\mxi\in D))$ and $p^D_0-\hat{p}^D_0$. Note that
\begin{align*}
&n^{-1}\od\bigoplus_{i=1}^n(\tilde{\mY}^i\od\I(\tilde{\mxi}^i\in D))\om\E(\mY\od\I(\mxi\in D))\\
&=n^{-1}\od\bigoplus_{i=1}^n((\tilde{\mY}^i\om\mY^i)\od\I(\mxi^i\in D))\op n^{-1}\od\bigoplus_{i=1}^n((\tilde{\mY}^i\om\mY^i)\od(\I(\tilde{\mxi}^i\in D)-\I(\mxi^i\in D)))\\
&\quad\op n^{-1}\od\bigoplus_{i=1}^n(\mY^i\od(\I(\tilde{\mxi}^i\in D)-\I(\mxi^i\in D)))\\
&\quad\op n^{-1}\od\bigoplus_{i=1}^n(\mY^i\od\I(\mxi^i\in D))\om\E(\mY\od\I(\xi\in D))\\
&=O_p\bigg(b_n+\sum_{j=1}^da_{nj}+n^{-1/2}\bigg).
\end{align*}
Also, $p^D_0-\hat{p}^D_0=O_p(\sum_{j=1}^da_{nj}+n^{-1/2})$ by Lemma \ref{p0 rate}. Hence,
\ba
\hat{\mf}_0\om\mf_0=O_p\bigg(b_n+\sum_{j=1}^da_{nj}+n^{-1/2}\bigg).
\ea

We now approximate $n^{-1}\od\bigoplus_{i=1}^n(\mf_j(\xi^i_j)\od\I(\tilde{\mxi}^i\in D))$. Note that
\ba
&\bigg\|n^{-1}\od\bigoplus_{i=1}^n(\mf_j(\xi^i_j)\od\I(\tilde{\mxi}^i\in D))\bigg\|\\
&=\bigg\|n^{-1}\od\bigoplus_{i=1}^n(\mf_j(\xi^i_j)\od(\I(\tilde{\mxi}^i\in D)-\I(\mxi^i\in D)))\op n^{-1}\od\bigoplus_{i=1}^n(\mf_j(\xi^i_j)\od\I(\mxi^i\in D))\bigg\|\\
&\leq n^{-1}\sum_{i=1}^n\|\mf_j(\xi^i_j)\|\I(\mxi^i\in D(n))\op O_p(n^{-1/2})\\
&=O_p\bigg(\sum_{j=1}^da_{nj}+n^{-1/2}\bigg),
\ea
where the inequality follows from the constraints (\ref{constraint 1}) and the last equality follows from the boundedness of $\mf_j$ on $D^+_j(\varepsilon)$ and the proof of Lemma \ref{p0 rate}.

We now approximate $\hat{\mdelta}_j^A$. We define
\ba
\check{\mdelta}^{A}_j(x_j)&=n^{-1}\odot\bigoplus_{i=1}^n((K_{h_j}(x_j,\xi^i_j)\I(\mxi^i\in D))\odot\mepsilon_i)
\ea
Note that
\begin{align}\label{decomposition A}
\begin{split}
&\sup_{x_j\in D_j}\|\hat{\mdelta}_j^A(x_j)\om\check{\mdelta}_j^A(x_j)\|\\
&\leq\sup_{x_j\in D_j}\bigg\|n^{-1}\od\bigoplus_{i=1}^n((K_{h_j}(x_j,\xi_j^i)\I(\mxi^i\in D)(\I(\tilde{\mxi}^i\in D)-\I(\mxi^i\in D)))\od\mepsilon^i)\bigg\|\\
&\quad+\sup_{x_j\in D_j}\bigg\|n^{-1}\od\bigoplus_{i=1}^n((K_{h_j}(x_j,\tilde{\xi}_j^i)\I(\tilde{\mxi}^i\in D)(\I(\tilde{\mxi}^i\in D)-\I(\mxi^i\in D)))\od\mepsilon^i)\bigg\|\\
&\quad+\sup_{x_j\in D_j}\bigg\|n^{-1}\od\bigoplus_{i=1}^n(((K_{h_j}(x_j,\tilde{\xi}_j^i)-K_{h_j}(x_j,\xi_j^i))\I(\mxi^i\in D)\I(\tilde{\mxi}^i\in D))\od\mepsilon^i)\bigg\|.
\end{split}
\end{align}
Define $V_{nj}=h_j^{-L_j}\cdot a_{nj}+\sum_{k\neq j}a_{nk}$. The first term on the right hand side of (\ref{decomposition A}) has the rate
\ba
O_p(n^{-1/2}\cdot(\log{n})^{1/2}\cdot V_{nj}^{1/2}\cdot h_j^{-L_j/2}),
\ea
provided that $n^{-1+2\beta+2/\alpha}\cdot h_j^{-L_j}=o(1)$ and $n^{-\beta}\cdot(\log{n})^{1/2}\cdot(\sum_{j=1}^da_{nj})^{1/\alpha}\cdot V_{nj}^{-1/2}=O(1)$ by Lemma \ref{uniform rates Dn2}.  Also, the second term on the right hand side of (\ref{decomposition A}) has the same rate. To see this, for sufficiently large $\kappa>0$, let $\{x_j^{(l)}:1\leq l\leq N(n^{-\kappa})\}\subset D_j$ be a set of points such that $N(n^{-\kappa})=O(n^{\kappa L_j})$ and $\{B_j(x_j^{(l)},n^{-\kappa}):1\leq l\leq N(n^{-\kappa})\}$ covers $D_j$. Define $R_n=\sum_{j=1}^da_{nj}$, $\mdelta_n^i=(\I(\tilde{\mxi}^i\in D)-\I(\mxi^i\in D))\od\mepsilon^i$ and
\ba
\mS_n(x_j)=n^{-1}\od\bigoplus_{i=1}^n((K_{h_j}(x_j,\tilde{\xi}_j^i)\I(\tilde{\mxi}^i\in D))\od\mdelta_n^i).
\ea
We also write $(\mxi^{\rm all},\tilde{\mxi}^{\rm all})=(\mxi^1,\ldots,\mxi^n,\tilde{\mxi}^1,\ldots,\tilde{\mxi}^n)$. By arguing as in the proof of Lemma \ref{uniform rates general}, we may show that
\begin{align*}
\sup_{x_j\in D_j}\|\mS_n(x_j)\om \E(\mS_n(x_j)|\mxi^{\rm all},\tilde{\mxi}^{\rm all})\|\leq&\max_{1\leq l\leq N(n^{-\kappa})}\|\mU_n(x_j^{(l)})\om \E(\mU_n(x_j^{(l)})|\mxi^{\rm all},\tilde{\mxi}^{\rm all})\|\\
&+o_p(n^{-1/2}\cdot(\log{n})^{1/2}\cdot V_{nj}^{1/2}\cdot h_j^{-L_j/2}),
\end{align*}
provided that $n^{-1+2\beta+2/\alpha}\cdot h_j^{-L_j}=o(1)$, where
\begin{align*}
\mU_n(x_j)=n^{-1}\od\bigoplus_{i=1}^n((K_{h_j}(x_j,\tilde{\xi}_j^i)\cdot\I(\tilde{\mxi}^i\in D)\cdot\I(\|\mdelta_n^i\|\leq n^{1/2-\beta}\cdot R_n^{1/\alpha}\cdot h_j^{L_j/2}))\od\mdelta_n^i).
\end{align*}
Define
\ba
\mZ_n^i(x_j)&=n^{-1}\od((K_{h_j}(x_j,\tilde{\xi}_j^i)\cdot\I(\tilde{\mxi}^i\in D)\cdot\I(\|\mdelta_n^i\|\leq n^{1/2-\beta}\cdot R_n^{1/\alpha}\cdot h_j^{L_j/2}))\od\mdelta_n^i)\\
&~~\om\E(n^{-1}\od((K_{h_j}(x_j,\tilde{\xi}_j^i)\cdot\I(\tilde{\mxi}^i\in D)\cdot\I(\|\mdelta_n^i\|\leq n^{1/2-\beta}\cdot R_n^{1/\alpha}\cdot h_j^{L_j/2}))\od\mdelta_n^i)|\mxi^{\rm all},\tilde{\mxi}^{\rm all}).
\ea
We note that $\E(\mZ_n^i(x_j)|\mxi^{\rm all},\tilde{\mxi}^{\rm all})=\bzero$ and $\|\mZ_n^i(x_j)\|\leq\const\cdot n^{-1/2-\beta}\cdot R_n^{1/\alpha}\cdot h_j^{-L_j/2}$. It holds that
\ba
&\E(\|\mZ_n^i(x_j)\|^2|\mxi^{\rm all},\tilde{\mxi}^{\rm all})\\
&\leq\const\cdot n^{-2}\cdot\bigg(h_j^{-L_j}K_j\bigg(\frac{\|x_j-\tilde{\xi}^i_j\|_j}{h_j}\bigg)\bigg)^2\cdot\E(\|\mdelta_n^i\|^2|\mxi^{\rm all},\tilde{\mxi}^{\rm all})\\
&\leq\const\cdot n^{-2}\cdot\bigg(h_j^{-L_j}K_j\bigg(\frac{\|x_j-\tilde{\xi}^i_j\|_j}{h_j}\bigg)-h_j^{-L_j}K_j\bigg(\frac{\|x_j-\xi^i_j\|_j}{h_j}\bigg)\bigg)^2\cdot\E(\|\mdelta_n^i\|^2|\mxi^{\rm all},\tilde{\mxi}^{\rm all})\\
&\quad+\const\cdot n^{-2}\cdot\bigg(h_j^{-L_j}K_j\bigg(\frac{\|x_j-\xi^i_j\|_j}{h_j}\bigg)\bigg)^2\cdot\E(\|\mdelta_n^i\|^2|\mxi^{\rm all},\tilde{\mxi}^{\rm all})\\
&\leq\const\cdot n^{-2}\cdot(h_j^{-L_j-1}a_{nj})^2\cdot\I(\mxi^i\in D(n))\cdot\E(\|\mepsilon^i\|^2|\mxi^{\rm all},\tilde{\mxi}^{\rm all})\\
&\quad+\const\cdot n^{-2}\cdot\bigg(h_j^{-L_j}K_j\bigg(\frac{\|x_j-\xi^i_j\|_j}{h_j}\bigg)\bigg)^2\cdot\I(\mxi^i\in D(n))\cdot\E(\|\mepsilon^i\|^2|\mxi^{\rm all},\tilde{\mxi}^{\rm all})\\
&\leq\const\cdot n^{-2}\cdot\bigg((h_j^{-L_j-1}a_{nj})^2+\bigg(h_j^{-L_j}K_j\bigg(\frac{\|x_j-\xi^i_j\|_j}{h_j}\bigg)\bigg)^2\bigg)\cdot\I(\mxi^i\in D(n))
\ea
almost surely. By arguing as in the proof of Lemma \ref{uniform rates Dn2}, we have
\begin{align*}
&\Prob\left(\bigg\|\bigoplus_{i=1}^n\mZ_n^i(\mx_J)\bigg\|>C_0\cdot n^{-1/2}\cdot(\log{n})^{1/2}\cdot V_{nj}^{1/2}\cdot h_j^{-L_j/2}\bigg|\mxi^{\rm all},\tilde{\mxi}^{\rm all}\right)\\
&\leq2n^{-C_0}\cdot\prod_{i=1}^n\big(1+\lambda_n^2\cdot\E(\|\mZ_n^i(\mx_J)\|^2|\mxi^{\rm all},\tilde{\mxi}^{\rm all})\cdot\exp\big(\const\cdot \lambda_n\cdot n^{-1/2-\beta}\cdot R_n^{1/\alpha}\cdot h_j^{-L_j/2}\big)\big)\\
&\leq2n^{-C_0}\cdot\prod_{i=1}^n\bigg(1+O(1)\cdot\lambda_n^2\cdot n^{-2}\cdot\bigg(h_j^{-2L_j-2}a^2_{nj}+h_j^{-2L_j}K^2_j\bigg(\frac{\|x_j-\xi^i_j\|_j}{h_j}\bigg)\bigg)\cdot\I(\mxi^i\in D(n))\bigg)
\end{align*}
almost surely, where $C_0>0$ is a constant, $\lambda_n=n^{1/2}\cdot(\log{n})^{1/2}\cdot V_{nj}^{-1/2}\cdot h_j^{L_j/2}$ and the last inequality follows from that $n^{-\beta}\cdot(\log{n})^{1/2}\cdot(\sum_{j=1}^da_{nj})^{1/\alpha}\cdot V_{nj}^{-1/2}=O(1)$.
Hence,
\begin{align*}
&\Prob\left(\bigg\|\bigoplus_{i=1}^n\mZ_n^i(\mx_J)\bigg\|>C_0\cdot n^{-1/2}\cdot(\log{n})^{1/2}\cdot V_{nj}^{1/2}\cdot h_j^{-L_j/2}\right)\\
&\leq2n^{-C_0}\\
&\quad\cdot\bigg(\E\bigg(1+O(1)\cdot\lambda_n^2\cdot n^{-2}\cdot\bigg(h_j^{-2L_j-2}a^2_{nj}+h_j^{-2L_j}K^2_j\bigg(\frac{\|x_j-\xi_j\|_j}{h_j}\bigg)\bigg)\cdot\I(\mxi\in D(n))\bigg)\bigg)^n.
\end{align*}
Since
\ba
&\E\bigg(\bigg(h_j^{-2L_j-2}a^2_{nj}+h_j^{-2L_j}K^2_j\bigg(\frac{\|x_j-\xi_j\|_j}{h_j}\bigg)\bigg)\cdot\I(\mxi\in D(n))\bigg)\\
&\leq h_j^{-2L_j-2}a^2_{nj}\cdot\Prob(\mxi\in D(n))+\int_{D(n)}h_j^{-2L_j}K^2_j\bigg(\frac{\|x_j-u_j\|_j}{h_j}\bigg)p(\bu)d\bu\\
&\leq\const\cdot h_j^{-2L_j-2}a^2_{nj}\cdot R_n+\const\cdot h_j^{-L_j}\cdot V_{nj}\\
&=O(h_j^{-L_j}\cdot V_{nj}),
\ea
we have
\ba
&\Prob\left(\bigg\|\bigoplus_{i=1}^n\mZ_n^i(\mx_J)\bigg\|>C\cdot n^{-1/2}\cdot(\log{n})^{1/2}\cdot V_{nj}^{1/2}\cdot h_j^{-L_j/2}\right)\\
&\leq2n^{-C_0}\cdot(1+O(1)\cdot\lambda_n^2\cdot n^{-2}\cdot O(h_j^{-L_j}\cdot V_{nj}))^n\\
&\leq2n^{-C_0}\cdot n^{O(1)}.
\ea
This with the fact that $N(n^{-\kappa})=O(n^{\kappa L_j})$ and $\E(\mepsilon^i|\mxi^{\rm all},\tilde{\mxi}^{\rm all})=\bzero$ for all $1\leq i\leq n$ gives the desired rate for the second term on the right hand side of (\ref{decomposition A}). We now prove that the third term on the right hand side of (\ref{decomposition A}) has the rate
\ba
O_p(n^{-1/2}\cdot(\log{n})^{1/2}\cdot h_j^{-L_j/2}\cdot h_j^{-1}\cdot a_{nj}).
\ea
To see this, note that the third term equals
\begin{align}\label{decomposition A2}
\begin{split}
&\sup_{x_j\in D_j}\bigg\|\frac{1}{n}\od\bigoplus_{i=1}^n\bigg(\bigg(\sum_{l=1}^{L_j}\frac{\partial}{\partial u_{jl}}K_{h_j}(x_j,u_j)\bigg|_{u_j=\bar{\xi}_j^i}(\tilde{\xi}_{jl}^i-\xi_{jl}^i)\I(\mxi^i\in D)\I(\tilde{\mxi}^i\in D)\bigg)\od\mepsilon^i\bigg)\bigg\|\\
&\leq\sum_{l=1}^{L_j}\sup_{x_j\in D_j}\bigg\|\frac{1}{n}\od\bigoplus_{i=1}^n\bigg(\bigg(\frac{\partial}{\partial u_{jl}}K_{h_j}(x_j,u_j)\bigg|_{u_j=\bar{\xi}_j^i}(\tilde{\xi}_{jl}^i-\xi_{jl}^i)\I(\mxi^i\in D)\I(\tilde{\mxi}^i\in D)\bigg)\od\mepsilon^i\bigg)\bigg\|\\
&\leq\sum_{l=1}^{L_j}\sup_{x_j\in D_j}\bigg\|\frac{1}{n}\od\bigoplus_{i=1}^n\bigg(\bigg(\bigg(\int_{D_j}\frac{1}{h_j^{L_j}}K_j\bigg(\frac{\|t_j-\bar{\xi}^i_j\|_j}{h_j}\bigg)dt_j\bigg)^{-1}\frac{1}{h_j^{L_j}}\frac{\partial}{\partial u_{jl}}K_j\bigg(\frac{\|x_j-u_j\|_j}{h_j}\bigg)\bigg|_{u_j=\bar{\xi}_j^i}\\
&\qquad\qquad\qquad\qquad\qquad\qquad\cdot(\tilde{\xi}_{jl}^i-\xi_{jl}^i)\I(\mxi^i\in D)\I(\tilde{\mxi}^i\in D)\bigg)\od\mepsilon^i\bigg)\bigg\|\\
&\quad+\sum_{l=1}^{L_j}\sup_{x_j\in D_j}\bigg\|\frac{1}{n}\od\bigoplus_{i=1}^n\bigg(\bigg(\bigg(\int_{D_j}\frac{1}{h_j^{L_j}}K_j\bigg(\frac{\|t_j-\bar{\xi}^i_j\|_j}{h_j}\bigg)dt_j\bigg)^{-2}\frac{1}{h_j^{L_j}}K_j\bigg(\frac{\|x_j-\bar{\xi}_j^i\|_j}{h_j}\bigg)\\
&\qquad\qquad\qquad\cdot\int_{D_j}\frac{1}{h_j^{L_j}}\frac{\partial}{\partial u_{jl}}K_j\bigg(\frac{\|t_j-u_j\|_j}{h_j}\bigg)\bigg|_{u_j=\bar{\xi}_j^i}dt_j(\tilde{\xi}_{jl}^i-\xi_{jl}^i)\I(\mxi^i\in D)\I(\tilde{\mxi}^i\in D)\bigg)\od\mepsilon^i\bigg)\bigg\|,
\end{split}
\end{align}
where $\bar{\xi}_j^i$ is a random vector lying on the line connecting $\tilde{\xi}_j^i$ and $\xi_j^i$. We first approximate the first term on the right hand side of (\ref{decomposition A2}). We define
\ba
\mS_n(x_j)=&n^{-1}\od\bigoplus_{i=1}^n\bigg(\bigg(\bigg(\int_{D_j}\frac{1}{h_j^{L_j}}K_j\bigg(\frac{\|t_j-\bar{\xi}^i_j\|_j}{h_j}\bigg)dt_j\bigg)^{-1}\frac{1}{h_j^{L_j}}\frac{\partial}{\partial u_{jl}}K_j\bigg(\frac{\|x_j-u_j\|_j}{h_j}\bigg)\bigg|_{u_j=\bar{\xi}_j^i}\\
&\qquad\qquad\qquad\qquad\cdot(\tilde{\xi}_{jl}^i-\xi_{jl}^i)\I(\mxi^i\in D)\I(\tilde{\mxi}^i\in D)\bigg)\od\mepsilon^i\bigg).
\ea
By arguing as in the proof of Lemma \ref{uniform rates general}, we may show that
\begin{align*}
\sup_{x_j\in D_j}\|\mS_n(x_j)\om \E(\mS_n(x_j)|\mxi^{\rm all},\tilde{\mxi}^{\rm all})\|\leq&\max_{1\leq l\leq N(n^{-\kappa})}\|\mU_n(x_j^{(l)})\om \E(\mU_n(x_j^{(l)})|\mxi^{\rm all},\tilde{\mxi}^{\rm all})\|\\
&+o_p(n^{-1/2}\cdot(\log{n})^{1/2}\cdot h_j^{-L_j/2}\cdot h_j^{-1}\cdot a_{nj})
\end{align*}
provided that $n^{-1+2\beta+2/\alpha}\cdot h_j^{-L_j}=o(1)$, where
\begin{align*}
\mU_n(x_j)=&n^{-1}\od\bigoplus_{i=1}^n\bigg(\bigg(\bigg(\int_{D_j}\frac{1}{h_j^{L_j}}K_j\bigg(\frac{\|t_j-\bar{\xi}^i_j\|_j}{h_j}\bigg)dt_j\bigg)^{-1}\frac{1}{h_j^{L_j}}\frac{\partial}{\partial u_{jl}}K_j\bigg(\frac{\|x_j-u_j\|_j}{h_j}\bigg)\bigg|_{u_j=\bar{\xi}_j^i}\\
&\qquad\qquad\qquad\cdot(\tilde{\xi}_{jl}^i-\xi_{jl}^i)\I(\mxi^i\in D)\I(\tilde{\mxi}^i\in D)\I(\|\mepsilon^i\|\leq n^{1/2-\beta}h_j^{L_j/2})\bigg)\od\mepsilon^i\bigg).
\end{align*}
Define
\ba
\mZ_n^i(x_j)=&n^{-1}\od\bigg(\bigg(\bigg(\int_{D_j}\frac{1}{h_j^{L_j}}K_j\bigg(\frac{\|t_j-\bar{\xi}^i_j\|_j}{h_j}\bigg)dt_j\bigg)^{-1}\frac{1}{h_j^{L_j}}\frac{\partial}{\partial u_{jl}}K_j\bigg(\frac{\|x_j-u_j\|_j}{h_j}\bigg)\bigg|_{u_j=\bar{\xi}_j^i}\\
&\qquad\qquad\qquad\cdot(\tilde{\xi}_{jl}^i-\xi_{jl}^i)\I(\mxi^i\in D)\I(\tilde{\mxi}^i\in D)\I(\|\mepsilon^i\|\leq n^{1/2-\beta}h_j^{L_j/2})\bigg)\od\mepsilon^i\bigg)\\
&\om\E\bigg(n^{-1}\od\bigg(\bigg(\bigg(\int_{D_j}\frac{1}{h_j^{L_j}}K_j\bigg(\frac{\|t_j-\bar{\xi}^i_j\|_j}{h_j}\bigg)dt_j\bigg)^{-1}\frac{1}{h_j^{L_j}}\frac{\partial}{\partial u_{jl}}K_j\bigg(\frac{\|x_j-u_j\|_j}{h_j}\bigg)\bigg|_{u_j=\bar{\xi}_j^i}\\
&\qquad\qquad\cdot(\tilde{\xi}_{jl}^i-\xi_{jl}^i)\I(\mxi^i\in D)\I(\tilde{\mxi}^i\in D)\I(\|\mepsilon^i\|\leq n^{1/2-\beta}h_j^{L_j/2})\bigg)\od\mepsilon^i\bigg)\bigg|\mxi^{\rm all},\tilde{\mxi}^{\rm all}\bigg).
\ea
Note that $\E(\mZ_n^i(x_j)|\mxi^{\rm all},\tilde{\mxi}^{\rm all})=\bzero$ and $\|\mZ_n^i(x_j)\|\leq\const\cdot n^{-1/2-\beta}\cdot h_j^{-L_j/2-1}\cdot a_{nj}$. It holds that
\ba
&\E(\|\mZ_n^i(x_j)\|^2|\mxi^{\rm all},\tilde{\mxi}^{\rm all})\\
&\leq\const\cdot n^{-2}\cdot a_{nj}^2\cdot\I(\mxi^i\in D)\cdot\bigg(\frac{1}{h_j^{L_j}}\frac{\partial}{\partial u_{jl}}K_j\bigg(\frac{\|x_j-u_j\|_j}{h_j}\bigg)\bigg|_{u_j=\bar{\xi}^i_j}\bigg)^2\cdot\E(\|\mepsilon^i\|^2|\mxi^{\rm all},\tilde{\mxi}^{\rm all})\\
&\leq\const\cdot n^{-2}\cdot a_{nj}^2\cdot\bigg(\frac{1}{h_j^{L_j}}\frac{\partial}{\partial u_{jl}}K_j\bigg(\frac{\|x_j-u_j\|_j}{h_j}\bigg)\bigg|_{u_j=\bar{\xi}^i_j}-\frac{1}{h_j^{L_j}}\frac{\partial}{\partial u_{jl}}K_j\bigg(\frac{\|x_j-u_j\|_j}{h_j}\bigg)\bigg|_{u_j=\xi^i_j}\bigg)^2\\
&\qquad\qquad\qquad\qquad\qquad\qquad\qquad\quad\cdot\E(\|\mepsilon^i\|^2|\mxi^{\rm all},\tilde{\mxi}^{\rm all})\\
&\quad+\const\cdot n^{-2}\cdot a_{nj}^2\cdot\I(\mxi^i\in D)\cdot\bigg(\frac{1}{h_j^{L_j}}\frac{\partial}{\partial u_{jl}}K_j\bigg(\frac{\|x_j-u_j\|_j}{h_j}\bigg)\bigg|_{u_j=\xi^i_j}\bigg)^2\cdot\E(\|\mepsilon^i\|^2|\mxi^{\rm all},\tilde{\mxi}^{\rm all})\\
&\leq\const\cdot n^{-2}\cdot(h_j^{-L_j-2}\cdot a_{nj}^2)^2\\
&\quad+\const\cdot n^{-2}\cdot a_{nj}^2\cdot\I(\mxi^i\in D)\cdot\bigg(\frac{1}{h_j^{L_j}}\frac{\partial}{\partial u_{jl}}K_j\bigg(\frac{\|x_j-u_j\|_j}{h_j}\bigg)\bigg|_{u_j=\xi^i_j}\bigg)^2
\ea
almost surely. Using this and by arguing as in the proof for the second term on the right hand side of (\ref{decomposition A}) with the facts that
\ba
&\E\bigg(\I(\mxi^i\in D)\cdot\bigg(\frac{1}{h_j^{L_j}}\frac{\partial}{\partial u_{jl}}K_j\bigg(\frac{\|x_j-u_j\|_j}{h_j}\bigg)\bigg|_{u_j=\xi^i_j}\bigg)^2\bigg)\\
&=\int_D\bigg(\frac{1}{h_j^{L_j}}\frac{\partial}{\partial u_{jl}}K_j\bigg(\frac{\|x_j-u_j\|_j}{h_j}\bigg)\bigg)^2p(\bu)d\bu\\
&=O(h_j^{-L_j-2})
\ea
and $n^{-\beta}\cdot(\log{n})^{1/2}=O(1)$, we get that the first term on the right hand side of (\ref{decomposition A2}) has the rate
\ba
O_p(n^{-1/2}\cdot(\log{n})^{1/2}\cdot h_j^{-L_j/2}\cdot h_j^{-1}\cdot a_{nj}).
\ea
We now approximate the second term on the right hand side of (\ref{decomposition A2}). Since
\ba
\bigg|\int_{D_j}\frac{1}{h_j^{L_j}}\frac{\partial}{\partial u_{jl}}K_j\bigg(\frac{\|t_j-u_j\|_j}{h_j}\bigg)\bigg|_{u_j=\bar{\xi}_j^i}dt_j\I(\mxi^i\in D)\I(\tilde{\mxi}^i\in D)\bigg|\leq\const\cdot h_j^{-1}
\ea
and
\ba
\bigg|\frac{1}{h_j^{L_j}}K_j\bigg(\frac{\|x_j-\bar{\xi}_j^i\|_j}{h_j}\bigg)-\frac{1}{h_j^{L_j}}K_j\bigg(\frac{\|x_j-\xi_j^i\|_j}{h_j}\bigg)\bigg|\leq\const\cdot h_j^{-L_j-1}\cdot a_{nj},
\ea
by arguing as above, we get the same rate for second term on the right hand side of (\ref{decomposition A2}). This gives the desired rate for the third term on the right hand side of (\ref{decomposition A}). Hence, we have
\begin{align}\label{A error}
\begin{split}
\sup_{x_j\in D_j}\|\hat{\mdelta}_j^A(x_j)\om\check{\mdelta}_j^A(x_j)\|&=O_p(n^{-1/2}\cdot(\log{n})^{1/2}\cdot h_j^{-L_j/2}\cdot(V_{nj}^{1/2}+h_j^{-1}\cdot a_{nj}))\\
&=O_p(n^{-1/2}\cdot(\log{n})^{1/2}\cdot h_j^{-L_j/2}\cdot V_{nj}^{1/2})
\end{split}
\end{align}
provided that $n^{-1+2\beta+2/\alpha}\cdot h_j^{-L_j}=o(1)$ and $n^{-\beta}\cdot(\log{n})^{1/2}\cdot(\sum_{j=1}^da_{nj})^{1/\alpha}\cdot V_{nj}^{-1/2}=O(1)$, where the last equality at (\ref{A error}) follows from the fact that $h_j^{-1}\cdot a_{nj}=O((h_j^{-L_j}\cdot a_{nj})^{1/2})$. Note that $\E(\check{\mdelta}_j^A(x_j))=\bzero$ and
\begin{align*}
\Var(\check{\mdelta}_j^A(x_j))=&n^{-1}\cdot\E(\|(K_{h_j}(x_j,\xi_j)\I(\mxi\in D))\od\mepsilon\|^2)\\
\leq&\const\cdot n^{-1}\cdot h_j^{-L_j}\cdot\int_D\E(\|\mepsilon\|^2|\mxi=\bu)\cdot p(\bu)\cdot h_j^{-L_j}K^2_j\bigg(\frac{\|x_j-u_j\|_j}{h_j}\bigg)d\bu\\
=&O(n^{-1}\cdot h_j^{-L_j})
\end{align*}
for any $x_j\in D_j$. Hence,
\begin{align*}
\check{\mdelta}_j^A(x_j)=O_p(n^{-1/2}\cdot h_j^{-L_j/2})
\end{align*}
for any $x_j\in D_j$ and thus,
\begin{align*}
\hat{\mdelta}_j^A(x_j)=O_p(n^{-1/2}\cdot h_j^{-L_j/2})
\end{align*}
for any $x_j\in D_j$ by (\ref{A error}). We also note that
\begin{align*}
&\int_{D_j}\E(\|\check{\mdelta}_j^A(x_j)\|^2)dx_j
\leq\const\cdot n^{-1}\cdot h_j^{-L_j}.
\end{align*}
Hence,
\ba
\bigg(\int_{D_j}\|\check{\mdelta}_j^A(x_j)\|^2dx_j\bigg)^{1/2}=O_p(n^{-1/2}\cdot h_j^{-L_j/2})
\ea
and thus,
\begin{align*}
\bigg(\int_{D_j}\|\hat{\mdelta}_j^A(x_j)\|^2dx_j\bigg)^{1/2}=O_p(n^{-1/2}\cdot h_j^{-L_j/2})
\end{align*}
by (\ref{A error}).
Also, Lemma \ref{uniform rates} implies that
\begin{align*}
\sup_{x_j\in D_j}\|\check{\mdelta}_j^A(x_j)\|=O_p(n^{-1/2}\cdot(\log{n})^{1/2}\cdot h_j^{-L_j/2}).
\end{align*}
Hence, we have
\begin{align*}
\sup_{x_j\in D_j}\|\hat{\mdelta}_j^A(x_j)\|=O_p(n^{-1/2}\cdot(\log{n})^{1/2}\cdot h_j^{-L_j/2})
\end{align*}
by (\ref{A error}).

We now approximate $\hat{\mdelta}_j^B$. We note that
\ba
&\|\hat{\mdelta}_j^B(x_j)\|\\
&\leq\const\frac{1}{n}\sum_{i=1}^n\frac{1}{h_j^{L_j}}K_j\bigg(\frac{\|x_j-\tilde{\xi}_j^i\|_j}{h_j}\bigg)\I(\tilde{\mxi}^i\in D)\|\mf_j(\xi_j^i)\om\mf_j(x_j)\|\\
&\leq\const\frac{1}{n}\sum_{i=1}^n\frac{1}{h_j^{L_j}}K_j\bigg(\frac{\|x_j-\tilde{\xi}_j^i\|_j}{h_j}\bigg)(\I(\mxi^i\in D(n))+\I(\mxi^i\in D))(\|\xi^i_j-\tilde{\xi}^i_j\|^{\nu_j}_j+\|\tilde{\xi}^i_j-x_j\|^{\nu_j}_j)\\
&\leq\const(a_{nj}^{\nu_j}+h_j^{\nu_j})\frac{1}{n}\sum_{i=1}^n\frac{1}{h_j^{L_j}}K_j\bigg(\frac{\|x_j-\tilde{\xi}_j^i\|_j}{h_j}\bigg)(\I(\mxi^i\in D(n))+\I(\mxi^i\in D))\\
&=O_p(h_j^{\nu_j})
\ea
uniformly for $x_j\in D_j$, where the equality follows by arguing as in the proof for the first part of Lemma \ref{marginal density approximation 1}. Hence,
\ba
&\hat{\mdelta}^{B}_j(x_j)=O_p(h_j^{\nu_j})\text{~for all~}x_j\in D_j,\quad\bigg(\int_{D_j}\|\hat{\mdelta}^{B}_j(x_j)\|^2dx_j\bigg)^{1/2}=O_p(h_j^{\nu_j}),\\
&\sup_{x_j\in D_j}\|\hat{\mdelta}^{B}_j(x_j)\|=O_p(h_j^{\nu_j}),\quad\sum_{k\neq j}\sup_{x_j\in D_j}\bigg\|\int_{D_k}\hat{\mdelta}^{B}_k(x_k)\odot
\frac{\hat{p}^D_{jk}(x_j,x_k)}{\hat{p}^D_k(x_k)}dx_k\bigg\|=O_p\bigg(\sum_{k\neq j}h_k^{\nu_k}\bigg).
\ea

We now approximate $\hat{\mdelta}_j^C$. Note that
\ba
&\|\hat{\mdelta}_j^C(x_j)\|\\
&\leq\const\sum_{k\neq j}\frac{1}{n}\sum_{i=1}^n\frac{1}{h_j^{L_j}}K_j\bigg(\frac{\|x_j-\tilde{\xi}_j^i\|_j}{h_j}\bigg)
\int_{D_k}\|\mf_k(\xi^i_k)\om\mf_k(x_k)\|\frac{1}{h_k^{L_k}}K_k\bigg(\frac{\|x_k-\tilde{\xi}_k^i\|_k}{h_k}\bigg)dx_k\\
&\qquad\qquad\qquad\qquad\qquad\cdot\I(\tilde{\mxi^i}\in D)\\
&\leq\const\sum_{k\neq j}\frac{1}{n}\sum_{i=1}^nK_{h_j}(x_j,\tilde{\xi}^i_j)
\int_{D_k}(\|\xi^i_k-\tilde{\xi}^i_k\|^{\nu_k}_k+\|\tilde{\xi}^i_k-x_k\|^{\nu_k}_k)\frac{1}{h_k^{L_k}}K_k\bigg(\frac{\|x_k-\tilde{\xi}_k^i\|_k}{h_k}\bigg)dx_k\\
&\qquad\qquad\qquad\qquad\qquad\cdot(\I(\mxi^i\in D(n))+\I(\mxi^i\in D))\\
&\leq\const\sum_{k\neq j}(a_{nk}^{\nu_k}+h_k^{\nu_k})\frac{1}{n}\sum_{i=1}^n\frac{1}{h_j^{L_j}}K_j\bigg(\frac{\|x_j-\tilde{\xi}_j^i\|_j}{h_j}\bigg)(\I(\mxi^i\in D(n))+\I(\mxi^i\in D))\\
&=O_p\bigg(\sum_{k\neq j}h_k^{\nu_k}\bigg)
\ea
uniformly for $x_k\in D_k$. Hence,
\ba
&\hat{\mdelta}^{C}_j(x_j)=O_p\bigg(\sum_{k\neq j}h_k^{\nu_k}\bigg)\text{~for all~}x_j\in D_j,\quad\bigg(\int_{D_j}\|\hat{\mdelta}^{C}_j(x_j)\|^2dx_j\bigg)^{1/2}=O_p\bigg(\sum_{k\neq j}h_k^{\nu_k}\bigg),\\
&\sup_{x_j\in D_j}\|\hat{\mdelta}^{C}_j(x_j)\|=O_p\bigg(\sum_{k\neq j}h_k^{\nu_k}\bigg).
\ea
This completes the proof. \qed

\subsection{Proof of Theorem \ref{differentiable rate}}

We apply Lemma \ref{general rate} for the proof. We approximated $D_n^{(1)}$, $D_{nj}^{(2)}$ and $\hat{\mdelta}^{A}$ in the proof of Theorem \ref{holder rate}. We approximate $\hat{\mdelta}_j^B$. We define
\ba
\check{\mdelta}_j^B(x_j)=n^{-1}\odot\bigoplus_{i=1}^n((K_{h_j}(x_j,\xi^i_j)\I(\mxi^i\in D))\odot(\mf_j(\xi^i_j)\om\mf_j(x_j))).
\ea
Note that
\begin{align}\label{B error}
\begin{split}
&\sup_{x_j\in D_j}\|\hat{\mdelta}_j^B(x_j)\om\check{\mdelta}_j^B(x_j)\|\\
&\leq\sup_{x_j\in D_j}\bigg\|n^{-1}\od\bigoplus_{i=1}^n((K_{h_j}(x_j,\xi_j^i)\I(\mxi^i\in D)(\I(\tilde{\mxi}^i\in D)-\I(\mxi^i\in D)))\od(\mf_j(\xi^i_j)\om\mf_j(x_j)))\bigg\|\\
&\quad+\sup_{x_j\in D_j}\bigg\|n^{-1}\od\bigoplus_{i=1}^n((K_{h_j}(x_j,\tilde{\xi}_j^i)\I(\tilde{\mxi}^i\in D)(\I(\tilde{\mxi}^i\in D)-\I(\mxi^i\in D)))\od(\mf_j(\xi^i_j)\om\mf_j(x_j)))\bigg\|\\
&\quad+\sup_{x_j\in D_j}\bigg\|n^{-1}\od\bigoplus_{i=1}^n(((K_{h_j}(x_j,\tilde{\xi}_j^i)-K_{h_j}(x_j,\xi_j^i))\I(\mxi^i\in D)\I(\tilde{\mxi}^i\in D))\od(\mf_j(\xi^i_j)\om\mf_j(x_j)))\bigg\|\\
&\leq\sup_{x_j\in D_j}n^{-1}\sum_{i=1}^n(K_{h_j}(x_j,\xi_j^i)\I(\mxi^i\in D\cap D(n))\|\mf_j(\xi^i_j)\om\mf_j(x_j)\|)\\
&\quad+\sup_{x_j\in D_j}n^{-1}\sum_{i=1}^n(K_{h_j}(x_j,\tilde{\xi}_j^i)\I(\tilde{\mxi}^i\in D)\I(\mxi^i\in D(n))\|\mf_j(\xi^i_j)\om\mf_j(x_j)\|)\\
&\quad+\sup_{x_j\in D_j}n^{-1}\sum_{i=1}^n(|K_{h_j}(x_j,\tilde{\xi}_j^i)-K_{h_j}(x_j,\xi_j^i)|\I(\mxi^i\in D)\I(\tilde{\mxi}^i\in D)\|\mf_j(\xi^i_j)\om\mf_j(x_j)\|)\\
&\leq h_j\sup_{x_j\in D_j}n^{-1}\sum_{i=1}^n(K_{h_j}(x_j,\xi_j^i)\I(\mxi^i\in D\cap D(n)))\\
&\quad+(a_{nj}+h_j)\sup_{x_j\in D_j}n^{-1}\sum_{i=1}^n(K_{h_j}(x_j,\tilde{\xi}_j^i)\I(\tilde{\mxi}^i\in D)\I(\mxi^i\in D(n)))\\
&\quad+\const\cdot\sum_{l=1}^{L_j}\max_{1\leq i\leq n}|\tilde{\xi}_{jl}^i-\xi_{jl}^i|\sup_{x_j\in D_j}\frac{1}{n}\sum_{i=1}^n\frac{1}{h_j^{L_j}}\bigg|\frac{\partial}{\partial u_{jl}}K_j\bigg(\frac{\|x_j-u_j\|_j}{h_j}\bigg)\bigg|_{u_j=\bar{\xi}_j^i}\bigg|\I(\mxi^i\in D)\\
&\qquad\qquad\qquad\qquad\qquad\qquad\qquad\qquad\qquad\qquad\cdot\|\mf_j(\xi^i_j)\om\mf_j(x_j)\|\\
&\quad+\const\cdot h_j^{-1}\max_{1\leq i\leq n,1\leq l\leq L_j}|\tilde{\xi}_{jl}^i-\xi_{jl}^i|\sup_{x_j\in D_j}\frac{1}{n}\sum_{i=1}^n\frac{1}{h_j^{L_j}}K_j\bigg(\frac{\|x_j-\xi_j^i\|_j}{h_j}\bigg)\I(\mxi^i\in D)\\
&\qquad\qquad\qquad\qquad\qquad\qquad\qquad\qquad\qquad\qquad\qquad\cdot\|\mf_j(\xi^i_j)\om\mf_j(x_j)\|\\
&=O_p(h_j\cdot(V_{nj}+h_j^{-L_j-2}\cdot a^2_{nj}+h_j^{-1}\cdot a_{nj}))\\
&=O_p(h_j\cdot V_{nj}+h_j^{-L_j-1}\cdot a^2_{nj})\\
&=:O_p(\mathcal{A}_{nj}^{(1)}),
\end{split}
\end{align}
where the last inequality and the first equality follow from the proof of Lemma \ref{marginal density approximation 1}.
By arguing as in the proof of Theorem 4.3 of Jeon et al. (2021), we may show that
\ba
\check{\mdelta}_j^B(x_j)&=O_p(h_j^2+n^{-1/2}\cdot h_j^{-L_j/2+1})\text{~for~}x_j\in D^-_j(2h_j),\\
\check{\mdelta}_j^B(x_j)&=O_p(h_j)\text{~for~}x_j\in D_j\setminus D^-_j(2h_j),
\ea
\ba
\bigg(\int_{D_j}\|\check{\mdelta}_j^B(x_j)\|^2dx_j\bigg)^{1/2}&=O_p(h_j^{3/2}+n^{-1/2}\cdot h_j^{-L_j/2+1}),\\ \bigg(\int_{D^-_j(2h_j)}\|\check{\mdelta}_j^B(x_j)\|^2dx_j\bigg)^{1/2}&=O_p(h_j^2+n^{-1/2}\cdot h_j^{-L_j/2+1}),
\ea
\ba
\sup_{x_j\in D_j}\|\check{\mdelta}_j^B(x_j)\|&=O_p(h_j),\\
\sup_{x_j\in D^-_j(2h_j)}\|\check{\mdelta}_j^B(x_j)\|&=O_p(h_j^2+n^{-1/2}\cdot (\log{n})^{1/2}\cdot h_j^{-L_j/2+1}),
\ea
\begin{align*}
\sum_{k\neq j}\sup_{x_j\in D_j}\bigg\|\int_{D_k}\check{\mdelta}^{B}_k(x_k)\odot
\frac{\check{p}^D_{jk}(x_j,x_k)}{\check{p}^D_k(x_k)}dx_k\bigg\|=O_p\bigg(\sum_{k\neq j}\bigg(h_k^2+n^{-1/2}\cdot h_k^{-L_k/2+1}\bigg)\bigg).
\end{align*}
These with (\ref{B error}) and Lemmas \ref{marginal density approximation 1} and \ref{marginal density approximation 2} entail that
\ba
\hat{\mdelta}_j^B(x_j)&=O_p(h_j^2+n^{-1/2}\cdot h_j^{-L_j/2+1}+\mathcal{A}_{nj}^{(1)})\text{~for~}x_j\in D^-_j(2h_j),\\
\hat{\mdelta}_j^B(x_j)&=O_p(h_j+\mathcal{A}_{nj}^{(1)})\text{~for~}x_j\in D_j\setminus D^-_j(2h_j),
\ea
\ba
\bigg(\int_{D_j}\|\hat{\mdelta}_j^B(x_j)\|^2dx_j\bigg)^{1/2}&=O_p(h_j^{3/2}+n^{-1/2}\cdot h_j^{-L_j/2+1}+\mathcal{A}_{nj}^{(1)}),\\ \bigg(\int_{D^-_j(2h_j)}\|\hat{\mdelta}_j^B(x_j)\|^2dx_j\bigg)^{1/2}&=O_p(h_j^2+n^{-1/2}\cdot h_j^{-L_j/2+1}+\mathcal{A}_{nj}^{(1)}),
\ea
\ba
\sup_{x_j\in D_j}\|\hat{\mdelta}_j^B(x_j)\|&=O_p(h_j+\mathcal{A}_{nj}^{(1)}),\\
\sup_{x_j\in D^-_j(2h_j)}\|\hat{\mdelta}_j^B(x_j)\|&=O_p(h_j^2+n^{-1/2}\cdot (\log{n})^{1/2}\cdot h_j^{-L_j/2+1}+\mathcal{A}_{nj}^{(1)}),
\ea
\begin{align*}
\sum_{k\neq j}\sup_{x_j\in D_j}\bigg\|\int_{D_k}\check{\mdelta}^{B}_k(x_k)\odot
\frac{\check{p}^D_{jk}(x_j,x_k)}{\check{p}^D_k(x_k)}dx_k\bigg\|=O_p\bigg(\sum_{k\neq j}\bigg(h_k^2+n^{-1/2}\cdot h_k^{-L_k/2+1}+\mathcal{A}_{nk}^{(1)}\bigg)\bigg).
\end{align*}

We now approximate $\hat{\mdelta}^{C}_j$. We define
\ba
\mQ_k^i(u_k)=\int_{D_k}(\mf_k(\xi^i_k)\om\mf_k(x_k))\odot K_{h_k}(x_k,u_k)dx_k
\ea
for $u_k\in\mathbb{R}^{L_k}$. Note that
\begin{align}\label{Q bound}
\I(\mxi^i\in D)\cdot\|\mQ_k^i(\xi^i_k)\|\leq \const\cdot h_k.
\end{align}
and
\begin{align}\label{Q bound2}
\begin{split}
\I(\tilde{\mxi^i}\in D)\cdot\|\mQ_k(\tilde{\xi}^i_k)\|&\leq\int_{D_k}(\|\mf_k(\xi^i_k)\om\mf_k(\tilde{\xi}^i_k)\|+\|\mf_k(\tilde{\xi}^i_k)\om\mf_k(x_k)\|)K_{h_k}(x_k,\tilde{\xi}^i_k)dx_k\\
&\leq\const\cdot (a_{nk}+h_k).
\end{split}
\end{align}
We also note that
\begin{align}\label{Q bound3}
\begin{split}
&\I(\mxi^i\in D)\I(\tilde{\mxi}^i\in D)\|\mQ_k^i(\tilde{\xi}^i_k)\om\mQ_k^i(\xi^i_k)\|\\
&\leq\I(\mxi^i\in D)\I(\tilde{\mxi}^i\in D)\int_{D_k}\|\mf_k(\xi_k^i)\om\mf_k(x_k)\|\cdot|K_{h_k}(x_k,\tilde{\xi}_k^i)-K_{h_k}(x_k,\xi_k^i)|dx_k\\
&\leq\const\cdot\sum_{l=1}^{L_k}|\tilde{\xi}_{kl}^i-\xi_{kl}^i|\int_{D_k}\|\mf_k(\xi_k^i)\om\mf_k(x_k)\|\frac{1}{h_k^{L_k}}\bigg|\frac{\partial}{\partial u_{kl}}K_k\bigg(\frac{\|x_k-u_k\|_k}{h_k}\bigg)\bigg|_{u_k=\bar{\xi}_k^i}\bigg|dx_k\\
&\quad\cdot\I(\mxi^i\in D)\I(\tilde{\mxi}^i\in D)\\
&\quad+\const\cdot h_k^{-1}\max_{1\leq l\leq L_k}|\tilde{\xi}_{kl}^i-\xi_{kl}^i|\int_{D_k}\|\mf_k(\xi_k^i)\om\mf_k(x_k)\|\frac{1}{h_k^{L_k}}K_k\bigg(\frac{\|x_k-\xi_k^i\|_k}{h_k}\bigg)dx_k\\
&\quad\cdot\I(\mxi^i\in D)\I(\tilde{\mxi}^i\in D)\\
&\leq\const\cdot h_k^{-1}\cdot a_{nk}\cdot(a_{nk}+h_k)\cdot\I(\mxi^i\in D)\I(\tilde{\mxi}^i\in D).
\end{split}
\end{align}
Define
\ba
\hat{\mdelta}^{C}_{jk}(x_j)&=n^{-1}\odot\bigoplus_{i=1}^n((K_{h_j}(x_j,\tilde{\xi}^i_j)\I(\tilde{\mxi^i}\in D))\od \mQ_k^i(\tilde{\xi}^i_k)),\quad\mbox{and}\\
\check{\mdelta}^{C}_{jk}(x_j)&=n^{-1}\odot\bigoplus_{i=1}^n((K_{h_j}(x_j,\xi^i_j)\I(\mxi^i\in D))\od\mQ_k^i(\xi^i_k)).
\ea
Note that
\begin{align*}
&\sup_{x_j\in D_j}\|\hat{\mdelta}_{jk}^C(x_j)\om\check{\mdelta}_{jk}^C(x_j)\|\\
&\leq\sup_{x_j\in D_j}\bigg\|n^{-1}\od\bigoplus_{i=1}^n((K_{h_j}(x_j,\xi_j^i)\I(\mxi^i\in D)(\I(\tilde{\mxi}^i\in D)-\I(\mxi^i\in D)))\od\mQ_k^i(\xi^i_k))\bigg\|\\
&\quad+\sup_{x_j\in D_j}\bigg\|n^{-1}\od\bigoplus_{i=1}^n((K_{h_j}(x_j,\tilde{\xi}_j^i)\I(\tilde{\mxi}^i\in D)(\I(\tilde{\mxi}^i\in D)-\I(\mxi^i\in D)))\od\mQ_k^i(\tilde{\xi}^i_k))\bigg\|\\
&\quad+\sup_{x_j\in D_j}\bigg\|n^{-1}\od\bigoplus_{i=1}^n((K_{h_j}(x_j,\tilde{\xi}_j^i)\od \mQ_k^i(\tilde{\xi}^i_k)\om K_{h_j}(x_j,\xi_j^i)\od \mQ_k^i(\xi^i_k))\\
&\qquad\qquad\qquad\qquad\qquad\od(\I(\mxi^i\in D)\I(\tilde{\mxi}^i\in D)))\bigg\|\\
&\leq h_k\sup_{x_j\in D_j}n^{-1}\sum_{i=1}^n(K_{h_j}(x_j,\xi_j^i)\I(\mxi^i\in D\cap D(n)))\\
&\quad+(a_{nk}+h_k)\sup_{x_j\in D_j}n^{-1}\sum_{i=1}^n(K_{h_j}(x_j,\tilde{\xi}_j^i)\I(\tilde{\mxi}^i\in D)\I(\mxi^i\in D(n)))\\
&\quad+\sup_{x_j\in D_j}\bigg\|n^{-1}\od\bigoplus_{i=1}^n((K_{h_j}(x_j,\tilde{\xi}_j^i)\I(\mxi^i\in D)\I(\tilde{\mxi}^i\in D))\od (\mQ_k^i(\tilde{\xi}^i_k)\om\mQ_k^i(\xi^i_k)))\bigg\|\\
&\quad+\sup_{x_j\in D_j}\bigg\|n^{-1}\od\bigoplus_{i=1}^n(((K_{h_j}(x_j,\tilde{\xi}_j^i)-K_{h_j}(x_j,\xi_j^i))\I(\mxi^i\in D)\I(\tilde{\mxi}^i\in D))\od \mQ_k^i(\xi^i_k))\bigg\|\\
&\leq O_p(h_k\cdot V_{nj})+\sup_{x_j\in D_j}n^{-1}\sum_{i=1}^nK_{h_j}(x_j,\tilde{\xi}_j^i)\I(\mxi^i\in D)\I(\tilde{\mxi}^i\in D)\|\mQ_k^i(\tilde{\xi}^i_k)\om\mQ_k^i(\xi^i_k)\|\\
&\quad+\const\cdot\sum_{l=1}^{L_j}\max_{1\leq i\leq n}|\tilde{\xi}_{jl}^i-\xi_{jl}^i|\sup_{x_j\in D_j}\frac{1}{n}\sum_{i=1}^n\frac{1}{h_j^{L_j}}\bigg|\frac{\partial}{\partial u_{jl}}K_j\bigg(\frac{\|x_j-u_j\|_j}{h_j}\bigg)\bigg|_{u_j=\bar{\xi}_j^i}\bigg|\I(\mxi^i\in D)\\
&\qquad\qquad\qquad\qquad\qquad\qquad\qquad\qquad\qquad\qquad\cdot\|\mQ_k^i(\xi^i_k)\|\\
&\quad+\const\cdot h_j^{-1}\max_{1\leq i\leq n,1\leq l\leq L_j}|\tilde{\xi}_{jl}^i-\xi_{jl}^i|\sup_{x_j\in D_j}\frac{1}{n}\sum_{i=1}^n\frac{1}{h_j^{L_j}}K_j\bigg(\frac{\|x_j-\xi_j^i\|_j}{h_j}\bigg)\I(\mxi^i\in D)\|\mQ_k^i(\xi^i_k)\|\\
&=O_p(h_k\cdot V_{nj})+O_p(a_{nk})+O_p(h_k\cdot(h_j^{-L_j-2}\cdot a^2_{nj}+h_j^{-1}\cdot a_{nj})),
\end{align*}
where the second inequality follows from (\ref{Q bound}) and (\ref{Q bound2}) and the last equality follows from (\ref{Q bound}) and (\ref{Q bound3}).
Thus, for $\check{\mdelta}_j^C(x_j)=\bigoplus_{k\neq j}\check{\mdelta}_{jk}^C(x_j)$, we have
\begin{align}\label{C error}
\begin{split}
\sup_{x_j\in D_j}\|\hat{\mdelta}_j^C(x_j)\om\check{\mdelta}_j^C(x_j)\|&=O_p\bigg((h_j^{-L_j}\cdot a_{nj}+h_j^{-L_j-2}\cdot a^2_{nj})\cdot\bigg(\sum_{k\neq j}h_k\bigg)+\sum_{k\neq j}a_{nk}\bigg)\\
&=:O_p(\mathcal{A}^{(2)}_{nj}).
\end{split}
\end{align}
By arguing as in the proof of Theorem 4.3 in Jeon et al. (2021), we may show that
\ba
\check{\mdelta}_j^C(x_j)=\sum_{k\neq j}O_p(h_k^2+h_k\cdot n^{-1/2}\cdot h_j^{-L_j/2})\text{~for all~}x_j\in D_j,
\ea
\ba
\bigg(\int_{D_j}\|\check{\mdelta}_j^C(x_j)\|^2dx_j\bigg)^{1/2}=\sum_{k\neq j}O_p(h_k^2+h_k\cdot n^{-1/2}\cdot h_j^{-L_j/2}),
\ea
\ba
\sup_{x_j\in D_j}\|\check{\mdelta}_j^C(x_j)\|=\sum_{k\neq j}O_p(h_k^2+h_k\cdot n^{-1/2}\cdot(\log{n})^{1/2}\cdot h_j^{-L_j/2}).
\ea
These with (\ref{C error}) entail that
\ba
\hat{\mdelta}_j^C(x_j)=\sum_{k\neq j}O_p(h_k^2+h_k\cdot n^{-1/2}\cdot h_j^{-L_j/2})+O_p(\mathcal{A}^{(2)}_{nj})\text{~for all~}x_j\in D_j,
\ea
\ba
\bigg(\int_{D_j}\|\hat{\mdelta}_j^C(x_j)\|^2dx_j\bigg)^{1/2}=\sum_{k\neq j}O_p(h_k^2+h_k\cdot n^{-1/2}\cdot h_j^{-L_j/2})+O_p(\mathcal{A}^{(2)}_{nj}),
\ea
\ba
\sup_{x_j\in D_j}\|\hat{\mdelta}_j^C(x_j)\|=\sum_{k\neq j}O_p(h_k^2+h_k\cdot n^{-1/2}\cdot(\log{n})^{1/2}\cdot h_j^{-L_j/2})+O_p(\mathcal{A}^{(2)}_{nj}).
\ea
This completes the proof. \qed

\subsection{Proof of Theorem \ref{Asymptotic distribution}}

Define
\begin{align*}
\tilde{\mDelta}_j^*(x_j)=\mdelta_j^*(x_j)\op\bigoplus_{k\neq j}\int_{ D_k}\mdelta_{jk}^*(x_j,x_k)
\od\frac{p^D_{jk}(x_j,x_k)}{p^D_j(x_j)}dx_k,
\end{align*}
where
\begin{align*}
\mdelta^*_j(x_j)&=\frac{h_j^2}{p^D_j(x_j)}\od\mathfrak{D}^1\mf_j(x_j)\bigg(\int_{B_j(\bzero_j,1)}t_j\cdot\mathfrak{D}^1p^D_j(x_j)(t_j)\cdot K_j(\|t_j\|_j)dt_j\bigg),\quad\mbox{and}\\
\mdelta^*_{jk}(x_j,x_k)
&=\frac{h_k^2}{p^D_{jk}(x_j,x_k)}\od\mathfrak{D}^1\mf_k(x_k)\bigg(\int_{B_k(\bzero_k,1)}t_k\cdot\mathfrak{D}^1_2\,p^D_{jk}(x_j,x_k)(t_k)\cdot K_k(\|t_k\|_k)dt_k\bigg).
\end{align*}
We claim that there exists a solution $(\mDelta^*_1,\ldots,\mDelta^*_d)$ of the system of equations
\begin{align}\label{Delta equation}
\mDelta^*_j(x_j)=\tilde{\mDelta}^*_j(x_j)\ominus\opluskj\int_{ D_k}
\mDelta^*_k(x_k)\odot\frac{p^D_{jk}(x_j,x_k)}{p^D_j(x_j)}dx_k, \quad \allj
\end{align}
subject to the constraints
\begin{align}\label{Delta constraint new}
\begin{split}
\int_{ D_j}\mDelta^*_j(x_j)\odot p^D_j(x_j)dx_j&=\int_{ D_j}\mdelta_j^*(x_j)\odot p^D_j(x_j)dx_j,\quad\allj,\\
\int_{ D_j}\|\mDelta^*_j(x_j)\|^2p^D_j(x_j)dx_j&<\infty,\quad\allj.
\end{split}
\end{align}
For this, define
\begin{align*}
\mA(\mx)=\bigoplus_{j=1}^{d}\frac{h_j^2}{p^D(\mx)}\od \mathfrak{D}^1\mf_j(x_j)\bigg(\int_{B_j(\bzero_j,1)}t_j\cdot\mathfrak{D}^1_j\,p^D(\mx)(t_j)\cdot K_j(\|t_j\|_j)dt_j\bigg),
\end{align*}
where $\mathfrak{D}^1_j\,p(\mx)$ is the partial Fr\'{e}chet derivative of $p$ at $\mx$ with respect to the $j$th argument. Denote by $\mA_j(\mx)$ the $j$th term of $\mA(\mx)$. Then,
\begin{align*}
&\sup_{\mx\in D}\|\mA_j(\mx)\|\\
&\leq h_j^2\cdot\sup_{\mx\in D}\frac{1}{p^D(\mx)}\cdot\sup_{x_j\in D_j}
\|\mathfrak{D}^1\mf_j(x_j)\|_{\rm op}\cdot\sup_{\mx\in D}\|\mathfrak{D}^1_j\,p^D(\mx)\|_{\rm op}\cdot\int_{B_j(\bzero_j,1)}\|t_j\|_j^2\cdot K_j(\|t_j\|_j)dt_j,
\end{align*}
where $\|\cdot\|_{\rm op}$ denotes the operator norms on the respective spaces of bounded linear operators. Hence, $\mA\in L_2^\mbH(p^D)$. We define $\mf^*=\bigoplus_{j=1}^d\mf^*_j$ as the minimizer of $F_{\mA}(\mg):=\int_D\|\mA(\mx)\om\mg(\mx)\|^2dP_{\mxi}^D(\mx)$ over $S^\mbH(p^D)$. Such $\mf^*$ exists by arguing as in the proof of Theorem 4.4 in Jeon et al. (2021a). By Theorem 5.3.19 in Atkinson and Han (2009), it holds that $\mathfrak{D}^1F_{\mA}(\mf^*)(\mk)=0$ for all $\mk\in S^\mbH(p^D)$, where $\mathfrak{D}^1F_{\mA}(\mf^*)(\mk)$ is the G\^{a}teaux derivative of $F_{\mA}$ at $\mf^*$ in the direction of $\mk$. This implies that
\begin{align*}
\bigoplus_{k=1}^d\int_{D_{-j}}(\mA_k(\mx)\om\mf^*_k(x_k))\od p^D(\mx)d\mx_{-j}=\bzero
\end{align*}
almost everywhere $x_j$ with respect to $\Leb_j$ for all $1\leq j\leq d$. The latter implies that
\begin{align*}
\mf_j^*(x_j)\od p^D_j(x_j)=&\int_{ D_{-j}}\mA_j(\mx)\od p^D(\mx)d\mx_{-j}\op\bigoplus_{k\neq j}\int_{ D_{-j}}\mA_k(\mx)\od p^D(\mx)d\mx_{-j}\\
&\om\bigoplus_{k\neq j}\int_{ D_k}\mf_k^*(x_k)\od p^D_{jk}(x_j,x_k)dx_k
\end{align*}
almost everywhere $x_j$ with respect to $\Leb_j$ for all $1\leq j\leq d$. One can show that
\begin{align*}
&\int_{ D_{-j}}\mA_j(\mx)\od p^D(\mx)d\mx_{-j}\\
&=h_j^2\od \mathfrak{D}^1\mf_j(x_j)\bigg(\int_{B_j(\bzero_j,1)}t_j\cdot\int_{D_{-j}}\mathfrak{D}^1_j\,p^D(\mx)(t_j)d\mx_{-j}\cdot K_j(\|t_j\|_j)dt_j\bigg)\\
&=\mdelta_j^*(x_j)\od p^D_j(x_j)
\end{align*}
for $1\leq j\leq d$. Similarly, one can show that
\begin{align*}
\int_{ D_{-j}}\mA_k(\mx)\od p^D(\mx)d\mx_{-j}=\int_{ D_k}\mdelta_{jk}^*(x_j,x_k)\od p^D_{jk}(x_j,x_k)dx_k
\end{align*}
for $1\leq j\neq k\leq d$. Hence, we get
\begin{align*}
\mf^*_j(x_j)=\tilde{\mDelta}^*_j(x_j)\ominus\opluskj\int_{ D_k}
\mf^*_k(x_k)\odot\frac{p^D_{jk}(x_j,x_k)}{p^D_j(x_j)}dx_k
\end{align*}
almost everywhere $x_j$ with respect to $\Leb_j$ for all $1\leq j\leq d$. We denote by $\mf^{**}_j(x_j)$ the right hand side of the above equality. Then, it holds that
\begin{align*}
\mf^{**}_j(x_j)=\tilde{\mDelta}^*_j(x_j)\ominus\opluskj\int_{ D_k}
\mf^{**}_k(x_k)\odot\frac{p^D_{jk}(x_j,x_k)}{p^D_j(x_j)}dx_k
\end{align*}
for all $x_j\in D_j$ and $1\leq j\leq d$. Since
\begin{align*}
\int_{ D_j}\mf^{**}_j(x_j)\od p^D_j(x_j)dx_j=\int_{ D_j}\tilde{\mDelta}^*_j(x_j)\od p^D_j(x_j)dx_j\om\bigoplus_{k\neq j}\int_{ D_k}\mf^{**}_k(x_k)\odot p^D_k(x_k)dx_k
\end{align*}
and
\begin{align*}
\int_{ D_j}\tilde{\mDelta}^*_j(x_j)\od p^D_j(x_j)dx_j=\bigoplus_{j=1}^d\int_{ D_j}\mdelta^*_j(x_j)\od p^D_j(x_j)dx_j,
\end{align*}
we get
\begin{align*}
\bigoplus_{j=1}^d\int_{ D_j}\mf^{**}_j(x_j)\od p^D_j(x_j)dx_j=\bigoplus_{j=1}^d\int_{ D_j}\mdelta^*_j(x_j)\od p^D_j(x_j)dx_j.
\end{align*}
This entails that $(\mDelta_j^*:1\leq j\leq d)$ with $\mDelta_j^*$ being defined by
\ba
\mDelta_j^*(x_j)=\mf^{**}_j(x_j)\om\int_{ D_j}\mf^{**}_j(x_j)\od p^D_j(x_j)dx_j\op\int_{ D_j}\mdelta^*_j(x_j)\od p^D_j(x_j)dx_j
\ea
satisfies (\ref{Delta equation}) and (\ref{Delta constraint new}).
Hence, the claim follows. We note that $\sup_{x_j\in D_j}\|\tilde{\mDelta}^*_j(x_j)\|=O(\sum_{k=1}^dh_k^2)$. Thus, by arguing as in the proof of Theorem 5.1 in Jeon et al. (2021a), we may show that $\|\mDelta^*_j\|_{2,D}=O(\sum_{k=1}^dh_k^2)$. This with (\ref{Delta equation}) gives that $\sup_{x_j\in D_j}\|\mDelta^*_j(x_j)\|=O(\sum_{k=1}^dh_k^2)$ for all $1\leq j\leq d$. Therefore, the bounded convergence theorem entails that
\begin{align*}
(\mDelta_1,\ldots,\mDelta_d):=\lim_{n\rightarrow\infty}(n^{2/(L_{\rm max}+4)}\od\mDelta^*_1,\ldots,n^{2/(L_{\rm max}+4)}\od\mDelta^*_d)
\end{align*}
satisfies (\ref{Delta}) and (\ref{Delta constraint}). The uniqueness of a solution of (\ref{Delta}) subject to (\ref{Delta constraint})
follows by arguing as in the proof of Theorem \ref{convergence2}. This completes the proof for the first part of the theorem.

Now, we prove the second part of the theorem. From (\ref{equation 3 new}), we get
\begin{align}\label{equation 3 new new}
\begin{split}
\hat{\mf}_j(x_j)\om\mf_j(x_j)=&\hat{\mf}_j^A(x_j)\op\hat{\mf}_j^B(x_j)\op\hat{\mf}_j^C(x_j)\op\mf_0\om\hat{\mf}_0\op\hat{\mm}_j(x_j)\om\check{\mm}_j(x_j)\\
&\om\bigoplus_{k\neq j}\int_{ D_k}(\hat{\mf}_k(x_k)\om\mf_k(x_k))\od\frac{\hat{p}^D_{jk}(x_j,x_k)}{\hat{p}^D_j(x_j)}dx_k,\quad1\leq j\leq d.
\end{split}
\end{align}
Define
\begin{align*}
\ma^*_k(x_k)=&\mathfrak{D}^1\mf_k(x_k)\bigg(\int_{ D_k}(u_k-x_k)K_{h_k}(x_k,u_k)du_k\bigg),\\
\mb^*_{jk}(x_j,x_k)=&\frac{1}{p^D_{jk}(x_j,x_k)}\od \mathfrak{D}^1\mf_k(x_k)\\
&\bigg(\int_{ D_k}(u_k-x_k)(\mathfrak{D}^1_2\,p^D_{jk}(x_j,x_k)(u_k-x_k)\cdot K_{h_k}(x_k,u_k))du_k\bigg),\\
\mc^*_k(x_k)=&\frac{1}{2}\od\int_{ D_k}\mathfrak{D}^2\mf_k(x_k)(u_k-x_k)(u_k-x_k)\od K_{h_k}(x_k,u_k)du_k,\\
\mu_{0,k}(x_k)=&\int_{D_k}K_{h_k}(x_k,u_k)du_k.
\end{align*}
It holds that
\begin{align}\label{a and c}
\begin{split}
&\sup_{x_k\in D^-_k(2h_k)}\|\ma^*_k(x_k)\|=0,\quad\sup_{x_k\in D_k\setminus D^-_k(2h_k)}\|\ma^*_k(x_k)\|=O(h_k),\\
&\sup_{(x_j,x_k)\in D_j\times D_k}\|\mb^*_{jk}(x_j,x_k)\|=O(h^2_k),\quad\sup_{x_k\in D_k}\|\mc^*_k(x_k)\|=O(h^2_k).
\end{split}
\end{align}
By arguing as in the proof of Theorem 4.4 in Jeon et al. (2021a) and using (\ref{a and c}) and Lemma \ref{uniform rates}, we can show that
\begin{align*}
\check{\mdelta}_j^C(x_j)=&\bigoplus_{k\neq j}\bigg(\int_{ D_k}(\mu_{0,k}(x_k)^{-1}\od\ma_k^*(x_k)\op\mc_k^*(x_k))\od(\check{p}^D_{jk}(x_j,x_k)\cdot \check{p}_0^D)dx_k\\
&\op(\check{p}_j^D(x_j)\cdot \check{p}_0^D)\od\int_{ D_k}\mb^*_{jk}(x_j,x_k)\od \frac{p^D_{jk}(x_j,x_k)}{p_j^D(x_j)}dx_k\op O_p\bigg(h_k\sqrt{\frac{\log n}{n\cdot h_j^{L_j}}}\bigg)\op o_p(h_k^2)\bigg)
\end{align*}
uniformly for $x_j\in D_j$. From this and using (\ref{a and c}) and Lemmas \ref{p0 rate}, \ref{marginal density approximation 1} and \ref{marginal density approximation 2}, we get
\begin{align*}
\check{\mdelta}_j^C(x_j)=&\bigoplus_{k\neq j}\bigg(\int_{ D_k}(\mu_{0,k}(x_k)^{-1}\od\ma_k^*(x_k)\op\mc_k^*(x_k))\od(\hat{p}^D_{jk}(x_j,x_k)\cdot \hat{p}_0^D)dx_k\\
&\op(\hat{p}_j^D(x_j)\cdot \hat{p}_0^D)\od\int_{ D_k}\mb^*_{jk}(x_j,x_k)\od \frac{p^D_{jk}(x_j,x_k)}{p_j^D(x_j)}dx_k\op O_p\bigg(h_k\sqrt{\frac{\log n}{n\cdot h_j^{L_j}}}\bigg)\op o_p(h_k^2)\bigg)
\end{align*}
uniformly for $x_j\in D_j$. Hence,
\ba
\hat{\mf}_j^C(x_j)&=(\hp^D_j(x_j)\cdot\hat{p}^D_0)^{-1}\od(\hat{\mdelta}_j^C(x_j)\om\check{\mdelta}_j^C(x_j))\op(\hp^D_j(x_j)\cdot\hat{p}^D_0)^{-1}\od\check{\mdelta}_j^C(x_j)\\
&=\bigoplus_{k\neq j}\bigg(\int_{ D_k}(\mu_{0,k}(x_k)^{-1}\od\ma_k^*(x_k)+\mc_k^*(x_k))\od\frac{\hat{p}^D_{jk}(x_j,x_k)}{\hat{p}^D_j(x_j)}dx_k\\
&\quad\op\int_{ D_k}\mb^*_{jk}(x_j,x_k)\od\frac{p^D_{jk}(x_j,x_k)}{p^D_j(x_j)}dx_k\op O_p\bigg(h_k\sqrt{\frac{\log n}{n\cdot h_j^{L_j}}}\bigg)\op o_p(h_k^2)\bigg)\op O_p(\mathcal{A}^{(2)}_{nj})
\ea
uniformly for $x_j\in D_j$.

By arguing as in the proof of Theorem 4.4 in Jeon et al. (2021a) and using (\ref{a and c}) and Lemma \ref{uniform rates}, we can also show that
\begin{align}\label{B check approximation}
\begin{split}
&\check{\mdelta}_j^B(x_j)\\
&=(\mu_{0,j}(x_j)\cdot p^D_j(x_j)\cdot p_0^D)\od(\mu_{0,j}(x_j)^{-1}\od\ma^*_j(x_j))\op(p^D_j(x_j)\cdot p_0^D)\od\mc_j^*(x_j)\op\tilde{\mb}_j(x_j)\\
&\quad\op\mr^{(1)}_j(x_j)\op O_p\bigg(h_j\sqrt{\frac{\log n}{n\cdot h_j^{L_j}}}\bigg)\\
&=(\check{p}^D_j(x_j)\cdot\check{p}_0^D)\od(\mu_{0,j}(x_j)^{-1}\od\ma^*_j(x_j)\op\mc_j^*(x_j))\op\tilde{\mb}_j(x_j)\\
&\quad\op\mr^{(1)}_j(x_j)\op O_p\bigg(h_j\sqrt{\frac{\log n}{n\cdot h_j^{L_j}}}\bigg)
\end{split}
\end{align}
uniformly for $x_j\in D_j$, where
\begin{align*}
\tilde{\mb}_j(x_j)=p_0^D\od\mathfrak{D}^1\mf_j(x_j)\bigg(\int_{ D_j}(u_j-x_j)(\mathfrak{D}^1p^D_j(x_j)(u_j-x_j)\cdot K_{h_j}(x_j,u_j))du_j\bigg)
\end{align*}
and $\mr^{(1)}_j: D_j\rightarrow\mbH$ is a generic stochastic map satisfying
\begin{align*}
\sup_{x_j\in D_j}\|\mr^{(1)}_j(x_j)\|=O_p(h_j^2),\quad\sup_{x_j\in D^-_j(2h_j)}\|\mr^{(1)}_j(x_j)\|=o_p(h_j^2).
\end{align*}
Note that
\begin{align}\label{b rate}
\sup_{x_j\in D_j}\|\tilde{\mb}_j(x_j)\|=O(h_j^2).
\end{align}
From (\ref{B check approximation}) and using (\ref{a and c}) and Lemmas \ref{p0 rate} and \ref{marginal density approximation 1}, we get
\ba
\check{\mdelta}_j^B(x_j)=&(\hat{p}^D_j(x_j)\cdot\hat{p}_0^D)\od(\mu_{0,j}(x_j)^{-1}\od\ma^*_j(x_j)\op\mc_j^*(x_j))\op\tilde{\mb}_j(x_j)\op\mr^{(1)}_j(x_j)\\
&\op O_p\bigg(h_j\sqrt{\frac{\log n}{n\cdot h_j^{L_j}}}\bigg)\op\mr^{(2)}_j(x_j),
\ea
where $\mr^{(2)}_j: D_j\rightarrow\mbH$ is a stochastic map satisfying
\ba
\sup_{x_j\in D_j}\|\mr^{(2)}_j(x_j)\|&=O_p\bigg(h_j\cdot\bigg(h_j^{-L_j}\cdot a_{nj}+\sum_{k\neq j}a_{nk}+h_j^{-L_j-2}\cdot a_{nj}^2\bigg)\bigg),\\
\sup_{x_j\in D^-_j(2h_j)}\|\mr^{(2)}_j(x_j)\|&=O_p\bigg(h^2_j\cdot\bigg(h_j^{-L_j}\cdot a_{nj}+\sum_{k\neq j}a_{nk}+h_j^{-L_j-2}\cdot a_{nj}^2\bigg)\bigg).
\ea
Hence,
\begin{align}\label{B result}
\begin{split}
\hat{\mf}_j^B(x_j)&=(\hp^D_j(x_j)\cdot\hat{p}^D_0)^{-1}\od(\hat{\mdelta}_j^B(x_j)\om\check{\mdelta}_j^B(x_j))\op(\hp^D_j(x_j)\cdot\hat{p}^D_0)^{-1}\od\check{\mdelta}_j^B(x_j)\\
&=\mu_{0,j}(x_j)^{-1}\od\ma^*_j(x_j)\op\mc_j^*(x_j)\op(\hp^D_j(x_j)\cdot\hat{p}^D_0)^{-1}\od\tilde{\mb}_j(x_j)\op\mr^{(1)}_j(x_j)\\
&\quad \op O_p\bigg(h_j\sqrt{\frac{\log n}{n\cdot h_j^{L_j}}}\bigg)\op\mr^{(2)}_j(x_j)\op O_p(\mathcal{A}^{(1)}_{nj})\\
&=\mu_{0,j}(x_j)^{-1}\od\ma^*_j(x_j)\op\mc_j^*(x_j)\op\mb^*_j(x_j)\op\mr^{(1)}_j(x_j)\\
&\quad \op O_p\bigg(h_j\sqrt{\frac{\log n}{n\cdot h_j^{L_j}}}\bigg)\op\mr^{(2)}_j(x_j)\op O_p(\mathcal{A}^{(1)}_{nj}),
\end{split}
\end{align}
where
\begin{align*}
\mb^*_j(x_j)=(p_j^D(x_j))^{-1}\od\mathfrak{D}^1\mf_j(x_j)\bigg(\int_{ D_j}(u_j-x_j)(\mathfrak{D}^1p^D_j(x_j)(u_j-x_j)\cdot K_{h_j}(x_j,u_j))du_j\bigg)
\end{align*}
and the last equality at (\ref{B result}) follows from (\ref{b rate}) and Lemmas \ref{p0 rate} and \ref{marginal density convergence 1}.

Define
\ba
\check{\mf}_j^A(x_j)=(\hp^D_j(x_j)\cdot\hat{p}^D_0)^{-1}\od\check{\mdelta}_j^A(x_j).
\ea
Note that we have $(\hp^D_j(x_j)\cdot\hat{p}^D_0)^{-1}$ instead of $(\check{p}^D_j(x_j)\cdot\check{p}^D_0)^{-1}$ in the definition of $\check{\mf}_j^A(x_j)$. It holds that
\ba
\sup_{x_j\in D_j}\|\hat{\mf}_j^A(x_j)\om\check{\mf}_j^A(x_j)\|&\leq\sup_{x_j\in D_j}(\hp^D_j(x_j)\cdot\hat{p}^D_0)^{-1}\cdot\sup_{x_j\in D_j}\|\hat{\mdelta}_j^A(x_j)\om\check{\mdelta}_j^A(x_j)\|\\
&=O_p(n^{-1/2}\cdot(\log{n})^{1/2}\cdot h_j^{-L_j/2}\cdot V_{nj}^{1/2}),
\ea
where the equality follows from (\ref{A error}) and Lemmas \ref{p0 rate} and \ref{marginal density convergence 1}.
We also claim that
\begin{align}\label{A claim}
\sup_{x_j\in D_j}\bigg\|\int_{ D_k}\check{\mf}_k^A(x_k)\od\frac{\hat{p}^D_{jk}(x_j,x_k)}{\hat{p}^D_j(x_j)}dx_k\bigg\|=O_p(n^{-1/2}\cdot(\log n)^{1/2}).
\end{align}
For this, it suffices to prove that
\ba
\sup_{x_j\in D_j}\bigg\|\int_{ D_k}\check{\mf}_k^A(x_k)\od\hat{p}^D_{jk}(x_j,x_k)dx_k\bigg\|=O_p(n^{-1/2}\cdot(\log n)^{1/2})
\ea
by Lemma \ref{marginal density convergence 1}. We define
\ba
A^i_k(x_j)=(\hat{p}_0^D)^{-1}\cdot\int_{D_k}K_{h_k}(x_k,\xi_k^i)\cdot\frac{\hat{p}^D_{jk}(x_j,x_k)}{\hat{p}^D_k(x_k)}dx_k\cdot\I(\mxi^i\in D).
\ea
Then,
\ba
\int_{ D_k}\check{\mf}_k^A(x_k)\od\hat{p}^D_{jk}(x_j,x_k)dx_k=n^{-1}\od\bigoplus_{i=1}^n(A^i_k(x_j)\od\mepsilon^i).
\ea
Since $\E(A^i_k(x_j)\od\mepsilon^i|\mxi^{\rm all},\tilde{\mxi}^{\rm all})=A^i_k(x_j)\od\E(\mepsilon^i|\mxi^{\rm all},\tilde{\mxi}^{\rm all})=\bzero$, it suffices to show that
\ba
&\sup_{x_j\in D_j}\bigg\|n^{-1}\od\bigoplus_{i=1}^n(A^i_k(x_j)\od\mepsilon^i)\om\E\bigg(n^{-1}\od\bigoplus_{i=1}^n(A^i_k(x_j)\od\mepsilon^i)\bigg|\mxi^{\rm all},\tilde{\mxi}^{\rm all}\bigg)\bigg\|\\
&=O_p(n^{-1/2}\cdot(\log n)^{1/2}).
\ea
We take any $\varsigma$ such that $1/\alpha<\varsigma<1/2$. Then,
\begin{align}\label{A claim 1}
\Prob(\|\mepsilon^i\|\leq n^{\varsigma}~\text{for all}~1\leq i\leq n)\geq1-\frac{n\cdot\E(\|\mepsilon\|^\alpha)}{n^{\alpha\cdot\varsigma}}\rightarrow1.
\end{align}
Also, it holds that
\begin{align}\label{A claim 2}
|A^i_k(x_j)-A^i_k(x^*_j)|\leq\const \cdot h_j^{-L_j-1}\cdot\|x_j-x_j^*\|_j
\end{align}
for any $x_j,x_j^*\in D_j$. By arguing as in the proof of Lemma \ref{uniform rates general} and using (\ref{A claim 1}) and (\ref{A claim 2}), one can check that proving
\begin{align}\label{A claim 3}
\sup_{x_j\in D_j}\Prob\bigg(\bigg\|\bigoplus_{i=1}^n\mZ_n^i(x_j)\bigg\|>C_0\cdot n^{-1/2}\cdot(\log{n})^{1/2}\bigg)\leq \const\cdot n^{-C_0}\cdot n^{O(1)}
\end{align}
gives (\ref{A claim}), where
\ba
\mZ_n^i(x_j)=&(n^{-1}\cdot A^i_k(x_j)\cdot\I(\|\mepsilon^i\|\leq n^{\varsigma}))\od\mepsilon^i\om\E((n^{-1}\cdot A^i_k(x_j)\cdot\I(\|\mepsilon^i\|\leq n^{\varsigma}))\od\mepsilon^i|\mxi^{\rm all},\tilde{\mxi}^{\rm all}).
\ea
By Lemmas \ref{p0 rate}, \ref{marginal density convergence 1} and \ref{marginal density convergence 2}, we may show that
\begin{align}\label{A bound}
\Prob\Big(\max_{1\leq i\leq n}\sup_{x_j\in D_j}|A^i_k(x_j)|\leq\const\Big)\rightarrow1.
\end{align}
Hence, we may assume that $\max_{1\leq i\leq n}\sup_{x_j\in D_j}|A^i_k(x_j)|\leq\const$. Then, by arguing as in the proof of (\ref{A error}) and using that
\ba
&\E(\mZ^i_n(x_j)|\mxi^{\rm all},\tilde{\mxi}^{\rm all})=\bzero,\quad\sup_{x_j\in D_j}\|\mZ_n^i(x_j)\|\leq\const\cdot n^{\varsigma-1},\quad\mbox{and}\\
&\sup_{x_j\in D_j}\E(\|\mZ^i_n(x_j)\|^2|\mxi^{\rm all},\tilde{\mxi}^{\rm all})\leq\const\cdot n^{-2}~~\text{almost surely},
\ea
we may prove that
\ba
&\sup_{x_j\in D_j}\Prob\bigg(\bigg\|\bigoplus_{i=1}^n\mZ_n^i(x_j)\bigg\|>C_0\cdot n^{-1/2}\cdot(\log{n})^{1/2}\bigg|\mxi^{\rm all},\tilde{\mxi}^{\rm all}\bigg)\\
&\leq\sup_{x_j\in D_j}2n^{-C_0}\cdot\prod_{i=1}^n\big(1+\lambda_n^2\cdot\E(\|\mZ^i_n(x_j)\|^2|\mxi^{\rm all},\tilde{\mxi}^{\rm all})\cdot\exp(\lambda_n\cdot n^{\varsigma-1})\big)\\
&\leq2n^{-C_0}\cdot(1+O(1)\cdot n^{-1}\cdot(\log{n}))^n\\
&\leq2n^{-C_0}\cdot n^{O(1)}
\ea
almost surely, where $\lambda_n=n^{1/2}\cdot(\log{n})^{1/2}$. This gives (\ref{A claim 3}). Hence, the claim (\ref{A claim}) follows.

Define
\begin{align*}
\hat{\mDelta}_j^*(x_j)=\hat{\mf}_j(x_j)\om\mf_j(x_j)\om\check{\mf}_j^A(x_j)\om\mu_{0,j}(x_j)^{-1}\od\ma^*_j(x_j)\om\mc_j^*(x_j)\om\mr^{(1)}_j(x_j)\om\mr^{(2)}_j(x_j)
\end{align*}
for $1\leq j\leq d$. Then, from (\ref{equation 3 new new}) and using all the above approximations with the properties
\begin{align*}
\mb^*_j(x_j)&=\mdelta^*_j(x_j)~\text{for}~x_j\in  D^-_j(2h_j),\\
\mb_{jk}^*(x_j,x_k)&=\mdelta^*_{jk}(x_j,x_k)~\text{for}~x_j\in D_j~\text{and}~x_k\in  D^-_k(2h_k),
\end{align*}
we get
\begin{align*}
\hat{\mDelta}_j^*(x_j)=\tilde{\mDelta}_j^*(x_j)\om\bigoplus_{k\neq j}\int_{ D_k}\hat{\mDelta}_k^*(x_k)\od
\frac{\hat{p}^D_{jk}(x_j,x_k)}{\hat{p}^D_j(x_j)}dx_k\op\mr_j^*(x_j),\quad1\leq j\leq d,
\end{align*}
where $\mr_j^*: D_j\rightarrow\mbH$ is a stochastic map satisfying
\begin{align*}
\sup_{x_j\in D_j}\|\mr_j^*(x_j)\|&=\sum_{k\neq j}o_p(h_k^2)+\sum_{k=1}^dO_p\bigg(h_k\sqrt{\frac{\log n}{n\cdot h_j^{L_j}}}\bigg)
+O_p(n^{-1/2}\cdot(\log{n})^{1/2}\cdot h_j^{-L_j/2}\cdot V_{nj}^{1/2})\\
&\quad+O_p(\mathcal{A}_{nj})+O_p(b_n)+O_p(n^{-1/2}\cdot(\log{n})^{1/2}).
\end{align*}
Now, we define $\mDelta_{\op}^*=\bigoplus_{j=1}^d\mDelta^*_j$ and $\hat{\mDelta}_{\op}^*=\bigoplus_{j=1}^d\hat{\mDelta}^*_j$.
Then, by arguing as in the proof of Lemma 4 in Jeon et al. (2021b), we may prove that
\begin{align}\label{expansion including boundary}
\esssup_{x_j\in D_j}\|\hat{\mDelta}^*_j(x_j)\om\mDelta^*_j(x_j)\|=o_p\bigg(\sum_{j=1}^dh_j^2\bigg)+O_p(\mathcal{R}^*_n),
\end{align}
where $\esssup$ denotes the essential supremum and
\begin{align*}
\mathcal{R}^*_n=&\sum_{k=1}^dh_k\cdot\sum_{j=1}^d(nh^{L_j}/\log{n})^{-1/2}+\sum_{j=1}^d(n^{-1/2}\cdot(\log{n})^{1/2}\cdot h_j^{-L_j/2}\cdot V_{nj}^{1/2})+\sum_{j=1}^d\mathcal{A}_{nj}+b_n\\
&+n^{-1/2}\cdot(\log{n})^{1/2}.
\end{align*}
Hence, from (\ref{a and c}) and (\ref{expansion including boundary}), we get
\begin{align*}
\hat{\mf}_j(x_j)=\mf_j(x_j)\op\check{\mf}_j^A(x_j)\op\mc_j^*(x_j)\op\mDelta_j^*(x_j)\op o_p\bigg(\sum_{j=1}^dh_j^2\bigg)
\op O_p(\mathcal{R}^*_n)
\end{align*}
for almost everywhere $x_j\in D^-_j(2h_j)$ with respect to $\Leb_j$. We note that
\begin{align*}
\lim_{n\rightarrow\infty}n^{2/(L_{\rm max}+4)}\od\mc_j^*(x_j)=\alpha_j^2\od\mc_j(x_j)
\end{align*}
for $x_j\in  D^-_j(2h_j)$. We also note that $n^{2/(L_{\rm max}+4)}\cdot \sum_{j=1}^dh_j^2=O(1)$ and $n^{2/(L_{\rm max}+4)}\cdot \mathcal{R}^*_n=o(1)$. Thus, the second part of the theorem follows from Lemma \ref{A.I.} and a version of Proposition 4.8 in Van Neerven (2008). This completes the proof. \qed

\subsection{Proof of Proposition \ref{iRSC response}}

The proposition follows from (\ref{parellel inequality2}), (\ref{log max order}) and Theorem 6 in Lin and Yao (2019) by replacing $(\mathcal{M},Z)$ to $(\mathcal{N},W)$. \qed

\subsection{Proof of Theorem \ref{extension theorem}}

The first part of the theorem follows by arguing as in the proofs of Propositions \ref{existence}-\ref{convergence1 individual}. For the second part, we note that $\Gamma_{\hat{\mu}_{W},\mu_{W}}$ is linear. Hence, $\hat{\mf}_j(x_j)$, $\hat\mf(\mx)$, $\hat{\mf}^{[r]}_j(x_j)$ and $\hat\mf^{[r]}(\mx)$ are indeed estimators defined as in Section \ref{estimation method} with $\tilde{\mY}^i$ being replaced by $\Gamma_{\hat{\mu}_{W},\mu_{W}}(\Log_{\hat{\mu}_{W}}W^i)$. Thus, the second part follows by arguing as in the proofs of Theorems \ref{convergence2}-\ref{Asymptotic distribution} with Proposition \ref{iRSC response}. \qed

\subsection{Hilbertian intrinsic Riemannian singular component scores and intrinsic Riemannian singular component scores}\label{SC scores Rimannian Y}

\begin{customexample}{S.1}\label{Hilbertian Riemannian singular fd}
(Hilbertian Intrinsic Riemannian singular component scores) Recall the defi\-nition of $\mX$ and $\mbX=\mX\om_*\E(\mX)$ in Example \ref{Fully observed hilbert}. For $\mY={\rm{Log}}_{\mu_{W}}W$, define the cross-covariance operators $C_{\mX\mY}$ and $C_{\mY\mX}$ and the operator $C_{\mX\mY\mX}$ in the same way as in Example \ref{singular hilbert} with $\mbH=\mathfrak{T}(\mu_{W})$. For the orthonormal eigenvectors $\mphi_r$ of $C_{\mX\mY\mX}$ induced from the spectral theorem on $C_{\mX\mY\mX}$, we call the mean 0 random variables $\lng \mbX,\mphi_r\rng_*$ the Hilbertian intrinsic Riemannian singular component (HiRSC) scores of $\mX$.
To estimate $\lng \mbX,\mphi_r\rng_*$, we estimate $\E(\mX)$, $\mu_{W}$ and $\mphi_r$ in the same manner to Example \ref{Riemannian Hilbertian fd}. For the estimator $\lng \mX\om_*\bar{\mX},\hat{\mphi}_r\rng_*$ of $\lng\mbX,\mphi_r\rng_*$ obtained in this way, we have the following new proposition.

\begin{customproposition}{S.1}\label{HiRSC}
Assume that the first $r$ eigenvalues of $C_{\mX\mY\mX}$ have multiplicity one, that condition \ref{L} holds for $\mathcal{M}^*=\mathcal{N}$ and $Z^*=W$, and that $\E(\|\mX\|_*^2\|\mY\|_{\mathfrak{T}(\mu_W)}^2)<\infty$. Then, it holds that $\|\hat{C}_{\mX\mY\mX}-C_{\mX\mY\mX}\|_{\rm HS}=O_p(n^{-1/2})$. Also, for $\mphi_r$ with $\lng\hat{\mphi}_r,\mphi_r\rng_*>0$, it holds that $\|\hat{\mphi}_r\om_*\mphi_r\|_*=O_p(n^{-1/2})$. Moreover, $\max_{1\leq i\leq n}|\lng\mX^i\om_*\bar{\mX},\hat{\mphi}_r\rng_*-\lng\mX^i\om_*\E(\mX),\mphi_r\rng_*|$ achieves the following rates:
\ba
\begin{cases}
O_p(n^{-1/2+1/\tau}), & \text{if}\ \E(\|\mX\|_*^{\tau})<\infty\ \text{for some}\ \tau\geq4, \\
O_p(n^{-1/2}\cdot\log{n}), & \text{if}\ \E\big(\exp(c\cdot\|\mX\|_*)\big)<\infty\ \text{for some}\ c>0, \\
O_p(n^{-1/2}), & \text{if}\ \|\mX\|_*<C\ \text{almost surely for some}\ C>0.
\end{cases}
\ea
\end{customproposition}

\begin{proof}
The proposition follows from the proof of Proposition \ref{iRHSC} by switching the roles of $X$ and $\mY$.
\end{proof}

\end{customexample}

\begin{customexample}{S.2}\label{Riemannian singular fd}
(Intrinsic Riemannian singular component scores) Recall the definition of $X={\rm{Log}}_{\mu_{Z}}Z$ in Example \ref{Riemannian fd}. For $\mY={\rm{Log}}_{\mu_{W}}W$, define the cross-covariance operators $C_{X\mY}$ and $C_{\mY X}$ and the operator $C_{X\mY X}$ in the same way as in Example \ref{singular hilbert} with $\mbH=\mathfrak{T}(\mu_{W})$ and $\mbH_*=\mathfrak{T}(\mu_{Z})$. For the orthonormal eigenfunctions $\mphi_r$ of $C_{X\mY X}$ induced from the spectral theorem on $C_{X\mY X}$, we call the mean 0 random variables $\lng X,\mphi_r\rng_{\mathfrak{T}(\mu_{Z})}$ the intrinsic Riemannian singular component (iRSC) scores of $X$.
We estimate $\mu_{Z}$, $\mu_{W}$ and $\mphi_r$ in the same manner to Example \ref{Riemannian Hilbertian fd} with $n-1$ in the cross-covariance estimators being replaced by $n$. For the estimator $\lng\Log_{\hat{\mu}_{Z}}Z,\hat{\mphi}_r\rng_{\mathfrak{T}(\hat{\mu}_{Z})}$ of $\lng X,\mphi_r\rng_{\mathfrak{T}(\mu_{Z})}$ obtained in this way, we have the following new proposition. Recall the definition of $\Phi_{\hat{\mu}_{Z},\mu_{Z}}:\HS(\mathfrak{T}(\hat{\mu}_{Z}))\rightarrow \HS(\mathfrak{T}(\mu_{Z}))$ given in Example \ref{Riemannian Hilbertian fd}

\begin{customproposition}{S.2}\label{iRSC}
Assume that the first $r$ eigenvalues of $C_{X\mY X}$ have multiplicity one, that condition \ref{L} holds for both pairs $(\mathcal{M}^*,Z^*)=(\mathcal{M},Z)$ and $(\mathcal{M}^*,Z^*)=(\mathcal{N},W)$, and that $\E(\|X\|_{\mathfrak{T}(\mu_{Z})}^2\|\mY\|_{\mathfrak{T}(\mu_{W})}^2)<\infty$. Then, it holds that $\|\Phi_{\hat{\mu}_{Z},\mu_{Z}}(\hat{C}_{X\mY X})-C_{X\mY X}\|_{\rm HS}=O_p(n^{-1/2})$. Also, for $\mphi_r$ with $\lng\Gamma_{\hat{\mu}_{Z},\mu_{Z}}(\hat{\mphi}_r),\mphi_r\rng_{\mathfrak{T}(\mu_{Z})}>0$, it holds that $\|\Gamma_{\hat{\mu}_{Z},\mu_{Z}}(\hat{\mphi}_r)-\mphi_r\|_{\mathfrak{T}(\mu_{Z})}=O_p(n^{-1/2})$. Moreover, $\max_{1\leq i\leq n}|\lng \Log_{\hat{\mu}_{Z}}Z^i,\hat{\mphi}_r\rng_{\mathfrak{T}(\hat{\mu}_{Z})}-\lng \Log_{\mu_{Z}}Z^i,\mphi_r\rng_{\mathfrak{T}(\mu_{Z})}|$ achieves the following rates:
\ba
\begin{cases}
O_p(n^{-1/2+1/\tau}), & \text{if}\ \E(\|X\|_{\mathfrak{T}(\mu_Z)}^{\tau})<\infty\ \text{for some}\ \tau\geq2, \\
O_p(n^{-1/2}\cdot\log{n}), & \text{if}\ \E\big(\exp(c\cdot\|X\|_{\mathfrak{T}(\mu_Z)})\big)<\infty\ \text{for some}\ c>0, \\
O_p(n^{-1/2}), & \text{if}\ \|X\|_{\mathfrak{T}(\mu_Z)}<C\ \text{almost surely for some}\ C>0.
\end{cases}
\ea
\end{customproposition}

\begin{proof}
Define $\hat{C}_{X\mY}^{\Gamma}:\mathfrak{T}(\mu_{W})\rightarrow\mathfrak{T}(\mu_{Z})$ and $\hat{C}_{\mY X}^{\Gamma}:\mathfrak{T}(\mu_{Z})\rightarrow\mathfrak{T}(\mu_{W})$ by
\ba
\hat{C}_{X\mY}^{\Gamma}(\cdot)&=\frac{1}{n}\sum_{i=1}^n(\lng\Gamma_{\hat{\mu}_{W},\mu_{W}}(\Log_{\hat{\mu}_{W}}W^i),\cdot\rng_{\mathfrak{T}(\mu_{W})}\cdot\Gamma_{\hat{\mu}_{Z},\mu_{Z}}(\Log_{\hat{\mu}_{Z}}Z^i)),\\
\hat{C}_{\mY X}^{\Gamma}(\cdot)&=\frac{1}{n}\sum_{i=1}^n(\lng\Gamma_{\hat{\mu}_{Z},\mu_{Z}}(\Log_{\hat{\mu}_{Z}}Z^i),\cdot\rng_{\mathfrak{T}(\mu_{Z})}\cdot\Gamma_{\hat{\mu}_{W},\mu_{W}}(\Log_{\hat{\mu}_{W}}W^i)).
\ea
and define $\hat{C}_{X\mY X}^{\Gamma}=\hat{C}_{X\mY}^{\Gamma}\circ\hat{C}_{\mY X}^{\Gamma}:\mathfrak{T}(\mu_{Z})\rightarrow\mathfrak{T}(\mu_{Z})$. A direct computation with Proposition 2 in Lin and Yao (2019) shows that $\Phi_{\hat{\mu}_{Z},\mu_{Z}}(\hat{C}_{X\mY X})=\hat{C}_{X\mY X}^{\Gamma}$.
By arguing as in (\ref{cross covariance decomposition}), we get
\begin{align}\label{cross covariance decomp3}
\begin{split}
&\|\Phi_{\hat{\mu}_{Z},\mu_{Z}}(\hat{C}_{X\mY X})-C_{X\mY X}\|_{\HS}\\
&=\|\hat{C}_{X\mY X}^{\Gamma}-C_{X\mY X}\|_{\HS}\\
&\leq\|\hat{C}^{\Gamma}_{X\mY}-C_{X\mY}\|_{\HS}\cdot\|\hat{C}^{\Gamma}_{\mY X}-C_{\mY X}\|_{\HS}+\|\hat{C}^{\Gamma}_{X\mY}-C_{X\mY}\|_{\HS}\cdot\|C_{\mY X}\|_{\HS}\\
&\quad+\|C_{X\mY}\|_{\HS}\cdot\|\hat{C}^{\Gamma}_{\mY X}-C_{\mY X}\|_{\HS}.
\end{split}
\end{align}
We first claim that $\|\hat{C}^{\Gamma}_{X\mY}-C_{X\mY}\|^2_{\HS}=O_p(n^{-1})$. Note that
\begin{align}\label{check c gamma decomp2}
\begin{split}
&\check{C}^{\Gamma}_{X\mY}(\cdot)\\
&=\frac{1}{n}\sum_{i=1}^n\lng\Gamma_{\hat{\mu}_{W},\mu_{W}}(\Log_{\hat{\mu}_{W}}W^i)-\Log_{\mu_{W}}W^i,\cdot\rng_{\mathfrak{T}(\mu_{W})}\cdot(\Gamma_{\hat{\mu}_{Z},\mu_{Z}}(\Log_{\hat{\mu}_{Z}}Z^i)-\Log_{\mu_{Z}}Z^i)\\
&\quad+\frac{1}{n}\sum_{i=1}^n\lng\Log_{\mu_{W}}W^i,\cdot\rng_{\mathfrak{T}(\mu_{W})}\cdot(\Gamma_{\hat{\mu}_{Z},\mu_{Z}}(\Log_{\hat{\mu}_{Z}}Z^i)-\Log_{\mu_{Z}}Z^i)\\
&\quad+\frac{1}{n}\sum_{i=1}^n\lng\Gamma_{\hat{\mu}_{W},\mu_{W}}(\Log_{\hat{\mu}_{W}}W^i)-\Log_{\mu_{W}}W^i,\cdot\rng_{\mathfrak{T}(\mu_{W})}\cdot\Log_{\mu_{Z}}Z^i\\
&\quad+\frac{1}{n}\sum_{i=1}^n\lng\Log_{\mu_{W}}W^i,\cdot\rng_{\mathfrak{T}(\mu_{W})}\cdot\Log_{\mu_{Z}}Z^i.
\end{split}
\end{align}
The sum of the squared Hilbert-Schmidt norms of the first three operators on the right hand side of (\ref{check c gamma decomp2}) achieves the rate $O_p(n^{-2})+O_p(n^{-1})+O_p(n^{-1})$. This follows by arguing as in the proof of Proposition \ref{iRHSC}. Also, the last operator on the right hand side of (\ref{check c gamma decomp2}) is an unbiased estimator of $C_{X\mY}(\cdot)$ and hence
\ba
\bigg\|\frac{1}{n}\sum_{i=1}^n\lng\Log_{\mu_{W}}W^i,\cdot\rng_{\mathfrak{T}(\mu_{W})}\cdot\Log_{\mu_{Z}}Z^i-C_{X\mY}(\cdot)\bigg\|^2_{\HS}=O_p(n^{-1}).
\ea
This proves the claim. Similarly, one can prove that $\|\hat{C}^{\Gamma}_{\mY X}-C_{\mY X}\|^2_{\HS}=O_p(n^{-1})$. Plugging those rates to (\ref{cross covariance decomp3}) gives that $\|\hat{C}_{X\mY X}^{\Gamma}-C_{X\mY X}\|_{\HS}=O_p(n^{-1/2})$. Now, using Proposition 2 in Lin and Yao (2019) and Lemma 2.3 in Horv\'{a}th and Kokoszka (2012), we get $\|\Gamma_{\hat{\mu}_{Z},\mu_{Z}}(\hat{\mphi}_r)-\mphi_r\|_{\mathfrak{T}(\mu_{Z})}=O_p(n^{-1/2})$.
The desired result follows by arguing as in the proof of Proposition \ref{iRFPC}.
\end{proof}

\end{customexample}

\bigskip
\noindent
\large{\textbf{References for Supplementary Material}

\bigskip
\normalsize
\noindent
1. Atkinson, K. and Han, W. (2009). {\it Theoretical Numerical Analysis}. Springer-Verlag New York.\\
2. Beltrami, E. J. (1967). On infinite-dimensional convex programs. {\it Journal of Computer and System Sciences}, \textbf{1}, 323-329.\\
3. Bosq, D. (2000). {\it Linear Processes in Function Spaces}. Springer-Verlag New York.\\
4. Horv\'{a}th, L. and Kokoszka, P. (2012). {\it Inference for Functional Data with Applications}. Springer New York.\\
5. Jeon, J. M. and Park, B. U. (2020). Additive regression with Hilbertian responses. {\it Annals of Statistics}, \textbf{48}, 2671-2697.\\
6. Jeon, J. M., Park, B. U. and Van Keilegom, I. (2021a). Additive regression for non-Euclidean responses and predictors. {\it Annals of Statistics}, \textbf{49}, 2611-2641.\\
7. Jeon, J. M., Park, B. U. and Van Keilegom, I. (2021b). Additive regression with predictors of variours natures and possibly incomplete Hilbertian responses. {\it Electronic Journal of Statistics}, \textbf{15}, 1473-1548.\\
8. Kallenberg, O. (2017). {\it Random Measures, Theory and Applications}. Springer International Publishing Switzerland.\\
9. Kendall, W. S. and Le, H. (2011). Limit theorems for empirical Fr\'{e}chet means of independent and non-identically distributed manifold-valued random variables. {\it Brazilian Journal of Probability and Statistics}, \textbf{25}, 323-352.\\
10. Kundu, S., Majumdar, S. and Mukherjee, K. (2000). Central limit theorems revisited. {\it Statistics and Probability Letters}, \textbf{47}, 265-275.\\
11. Lee, J. M. (2018). {\it Introduction to Riemannian Manifolds}. Springer International Publishing AG.\\
12. Lin, Z., M\"uller, H.-G. and Park, B. U. (2022). Additive models for symmetric positive-definite matrices and Lie groups. To appear in {\it Biometrika}.\\
13. Lin, Z. and Yao, F. (2019). Intrinsic Riemannian functional data analysis. {\it Annals of Statistics}, \textbf{47}, 3533-3577.\\
14. Sacks, P. (2017). {\it Techniques of functional analysis for differential and integral equations}. Academic Press, London.\\
15. Tu, L. W. (2017). {\it Differential Geometry}. Springer International Publishing AG.\\
16. Van Neerven, J. (2008). {\it Stochastic evolution equations}. Lecture Notes of the 11th Internet Seminar, TU Delft OpenCourseWare, http://ocw.tudelft.nl.


\begin{thebibliography}{9}
\bibitem[Abdulaziz(2021)]{Abdulaziz (2021)} Abdulaziz, A. (2021). A review of compositional data analysis and recent advances. {\it Communications in Statistics - Theory and Methods}, DOI: 10.1080/03610926.2021.2014890.
\bibitem[Bertrand et al.(2019)]{Bertrand et al. (2019)} Bertrand, A., Van Keilegom, I. and Legrand, C. (2019). Flexible parametric approach to classical measurement error variance estimation without auxiliary data. {\it Biometrics}, \textbf{75}, 297-307.
\bibitem[Bunkure(2019)]{Bunkure (2019)} Bunkure, J. K. (2019). Lebegue-Bochner spaces and evolution triples. {\it International Journal of Mathematics And its Applications}, \textbf{7}, 41-52.
\bibitem[Cai and Hall(2006)]{Cai and Hall (2006)} Cai, T. T. and Hall, P. (2006). Prediction in functional linear regression. {\it Annals of Statistics}, \textbf{34}, 2159-2179.
\bibitem[Cheng and Wu(2013)]{Cheng and Wu (2013)} Cheng, M.-Y. and Wu, H.-T. (2013). Local linear regression on manifolds and its geometric interpretation. {\it Journal of the American Statistical Association}, \textbf{108}, 1421-1434.
\bibitem[Cohn(2013)]{Cohn (2013)} Cohn, D. L. (2013). {\it Measure Theory}. Birkh\"{a}user Basel.
\bibitem[Dai and M\"{u}ller(2018)]{Dai and Muller (2018)} Dai, X. and M\"{u}ller, H.-G. (2018). Principal component analysis for functional data on Riemannian manifolds and spheres. {\it Annals of Statistics}, \textbf{46}, 3334-3361.
\bibitem[Delaigle(2008)]{Delaigle (2008)} Delaigle, A. (2008). An alternative view of the deconvolution problem. {\it Statistica Sinica}, \textbf{18}, 1025-1045.
\bibitem[Delaigle et al.(2008)]{Delaigle et al. (2008)} Delaigle, A., Hall, P. and Meister, A. (2008). On deconvolution with repeated measurements. {\it Annals of Statistics}, \textbf{36}, 665-685.
\bibitem[Fan(1991)]{Fan (1991)} Fan, J. (1991). On the optimal rates of convergence for nonparametric deconvolution problems. {\it Annals of Statistics}, \textbf{19}, 1257-1272.
\bibitem[Fan(1992)]{Fan (1992)} Fan, J. (1992). Deconvolution with supersmooth distributions. {\it The Canadian Journal of Statistics}, \textbf{20}, 155-169.
\bibitem[Hall and Horowitz(2007)]{Hall and Horowitz (2007)} Hall, P. and Horowitz, J.~L. (2007). Methodology and convergence rates for functional linear regression. {\it Annals of Statistics}, \textbf{35}, 70-91.
\bibitem[Hall and Vial(2006)]{Hall and Vial (2006)} Hall, P. and Vial, C. (2006). Assessing the finite dimensionality of functional data. {\it Journal of the Royal Statistical Society: Series B (Statistical Methodology)}, \textbf{68}, 689-705.
\bibitem[Han et al.(2018)]{Han et al. (2018)} Han, K., M\"uller, H.-G. and Park, B. U. (2018). Smooth backfitting for additive modeling with small errors-in-variables, with an application to additive functional regression for multiple predictor functions. {\it Bernoulli}, \textbf{24}, 1233-1265.
\bibitem[Han et al.(2020)]{Han et al. (2020)} Han, K., M\"uller, H.-G. and Park, B. U. (2020). Additive functional regression for densities as responses. {\it Journal of the American Statistical Association}, \textbf{115}, 997-1010.
\bibitem[Han and Park(2018)]{Han and Park (2018)} Han, K. and Park, B. U. (2018). Smooth backfitting for error-in-variables additive models. {\it Annals of Statistics}, \textbf{46}, 2216-2250.
\bibitem[Hron et al.(2016)]{Hron et al. (2016)} Hron, K., Menafoglio, A., Templ, M., Hr\.{u}zov\'{a} and Filzmoser, P. (2016). Simplicial principal component analysis for density functions in Bayes spaces. {\it Computational Statistics and Data Analysis}, \textbf{94}, 330-350.
\bibitem[Hsing and Eubank(2015)]{Hsing and Eubank (2015)} Hsing, T. and Eubank, R. (2015). {\it Theoretical Foundations of Functional Data Analysis, with an Introduction to Linear Operators}. John Wiley \& Sons.
\bibitem[Jeon and Park(2020)]{Jeon and Park (2020)} Jeon, J. M. and Park, B. U. (2020). Additive regression with Hilbertian responses. {\it Annals of Statistics}, \textbf{48}, 2671-2697.
\bibitem[Jeon et al.(2021a)]{Jeon et al. (2021a)} Jeon, J. M., Park, B. U. and Van Keilegom, I. (2021a). Additive regression for non-Euclidean responses and predictors. {\it Annals of Statistics}, \textbf{49}, 2611-2641.
\bibitem[Jeon et al.(2021b)]{Jeon et al. (2021b)} Jeon, J. M., Park, B. U. and Van Keilegom, I. (2021b). Additive regression with predictors of variours natures and possibly incomplete Hilbertian responses. {\it Electronic Journal of Statistics}, \textbf{15}, 1473-1548.
\bibitem[Jeon et al.(2022)]{Jeon et al. (2022)} Jeon, J. M., Park, B. U. and Van Keilegom, I. (2022). Nonparametric regression on Lie groups with measurement errors. {\it Annals of Statistics}, \textbf{50}, 2973-3008.
\bibitem[Johannes(2009)]{Johannes (2009)} Johannes, J. (2009). Deconvolution with unknown measurement error distribution. {\it Annals of Statistics}, \textbf{37}, 2301-2323.
\bibitem[Lee et al.(2010)]{Lee et al. (2010)} Lee, Y. K., Mammen, E. and Park, B. U. (2010). Backfitting and smooth backfitting for additive quantile models. {\it Annals of Statistics}, \textbf{38}, 2857-2883.
\bibitem[Lee et al.(2012)]{Lee et al. (2012)} Lee, Y. K., Mammen, E. and Park, B. U. (2012). Flexible generalized varying coefficient regression models. {\it Annals of Statistics}, \textbf{40}, 1906-1933.
\bibitem[Lin et al.(2022)]{Lin et al. (2022)} Lin, Z., M\"uller, H.-G. and Park, B. U. (2022). Additive models for symmetric positive-definite matrices and Lie groups. To appear in {\it Biometrika}.
\bibitem[Lin and Yao(2019)]{Lin and Yao (2019)} Lin, Z. and Yao, F. (2019). Intrinsic Riemannian functional data analysis. {\it Annals of Statistics}, \textbf{47}, 3533-3577.
\bibitem[Linton et al. (2008)]{Linton et al. (2008)} Linton, O., Sperlich, S. and Van Keilegom, I. (2008). Estimation of a semiparametric transformation model. {\it Annals of Statistics}. \textbf{36}, 686-718.
\bibitem[Maier et al.(2021)]{Maier et al. (2021)} Maier, E.-M., St\"{o}cker, A., Fitzenberger, B. and Greven, S. (2021). Additive density-on-scalar regression in Bayes
Hilbert spaces with an application to gender economics. arXiv:2110.11771v1.
\bibitem[Mammen et al.(1999)]{Mammen et al. (1999)} Mammen, E., Linton, O. B. and Nielsen, J. P. (1999). The existence and asymptotic properties of a backfitting projection algorithm under weak conditions. {\it Annals of Statistics}, \textbf{27}, 1443-1490.
\bibitem[Mammen and Park(2005)]{Mammen and Park (2005)} Mammen, E. and Park, B. U. (2005). Bandwidth selection for smooth backfitting in additive models. {\it Annals of Statistics}, \textbf{33}, 1260-1294.
\bibitem[Mammen and Park(2006)]{Mammen and Park (2006)} Mammen, E. and Park, B. U. (2006). A simple smooth backfitting method for additive models. {\it Annals of Statistics}, \textbf{34}, 2252-2271.
\bibitem[Lafaye de Micheaux et al.(2019)]{Lafaye de Micheaux et al. (2019)} Lafaye de Micheaux, P., Liquet, B. and Sutton, M. (2019). PLS for Big Data: A unified parallel
algorithm for regularised group PLS. {\it Statistics Surveys}, \textbf{13}, 119-149.
\bibitem[Park et al.(2018)]{Park et al. (2018)} Park, B. U., Chen, C.-J., Tao, W. and M\"uller, H.-G. (2018). Singular additive models for function to function regression. {\it Statistica Sinica}, \textbf{28}, 2497-2520.
\bibitem[Pelletier(2005)]{Pelletier (2005)} Pelletier, B. (2005). Kernel density estimation on Riemannian manifolds. {\it Statistics and Probability Letters}, \textbf{73}, 297-304.
\bibitem[Petersen and M\"{u}ller(2016)]{Petersen and Muller (2016)} Petersen, A. and M\"{u}ller, H.-G. (2016). Functional data analysis for density functions by transformation to a Hilbert space. {\it Annals of Statistics}, \textbf{44}, 183-218.
\bibitem[Petersen and M\"{u}ller(2019)]{Petersen and Muller (2019)} Petersen, A. and M\"{u}ller, H.-G. (2019). Fr\'{e}chet regression for random objects with Euclidean predictors. {\it Annals of Statistics}, \textbf{47}, 691-719.
\bibitem[Preston(2008)]{Preston (2008)} Preston, C. (2008). Some notes on standard Borel and related spaces. arXiv:0809.3066.
\bibitem[Ramsay and Silverman(2005)]{Ramsay and Silverman (2005)} Ramsay, J. and Silverman, B. (2005). {\it Functional Data Analysis}. Springer New York.
\bibitem[Talsk\'{a} et al.(2018)]{Talska et al. (2018)} Talsk\'{a}, R., Menafoglio, A., Machalov\'{a}, J., Hron, K. and Fi\v{s}erov\'{a}, E. (2018). Compositional regression with functional response. {\it Computational Statistics and Data Analysis}, \textbf{123}, 66-85.
\bibitem[van Es and Gugushvili(2010)]{van Es and Gugushvili (2010)} van Es, B. and Gugushvili, S. (2010). Asymptotic normality of the deconvolution kernel density estimator under the vanishing error variance. {\it Journal of the Korean Statistical Society}, \textbf{39}, 103-115.
\bibitem[Wang et al.(2015)]{Wang et al. (2015)} Wang, H., Shangguan, L., Guan, R. and Billard, L. (2015). Principal component analysis for compositional data vectors. {\it Computational Statistics}, \textbf{30}, 1079-1096.
\bibitem[Wang et al.(2016)]{Wang et al. (2016)} Wang, J.-L., Chiou, J.-M. and M\"{u}ller, H.-G. (2016). Functional data analysis. {\it Annual Review of Statistics and Its Application}, \textbf{3}, 257-295.
\bibitem[Yang et al.(2011)]{Yang et al. (2011)} Yang, W., M\"uller, H.-G. and Stadtm\"uller, U. (2011). Functional singular component analysis. {\it Journal of the Royal Statistical Society: Series B (Statistical Methodology)}, \textbf{73}, 303-324.
\bibitem[Yu et al.(2011)]{Yu et al. (2011)} Yu, K., Mammen, E. and Park, B. U. (2011). Semi-parametric regression: Efficiency gains from modeling the nonparametric part.
{\it Bernoulli} \textbf{17}, 736-748.
\bibitem[Yu et al.(2008)]{Yu et al. (2008)} Yu, K., Park, B. U. and Mammen, E. (2008). Smooth backfitting in generalized additive models. {\it Annals of Statistics}, \textbf{36}, 228-260.
\end{thebibliography}
\end{document}